\documentclass[a4paper, 12pt]{amsart}
\usepackage[utf8]{inputenc}
\usepackage[OT2, T1]{fontenc}

\usepackage{fancyhdr}

\usepackage{subfiles} 

\usepackage{amsthm}
\usepackage{amssymb}
\usepackage{mathtools}
 \usepackage{stmaryrd}	
\usepackage{mathrsfs}
\usepackage{bm}
\usepackage{tikz-cd}
\usepackage{array} 
\usetikzlibrary{positioning, shapes.geometric, patterns, graphs, decorations.pathmorphing}

\usepackage{paralist}

\usepackage[colorlinks]{hyperref}
 \hypersetup{
	linkcolor=blue,
}

\newcommand{\Z}{\ensuremath{\mathbb{Z}}} 
\newcommand{\bZ}{\ensuremath{\breve{\Z}}} 
\newcommand{\bF}{\ensuremath{\breve{F}}} 
\newcommand{\Q}{\ensuremath{\mathbb{Q}}} 
\newcommand{\bQ}{\ensuremath{\breve{\Q}_p}} 

\newcommand{\R}{\ensuremath{\mathbb{R}}} 
\newcommand{\CC}{\ensuremath{\mathbb{C}}} 
\newcommand{\A}{\ensuremath{\mathbb{A}}} 

\newcommand{\gr}{\operatorname{gr}} 
\newcommand{\Tor}{\operatorname{Tor}} 

\newcommand{\Stab}{\ensuremath{\mathrm{Stab}}} 

\newcommand{\Gal}{\operatorname{Gal}}  

\newcommand{\lrangle}[1]{\ensuremath{\langle #1 \rangle}}
\newcommand{\lrbracket}[1]{\ensuremath{\{ #1 \}}}

\newcommand{\wdt}[1]{\ensuremath{\widetilde{#1}}} 
\newcommand{\wdh}[1]{\ensuremath{\widehat{#1}}} 
\newcommand{\ovl}[1]{\ensuremath{\overline{#1}}} 

\newcommand{\pr}{\operatorname{pr}} 

\newcommand{\identity}{\ensuremath{\mathrm{id}}}
\newcommand{\prolim}{\ensuremath{\underleftarrow{\lim}}}

\newcommand{\Hom}{\operatorname{Hom}}

\newcommand{\rightiso}{\ensuremath{\stackrel{\sim}{\rightarrow}}}

\newcommand{\Ker}{\operatorname{ker}}

\newcommand{\MH}{\operatorname{\mathrm{MH}}}

\newcommand{\Lie}{\operatorname{Lie}}
\newcommand{\Ad}{\operatorname{Ad}}
\newcommand{\ad}{\ensuremath{\mathrm{ad}}}  
\newcommand{\der}{\ensuremath{\mathrm{der}}} 

\newcommand{\Spec}{\operatorname{Spec}}
\newcommand{\Proj}{\operatorname{Proj}}
\newcommand{\Gm}{\ensuremath{\mathbb{G}_\mathrm{m}}}

\newcommand{\Res}{\operatorname{Res}}

\newcommand{\et}{\textit{\'et}} 

\newcommand{\Spf}{\operatorname{Spf}}

\newcommand{\Spd}{\operatorname{Spd}}

\newcommand{\loc}{\ensuremath{\mathrm{loc}}} 
\newcommand{\XX}{\ensuremath{\mathcal{X}}} 

\newcommand{\GL}{\operatorname{GL}}

\newcommand{\GSp}{\operatorname{GSp}}
\newcommand{\GSP}{\operatorname{\mathcal{GSP}}}
\newcommand{\GLL}{\operatorname{\mathcal{GL}}}

\newcommand{\cl}{\ensuremath{\mathrm{cl}}} 

\newcommand{\lieg}{\ensuremath{\mathfrak{g}}} 

\newcommand{\OO}{\ensuremath{\mathcal{O}}} 

\newcommand{\DS}{\ensuremath{\mathbb{S}}} 
\newcommand{\tor}{\ensuremath{\mathrm{tor}}} 
\newcommand{\Shum}[1]{\ensuremath{\mathscr{S}_{#1}}} 
\newcommand{\shu}[1]{\ensuremath{\mathrm{Sh}_{#1}}} 
\newcommand{\shuc}[1]{\ensuremath{\mathrm{Sh}_{#1}^{\tor}}}
\newcommand{\shum}[1]{\ensuremath{\mathrm{Sh}_{#1}^{\min}}}
\newcommand{\Shumc}[1]{\ensuremath{\mathscr{S}_{#1}^{\tor}}}
\newcommand{\Shumm}[1]{\ensuremath{\mathscr{S}_{#1}^{\min}}}

\newcommand{\Cusp}{\operatorname{Cusp}} 

\newcommand{\bigsur}{\ensuremath{\star}} 

\newcommand{\imin}[1]{\ensuremath{i_{#1}^{\mathrm{min}}}}
\newcommand{\itor}[1]{\ensuremath{i_{#1}^{\tor}}}
\newcommand{\Zb}{\ensuremath{\mathrm{Z}}} 

\newcommand{\GG}{\ensuremath{\mathcal{G}}}

\newcommand{\GGc}{\ensuremath{\mathcal{G}^{\circ}}}

\newcommand{\KK}{\ensuremath{\Breve{K}}} 
\newcommand{\KKc}{\ensuremath{\Breve{K}^{\circ}}} 
\newcommand{\KKf}{\Breve{K}^{\flat}} 
\newcommand{\KKn}{\Breve{K}^{\natural}} 

\newcommand{\Kc}{\ensuremath{K}^{\circ}} 
\newcommand{\Kf}{\ensuremath{K}^{\flat}} 
\newcommand{\Kn}{\ensuremath{K}^{\natural}} 

\newcommand{\red}{\ensuremath{\textrm{red}}} 
\newcommand{\ext}{\ensuremath{\textrm{ext}}} 

\newcommand{\spl}{\ensuremath{\textrm{sp}}} 

\newcommand{\Gg}{\ensuremath{\mathfrak{G}}} 
\newcommand{\Ggh}{\ensuremath{\hat{\mathfrak{G}}}} 
\newcommand{\Ggc}{\ensuremath{\mathfrak{G}^{\circ}}} 
\newcommand{\Ggf}{\ensuremath{\mathfrak{G}^{\flat}}}

\newcommand{\bc}{\ensuremath{\mathfrak{C}}} 

\newcommand{\MPF}{\ensuremath{\mathrm{MP}}} 

\newcommand{\Bui}{\operatorname{\mathcal{B}}} 
\newcommand{\Aff}{\operatorname{Aff}} 
\newcommand{\alc}{\ensuremath{\mathfrak{a}}} 
\newcommand{\aff}{\ensuremath{\textrm{aff}}} 
\newcommand{\FF}{\ensuremath{\mathcal{F}}} 

\newcommand{\ab}{\ensuremath{\mathcal{A}}} 

\newcommand{\ZZ}{\ensuremath{\mathfrak{Z}}} 
\newcommand{\TT}{\ensuremath{\mathfrak{T}}} 
\newcommand{\UU}{\ensuremath{\mathfrak{U}}} 

\newcommand{\dd}{\ensuremath{\mathrm{\ddagger}}} 

\theoremstyle{plain}
\newtheorem{proposition}{Proposition}
\newtheorem{lemma}[proposition]{Lemma}
\newtheorem{theorem}[proposition]{Theorem}
\newtheorem{corollary}[proposition]{Corollary}

\theoremstyle{definition}
\newtheorem{definition}[proposition]{Definition}
\newtheorem{definition-theorem}[proposition]{Definition-Theorem}
\newtheorem{definition-proposition}[proposition]{Definition-Proposition}
\newtheorem{remark}[proposition]{Remark}

\newtheorem{example}[proposition]{Example}

\theoremstyle{definition}

\theoremstyle{plain}

\numberwithin{equation}{section}
\numberwithin{table}{section}
\numberwithin{proposition}{section}
\numberwithin{conj}{section}	

\usepackage{amsmath}

\usepackage{geometry}
\usepackage{lmodern}

\title[Boundary structures]{Boundary structures of integral models of Hodge-type Shimura Varieties}
\author[Shengkai Mao]{Shengkai Mao}
\address{Morningside Center of Mathematics}
\email{maoshengkai@amss.ac.cn}

\date{}

\begin{document}

\begin{abstract}

  We compute the level groups associated with mixed Shimura varieties that appear at the boundaries of compactifications of Shimura varieties and show that the boundaries of minimal compactifications of Pappas-Rapoport integral models are finite quotients of smaller Pappas-Rapoport integral models. Additionally, we prove that the compactifications of integral models of Hodge-type Shimura varieties with quasi-parahoric level structures are independent of the choice of Siegel embedding, and use this to construct and analyze the change-of-parahoric morphisms on these compactifications.
\end{abstract}

\maketitle

\tableofcontents

\section{Introduction}


   Shimura varieties play a central role in number theory, and their compactifications are an important topic of study.

    Let $(G, X)$ be a Shimura datum, $\shu{K}:=\shu{K}(G, X)$, we focus on the \emph{minimal compactification} $\shum{K}$ and the \emph{toroidal compactification} $\shuc{K, \Sigma}$. See \cite{baily1966compactification} (resp. \cite{ash2010smooth}) for constructions of the minimal (resp. toroidal) compactifications, and see \cite{pink1989arithmetical} for compactifications of \emph{mixed} Shimura varieties. The constructions are initially carried out over $\CC$ and subsequently descended to $E$ under certain conditions.
     
     Roughly speaking, the \emph{minimal compactification} $\shum{K}$ of $\shu{K}$ is a normal projective variety over $E=E(G,X)$. Moreover, it is \emph{canonical}, i.e., it does not depend on the choice of extra data. Its boundary strata are themselves quotients of Shimura varieties. Unfortunately, $\shum{K}$ might have uncontrolled singularities along the boundaries. The \emph{toroidal compactification} $\shuc{K, \Sigma}$ is a normal algebraic space over $\CC$ which dominants $\shum{K}$, it depends on the choice of a \emph{cone decomposition} $\Sigma$. If $\Sigma$ is smooth or projective, then $\shuc{K, \Sigma}$ is an algebraic variety. We can furthur refine $\Sigma$ such that $\shuc{K, \Sigma}$ is a smooth, projective algebraic variety over $E$, and the boundaries of $\shuc{K, \Sigma}$ are normal crossings divisors. The boundary charts of $\shuc{K, \Sigma}$ are toroidal embeddings of torus-torsors over families of abelian varieties over the boundary strata of $\shum{K}$.

    \subsubsection{}
    
     Our first main result is the study of level groups of mixed Shimura varieties appear at boundary. 
     
     First, we work in an absolute setting. This is a pure group-theoretical problem and can be separated from the theory of Shimura varieties. Let $G$ be a reductive group over $F$, where $F$ is a finite extension of $\Q_p$ or the maximal unramified extension of some finite extension of $\Q_p$, $Q \subset G$ be a parabolic subgroup, $\pi: Q \to L$ be the Levi quotient, $G_h \subset L$ be a normal subgroup of $L$. Let $K \subset G(F)$ be an open subgroup, we determine $K_L:=\pi(K \cap Q(F)) \subset L(F)$ and $K_h:= G_h(F) \cap K_L \subset G_h(F)$.

     \begin{proposition}[{Proposition \ref{prop: main prop for section BT theory}}]\label{prop: main prop for BT theory, intro}\leavevmode
        \begin{enumerate}
            \item When $K \subset G(F)$ is a stablizer quasi-parahoric subgroup (resp. quasi-parahoric subgroup, neutral quasi-parahoric subgroup, Moy-Prasad subgroup with level $r > 0$), then $K_L \subset L(F)$ and $K_h \subset G_h(F)$ are stablizer quasi-parahoric subgroups (resp. quasi-parahoric subgroups, neutral quasi-parahoric subgroups, Moy-Prasad subgroups with level $r > 0$) respectively.
            \item When $K \subset G(F)$ is a parahoric subgroup, then $K_L$ is a parahoric subgroup, but $K_h$ might not be a parahoric subgroup. $K_h$ is a parahoric subgroup in some special cases, for example, when $G$ is unramified or when $G_h \to L$ has a section.
        \end{enumerate}
     \end{proposition}

     Second, we work in a relative setting. Fix a Shimura datum $(G, X)$ and a cusp label representative $[\Phi] \in \Cusp_K(G, X)$, we have an admissable parabolic subgroup $Q_{\Phi}$ and a specific normal subgroup $P_{\Phi}$ of $Q_{\Phi}$. Let $W_{\Phi} \subset P_{\Phi}$ be the unipotent radical, it is also the unipotent radical of $Q_{\Phi}$. Let $U_{\Phi}$ be the center of $W_{\Phi}$ and $V_{\Phi}$ be the quotient $W_{\Phi}/U_{\Phi}$. Let $\pi_h: P_{\Phi} \to G_{h, \Phi}$ and $\pi_L: Q_{\Phi} \to L_{\Phi}$ be Levi quotients, and $G_{\Phi, l} = Q_{\Phi}/P_{\Phi}$. Then one has an exact sequence of unipotent groups (resp. reductive groups):
     \begin{equation}
         1 \to U_{\Phi} \to W_{\Phi} \stackrel{\pi_V}{\to} V_{\Phi} \to 1,\quad (\textit{resp.}\ 1 \to G_{h, \Phi} \to L_{\Phi} \to G_{l, \Phi} \to 1).
     \end{equation}
     A cusp label representative $\Phi$ determines an element $g \in G(\A_f)$. Given an open compact subgroup $K \subset G(\A_f)$, we denote $K_{\Phi, h} = \pi_h(gKg^{-1} \cap P_{\Phi}(\A_f))$, $K_{\Phi, W} = gKg^{-1} \cap W_{\Phi}(\A_f)$, $K_{\Phi, V} = \pi_{V}(K_{\Phi, W})$. 
     
     Recall that the boundary of minimal compactification $\shum{K}(G, X)$ of $\shu{K}(G, X)$ is stratified by finite quotients of Shimura varieties $\shu{K_{\Phi, h}}(G_{\Phi, h}, X_{\Phi, h})$, indexed by cusp labels representatives $[\Phi] \in \Cusp_K(G, X)$. Fix an embedding of Shimura data $(G_1, X_1) \hookrightarrow (G_2, X_2)$, and neat level groups $K_1 \subset K_2$, the induced finite morphism $\iota: \shu{K_1}(G_1, X_1) \to  \shu{K_2}(G_2, X_2) $ extends to a finite morphism $\iota^{\min}: \shum{K_1}(G_1, X_1) \to  \shum{K_2}(G_2, X_2)$ that is compatible with the induced morphism of Shimura varieties at boundary 
     \begin{equation}
         \iota_{\Phi}: \shu{K_{1, \Phi_1, h}}(G_{1, \Phi_1, h}, X_{1, \Phi_1, h}) \to \shu{K_{2, \Phi_2, h}}(G_{2, \Phi_2, h}, X_{2, \Phi_2, h}),
     \end{equation}
     for all $[\Phi_1] \in \Cusp_{K_1}(G_1, X_1)$ and $[\Phi_2] = [\iota_*\Phi_1] \in \Cusp_{K_2}(G_2, X_2)$. The following two conditions are important when we understand the boundary structures of Shimura varieties.
     \begin{equation}\label{eq: G_h}
         K_{1, \Phi_1, h, p} = G_{1, \Phi_1, h}(\Q_p) \cap K_{2, \Phi_2, h, p},
     \end{equation}
     \begin{equation}\label{eq: V}
         K_{1, \Phi_1, V_1, p} = V_{1, \Phi_1}(\Q_p) \cap K_{2, \Phi_2, V_2, p}
     \end{equation}
     
     \begin{proposition}[{Corollary \ref{cor: very nice emb implies very nice emb}, \ref{corollary: Bruhat Tits group are still intersection from GSp side, 2}, \ref{corollary: equality of cpt groups on V}}]\label{prop: main prop in intro}
         Fix any prime $p$. 
         \begin{enumerate}
             \item If $(G_1, X_1, K_{1, p}) \hookrightarrow (G_2, X_2, K_{2, p})$ is a very nice embedding (Definition \ref{def: nice embedding}), let $K_{1, p}$ and $K_{2, p}$ be the corresponding stablizer quasi-parahoric subgroups (resp. quasi-parahoric subgroups, neutral quasi-parahoric subgroups, Moy-Prasad subgroups with level $r > 0$). Assume $K_{1, p} = K_{2, p} \cap G_1(\Q_p)$, then equations \ref{eq: G_h} and \ref{eq: V} hold for all $\Phi_1$ and $\Phi_2 = \iota_*\Phi_1$.
             \item If $(G_1, X_1, K_{1, p}) \hookrightarrow (G_2, X_2, K_{2, p})$ is an adjusted Siegel embedding (Definition \ref{def: adjusted embedding}), let $K_{2, p} = \Ker(\GSP(\Z_p) \to \GSP(\Z_p/p^n\Z_p))$ for some $n\geq 0$ (we denote $K_{2, p} = \GSP(\Z_p)$ when $n = 0$) for the corresponding integral model $\GSP$ of $G_2 = \GSp$. Assume $K_{1, p} = K_{2, p} \cap G_1(\Q_p)$, then equations \ref{eq: G_h} and \ref{eq: V} hold for all $\Phi_1$ and $\Phi_2 = \iota_*\Phi_1$.
         \end{enumerate}
     \end{proposition}

     This Proposition implies that equations \ref{eq: G_h} and \ref{eq: V} hold for all $\Phi_1$ and $\Phi_2 = \iota_*\Phi_1$ for the Siegel embeddings in \cite{kisin2018integral}, \cite{kisin2024independence}, \cite{kisin2024integral}, \cite{pappas2024p} with appropriate chosen level groups.
     
\subsubsection{}
       In many mixed characteristic $(0, p)$ scenarios, we also have good minimal compactifications denoted as $\Shumm{K}$ and toroidal compactifications denoted as $\Shumc{K, \Sigma}$ for the integral model $\Shum{K}=\Shum{K}(G, X)$ of Hodge type. These compactifications extend the compactifications provided over the generic fiber. In the Siegel case, Chai and Faltings, as detailed in \cite{faltings2013degeneration}, constructed toroidal compactifications for smooth Siegel moduli problems by investigating the degenerations of abelian varieties. Kai-Wen Lan, in \cite{lan2013arithmetic}, extended this construction to compactify smooth PEL moduli problems. Later, Lan furthur extended it to cover cases in arbitary ramified characteristics and arbitary levels in \cite{lan2016compactifications} and \cite{lan2017integral}. 
       
       In the Hodge-type case, the absence of a modular interpretation posed challenges. Keerthi Madapusi Pera constructed compactifications of integral models of Hodge-type Shimura varieties in \cite{pera2019toroidal}. This approach relies on relative normalizations and extends a group-theoretical approach which was also introduced in \cite{lan2016compactifications} within the context of PEL-type cases. In the construction \cite{pera2019toroidal}, the author used the comparison of analytic construction and algebraic construction of toroidal compactifications of PEL-type Shimura varieties, and this comparison is the main topic of \cite{lan2012comparison}. Recently, the compactifications of integral models of abelian-type Shimura varieties have been worked out in \cite{wu2025arith}.

       Let us recall very briefly how these level groups are involved. By the main results of \cite{pera2019toroidal}, 
       \begin{enumerate}
           \item The toroidal compactification $\Shumc{K, \Sigma}$ of $\Shum{K}$ \'etale locally at boundary is a toroidal embedding $\Shum{K_{\Phi, P}}(P, X_{\Phi}) \hookrightarrow \Shum{K_{\Phi, P}}(P, X_{\Phi})(\sigma)$ for some $\sigma \in \Sigma(\Phi)$, this toroidal embedding is over $\Shum{K_{\Phi, \ovl{P}}}(\ovl{P}, X_{\Phi, \ovl{P}})$. $\Shum{K_{\Phi, P}}(P, X_{\Phi}) \to \Shum{K_{\Phi, \ovl{P}}}(\ovl{P}, X_{\Phi, \ovl{P}})$ is a torsor under a split torus $\mathbf{E}_{K}(\Phi)$ with cocharacter group $\mathbf{B}_{K}(\Phi)$ which spans $\mathbf{B}_K(\Phi)_{\R} = U_{\Phi}(\R)(-1)$, and collection of cone decompositions happens in $U_{\Phi}(\R)(-1)$. 
           \item There is a proper surjection $\Shum{K_{\Phi, \ovl{P}}}(\ovl{P}, X_{\Phi, \ovl{P}}) \to \Shum{K_{\Phi, h}}(G_h, X_{\Phi, h})$, which is a torosr under abelian schemes when the equation \ref{eq: V} holds at prime $p$. Here $\Shum{K_{\Phi, h}}(G_h, X_{\Phi, h})$ is an integral model of $\shu{K_{\Phi, h}}(G_h, X_{\Phi, h})$ defined via relative normalizations.
           \item There is a $\Lambda_{\Phi, K}$-action on $\Shum{K_{\Phi, h}}(G_h, X_{\Phi, h})$ which factors through a finite quotient, and $\Lambda_{\Phi, K}\backslash \Shum{K_{\Phi, h}}(G_h, X_{\Phi, h})$ form a stratification of the boundary of $\Shumm{K}(G, X)$ when the cusp label representatives $[\Phi] \in \Cusp_K(G, X)$ vary.
       \end{enumerate}

       We specify the construction in \cite{pera2019toroidal} in the case when $K_p$ is quasi-parahoric, where we have good integral models of $\shu{K}(G, X)$. In Kisin-Pappas integral models \cite{kisin2018integral} and Kisin-Pappas-Zhou integral models \cite{kisin2024independence}, \cite{kisin2024integral}, given certain conditions on $(p, G)$, by choosing a Siegel embedding carefully, the integral model $\Shum{K}(G, X)$ of $\shu{K}(G, X)$ defined via the relative normalization supports a local model diagram, the local properties of $\Shum{K}(G, X)$ can be understood using the deformation of $p$-divisible groups together with Hodge tensors. In Pappas-Rapoport integral models \cite{pappas2024p} (when $K_p$ is stablizer quasi-parahoric), \cite{daniels2024conjecture} (when $K_p$ is a general quasi-parahoric), without conditions on $(p, G)$, the integral model $\Shum{K}(G, X)$ of $\shu{K}(G, X)$ defined via the relative normalization can be characterized by the theory of $\GG$-Shtukas and integral local Shimura varieties, it has good functoriaty and is independent of the choice of Siegel embedding.
       
       The Proposition \ref{prop: main prop for BT theory, intro} and \ref{prop: main prop in intro} show that the integral model $\Shum{K_{\Phi, h}}(G_h, X_{\Phi, h})$ has correct quasi-parahoric level group at $p$, which says that it is defined via the relative normalization of a correct Siegel embedding. With some additional work, we show that
       \begin{theorem}[{Theorem \ref{thm: boundary of kp is again kp}}]\label{thm: main theorem intro}
           Minimal compactification of Pappas-Rapoport integral model (resp. Kisin-Pappas integral models, Kisin-Pappas-Zhou integral models) of Shimura variety of Hodge type with quasi-parahoric level structure is stratified by finite \'etale quotients of Pappas-Rapoport integral models (resp. Kisin-Pappas integral model, Kisin-Pappas-Zhou integral models) of smaller Shimura variety of Hodge type with quasi-parahoric level structure (resp. under some extra conditions).
       \end{theorem}
          Let us explain the \emph{extra conditions} here. In \cite{kisin2024integral}, the authors showed that the construction in \cite{kisin2018integral} needs an extra condition, which says the good embedding there should be replaced with a very good embedding. The existence of a very good embedding for a Hodge-type Shimura datum holds for most of the cases. In this Theorem, when we consider Kisin-Pappas integral models or Kisin-Pappas-Zhou integral models, we show that the morphisms at boundary are good, but need to assume these good embeddings are all very good embeddings, this can be easily verified in many cases using \cite[Theorem 4.4.3, Proposition 5.3.10]{kisin2024integral}. Moreover, in the case of Kisin-Pappas-Zhou integral models, we need to assume that for all boundary Shimura data $(G_{\Phi, h}, X_{\Phi, h})$ of $(G, X)$, $G_{\Phi, h}$ are $R$-smooth.
       
\subsubsection{}
       When $K_p$ is quasi-parahoric, the Pappas-Rapoport integral models $\Shum{K}(G, X)$ as well as $\lrbracket{\Shum{K_{\Phi, h}}(G_h, X_{\Phi, h})}_{[\Phi] \in \Cusp_K(G, X)}$ are independent of the choice of the Siegel embeddings. Moreover, $\Shum{K_{\Phi, \ovl{P}}}(\ovl{P}, X_{\Phi, \ovl{P}}) \to \Shum{K_{\Phi, h}}(G_h, X_{\Phi, h})$ is a torsor under an abelian scheme under the equation \ref{eq: V}, and $\Shum{K_{\Phi, P}}(P, X_{\Phi}) \to \Shum{K_{\Phi, \ovl{P}}}(\ovl{P}, X_{\Phi, \ovl{P}})$ is a torsor under a split torus over $\Spec \Z$, both parts are more or less determined by the structures on their generic fiber.
       \begin{proposition}{(Proposition \ref{prop: uniqueness of toroidal compactification})}\label{thm: uniqueness, intro}
           Minimal and toroidal compactifications of integral models of Hodge type Shimura varieties with quasi-parahoric level structres are independent of the choice of (adjusted) Siegel embeddings.
       \end{proposition}

       On toroidal compactifications, similar uniqueness properties have also been studied independently in \cite{wu2025arith}, where the author examined the functoriality of the conjugation action by $G^{\ad}(\Q)$ and established uniqueness along the way. His proof also replies on the calculation in Proposition \ref{prop: main prop in intro}. 
       
       With these uniqueness properties and their proofs, we further construct the change-of-parahoric morphisms on compactifications of integral models of Hodge-type Shimura varieties with quasi-parahoric level structures.

\subsubsection{}
       Given $K_1 \subset K_2$ two quasi-parahoric level groups, then the natural morphism $\shu{K_1} \to \shu{K_2}$ extends to a morphism between the Pappas-Rapoport integral models $\Shum{K_1} \to \Shum{K_2}$, due to \cite[Corollary 4.3.2]{pappas2024p} and \cite[Corollary 4.1.10]{daniels2024conjecture}. We show that
       \begin{theorem}{(Theorem \ref{Theorem: change of parahoric})}.
         The change-of-parahoric morphism $\pi: \Shum{K_1} \to \Shum{K_2}$ extends (uniquely) to the proper morphisms $\pi^{\min}: \Shumm{K_1} \to \Shumm{K_2}$, $\pi^{\tor}: \Shumc{K_1, \Sigma_1} \to \Shumc{K_2, \Sigma_2}$ for some collection of cone decompositions $\Sigma_1$ refine $\Sigma_2$. The properties of $\pi^{\min}$ and $\pi^{\tor}$ are characterized in Proposition \ref{Proposition: change of parahoric, prop}.
       \end{theorem}
      As an immediate corollary, $\Shum{K_1} \to \Shum{K_2}$ is a proper morphism.

\subsection*{Acknowledgement}

This article is based on the first half of the author's thesis, with many revisions and strengthened results. The author thanks his advisor, Kai-Wen Lan, for introducing the topic and providing many valuable insights. The author also thanks Peihang Wu for many helpful discussions.

\section{Level groups at boundary}\label{sec: Bruhat-Tits buildings}

 As the boundaries of the minimal compactifications of Hodge-type Shimura varieties consist of finite quotients of Hodge-type Shimura varieties, it is necessary to determine their level groups at a prime $p$ and related functoriality. Here we do not restrict ourselves to parahoric subgroups, we consider any deep level.

In the first subsection (\ref{subsec: levi and normal subgroups}), we present an axiomatic framework that enables us to simultaneously handle (quasi-)parahoric subgroups and Moy-Prasad subgroups. In the subsequent subsections (\ref{subsec: bruhat tits theory} and \ref{subsec: parahoric group schemes}), we verify that the level groups under consideration satisfy the outlined axiomatic framework. 

In the fourth subsection (\ref{subsec: compute levels at boundary}), we compute the level groups at the boundary. The main result of this section is Proposition \ref{prop: main prop for section BT theory}. 

In the fifth subsection (\ref{subsec: general facts}), we establish that certain properties of the group $G$ are inherited by its Levi subgroups and normal subgroups. 

In the sixth subsection (\ref{subsec: functoriality of levis}), we introduce functoriality results that play a crucial role in computing level groups associated with Siegel embeddings. 

Finally, in the last subsection (\ref{subsec: quasi-parahoric subgroups}), we provide additional remarks on quasi-parahoric subgroups.

This section is purely group theoretical. Many properties presented in this section can be extended to more general fields without much difficult, but for the sake of relevance to Shimura varieties, we will maintain our focus on the case when $F$ is a local field of mixed characteristics.

\subsection{Levi subgroups and normal subgroups}\label{subsec: levi and normal subgroups}

Let $G$ be a $F$-reductive group, $S$ be a maximal $F$-split torus of $G$, $Z = Z_G(S)$ be its centralizer, which forms a minimal Levi subgroup of $G$. In the case where $G$ is quasi-split over $F$, $Z$ is contained in an $F$-Borel subgroup, $Z$ itself is an $F$-torus. Let $Z_c$ be the maximal central $F$-torus of $G$, $Z_{c, \spl}$ be the maximal $F$-split $F$-torus contained in $Z_c$. 

Let $N = N_G(S)$, $\Phi = \Phi(G, S)$, $W=N_G(S)(F)/Z_G(S)(F)=N(F)/Z(F)$, $\Phi_{\red} = \lrbracket{\alpha \in \Phi|\ \frac{1}{2}\alpha \not\in \Phi}$. Fix a system of positive roots $\Phi^+ \subset \Phi$, and $\Delta$ the set of simple roots. Let $(\alpha, \beta) = \lrbracket{p\alpha + q\beta |\ p, q \in \Z_{>0}} \cap \Phi$. $\Psi \subset \Phi$ is called \emph{closed} if $(\alpha, \beta) \subset \Psi$ for all $\alpha, \beta \in \Psi$, and is called \emph{positively closed} if $\Psi \subset \Phi^+$. An \emph{ordering} refers to an arbitrary total ordering on either $\Psi_{\red}$ or $\Phi_{\red}$.

For each $\alpha \in \Phi$ (resp. positively closed set $\Psi \subset \Phi$), let $U_{\alpha}$ (resp. $U_{\Psi}$) be the unique closed connected unipotent $F$-subgroup of $G$ that is normalized by $Z$ and has Lie algebra $\lieg_{\alpha} + \lieg_{2\alpha}$ (set $\lieg_{2\alpha} = 0$ if $2\alpha \not\in \Phi$) (resp. $\Sigma_{\alpha \in \Psi} \lieg_{\alpha}$). There exists an isomorphism of $F$-varieties $\prod_{\alpha \in \Psi_{\red}} U_{\alpha} \to U_{\Psi}$ in any ordering. Let $U^{\pm}$ denote $U_{\Phi^{\pm}}$ respectively.

We are considering the following two conditions on a subgroup $K \subset G(F)$. For any algebraic subgroup $(\ast) \subset G$, $(\ast)_K$ denotes $K \cap (\ast)(F)$.
\begin{equation}\label{eq: U^-U^+N}
     K = U^-_{K}U^+_{K}N_{K} = U^+_{K}U^-_{K}N_{K},\quad \prod_{\alpha\in\Phi_{\red}^{\pm}}U_{\alpha, K} \rightiso U^{\pm}_{K}\ \textit{is bijective for any ordering}.
\end{equation}
\begin{equation}\label{eq: U^+U^-U^-Z}
     K = U_K^+U_K^-U_K^+Z_K.
\end{equation}

\subsubsection{Levi subgroups}

Let $J \subset \Delta$ be a subset, $\Phi_J = \Z J \cap \Phi$, let $Q = Q_J = \lrangle{Z, U_{\alpha}|\ \alpha\in \Phi_J \cup \Phi^+}$, $S_J = (\cap_{\alpha\in J}\Ker \alpha)^{\circ} \subset S$, $L=L_J = Z_G(S_J) = \lrangle{Z, U_{\alpha}|\ \alpha \in \Phi_J}$, $R_uQ = \lrangle{U_{\alpha}|\ \alpha \in \Phi^+\setminus \Phi^+_J}$, $U^{\pm}_J := U_{\Phi_J^{\pm}}$. Such pair $(Q, L)$ is called \emph{at standard position} with respect to $(S, J \subset \Delta)$. In this section, we always assume that $(Q, L)$ is standard with respect to $(S, J \subset \Delta)$ for some $J \subset \Delta$.

\begin{definition}
    We say that $K_L \subset L(F)$ satisfies \ref{eq: U^-U^+N} and \ref{eq: U^+U^-U^-Z} when these equations or bijections hold with $U^{\pm}$, $N$, $Z$ replaced with $U^{\pm}_J = L \cap U^{\pm}$, $N_L:=N_L(S) = N \cap L$, $Z_L(S) = Z_G(S) = Z$ respectively.
\end{definition}

\begin{proposition}\label{proposition: henniart factors}
Assume $K \subset G(F)$ satisfies the condition \ref{eq: U^-U^+N}, then
    \begin{enumerate}
    \item $Q(F) \cap K = ((R_uQ)(F) \cap K)  \rtimes (L(F) \cap K)$, in particular, $\pi(Q(F) \cap K) = L(F) \cap K$, where $\pi: Q \to L$ is the projection.
    \item $L(F) \cap K = (U^+(F) \cap L(F) \cap K)(U^-(F) \cap L(F) \cap K)(N(F) \cap L(F) \cap K)$.
    \item $L(F) \cap K$ satisfies \ref{eq: U^-U^+N}.
\end{enumerate}
\end{proposition}
The formulation of the Proposition and its proof are quite similar to \cite[\S 6.6, 6.7, 6.8]{henniart2015satake}. Since the proof is short, we write it out and provide more details here.
\begin{lemma}\label{lemma: a lemma for Bruhat decomposition}
    $R_uQ^{-}(F)N(F) \cap Q(F) = L(F) \cap N(F)$.
\end{lemma}
\begin{proof}
It suffices to prove $R_uQ^{-}(F)N(F) \cap Q(F) \subset L(F) \cap N(F)$. Recall the Bruhat decomposition, see \cite[Lemma 21.14, Theorem 21.15]{borel2012linear}, where we conjugate $U^+$ to $U^-$ using a representative $\wdt{w}_0 \in N(F)$ of the longest element in the Weyl group:
\begin{equation}\label{eq: bruhat decomposition}
    G(F) = \wdt{w}_0G(F) = U^-(F)N(F)U^+(F)
\end{equation}

    Let $ u^-n \in R_uQ^{-}(F)N(F) \cap Q(F)$, $u^{-}n = lu$ for some $u^- \in R_uQ^{-}(F) \subset U^-(F)$, $n \in N(F)$, $l \in L(F)$, $u \in R_uQ(F) \subset U^+(F)$. The decomposition \ref{eq: bruhat decomposition} implies $l = u_l^-n_lu_l^+$, $u_l^- \in U_J^-(F) \subset U^-(F)$, $n_l \in N_L(S)(F) \subset N(F)$, $u_l^+ \in U_J^+(F) \subset U^+(F)$, then $n = (u^-)^{-1}u_l^-n_lu_l^+u \in U^-(F)N(F)U^+(F)$. \cite[Lemma 21.14]{borel2012linear} also implies that given $n_1, n_2 \in N(F)$, 
    \[U^-(F)n_1U^+(F) = U^-(F)n_2U^+(F)\ \textit{if and only if}\ n_1 = n_2,\] 
    thus $n = n_l \in N_L(S)(F)$, $u^- = (lul^{-1})(ln^{-1})$. Since $L(F)$ normalizes $R_uQ^+(F)$, then $u^{-1} = R_uQ^{-}(F) \cap Q(F)$, which forces $u^- = 1$, therefore $R_uQ^{-}(F)N(F) \cap Q(F) \subset N(F) \cap L(F)$.
\end{proof}
\begin{proof}{Proof of Proposition \ref{proposition: henniart factors}}
 \begin{enumerate}
     \item Let $q \in Q(F) \cap K$, then $q = u^+u^-n$, $u^{\pm}\in U^{\pm}_K$, $n \in N_K$. Since $u^-$ is a product of $u_{\alpha} \in U_{\alpha, K}$ for $\alpha \in \Phi^{-}_{\red}$ in any ordering, we decompose $u^- = u^-_1u^-_2$, $u^-_1 \in L(F) \cap U_K^-$, $u^-_2 \in R_uQ^-(F) \cap U_K^-$. Then $u^-_2n = (u^-_1)^{-1}(u^+)^{-1}q \in  R_uQ^{-}(F)N(F)\cap Q(F)$, $u^-_2 = 1$ and $n \in L(F)$ due to the proof of Lemma \ref{lemma: a lemma for Bruhat decomposition}, $q = u_1^+u_2^+u^-_1n$, $u^+_1 \in R_uQ(F) \cap U_K^+$, $u_2^+ \in L(F) \cap U_K^+$ (thus $u_2^+u^-_1n \in L(F)\cap K$, we get part $(1)$.
     \item In the proof of part $(1)$, let $q \in L(F)\cap K$, then $u^+ = q n^{-1}(u_1^-)^{-1} \in L(F) \cap U_K^+$, $q = u^+ u^-_1 n$ is the factorization we need.
     \item Let $K_L = L(F) \cap K$, part $(2)$ means $K_L = U^+_{J, K_L}U^-_{J, K_L}N_{L, K_L}$. Similarly, $K_L = U^-_{J, K_L}U^+_{J, K_L}N_{L, K_L}$. Since $K \cap U^{\pm} = \prod_{\alpha \in \Phi^{\pm}_{\red}} U_{\alpha, K}$ in any ordering and $U^{\pm} = U^{\pm}_J \times U^{\pm}_{\Delta \setminus J}$, then $K \cap U^{\pm}_J = \prod_{\alpha \in \Phi^{\pm}_{J, \red}} U_{\alpha, K}$ in any ordering.
 \end{enumerate} 
\end{proof}

\begin{lemma}\label{lemma: open cell decomp into levi}
    Assume $K$ satisfies both conditions \ref{eq: U^-U^+N} and \ref{eq: U^+U^-U^-Z} , then $K \cap L(F)$ also satisfies both conditions \ref{eq: U^-U^+N} and \ref{eq: U^+U^-U^-Z} .
\end{lemma}
\begin{proof}
    The condition \ref{eq: U^-U^+N} follows from Proposition \ref{proposition: henniart factors}. To show condition \ref{eq: U^+U^-U^-Z} , it suffices to show 
        \begin{equation}\label{eq: open cell decomp into levi}
            L(F) \cap K = (U_K^+ \cap L(F))(U_K^- \cap L(F))(U_K^+ \cap L(F))Z_K.
        \end{equation}

    Let $p \in Q(F) \cap K$, $ p = u^+u^-v^+z$, $u^+, v^+ \in U_K^+$, $u^- \in U_K^-$, $z\in Z_K$. Decompose $u^- = u_1^- u_2^-$, $u_1^- \in L(F) \cap U_K^-$, $u_2^- \in R_uQ^-(F) \cap U_K^-$, then $u_2^- = (u_1^{-})^{-1} (u^{+})^{-1}pz^{-1}(v^{+})^{-1} \in Q(F) \cap R_uQ^-(F)$, $u_2^- = 1$. Let $u^+ = u^+_1 u^+_2$, $v^+ = v^+_1 v^+_2$, $u^+_1, v^+_1 \in L(F) \cap U_K^+$, $u^+_2, v^+_2 \in R_uQ(F) \cap U_K^+$, then $\pi(p) =  u^+_1u_1^-v^+_1z$. Proposition \ref{proposition: henniart factors} implies $L(F) \cap K = \pi(Q(F) \cap K)$, we are done.
\end{proof}
\begin{corollary}\label{cor: when K cap L is K_L}
    Assume $K \subset G(F)$ and $K_L \subset L(F)$ be two groups that satisfy both conditions \ref{eq: U^-U^+N} and \ref{eq: U^+U^-U^-Z} with respect to $G$ and to $L$ respectively, then $L(F) \cap K = K_L$ if and only if
    \begin{equation}
        \forall\alpha \in \Phi_{J, \red},\ K_L \cap U_{\alpha}(F) = K \cap  U_{\alpha}(F), \quad Z_{K_L} = Z_K.
    \end{equation}
\end{corollary}
\begin{proof}
    The \emph{only if} part is trivial. We show the \emph{if} part as follows. Since $U^{\pm}_J = L \cap U^{\pm}$, and the comparison between the decomposition
       \begin{equation}
         K_L= (K_L \cap U^+_J(F)) (K_L \cap U^-_J(F)) (K_L \cap U^+_J(F)) Z_{K_L}.
     \end{equation}
     and decomposition \ref{eq: open cell decomp into levi}, then $L(F) \cap K = K_L$ follows from $Z_{K_L} = Z_K$ and $K \cap U^{\pm}_J = K_L \cap U^{\pm}_J$. Due to Proposition \ref{proposition: henniart factors}, $K \cap U^{\pm}_J = \prod_{\alpha \in \Phi_{J, \red}^{\pm}} U_{\alpha, K}$, compare it with the decomposition $K_L \cap U^{\pm}_J = \prod_{\alpha \in \Phi_{J, \red}^{\pm}} U_{\alpha, K_L}$, $K \cap U^{\pm}_J = K_L \cap U^{\pm}_J$ follows from $K_L \cap U_{\alpha}(F) = K \cap  U_{\alpha}(F)$ for all $\alpha \in \Phi_{J, \red}$.
\end{proof}

\subsubsection{Normal subgroups}

We start with a general fact, it follows easily from the structure theory of reductive groups:
 \begin{proposition}\label{proposition: some props about normal reductive groups and its quotients}
          Let $1 \to G_h \to L \to G_l \to 1$ be an exact sequence of reductive groups over $F$, then:
          \begin{enumerate}
          \item There exists a normal reductive $F$-subgroup $G_l'$ of $L$ such that $L$ is an almost product of $G_h$ and $G_l'$.
          \item Let $L^{\ad}$ (resp. $L^{\der}$) be the adjoint (resp. derived) group of $L$, $\wdt{L}\to L^{\der}$ be the simply connected covering of $L^{\der}$, then one has induced exact sequences $1 \to G_h^{\ad} \to L^{\ad} \to G_l^{\ad} \to 1$ and $1 \to \wdt{G}_h \to \wdt{L} \to \wdt{G}_l \to 1$. These exact sequences split in a canonical way, $L^{\ad} \cong G_h^{\ad} \times G_l^{\ad}$, $\wdt{L} = \wdt{G}_h \times \wdt{G}_l$, and $G_l^{\prime\ad} \cong G_l^{\ad}$, $\wdt{G}_l' \cong \wdt{G}_l$.
   
          \end{enumerate}
      \end{proposition}

Let $G_h$ be a normal reductive $F$-subgroup of $L$, $F_1$ be a field extension of $F$. The following lemma is useful and standard, and we provide some proofs from the literature while also supplying the remaining. We have following standard facts in algebraic groups.
\begin{lemma}\label{lemma: good sub and quotient for centralizers and normalizers}
    Let $S \subset L$ be a $F$-torus , $S_h:=(S \cap G_h)^{\circ}$, $S_{l'}:=(S \cap G_l')^{\circ}$, $S_l = \pi(S)$, then
    \begin{enumerate}
        \item $S_h \times S_{l'} \to S$ is an almost product, $S_{l'} \to S_l$ is an isogeny. When $S$ is a maximal torus (resp. maximal $F$-split torus, maximal $F_1$-split torus defined over $F$, maximal central torus, maximal central $F$-split torus), then so are $S_h$ and $S_l$.
        \item Let $Z$ (resp. $Z_h, Z_{l'}, Z_l$) be the centralizer of $S$ in $L$ (resp. $S_h$ in $G_h$, $S_{l'}$ in $G_l'$, $S_l$ in $G_l$), then $Z_h = Z \cap G_h$, $Z_{l'} = Z \cap G_l'$, $Z_l = \pi(Z)$, $Z_h \times Z_{l'} \to Z$ is an almost product. The same statements also hold when the centralizers $Z$, $Z_h, Z_{l'}, Z_l$ are replaced with normalizers $N$, $N_h$, $N_{l'}$, $N_l$.
    \end{enumerate}.
\end{lemma}

Let $S$ be a maximal $F$-split torus of $L$, due to Lemma \ref{lemma: good sub and quotient for centralizers and normalizers}, $Z_h \times Z_{l'} \to Z$ is an almost product, and $Z_{l'} \to Z_l$ is an isogeny. Note that the isogeny $G_l' \to G_l$ identifies $\Phi_{l'}$ and $\Phi_l$ via $X^*(S_l)_{\Q} \cong X^*(S_{l'})_{\Q}$ canonically, $\Phi:=\Phi_L \cong \Phi_h \sqcup \Phi_{l'}$, $J = \Delta_h \sqcup \Delta_{l'}$, and these decompositions are orthogonal. $G_h = \lrangle{Z_{G_h}(S_h), U_{\alpha}|\ \alpha \in\Phi_h}$, $G_l' = \lrangle{Z_{G_l'}(S_{l'}), U_{\alpha}|\ \alpha \in\Phi_{l'}}$, $G_{l} = \lrangle{Z_{G_l}(S_l), U_{\alpha}|\ \alpha \in\Phi_l}$. The isogeny $G_l' \to G_l$ maps $U_{\alpha}$ ($\alpha\in \Phi_{l'}$) isomorphically onto $U_{\alpha}$ ($\alpha\in\Phi_l$) when we identify $\alpha\in \Phi_{l'}$ with $\alpha\in\Phi_l$. Similarly, replace $S$ with $Z_c$ (resp. $Z_{c, \spl}$, $T$, $Z':=Z_T(L)$), the maximal central $F$-torus (resp. the maximal central split $F$-torus, a maximal split $F_1$-torus defined over $F$, the centralizer of $T$) of $L$, with Lemma \ref{lemma: good sub and quotient for centralizers and normalizers}, we have an almost product $Z_{h, c} \times Z_{l', c} \to Z_{c}$ (resp. $Z_{h, c, \spl} \times Z_{l', c, \spl} \to Z_{c, \spl}$, $T_h \times T_l \to T$, $Z'_h \times Z'_l \to Z'$) and an isogeny $Z_{l', c} \to Z_{l, c}$ (resp. $Z_{l', c, \spl} \to Z_{l, c, \spl}$, $T_{l'} \to T_l$, $Z'_{l'} \to Z'_l$).

\begin{definition-proposition}
    We say that $K_h \subset G_h(F)$ satisfies conditions \ref{eq: U^-U^+N} and \ref{eq: U^+U^-U^-Z} when these equations or bijections hold with $U^{\pm}$, $N$, $Z$ replaced with $U^{\pm}_h := U_{\Phi_h^{\pm}} = U^{\pm}_J \cap G_h$, $N_h:= N_{G_h}(S_h) = N_L(S) \cap G_h$, $Z_h = Z_{G_h}(S_h) = Z \cap G_h$ respectively.
\end{definition-proposition}

\begin{lemma}\label{lem: K_L cap G_h is good as K_L}
    When $K_L \subset L(F)$ satisfies condition \ref{eq: U^-U^+N}, then $K_L \cap G_h(F)$ satisfies condition \ref{eq: U^-U^+N}. If moreover $K_L \subset L(F)$ satisfies condition \ref{eq: U^+U^-U^-Z} , then  $K_L \cap G_h(F)$ satisfies condition \ref{eq: U^+U^-U^-Z} . 
\end{lemma}
\begin{proof}
    Since $U_J^{\pm} = U_h^{\pm} \times U_{l'}^{\pm} = U_{l'}^{\pm} \times U_h^{\pm}$ and $U_h^{\pm}$ commute with $U_{l'}^{\pm}$, the decompositions $\prod_{\alpha \in \Phi_{J, \red}^{\pm}} U_{\alpha, K_L} = U_{J, K_L}^{\pm}$ factor as $\prod_{\alpha \in \Phi_{h, \red}^{\pm}} U_{\alpha, K_L} = U_{h, K_L}^{\pm}$ and $\prod_{\alpha \in \Phi_{l', \red}^{\pm}} U_{\alpha, K_L} = U_{l', K_L}^{\pm}$.
    
    Let $g \in K_L \cap G_h(F)$, $g = u^+u^-n$, $u^+ \in U_{J, K_L}^+$, $u^- \in U_{J, K_L}^-$, $n \in N_L(F) \cap K_L$. Decompose $u^+ = u_h^+u_l^+$, $u^-=u_h^-u_l^-$, $u_h^{\pm} \in U_{h, K_L}^{\pm}$, $u_l^{\pm} \in U_{l', K_L}^{\pm}$. The projection $\pi_l: L \to G_l$ kills the $G_h$-part: $\pi_l(g) = u_l^+u_l^-\pi_l(n)$. Since $\pi_l(n) \in N_l(F)$ and $N_l(F) \cap U_l^+(F)U_l^-(F) = \lrbracket{1}$, then $u_l^+ = u_l^- = 1$, $n \in N_h(F) \cap K_L$, thus $K_L \cap G_h(F) = (U_h^+(F) \cap K_L)(U_h^-(F) \cap K_L)(N_h(F) \cap K_L)$, $K_L \cap G_h(F)$ satisfies condition \ref{eq: U^-U^+N} when $K_L \subset L(F)$ does.

    Similar arguments also show $K_L \cap G_h(F)$ satisfies condition \ref{eq: U^+U^-U^-Z}  when $K_L \subset L(F)$ does: Let $g \in K_L \cap G_h(F)$, $g = u^+u^-v^+z$, $u^+, v^+ \in U_{J, K_L}^+$, $u^- \in U_{J, K_L}^-$, $z \in Z_{K_L}$. Decompose $u^+ = u_h^+ u_l^+$, $u^- = u_h^- u_l^-$, $v^+ = v_h^+ v_l^+$, $u_h^+, v_h^+ \in U_{h, K_L}^+$, $u_l^+, v_l^+ \in U_{l, K_L}^+$, $u_h^- \in U_{h, K_L}^-$, we have $ u_l^+u_l^-v_l^+ \pi_l(z) = 1$. Since $Z_l$ and $Z'_l$ normalizes $U_l^{\pm}$, and $N_l(F) \cap U_l^+(F)U_l^-(F) = \lrbracket{1}$, we have $u_l^+ v_l^+ = u_l^- = 1$, $z \in K_L \cap Z_h(F)$.
\end{proof}

Using similar arguments as in Corollary \ref{cor: when K cap L is K_L}, but with the input replaced by Lemma \ref{lem: K_L cap G_h is good as K_L}, we can conclude that:
\begin{corollary}\label{cor: when K_L cap G_h is K_h}
    Assume $K_L \subset L(F)$ and $K_h \subset G_h(F)$ satisfy both conditions \ref{eq: U^-U^+N} and \ref{eq: U^+U^-U^-Z} with respect to $L$ and to $G_h$ respectively, then $G_h(F) \cap K_L = K_h$ if and only if for all $\alpha \in \Phi_{h, \red}$, $K_L \cap U_{\alpha}(F) = K_h \cap  U_{\alpha}(F)$ and
    \begin{equation}
        Z_h(F) \cap K_L = Z_h(F) \cap K_h.
    \end{equation}
\end{corollary}

\subsection{Yu's work on open big cells}\label{subsec: bruhat tits theory}

Let $\bF$ be the completion of the maximal unramified extension of $F$. One of the most important proposition in the construction of Bruhat-Tits building over $F$ descent from $\bF$ is the following:

\begin{definition-proposition}{{\cite[\S 5.1.10]{bruhat1984groupes}}}\label{def: T}
     There exists a maximal $\bF$-split torus $T$ containing $S$ and defined over $F$. $Z' = Z_G(T)$ is a maximal torus of $G$ defined over $F$ (since $G_{\bF}$ is quasi-split), and $Z' \subset Z$.
\end{definition-proposition}

In this and the following sections, we will recall the constructions of the following list of groups:
\begin{definition}\label{definition: type ast}
We say $K$ is of type $(\ast)$ when $K$ belongs to one of the following types:
    \begin{enumerate}
    \item A neutral quasi-parahoric group
    \item A stablizer quasi-parahoric group.
    \item A parahoric group.
    \item An absolutely special/hyperspecial/superspecial/special parahoric group.
    \item A Moy-Prasad group of level $r \in \R_{>0}$.
\end{enumerate}
   We recall the definitions in \ref{def: parahoric subgroup}, \ref{def: type ast groups, in an appartment}, \ref{def: type ast groups, in a building} and \ref{def: various special points}.
\end{definition}

The constructions involve Bruhat-Tits theory. We provide a brief overview of some basic definitions and propositions. Here we only recall a very small portion of the theory. For a comprehensive understanding and detailed proofs, we recommend referring to \cite{bruhat1972groupes}, \cite{bruhat1984groupes}, and \cite{tits1979reductive}. Additional insights and approaches can also be found in \cite{landvogt2000some}, \cite{landvogt2006compactification}, and \cite{yu2015smooth}. The main result of this subsection is Proposition \ref{proposition: MP-group satisfies two conditions}.

\subsubsection{}

Let $G$ be a reductive group over $F$, Kottwitz constructed the following morphisms in \cite{kottwitz1997isocrystals} (for $p$-adic fields $F$), where $I = \Gal(\ovl{F}|\bF)$, $\Sigma = \Gal(\ovl{F}|F)$, $\Sigma_0 = \Gal(\bF|F)$:
\begin{equation}\label{eq: Kottwitz, over bF}
  \nu_G: G(\bF) \stackrel{\kappa_G}{\twoheadrightarrow} X^*(Z(\hat{G})^I) \stackrel{q_G}{\twoheadrightarrow} \Hom(X_*(Z(\hat{G}))^I, \Z),
\end{equation}
here $\nu_G$ is the composition. By taking $\Sigma_0$-invariants, we have
\begin{equation}\label{eq: Kottwitz, over F}
  \nu_G: G(F) \stackrel{\kappa_G}{\twoheadrightarrow} X^*(Z(\hat{G})^I)^{\Sigma_0} \stackrel{q_G}{\longrightarrow} \Hom(X_*(Z(\hat{G}))^{\Sigma}, \Z),
\end{equation}
These morphisms can be similarly defined over any henselian discrete valued field. 

Let $G(\bF)^1 = \Ker \nu_G$, $G(F)^1 = \Ker \nu_G|_{G(F)}$ (use \ref{eq: Kottwitz, over F}), $G(\bF)_0 = \Ker \kappa_G$, $G(F)_0 = \Ker \kappa_G|_{G(F)}$ (use \ref{eq: Kottwitz, over F}). Since taking $\Sigma_0$-invariants is a left exact functor, $G(\bF)^1 \cap G(F) = G(F)^1$, $G(\bF)_0 \cap G(F) = G(F)_0$\footnote{ Here, $G(\bF)_0$ and $G(\bF)^1$ are distinct from the notions defined in \cite{kottwitz1997isocrystals}, where $G(\bF)_b = \Ker \nu_G$ and $G(\bF)_1 = \Ker \kappa_G$. We have chosen to use different notation to avoid confusion later when we discuss the Moy-Prasad filtrations, as the subscript $r \geq 0$ in $G(\bF)_r$ will refer to the levels.}. Let us denote $\Lambda_G = X^*(Z(\hat{G})^I) = \pi_1(G)_I$.

When we replace $G$ with $Z = Z_G(S)$, and when $F$ is locally compact (meaning that $F$ is complete and $k_F$ is finite), \cite[Proposition 1.2]{landvogt2006compactification} implies that $Z(F)^1$ is the maximal open compact subgroup of $Z(F)$, and $Z(F)/Z(F)^1$ is a free abelian group.

\subsubsection{}

We will continue to use the notations from subsection \ref{subsec: levi and normal subgroups} and follow the presentations in \cite{yu2015smooth}. The work there was primarily written out over $\bF$. We can state similar results (for example, the open cell embedding) over $F$ using étale descent, as described in \cite[\S 5]{bruhat1984groupes}.

\begin{definition}
    Given a scheme $X$ over $\Spec F$, we say $\mathcal{X}$ is a model of $X$ when $\mathcal{X}$ is a scheme over $\Spec \OO_F$ with the generic fiber $\mathcal{X} \times_{\Spec \OO_F} F = X$.
\end{definition}
Note that when $\XX$ is a smooth scheme over $\OO_F$, then $\XX(\OO_F) \to \XX(k_F)$ is surjective.

Let $\wdh{\Phi}:= \Phi \cup \lrbracket{0}$, $\wdh{\Phi}_{\red} = \Phi_{\red} \cup \lrbracket{0}$, $U_0:= Z':= Z_G(T)$, $\wdt{\R}$ be the totally ordered monoid $\R \cup \lrbracket{r+ | r \in \R} \cup \lrbracket{\infty}$ defined in \cite[6.4.1]{bruhat1972groupes}. 

Let $A_{\red}(G, S)$ be the reduced apartment, and fix a point $x \in A_{\red}(G, S)$. We identify $x$ with its image in $A_{\red}(G_{\bF}, T_{\bF})$ under the identification $A_{\red}(G, S) = A_{\red}(G_{\bF}, T_{\bF})^{\Sigma_0}$. To distinguish objects defined over $\bF$ from those defined over $F$, we add a superscript $\prime$ to those data defined related with $\bF$, such as $Z':= Z_G(T)$.

The Moy-Prasad groups of $G(\bF)$ depend on certain filtrations on unipotent groups $\lrbracket{U_{\alpha}(\bF)}$ and on the maximal torus $U_0(\bF) = Z'(\bF)$. The former was established in \cite{bruhat1972groupes} and \cite{bruhat1984groupes}, while the latter was studied in \cite{moy1996jacquet}. In addition, \cite{yu2015smooth} provides a comprehensive description of the filtrations on tori:

\subsubsection{Tori} Let $T$ be a torus over $F$, and let $\TT$ be a smooth model contained in the lft-N\'eron model of $T$ over $\Spec \OO_F$. We choose $\TT$ to be the finite-type N\'eron model such that $\TT(\OO_{\bF}) = T(\bF)^1$. Additionally, let $\TT^{\circ}$ be the connected component of $\TT$ (which is also the connected component of the lft-N\'eron model and also the finite-type N\'eron model). It is a well-known fact that $\TT^{\circ}(\OO_{\bF}) = T(\bF)_0$, for example, see \cite{haines2008parahoric}.
\begin{definition}{{\cite[\S 4.2]{yu2015smooth}}}\label{def: Yu MP}
    The Moy-Prasad filtration on $T(F)$ is defined as follows: $T(F)_0:= \Ker \kappa_T$. When $r > 0$,
    \begin{enumerate}
        \item When $T = R = \prod_{i=1}^m \Res_{F_i|F} \Gm$ (we call such tori induced tori), let 
        \begin{equation}
            R(F)^{\MPF}_r:= \lrbracket{ (x_i) \in \prod_{i = 1}^m F_i^{\times} |\ \omega_F(x_i - 1) \geq r\ \textit{for each}\ i} \subset R(F) = \prod_{i = 1}^m F_i^{\times}
        \end{equation}
        \item In general, let $R$ be an induced torus containing $T$, let
        \begin{equation}\label{eq: Yu MP}
            T(F)^{\MPF}_r := R(F)^{\MPF}_r \cap T(F)_0.
        \end{equation}
    \end{enumerate}
\end{definition}
\begin{remark}\label{remark: MP independent of embedding}
    Since $\omega_F$ is a discrete valuation, this definition is equivalent to the one defined by Moy and Prasad in \cite[\S 3.2]{moy1996jacquet}:
    \begin{equation}\label{eq: MP original}
        T(F)^{\MPF}_n = \lrbracket{ x \in T(F)_0 |\ \omega_F(\chi(x) -1) \geq n,\ \forall \chi \in X^*(T) }
    \end{equation}
    Note that they use $\omega_L$ ($L$ is the splitting field of $T$ and $\omega_L(L^{\times}) = \Z$) instead of $\omega_F$. The latter definition \ref{eq: MP original} implies that \ref{eq: Yu MP} is independent of the choice of $R$. 
\end{remark} 

The author also considered more general filtrations on $T$ in \cite[\S 4]{yu2015smooth}, for instance, an \emph{admissible}, \emph{schematic} filtration. In fact, the author introduced a functor denoted by $\dagger$, which maps tori over henselian local fields (not limited to $F$) to sequences of groups. 
\begin{equation}\label{eq: adm filtrations}
    \dagger: T/F \mapsto \lrbracket{T(F)_r^{\dagger}}_{r \geq 0}
\end{equation}

The filtration is considered \emph{admissable} if it satisfies certain conditions in comparison to the Moy-Prasad filtration, as detailed in \cite[\S 4.3]{yu2015smooth}. The condition \cite[\S 4.4 (S1)]{yu2015smooth} for $\dagger$ to be schematic is that for each $r \geq 0$, there exists a unique smooth group scheme model $\TT^{\dagger}_r$ of $T$ over $\OO_F$ such that $\TT^{\dagger}_r(\OO_{\bF}) = T(\bF)^{\dagger}_r$. It is worth noting that the Moy-Prasad filtration $\MPF: T/F \mapsto \lrbracket{T(F)^{\MPF}_r}$ is both admissible and schematic, as stated in \cite[Proposition 4.1]{yu2015smooth}.

\subsubsection{Unipotent groups}\label{subsubsec: unipotent groups}

Let us recall the filtrations on unipotent groups. Given any $x \in A_{\red}(G_{\bF}, T_{\bF})$, Bruhat and Tits defined a filtration $\lrbracket{U_{\alpha'}(\bF)_{x, r}}_{\alpha' \in \Phi', r \in \R}$ on $\lrbracket{U_{\alpha'}(\bF)}_{\alpha'\in \Phi'}$ characterized in \cite[\S 6.2.1]{bruhat1972groupes}. In fact, $x \in A_{\red}(G_{\bF}, T_{\bF})$ determines and is determined by a valuation in \cite[\S 6.2.1]{bruhat1972groupes}. In the following, we assume $x \in A_{\red}(G, S) = A_{\red}(G_{\bF}, T_{\bF})^{\Sigma_0}$ as usual.

\begin{proposition}{{\cite[\S 4.3, 5.2]{bruhat1984groupes}}}\label{prop: etale descent for unipotent groups}\leavevmode
    \begin{enumerate}
        \item For any $\alpha' \in \Phi'$, $r \in \R$, there exists a unique smooth group scheme model $\UU_{\alpha', x, r}$ of $U_{\alpha'}$ such that $\UU_{\alpha', x, r}(\OO_{\bF}) = U_{\alpha'}(\bF)_{x, r}$. Moreover, $\UU_{\alpha', x, r}$ is connected and its special fiber $\UU_{\alpha', x, r} \otimes k_F$ is unipotent.
        \item For any $\alpha' \in \Phi'$ with $2\alpha' \in \Phi'$, $r, s \in \R$ such that $2r \geq s$, there exists a unique smooth group scheme model $\UU_{\alpha', x, r, s}$ of $U_{\alpha'}$ such that 
        \[\UU_{\alpha', x, r, s}(\OO_{\bF}) = U_{\alpha'}(\bF)_{x, r}U_{2\alpha'}(\bF)_{x, s}.\] Moreover, $\UU_{\alpha', x, r, s}$ is connected and its special fiber $\UU_{\alpha', x, r, s} \otimes k_{\bF}$ is unipotent.
        \item For any $\alpha \in \Phi$, $r \in \R$, let
\begin{equation}
  \wdt{U}_{\alpha}(\bF)_{x, r} := \prod_{\alpha' \in \Phi', \alpha'|_S = \alpha} U_{\alpha'}(\bF)_{x, r} \prod_{\alpha' \in \Phi'_{\red}, \alpha'|_S = 2\alpha}U_{\alpha'}(\bF)_{x, 2r} \subset U_{\alpha}(\bF).
\end{equation}
and $U_{\alpha}(F)_{x, r} = U_{\alpha}(F) \cap \wdt{U}_{\alpha}(\bF)_{x, r}$. There exists a unique smooth group scheme model $\UU_{\alpha, x, r}$ of $U_{\alpha}$ over $\OO_F$ such that $\UU_{\alpha, x, r}(\OO_{\bF}) = \wdt{U}_{\alpha}(\bF)_{x, r}$, $\UU_{\alpha, x, r}(\OO_{F}) = U_{\alpha}(F)_{x, r}$, and the multiplicative morphism is an isomorphism over $\OO_{\bF}$:
     \begin{equation}
         \prod_{\alpha' \in \Phi', \alpha'|_S = \alpha} \UU_{\alpha', x, r} \times \prod_{\alpha' \in \Phi'_{\red}, \alpha'|_S = 2\alpha} \UU_{\alpha', x, 2r} \to \UU_{\alpha, x, r}.
     \end{equation}
\item For any $\alpha \in \Phi$ with $2\alpha \in \Phi$, $r, s \in \R$ such that $2r \geq s$, there exists a unique smooth group scheme model $\UU_{\alpha, x, r, s}$ of $U_{\alpha}$ over $\OO_F$ such that $\UU_{\alpha, x, r, s}(\OO_{\bF}) = \wdt{U}_{\alpha}(\bF)_{x, r}\wdt{U}_{2\alpha}(\bF)_{x, s}$, $\UU_{\alpha, x, r, s}(\OO_{F}) = U_{\alpha}(F)_{x, r}U_{2\alpha}(\bF)_{x, s}$, and the multiplicative morphism is an isomorphism over $\OO_{\bF}$:
     \begin{equation}
         \prod_{\alpha' \in \Phi', \alpha'|_S = \alpha} \UU_{\alpha', x, r, s} \times \prod_{\alpha' \in \Phi'_{\red}, \alpha'|_S = 2\alpha} \UU_{\alpha', x, s} \to \UU_{\alpha, x, r, s}.
     \end{equation}
    \end{enumerate}
\end{proposition}

\begin{definition}\label{def: concave function}
    A function $f: \wdh{\Phi} \to \wdt{\R}$ is a \emph{concave} function if for any non-empty finite set $\lrbracket{\alpha_i}_{i \in I}$ of $\wdh{\Phi}$, $\Sigma_{i \in I} f(\alpha_i) \geq f(\Sigma_{i \in I} \alpha_i)$.
\end{definition}
\begin{remark}\label{rmk: concave and quasi-concave}
    Since $0 \in \wdh{\Phi}$, $f(0) \geq 0$ when $f$ is concave. $f$ is concave in the sense of \ref{def: concave function} implies that $f$ is concave in the sense of \cite[\S 4.5.3]{bruhat1984groupes}, thus $f$ is quasi-concave in the sense of \cite[\S 4.5.3]{bruhat1984groupes}.
\end{remark}
\begin{remark}\label{remark: 2f(alpha) notation 1}
    For convenience, if $\alpha \in \Phi$, but $2\alpha \not\in \Phi$, we define $f(2\alpha) = 2f(\alpha)$ and set $U_{2\alpha}(F)_{x, r} = \lrbracket{1}$ for any $r \in \R$.
\end{remark}

\subsubsection{Nice decompositions}
    Let $f$ be a concave function, let $U_{x, f} \subset G(F)$ be the subgroup generated by $U_{\alpha}(F)_{x, f(\alpha)}$ for all $\alpha \in \Phi$, $U_{x, f}^{(\alpha)} \subset U_{x, f}$ be the subgroup generated by $U_{i\alpha}(F)_{x, f(i\alpha)}$ for $i \in \lrbracket{\pm 1, \pm 2}$, $N_{x, f}^{(\alpha)}:= U_{x, f}^{(\alpha)} \cap N(F)$, $U_{x, f}^{\pm} := U_{x, f} \cap U^{\pm}(F)$, then 
    \begin{proposition}{{\cite[Proposition 6.4.9]{bruhat1972groupes}}}\label{prop: U x, f good factorization}
        \begin{enumerate}
        \item $U_{\alpha, x, f}:= U_{x, f} \cap U_{\alpha}(F) = U_{\alpha}(F)_{x, f(\alpha)}U_{2\alpha}(F)_{x, f(2\alpha)}$.
        \item $N_{x, f}:= U_{x, f} \cap N(F)$ is generated by $N_{x, f}^{(\alpha)}$ for all $\alpha \in \Phi$.
        \item $\prod_{\alpha \in \Phi_{\red}^{\pm}} U_{\alpha, x, f} \to U_{x, f}^{\pm}$ is a bijection in any ordering.
        \item $U_{x, f} = U_{x, f}^-U_{x, f}^+N_{x, f} = U_{x, f}^+U_{x, f}^-N_{x, f}$
    \end{enumerate}
    In particular, $U_{x, f}$ satisfies the nice decomposition \ref{eq: U^-U^+N}.
    \end{proposition}

  Because $Z(F)^1$ has a trivial action on $A_{\red}(G, S)$ (we briefly recall these properties in subsection \ref{subsec: parahoric group schemes}), it normalizes $U_{\alpha}(F)_{x, f(\alpha)}$, which implies its normalization of $U_{\alpha, x, f}$, $U_{x, f}^{(\alpha)}$, $N_{x, f}^{(\alpha)}$, $U^{\pm}_{x, f}$, $N_{x, f}$, and $U_{x, f}$. For any subgroup $X \subset Z(F)^1$, we define $X_{x, f}$ as the subgroup generated by $X$ and $U_{x, f}$. Due to the factorization on Proposition \ref{prop: U x, f good factorization}, we have 
    \begin{equation}
        X_{x, f} = U_{x, f}X_{x, f} = U_{x, f}^-U_{x, f}^+(N(F) \cap X_{x, f}),\ N(F) \cap X_{x, f} = N_{x, f}X.
    \end{equation}
    and $X_{x, f} \cap U_{\alpha}(F) = U_{\alpha, x, f}$, $X_{x, f} \cap U^{\pm}(F) = U_{x, f}^{\pm}$. In particular,
    \begin{corollary}\label{cor: U, X, f, good factorization}
         $X_{x, f}$ satisfies the condition \ref{eq: U^-U^+N}.
    \end{corollary}

    These results with $F$ are deduced from the results with $\bF$ using \'etale descent.

\subsubsection{Yu's theorem}
 Let $U_0/F \mapsto \lrbracket{U_{0}(F)_{x, r}}_{r\in \R}$ be a schematic admissible filtration, let $\UU_{0, x, r}$ be the smooth model of $U_0$ with the set of integral points $\UU_{0, x, r}(\OO_{\bF}) = U_{0}(\bF)_{x, r}$. Let $f': \wdh{\Phi}' \to \wdt{\R}\setminus\lrbracket{\infty}$ be a concave function. Let $G(\bF)_{x, f'} \subset G(\bF)$ be the subgroup generated by $U_{\alpha'}(\bF)_{x, f'(\alpha')}$ for all $\alpha' \in \wdh{\Phi}'$. Let $\UU_{\alpha', x, f'} = \UU_{\alpha', x, f'(\alpha')}$ when $\alpha' \in \Phi'$ but $2\alpha' \not\in \Phi'$, and $\UU_{\alpha', x, f'} = \UU_{\alpha', x, f'(\alpha'), 2f'(\alpha')}$ when $\alpha', 2\alpha' \in \Phi'$.

 \begin{theorem}{{\cite[Theorem 8.1]{yu2015smooth}}}\label{thm: yu main theorem}\leavevmode
     \begin{enumerate}
         \item There is a unique smooth model $\Gg_{x, f'}$ of $G$ over $\OO_{\bF}$ with $\Gg_{x, f'}(\OO_{\bF}) = G(\bF)_{x, f'}$.
         \item The schematic closure of $U_{\alpha'}$ in $\Gg_{x, f'}$ is $\UU_{\alpha', x, f'}$ for any $\alpha' \in \wdh{\Phi}'_{\red}$.
         \item The multiplication morphism
         \begin{equation}\label{eq: open big cell}
             \prod_{\alpha' \in \Phi^{\prime, -}_{\red}} \UU_{\alpha', x, f'} \times \UU_{0, x, f'} \times \prod_{\alpha' \in \Phi^{\prime, +}_{\red}} \UU_{\alpha', x, f'} \to \Gg_{x, f'}
         \end{equation}
         is an open immersion and the products $\prod_{\alpha' \in \Phi^{\prime, \pm}_{\red}}$ can be arranged in any ordering. Moreover, when $f'(0) > 0$, \ref{eq: open big cell} is an isomorphism on the special fiber. 
     \end{enumerate}
 \end{theorem}
 \begin{remark}\label{rmk: combine to positively closed set}
    \cite[\S 4.6.2]{bruhat1984groupes} says that, let $\Psi' \subset \Phi'$ be a positively closed subset, there exists a unique smooth group scheme model $\UU_{\Psi', x, f'}$ of $U_{\Psi'}$ such that the $\bF$-isomorphism $\prod_{\alpha' \in \Psi'_{\red}} U_{\alpha'} \to U_{\Psi'}$ extends to an $\OO_{\bF}$-isomorphism $\prod_{\alpha' \in \Psi'_{\red}} \UU_{\alpha', x, f'} \to \UU_{\Psi', x, f'}$.
    
    The proof for part $(2)$ in Theorem \ref{thm: yu main theorem} also imply that the closure of $U_{\Psi'}$ in $\Gg_{x, f'}$ is $\UU_{\Psi', x, f'}$. The open big cell in part $(4)$ can also be written as
    \begin{equation}
        \UU_{x, f'}^+ \times \UU_{0, x, f'} \times \UU_{x, f'}^- \to \Gg_{x, f'},
    \end{equation}
    here $\UU^{\pm}_{x, f'}$ denotes $\UU_{\Phi^{\prime, \pm}, x, f'}$.
\end{remark}

 \subsubsection{\'Etale descent}
 Due to \cite[\S 5.15]{bruhat1984groupes}, when $x \in B_{\red}(G, F)$, $U_{\sigma(\alpha')}(\bF)_{x, r} = \sigma(U_{\alpha'}(\bF)_{x, r})$ for any $\sigma \in \Sigma_0$. The relations between $\Phi$ and $\Phi'$ can be explained as follows: Let $S_a$ be the maximal $F$-anisotropic torus inside $T$, then $S_a \times S \to T$ is an almost product. Let $\alpha' \in \Phi'$, if $\alpha'|_S$ is trivial, then $\alpha'|_{S_a}$ is nontrivial, thus it belongs to $\Phi_a':=\Phi(S_a, Z)$. We separate $\Phi'$ into two parts: $\Phi_s'$ and $\Phi_a'$, $\Phi_s'$ restricts to nontrivial elements in $\Phi$, and $\Phi_a'$ restricts to trivial elements in $\Phi$, for example, see \cite[\S 21.8]{borel2012linear}.
 \begin{definition}\label{def: sigma invariant f}\leavevmode
 \begin{enumerate}
     \item A concave function $f': \wdh{\Phi}' \to \wdt{\R}$ is $\Sigma_0$-\emph{invariant} if $f'(\alpha') = f'(\sigma(\alpha'))$ for any $\alpha' \in \wdt{\Phi}'$ and any $\sigma \in \Sigma_0$.
     \item A concave function $f': \wdh{\Phi}' \to \wdt{\R}$ is \emph{induced} by $f: \wdh{\Phi} \to \wdt{\R}$ if $f'|_{\wdh{\Phi}} = f$.
 \end{enumerate}
 \end{definition}
 In particular, when $f'$ is $\Sigma_0$-invariant, let $\alpha \in \Phi$, we can define $f'(\alpha):= f'(\alpha')$ for any $\alpha' \in \Phi'$ (all such $\alpha'$ form a $\Sigma_0$-conjugacy class).

 \begin{remark}\label{rmk: G x,f descends over F}
     Let $f'$ be a $\Sigma_0$-invariant concave function, then $\UU_{\sigma(\alpha'), x, f'} = \sigma(\UU_{\alpha', x, f'})$. Given the factorization described in Proposition \ref{prop: U x, f good factorization} and the admissible filtration on $U_0$, it follows that $G(\bF)_{x, f'}$ is $\Sigma_0$-invariant. Consequently, $\Gg_{x, f'}$ descends over $\OO_F$, as demonstrated in \cite[Proposition 10.9]{landvogt2006compactification}.
 \end{remark}

  \begin{remark}
     Using the notations in Remark \ref{remark: 2f(alpha) notation 1}, we have for all $\alpha' \in \Phi'$
     \[ \UU_{\alpha', x, f'}(\OO_{\bF}) = U_{\alpha'}(\bF)_{x, f'(\alpha')}U_{2\alpha'}(\bF)_{x, 2f'(\alpha')}.\] 
     When $f'$ is $\Sigma_0$-invariant, by denoting $\UU_{\alpha, x, f} = \UU_{\alpha, x, f(\alpha)}$ when $\alpha \in \Phi$ but $2\alpha \not\in \Phi$, and $\UU_{\alpha, x, f} = \UU_{\alpha, x, f(\alpha), 2f(\alpha)}$ when $\alpha, 2\alpha \in \Phi$, then for all $\alpha \in \Phi$, we have
     \[\UU_{\alpha, x, f}(\OO_{\bF}) = \wdt{U}_{\alpha}(\bF)_{x, f(\alpha)}\wdt{U}_{2\alpha}(\bF)_{x, 2f(\alpha)}.\]
     If $2\alpha \not\in \Phi$, then $\lrbracket{\alpha' \in \Phi', \alpha'|_S = 2\alpha}$ is empty, the notation here is well-defined.
 \end{remark}
 
  While the \'etale descent process for unipotent groups is standard, and could be found in \cite[\S 5]{bruhat1984groupes}, we provide some detailed calculation specified in our situation. By arranging the ordering, the big cell \ref{eq: open big cell} is:
  \begin{equation}
      (\prod_{\alpha' \in \Phi^{\prime, -}_{s, \red}} \UU_{\alpha', x, f'})\times (\prod_{\alpha' \in \Phi^{\prime, -}_{a, \red}} \UU_{\alpha', x, f'}) \times \UU_{0, x, f'} \times (\prod_{\alpha' \in \Phi^{\prime, +}_{a, \red}} \UU_{\alpha', x, f'}) \times (\prod_{\alpha' \in \Phi^{\prime, +}_{s, \red}} \UU_{\alpha', x, f'}) \to \Gg_{x, f'}
  \end{equation}

  \begin{lemma}\label{lem: Sigma invariance of unipotent part}
   Let $f'$ be $\Sigma_0$-invariant, then
   \begin{enumerate}
       \item $\prod_{\alpha' \in \Phi^{\prime, \pm}_{a, \red}} \UU_{\alpha', x, f'}$ are $\Sigma_0$-invariant,
       \item $\prod_{\alpha' \in \Phi^{\prime, \pm}_{s, \red}} \UU_{\alpha', x, f'}$ are $\Sigma_0$-invariant and identify with $\prod_{\alpha \in \Phi^{\pm}_{\red}} \UU_{\alpha, x, f} \otimes \OO_{\bF}$.
   \end{enumerate}
  \end{lemma}
  \begin{proof}
      These two products are $\Sigma_0$-invariant since $\UU_{\sigma(\alpha'), x, f'} = \sigma(\UU_{\alpha', x, f'})$ and both sets $\Phi_s'$ and $\Phi_a'$ are $\Sigma_0$-invariant. To show $\prod_{\alpha' \in \Phi^{\prime, \pm}_{s, \red}} \UU_{\alpha', x, f'} = \prod_{\alpha \in \Phi^{\pm}_{\red}} \UU_{\alpha, x, f} \otimes \OO_{\bF}$, it suffices to check $\OO_{\bF}$-points since both sides are smooth models for the same group.
      
      Fix $\alpha \in \Phi_{\red}$. When $\alpha'|_S = \alpha$, then $\alpha'$ is also non-reducible. Given the index set $\lrbracket{\alpha' \in \Phi', \alpha'|_S = 2\alpha}$, we separate it into $\lrbracket{\alpha' \in \Phi'_{\red}, \alpha'|_S = 2\alpha}$ and $\lrbracket{\alpha' \in \Phi'\setminus \Phi'_{\red}, \alpha'|_S = 2\alpha}$, the latter set bijects to $\lrbracket{\alpha' \in \Phi'_{\red}, \alpha'|_S = \alpha}$ by dividing $2$. Thus
      \begin{equation*}
          \prod_{\alpha' \in \Phi', \alpha'|_S = 2\alpha} U_{\alpha'}(\bF)_{x, 2f(\alpha)} = \prod_{\alpha' \in \Phi'_{\red}, \alpha'|_S = 2\alpha} U_{\alpha'}(\bF)_{x, 2f(\alpha)} \prod_{\alpha' \in \Phi'_{\red}, \alpha'|_S = \alpha} U_{2\alpha'}(\bF)_{x, 2f(\alpha)}.
      \end{equation*}

      \begin{align*}
          &\UU_{\alpha, x, f}(\OO_{\bF}) = \wdt{U}_{\alpha}(\bF)_{x, f(\alpha)}\wdt{U}_{2\alpha}(\bF)_{x, 2f(\alpha)} \\
          &= \prod_{\alpha' \in \Phi'_{\red}, \alpha'|_S = \alpha} U_{\alpha'}(\bF)_{x, f(\alpha)} \prod_{\alpha' \in \Phi'_{\red}, \alpha'|_S = 2\alpha}U_{\alpha'}(\bF)_{x, 2f(\alpha)} \\
           & \quad \quad \quad  \prod_{\alpha' \in \Phi', \alpha'|_S = 2\alpha} U_{\alpha'}(\bF)_{x, 2f(\alpha)} \prod_{\alpha' \in \Phi'_{\red}, \alpha'|_S = 4\alpha} U_{\alpha'}(\bF)_{x, 4f(\alpha)} \\
          & = (\prod_{\alpha' \in \Phi'_{\red}, \alpha'|_S = \alpha} U_{\alpha'}(\bF)_{x, f(\alpha)}\prod_{\alpha' \in \Phi'_{\red}, \alpha'|_S = \alpha} U_{2\alpha'}(\bF)_{x, 2f(\alpha)}) \\
          & \quad \quad (\prod_{\alpha' \in \Phi'_{\red}, \alpha'|_S = 2\alpha}U_{\alpha'}(\bF)_{x, 2f(\alpha)} \prod_{\alpha' \in \Phi'_{\red}, \alpha'|_S = 2\alpha} U_{\alpha'}(\bF)_{x, 2f(\alpha)} \prod_{\alpha' \in \Phi'_{\red}, \alpha'|_S = 4\alpha} U_{\alpha'}(\bF)_{x, 4f(\alpha)}) \\
          & = \prod_{\alpha' \in \Phi'_{\red}, \alpha'|_S = \alpha} \UU_{\alpha', x, f'}(\OO_{\bF}) \prod_{\alpha' \in \Phi'_{\red}, \alpha'|_S = 2\alpha} \UU_{\alpha', x, f'}(\OO_{\bF}).
      \end{align*}
  \end{proof}

   Recall that $Z$ is a Levi subgroup of $G$. Thanks to \cite[\S 7.6.4]{bruhat1972groupes} or \cite[\S 2.1.4]{landvogt2000some}, there exists a unique $N'(\bF)$-equivariant morphism $\pr: A_{\red}(G_{\bF}, T_{\bF}) \to A_{\red}(Z_{\bF}, T_{\bF})$ of affine spaces. Since $T$ is defined over $F$, $\pr$ is also $\Sigma_0$-invariant. We are using a slight abuse of notation here, where we let $x$ represent the image of $x$ under the composition
   \begin{equation}
       A_{\red}(G, S) \to A_{\red}(G_{\bF}, T_{\bF}) \to A_{\red}(Z_{\bF}, T_{\bF}),
   \end{equation}

   Apply the theorem \ref{thm: yu main theorem} to $Z$ with respect to $x \in A_{\red}(Z_{\bF}, T_{\bF})^{\Sigma_0} = A_{\red}(Z, S)$ and to the restriction of a $\Sigma_0$-invariant concave function $f'|_{\wdh{\Phi}_a'}: \wdh{\Phi}'_a \to \wdt{\R}$, we get an open big cell:
   \begin{equation}\label{eq: apply theory to Z}
       \prod_{\alpha' \in \Phi^{\prime, \pm}_{a, \red}} \UU_{\alpha', x, f'} \times \UU_{0, x, f'} \times \prod_{\alpha' \in \Phi^{\prime, \pm}_{a, \red}} \UU_{\alpha', x, f'} \to \ZZ_{x, f'}
   \end{equation}
   Each part of \ref{eq: apply theory to Z} descends to $F$ due to Remark \ref{rmk: G x,f descends over F} and Lemma \ref{lem: Sigma invariance of unipotent part}, therefore:
   \begin{corollary}\label{cor: open big cell, over F}
      Let $f'$ be a $\Sigma_0$-invariant concave function, then
       \begin{enumerate}
           \item The schematic closure of $U_{\alpha}$ in $\Gg_{x, f'}$ is $\UU_{\alpha, x, f'}$ for any $\alpha \in \Phi_{\red}$.
         \item The multiplication morphism
         \begin{equation}\label{eq: open big cell, over F}
             \prod_{\alpha \in \Phi^{-}_{\red}} \UU_{\alpha, x, f'} \times \ZZ_{x, f'} \times  \prod_{\alpha \in \Phi^{+}_{\red}} \UU_{\alpha, x, f'} \to \Gg_{x, f'}
         \end{equation}
         is an open immersion of smooth models over $\OO_F$ and the product $\prod_{\alpha \in \Phi^{\pm}_{\red}}$ can be taken in any ordering. Moreover, when $f'(0) > 0$, \ref{eq: open big cell, over F} is an isomorphism on the special fiber. 
       \end{enumerate}
   \end{corollary}
   \begin{proof}
       It suffices to check $\OO_{\bF}$-points. The first part comes from the part $(2)$ of Proposition \ref{thm: yu main theorem} and Lemma \ref{lem: Sigma invariance of unipotent part}. The second part comes from the open big cell \ref{eq: apply theory to Z} and the part $(3)$ of Proposition \ref{thm: yu main theorem}.
   \end{proof}

\begin{remark}\label{rmk: combine to positively closed set, over F}
    Remark \ref{rmk: combine to positively closed set} can also be made over $F$. According to \cite[\S 4.6.2, 5.2.3]{bruhat1984groupes}, let $\Psi \subset \Phi$ be a positively closed subset, there exists a unique smooth group scheme model $\UU_{\Psi, x, f}$ of $U_{\Psi}$ such that the $F$-isomorphism $\prod_{\alpha \in \Psi_{\red}} U_{\alpha} \to U_{\Psi}$ extends to an $\OO_F$-isomorphism $\prod_{\alpha \in \Psi_{\red}} \UU_{\alpha, x, f} \to \UU_{\Psi, x, f}$. $\UU^{\pm}_{x, f}$ denotes $\UU_{\Phi^{\pm}, x, f}$.
    
    The proof for part $(2)$ in Theorem \ref{thm: yu main theorem} also implies that the closure of $U_{\Psi}$ in $\Gg_{x, f}$ is $\UU_{\Psi, x, f}$, note that it does not rely on the assumption that $F$ is strictly henselian since it only uses \cite[Lemma 7.1]{yu2015smooth}. The open big cell expression in \ref{eq: open big cell} (resp. \ref{eq: open big cell, over F}) can also be alternatively written as
    \begin{equation}
        \UU_{\Phi^+, x, f'} \times \UU_{0, x, f'} \times \UU_{\Phi^-, x, f'} \to \Gg_{x, f'}.\quad (\textit{resp.}\ \UU_{x, f'}^+ \times \ZZ_{x, f'} \times \UU_{x, f'}^- \to \Gg_{x, f'}),
    \end{equation}
    where $\UU^{\pm}_{x, f}$ denotes $\UU_{\Phi^{\pm}, x, f}$.
\end{remark}

 \begin{remark}
     When considering the open big cell \ref{eq: open big cell, over F}, we typically use $f'$ induced by $f$ in practice. In this case, we can omit the subscript $\prime$ on $f'$.
 \end{remark}
 \begin{remark}
     Let $f(\alpha) = r$ for all $\alpha \in \wdh{\Phi}$, $f$ is a concave function.
     \begin{enumerate}
         \item when $r = 0$, then $G(F)_{x, 0}:= \Gg_{x, f}(\OO_F)$ is the \emph{parahoric subgroup} of $G(F)$ with respect to $x$, see \cite[\S 3.2]{yu2015smooth}.
         \item when $r > 0$, then $G(F)_{x, r}:= \Gg_{x, f}(\OO_F)$ is the \emph{Moy-Prasad group} with respect to $x$ and level $r$, see \cite[\S 8.2]{yu2015smooth}.
     \end{enumerate}
     Also, see Definition \ref{def: parahoric subgroup}.
 \end{remark}
\subsubsection{}
\begin{proposition}\label{proposition: MP-group satisfies two conditions}
   Let $x \in A_{\red}(G, S)$, $f': \wdh{\Phi}' \to \wdh{\R}$ be a $\Sigma_0$-invariant concave function, $K = G(F)_{x, f'}:= \Gg_{x, f'}(\OO_F)$, then $K$ satisfies Conditions \ref{eq: U^-U^+N} and \ref{eq: U^+U^-U^-Z}.
\end{proposition}
\begin{proof}
   The condition \ref{eq: U^+U^-U^-Z} is implied by the open big cell \ref{eq: open big cell, over F}, the third part of Lemma \ref{lem: special fiber and integral points} when $f'(0) > 0$, and the first two parts of Lemma \ref{lem: special fiber and integral points} when $f'(0) = 0$. This holds because, in the latter case, $\Gg_{x, f}$ becomes the parahoric subgroup scheme and is connected.

   The decomposition \ref{eq: U^+U^-U^-Z} implies that $G(F)_{x, f'}$ is generated by $U_{x, f'}$ (defined in Proposition \ref{prop: U x, f good factorization}) and $Z(F)_{x, f'}$. Consequently, Condition \ref{eq: U^-U^+N} follows from Corollary \ref{cor: U, X, f, good factorization}.
\end{proof}

\begin{lemma}\label{lem: special fiber and integral points}
    Let $\bc$ be a smooth scheme over $\OO_F$, $\Gg$ be a smooth group scheme over $\OO_F$, assume there exists an open immersion $\bc \to \Gg$.
    \begin{enumerate}
        \item If $\Gg\otimes k_F$ is connected, then $\Gg(k_F) = \bc(k_F)^{-1}\bc(k_F)$.
        \item If $\Gg(k_F) = \bc(k_F)^{-1}\bc(k_F)$, then $\Gg(\OO_F) = \bc(\OO_F)^{-1}\bc(\OO_F)$.
        \item If $\Gg(k_F) = \bc(k_F)$, then $\Gg(\OO_F) = \bc(\OO_F)$.
    \end{enumerate}
\end{lemma}
\begin{proof}
    \begin{enumerate}
        \item The inverse group law on $\Gg$ gives rise to a smooth scheme $\bc^{-1}$ over $\OO_F$ along with an open immersion $\bc^{-1} \to \Gg$. When $\Gg\otimes k_F$ is connected, both $\bc \otimes k_F$ and $\bc^{-1} \otimes k_F$ are open and Zariski dense. Consequently, we have $\Gg\otimes k_F = (\bc^{-1} \otimes k_F)(\bc \otimes k_F)$. Taking $k_F$-points into account, we obtain the equation $\Gg(k_F) = \bc^{-1}(k_F)\bc(k_F) = \bc(k_F)^{-1}\bc(k_F)$.
        \item Let $g \in \Gg(\OO_F)$, and $\Bar{g} \in \Gg(k_F)$. Given that $\bc^{-1}$ is smooth over $\OO_F$, we have $\bc^{-1}(\OO_F) = \bc(\OO_F)^{-1} \to \bc^{-1}(k_F)$ being surjective. Hence, we can choose $g_1^{-1} \in \bc(\OO_F)^{-1}$ such that $\Bar{g}_1\Bar{g} \in \bc(k_F)$. Now, we claim that $g_1g \in \bc(\OO_F)$.

The proof is standard, and it follows a similar argument as presented in \cite[Proof of the first Lemma in \S 7.4]{yu2015smooth}. We need to demonstrate that $g_1g \in \bc(F)$. If, by way of contradiction, we assume that $g_1g \not\in \bc(F)$, let $\mathfrak{Z}$ be the complement of $\bc$ in $\Gg$, equipped with the induced reduced closed subscheme structure. Also, let $Z$ be the generic fiber of $\mathfrak{Z}$. Then, we would have $g_1g \in Z(F)$.

Since the closure of $Z$ in $\Gg$ is contained within $\mathfrak{Z}$, we deduce that $g_1g \in \mathfrak{Z}$, which in turn implies that $\Bar{g}_1\Bar{g} \in \mathfrak{Z}(k_F)$. However, this contradicts the fact that $\Gg(k_F) = \mathfrak{Z}(k_F) \sqcup \bc(k_F)$, where the two sets are disjoint.
        \item It follows from the proof of the second part.
    \end{enumerate}
\end{proof}

\subsection{Quasi-parahoric subgroups and Bruhat-Tits theory}\label{subsec: parahoric group schemes}

In this subsection, we will recall other groups defined using Bruhat-Tits theory and examine their intersections with Levi subgroups and normal subgroups. We will continue to use the notations introduced in the previous subsection. The main definitions are \ref{def: type ast groups, in an appartment} and \ref{def: type ast groups, in a building}, the main proposition is \ref{lem: parahoric subgroup satisfies the condition}.

\subsubsection{}
Let $V_1 = X_*(S)$, $V_0 = \lrbracket{v \in V_1|\ \lrangle{v, \alpha} = 0, \forall \alpha\in \Phi}$, $V = X_*(S/Z_{c, sp}) \otimes \R \cong V_1/V_0$, let
\begin{equation}
    \nu_Z: Z(F) \to \Hom_{\Z}(X^*(Z)_F, \Z) \otimes \R = \Hom_{\Z}(X^*(S), \Z)\otimes\R = V_1.
\end{equation}
Let $\nu: Z(F) \to V$ be the composition of $\nu_Z$ and the projection $V_1 \to V$. $A_{\red}(G, S)$ is the unique affine space under $V$ unique up to a unique isomorphism such that $\nu$ extends to morphism $\nu_N: N(F) \to \Aff(A_{\red}(G, S))$. 
\begin{equation*}
    Z(F)\to V \to \Aff(A_{\red}(G, S)): v\mapsto [w \mapsto v + w]
\end{equation*}
acts as translation, $N(F)/Z(F) = W$ acts as $W \to \GL(V)$. Let
\begin{equation}\label{equation: nuG another definition}
    \nu_G: G(F) \to \Hom_{\Z}(X^*(G)_F, \Z) \otimes \R = X_*(Z_{c, sp})\otimes \R = V_0.
\end{equation}
It defines an affine action (translation) on $V_0$:
\[\nu_G: G(F) \to \Aff(V_0), \quad \nu_G(g)(v) = v + \nu_G(g)    \]
The extended affine apartment $A_{\ext}(G, S) = A_{\red}(G, S) \times V_0$ associated with $(G, S)$ is endowed with an affine action:
\[ \nu_{N, \ext}: N(F) \stackrel{(\nu_N, \nu_G)}{\longrightarrow} \Aff(A_{\red}(G, S)) \times \Aff(V_0) \subset \Aff(A_{\ext}(G, S)).  \]
Note that this $\nu_G$ is the same as the one defined in \ref{eq: Kottwitz, over F}, thus $\Ker \nu_G = G(F)^1$. If $g \in G(F)$ fixes $x_0$, then $g$ fixes $V_0$, and $g \in G(F)^1$, vice versa.

Given any point $x \in A_{\red}(G, S)$, one can associate with it a family of filtrations $\varphi = \lrbracket{\varphi_{\alpha}}_{\alpha \in \Phi}$ on $\lrbracket{U_{\alpha}(F)}$ in the sense of \cite[Definition 6.2.1]{bruhat1972groupes} (these define the filtrations we used in Subsection \ref{subsubsec: unipotent groups}). Let 
\begin{enumerate}
    \item $\Gamma_{\alpha} = \varphi_{\alpha}(U_{\alpha}(F)\setminus\lrbracket{1})$,
    \item $\Gamma_{\alpha}' = \lrbracket{\varphi_{\alpha}(u)|\ u \in U_{\alpha}(F)\setminus\lrbracket{1}, \varphi_{\alpha}(u) = \sup \varphi_{\alpha}(uU_{2\alpha}(F))}$,
\end{enumerate}
then $\Gamma_{-\alpha} = - \Gamma_{\alpha}$, $\Gamma_{\alpha} = \Gamma_{\alpha}'$ (when $2\alpha \not\in \Phi$), $\Gamma_{\alpha} = \Gamma_{\alpha}' \cup \frac{1}{2}\Gamma_{2\alpha}$. See \cite[\S 4.2.2]{bruhat1984groupes} for more properties about these two sets. Moreover, the set of affine roots $\Phi_{\aff}$ is the set $\lrbracket{\alpha + k|\ \alpha \in \Phi, k \in \Gamma_{\alpha}'}$. Here $\Gamma_{\alpha}$ is a discrete subset of $\R$ since we start with a discrete valuation on $F$.

\subsubsection{Concave functions}
Fix any point $x \in A_{\red}(G, S)$ (thus $\Gamma_{\alpha}, \Gamma_{\alpha}'$, etc. are fixed), let $\Omega \subset A_{\red}(G, S)$ be a non-empty set, let
\begin{equation}
    f_{\Omega}(\alpha) = \inf \lrbracket{ r \in \R | \ \alpha(a - x) + r \geq 0, \forall a\in \Omega},
\end{equation}
\begin{equation}
    f_{\Omega}'(\alpha) = \inf \lrbracket{ r \in \Gamma_{\alpha}' | \ \alpha(a - x) + r \geq 0, \forall a\in \Omega},
\end{equation}
here $f'_{\Omega}$ is the optimazation of $f_{\Omega}$ in the sense of \cite[Proposition 6.4.10]{bruhat1972groupes}, also, see \cite[Remark 6.4.12]{bruhat1972groupes}. The function $f'_{\Omega}(\alpha) = r$ represents the minimal value of $r$ such that the semi-space $H_{\alpha, r} = \lrbracket{a \in A_{\red}|\ \alpha(a - x) + r \geq 0}$ defined by the affine root $\alpha + r$ contains $\Omega$. When $\Omega = \lrbracket{a}$ is a point, then $f_{a}(\alpha) = -\alpha(a - x)$, and $f_{a}'(\alpha) = r$ is the minimal $r \geq -\alpha(a - x)$ such that $r \in \Gamma_{\alpha}'$.

\begin{remark}\label{remark: parahoric subgroup, optimal function}
    It follows from the definition that $f_{\Omega}$ is concave, but $f_{\Omega}'$ might not be concave. In \cite[\S 4.6.26]{bruhat1984groupes}, the authors defined parahoric subgroups using $f'_{\Omega}$ instead of $f_{\Omega}$. $f'_{\Omega}$ is quasi-concave (a condition weaker than being concave) and optimal ($f_{\Omega}' = (f_{\Omega}')'$), allowing for the full application of the theory.
    
    Note that $f_{\Omega}'$ only depends on the simplicial closure $\cl(\Omega)$ of $\Omega$, which is the intersection of all $H_{\alpha, f_{\Omega}'(\alpha)}$, see \cite[\S 4.6.27]{bruhat1984groupes}. However, $f_{\Omega}$ might depend on $\Omega$ itself, not just $\cl(\Omega)$.

    Nevertheless, due to \cite[\S 4.5.2, 4.6.23]{bruhat1984groupes}, since the valuation $\omega$ on $F$ is discrete, $U_{\alpha, x, f} = U_{\alpha, x, f'}$, $U_{x, f'} = U_{x, f}$. Therefore, it does not matter whether to use $f$ or $f'$ when considering groups that only involve $U_{x, f}$ and $Z(F)^1$.
    
\end{remark}

\begin{remark}\label{rmk: change point, change valuation}
    In the last section, we were given a pair $(x, f)$, where $x \in A_{\red}(G, S)$ and $f: \wdh{\Phi} \to \wdh{\R}$ is a concave function. The choice of $x$ is not important in the following sense: if we have another point $x' \in A_{\red}(G, S)$ and define $f'(\alpha) = f(\alpha) + \alpha(x' - x)$ for all $\alpha \in \wdh{\Phi}$, then $G(F)_{x, f} = G(F)_{x', f'}$, as shown in \cite[\S 4.6.16]{bruhat1984groupes}. When $f$ is a concave function, $f'$ is also a concave function since $f' - f$ is a linear function.

Let $\Omega = \lrbracket{x}$, $f_x := f_{\lrbracket{x}}$, then $f_x(\alpha) = 0$ for any $\alpha \in \Phi$. Let $f_{x, r}:= f_x + r : \wdh{\Phi} \to \wdh{\R}$, it is a constant function on $\wdh{\Phi}$ with value $r$. We will denote this function by $r$. Let us pick another point $y \in A_{\red}$ and let $\Omega = \lrbracket{y}$ (note that we fix $x$ as the base point), and set $f_y(0) = 0$, $f_{y, r} = f_y + r$. Then, $f_{y, r}$ is a concave function, and the transformation $f_{x, r} \mapsto f_{y, r}: f_{y,r}(\alpha) = f_{x, r}(\alpha) + \alpha(y - x)$ occurs when we change the valuation from $x$ to $y$. These pairs $(x, f_{x, r})$ and $(y, f_{y, r})$ define same filtrations on $U_{\alpha}(F)$ for all $\alpha \in \wdh{\Phi}$.
\end{remark}

From now on, we fix a base special point $x_o \in A_{\red}$, and all $f_{\Omega}$ are defined relative to this point $x_o$.

\subsubsection{}

We add superscript $\prime$ to those objects defined over $\bF$. Let $(\ast)$ be either $\red$ or $\ext$. There is a canonical $G(F)$-equivariant embedding
\[ \iota: B_{(\ast)}(G, F) \hookrightarrow B_{(\ast)}(G, \bF),\quad  B_{(\ast)}(G, F) =  B_{(\ast)}(G, \bF)^{\Sigma_0} \]
which maps bijectively the facets $\FF$ of $B_{(\ast)}(G, F)$ to $\Sigma_0$-invariant $\FF'$ facets of $B_{(\ast)}(G, \bF)$, the inverse of this bijection is $\FF' \mapsto \FF^{\prime \Sigma_0} = \FF $. In general, there are more affine roots over $\bF$ than over $F$, thus $A_{\red}'$ has finer simplicial structures.

Fix an embedding $A_{\red}(G, S) \to A_{\ext}(G, S)$, $B_{\red}(G, F) \to B_{\ext}(G, F)$. Let $x_e \in B_{\ext}(G, F)$, $x_e = (x, x_0)$, $x \in B_{\red}(G, F)$, $x_0 \in V_0$. Consider the following groups:
\begin{equation}
   \hat{P}_x^1:= \hat{P}_{x_e} :== \lrbracket{ g \in G(F)|\ g x_e = x_e}, \quad \hat{P}_{x} = \lrbracket{ g \in G(F)|\ g x = x},
\end{equation}
then $\hat{P}_{x_e} = \hat{P}_{x} \cap G(F)^1$. In \cite[\S 4.6]{bruhat1984groupes}, $\hat{P}_{x_e}$ is denoted by $\hat{P}_x^1$.

Also, consider the following groups as in Proposition \ref{prop: U x, f good factorization}:
\begin{itemize}
    \item $\hat{N}_{x}=\lrbracket{n\in N(F)|\ n\cdot x = x}$,
    \item $U_{\alpha, x}= U_{\alpha, f_{x}(\alpha)}$, $U_{x} = U_{x, f_x}:= \lrangle{U_{\alpha, x}: \alpha \in \Phi}$, $U^{\pm}_{x}= U^{\pm}(F) \cap U_{x}$.
\end{itemize}

The identification follows from \cite[Proposition 7.4.4]{bruhat1972groupes}, \cite[Proposition 13.3]{landvogt2006compactification}:
\begin{equation}
    \hat{P}_{x} = \lrangle{U_{x}, \hat{N}_{x}} = \lrbracket{ g \in G(F)|\ g x = x}.
\end{equation}

We similarly consider following groups (see \cite[Proposition 4.6.3, 4.6.28]{bruhat1984groupes}):
\begin{itemize}
    \item $\hat{N}_{x}^1:= \hat{N}_x \cap G(F)^1$, $N_x = U_x \cap N(F)$ (as in Proposition \ref{prop: U x, f good factorization}),
    \item $P^1_x = Z(F)^1U_x$, $P^0_x = Z(F)_0U_x$,
    \item $N_x^1 = P_x^1 \cap N(F) = Z(F)^1N_x$, $N_x^0 = P^0_x \cap N(F) = Z(F)_0N_x$.
\end{itemize}

Similar groups are also defined over $\bF$. There exist smooth models of $G$ defined over $\OO_{F}$:
\begin{equation}
   \Ggc_x \subset \Gg_x \subset \Ggh_x, 
\end{equation}
such that
\[\Ggc_x(\OO_{\bF}) = P_{\iota(x)}^0,\quad \Gg_x(\OO_{\bF}) = P_{\iota(x)}^1,\quad \Ggh_x(\OO_{\bF}) = \hat{P}_{\iota(x)}^1,\] 
and $\Ggc_x$ is the neutral connected component of both $\Gg_x$ and $\Ggh_x$.

Since the canonical embeddings $B_{\red}(G, F) \to B_{\red}(G, \bF)$, $B_{\ext}(G, F) \to B_{\ext}(G, \bF)$ are $G(F)$-invariant, and since $G(F)^1 = G(\bF)^1 \cap G(F)$:
\begin{equation}\label{equation: factorization of Ggh}
    \Ggh_{x}(\OO_F) = \hat{P}_{\iota(x)} \cap G(\bF)^1 \cap G(F) = \hat{P}_{x} \cap G(F)^1.
\end{equation}
Moreover, $\Ggc_x(\OO_F) = P^0_{x}$ and $\Gg_x(\OO_F) = P^1_{x}$ by \'etale descent, see \cite[\S 5]{bruhat1984groupes} or \cite[Chapter IV]{landvogt2006compactification}. Therefore, we omit $\iota$ in the subscript if there is no confusion.

The above groups $P^0_{x}$, $P^1_{x}$, $\hat{P}_{x}^1$ and $\hat{P}_{x}$ are the main objects we are considering in this section.

\subsubsection{}
We list some good properties for $\hat{P}_x$, see \cite[\S 5]{bruhat1972groupes} or \cite[\S Chapter IV]{landvogt2006compactification}:
\begin{proposition}\label{proposition: basic properties about groups omega}\leavevmode
    \begin{enumerate}
\item $U_x$ satisfies Proposition \ref{prop: U x, f good factorization},
     \item $\hat{P}_{x} = U^-_{x}U^+_{x}\hat{N}_{x} = U^+_{x}U^-_{x}\hat{N}_{x}$, $\hat{P}_{x} \cap U^{\pm}(F) = U^{\pm}_{x}$, $\hat{P}_{x} \cap N(F) = \hat{N}_{x}$
      \item $nU_{\alpha, f_{x}(\alpha)}n^{-1} = U_{\jmath_W(n)(\alpha), f_{\nu_N(n)(x)}(\jmath_W(n)(\alpha))}$, $n\hat{P}_{x}n^{-1} = \hat{P}_{\nu_N(n)(x)}$ for all $n \in N(F)$, where $\jmath_W: N(F) \to W$ is the natural projection.
     \item $G(F) = U^+(F)N(F)\hat{P}_x$ for any $x \in A_{\red}$ (Iwasawa decomposition).
        \item $G(F) = \hat{P}_xN(F)\hat{P}_y$ for any $x, y \in A_{\red}$ (mixed Bruhat decomposition).
     \item $N(F)$ is the stablizer of $A_{\red}$.
        \item $U_{x}$ acts transitively on the set of apartments containing $x$.
\end{enumerate}
\end{proposition}

Compare to $\hat{P}_x$, it follows that $\hat{P}_x^1$ satisfies the condition \ref{eq: U^-U^+N}: $\prod_{\alpha \in \Phi^{\pm}_{\red}} U_{\alpha, x} \to U_{x}^{\pm}$ is a bijection for any ordering, and
\begin{equation}
    \hat{P}_{x}^1 = U^-_{x}U^+_{x}\hat{N}_{x}^1 = U^+_{x}U^-_{x}\hat{N}_{x}^1.
\end{equation}

Replace $G$ with the reductive group $Z$, $A_{\red}(Z, S)$ consists of a single point $\lrbracket{z}$. Let $\ZZ = \ZZ_{z}$, then it is the schematic closure of $Z$ in any $\Ggh_x$, $x \in A_{\red}(G, S)$. $\ZZ(\OO_{\bF}) = Z(\bF)^1$, see \cite[\S 5.2.1]{bruhat1984groupes}. Let $\ZZ^{\circ}$ be its connected component, then $\ZZ^{\circ}(\OO_{\bF}) = Z(\bF)_0$. 

Also, due to \cite{haines2008parahoric} (over $\bF$) and \cite{richarz2016iwahori} (over $F$), we always have:
\begin{equation}\label{eq: HR appendix}
    \Ggc_x(\OO_F) = \Ggh_x(\OO_F) \cap G(F)_0 = \Ggh_x(\OO_F) \cap \Ker \kappa_G.
\end{equation}

Referring to \cite[\S 5.2.4]{bruhat1984groupes}, as noted in \cite{haines2009corrigendum}, it should be mentioned that all the $\hat{(\ast)}$ symbols should be omitted, one has open big cells:
\begin{equation}\label{eq: big cell for parahoric}
    \UU_{x}^+ \times \ZZ \times \UU_x^- \to \Gg_x,\quad \UU_{x}^+ \times \ZZ^{\circ} \times \UU_x^- \to \Ggc_x.
\end{equation}

In conjunction with Corollary \ref{cor: U, X, f, good factorization} and Lemma \ref{lem: special fiber and integral points}, we obtain:

\begin{proposition}\label{lem: parahoric subgroup satisfies the condition}\leavevmode
\begin{enumerate}
    \item $\Ggh_x(\OO_F) = \hat{P}^1_x$ satisfies condition \ref{eq: U^-U^+N},
    \item $\Gg_x(\OO_F) = P^1_x$, $\Ggc_x(\OO_F) = P^0_x$ satisfy conditions \ref{eq: U^-U^+N} and \ref{eq: U^+U^-U^-Z}.
\end{enumerate}
\end{proposition}

\subsubsection{}

Let us recall the basic definitions:
\begin{definition}\label{def: parahoric subgroup}\leavevmode
   \begin{enumerate}
       \item A \emph{parahoric subgroup} (resp. a \emph{neutral quasi-parahoric subgroup}, \emph{stablizer quasi-parahoric subgroup}) of $G(F)$ with respect to $x \in A_{\red}(G, S)$ is defined as $\Ggc_x(\OO_{F})$ (resp. $\Gg_x(\OO_{F})$, $\Ggh_{x}(\OO_F)$).
       \item A \emph{quasi-parahoric subgroup} of $G(F)$ with respect to $x \in A_{\red}(G, S)$ is defined as $\Kf_p = (\KKf_p)^{\Sigma_0}$, for some $\Sigma_0$-stable open subgroup $\KKf_p \subset G(\bF)$ controlled by $\Ggc_x(\OO_{\bF}) \subset \KKf_p \subset \Ggh_x(\OO_{\bF})$.
       \item Let $r \in \R_{\geq 0}$, the \emph{Moy-Prasad} subgroup $G_{x, r}$ of $G(F)$ with respect to $x \in A_{\red}(G, S)$ and $r$ is defined as $G(F)_{x, r} := \Gg_{x, r}(\OO_F)$, where $r: \wdh{\Phi} \to \R$ is the constant function with value $r$. When $r = 0$, $G_{x, 0}$ is the parahoric subgroup $\Ggc_x(\OO_{F})$. Also, let $G_{x, r +} = \cup_{\epsilon > 0} G_{x, r + \epsilon}$.
   \end{enumerate}
\end{definition}

\begin{remark}\label{rmk: simplicial closure of x}
    Due to Remark \ref{remark: parahoric subgroup, optimal function}, $\Gg_x$ and $\Ggc_x$ only depend on the simplicial closure $\cl(x)$ of $x$. Therefore, we write $\Gg_{x} = \Gg_{\FF}$, $\Ggc_{x} = \Ggc_{\FF}$ for a facet $\FF \subset A_{\red}$ when $\cl(x) = \cl(\FF)$.
\end{remark}

\subsubsection{}
To define the Moy-Prasad groups with respect to a general point $x \in B_{\ext}(G, F)$, we introduce the following two lemmas, which are also useful in other contexts:

\begin{lemma}\label{lemma: enhanced Iwasawa decomposition}
    Let $x \in A_{\red}(G, S)$, $Q$ be any parabolic $F$-subgroup of $G$ containing $S$. Then
    \[  G(F) = Q(F)N(F)P_{x}^{0}, \quad G(F) = P_{x}^{0}N(F)P_{x}^{0} \]
    where $P_{x}^{0} = \Ggc_{x}(\OO_F)$.
\end{lemma}
\begin{proof}
    This is a slightly different version of the Iwasawa decomposition presented in Proposition \ref{proposition: basic properties about groups omega}. Since $G(F) = U^+(F)N(F)\hat{P}_{x}$ holds for any ordering on $\Phi$, it follows that $G(F) = Q(F)N(F)\hat{P}_{x}$. Since $\hat{P}_{x} \cap G(F)_0 = P_{x}^{0}$ (see \ref{eq: HR appendix}) and since $G(F)_0$ is a normal subgroup of $G(F)$, then $\hat{P}_{x}$ normalizes $P_{x}^{0}$. Since $\hat{P}_x = U_x^-U_x^+(\hat{P}_x \cap N(F))$ and $P_x^0 = U_x^-U_x^+(P_x^0 \cap N(F))$, then $\hat{P}_{x} \subset P_{x}^{0}N(F)$. Similarly,  $\hat{P}_{x} \subset N(F)P_{x}^{0}$. Consequently, the mixed Iwasawa decomposition $G(F) = Q(F)N(F)\hat{P}_{x}$ implies $G(F) = Q(F)N(F)P_{x}^{0}$. The mixed Bruhat-decomposition $G(F) = \hat{P}_{x}N(F)\hat{P}_{x}$ implies $G(F) = P_{x}^{0}N(F)P_{x}^{0}$.
\end{proof}

\begin{remark}\label{remark: bon group}
    When $x$ is a special point, $\hat{P}_x^1=\Ggh(\OO_F)$ is a \emph{bon} group in the sense of \cite[\S 4.4]{bruhat1972groupes}, and as a result, we have $G(F) = Q(F)\hat{P}_x^1$. Furthermore, this relationship is extended to $G(F)=Q(F)P_x^0$ thanks to \cite[Corollary 9.1.2]{haines2010satake}.
\end{remark}

\begin{definition}\label{def: type ast groups, in an appartment}
    Let $x \in A_{\red}(G, S)$, $r \in \R_{\geq 0}$. A subgroup $K_{x, (r)} \subset G(F)$ is said to be of type $(\ast)$ (see also \ref{definition: type ast}) if $K_{x, (r)}$ is one of the following groups:
    \begin{equation}
        \hat{P}_x^1=\Ggh_{x}(\OO_F), \quad P_x^0=\Ggc_{x}(\OO_F), \quad P_x^1= \Gg_{x}(\OO_F), \quad G_{x, r} = G(F)_{x, r}.
    \end{equation}
\end{definition}

\begin{lemma}\label{lemma: N action keeps the type}
    Let $K_{x, (r)}$ be a subgroup of $G(F)$ of type $(\ast)$, then for any $n \in N(F)$, $nK_{x, (r)}n^{-1} = K_{\nu_N(n)x, (r)}$.
\end{lemma}
\begin{proof}
    When $K_{x, (r)}$ is $\hat{P}^1_x = \hat{P}_{x} \cap G(F)^1$ or $P_x^0 = \hat{P}_{x} \cap G(F)_0$ or $P_x^1 = P_x^0\ZZ(\OO_F)$, $nK_{x, (r)}n^{-1} = K_{\nu_N(n)x, (r)}$ follows from: 
    \begin{itemize}
        \item $n\hat{P}_{x}n^{-1} = \hat{P}_{\nu_N(n)x}$ (see Proposition \ref{proposition: basic properties about groups omega}),
        \item $G(F)^1, G(F)_0 \subset G(F)$ are normal subgroups since they are the kernels of $\nu_G$ and $\kappa_G$ respectively,
        \item $N$ normalizes the open subgroup $\ZZ(\OO_F) \subset Z(F)$ since the apartment of $Z$ consists of a single point.
    \end{itemize}

    Let $K_{x, r} = G_{x, r}$ be a Moy-Prasad group of level $r$. $G_{x, r}$ is generated by $U_{\alpha, x, r}$ for $\alpha \in \Phi$ and $Z(F)_{x, r}:= \ZZ_{x, r}(\OO_F)$ (here $r: \wdh{\Phi} \to \R$ is the constant function with value $r$). Since $Z(F)_0$ normalizes $Z(F)_{x, r}$ and $U_{\alpha, x, r}$, and there is a natural isomorphism $N'(F)/Z'(F)_0 \rightiso N(F)/Z(F)_0$ due to \cite[Lemma 1.6]{richarz2016iwahori}, considering the action of $n \in N(F)$ on $G_{x, r}$ is equivalent to considering the action of $n' \in N'(F)$. It follows from the definition that $N'(F)$ normalizes each $Z'(F)_{r}^{\MPF}$. Due to Proposition \ref{proposition: basic properties about groups omega}, for any $\alpha' \in \Phi'$, one has:
    \[n'U_{\alpha', f_{x}(\alpha')}n^{\prime, -1} = U_{\jmath_{W'}(n')(\alpha'), f_{\nu_N(n)(x)}(\jmath_{W'}(n')(\alpha'))}.\] 
    Therefore, $nK_{x, r}n^{-1} = K_{\nu_N(n)(x), r}$.
\end{proof}

\begin{definition}\label{def: type ast groups, in a building}
    Recall that $B_{\red}(G, F) = \cup_{g \in G(F)}gA_{\red}$, let $x_1 \in B_{\red}(G, F)$, then $x_1 = g x$ for some $x \in A_{\red}(G, S)$, $g \in G(F)$. Following the Definition \ref{def: type ast groups, in an appartment}, we define the type $(\ast)$ groups as $K_{x_1, (r)} = gK_{x, (r)}g^{-1}$, $K_{x_1, (r+)} = gK_{x, (r+)}g^{-1}$. When $x_1 \in B_{\ext}(G, F)$, we define $K_{x_1, r}$ as $K_{\pr(x_1), r}$.
\end{definition}
This is well-defined since it does not depend on the apartment containing $x$, see \cite[\S 4.6.30]{bruhat1984groupes} and the following calculation: when $x_1 = g_1x = g_2x$ for $g_1, g_2 \in G(F)$, then $g_1^{-1}g_2 \in \hat{P}_{x}$. $\hat{P}_x \subset N(F)P_x^0$ implies that $g_1^{-1}g_2 = nk$ for some $n \in N(F) \cap \hat{P}_{x}$ and $k\in P_x^0$. Due to Lemma \ref{lemma: N action keeps the type}, 
\[g_1^{-1}g_2 K_{x, (r)} g_2^{-1}g_1= nkK_{x, (r)}k^{-1}n^{-1} = nK_{x, (r)}n^{-1} = K_{x, (r)} \Longrightarrow g_1 K_{x, (r)} g_1^{-1}= g_2 K_{x, (r)} g_2^{-1}.\] 

\subsection{Computing levels at boundary}\label{subsec: compute levels at boundary}

In this subsection, we utilize the axiomatic descriptions provided in Subsection \ref{subsec: levi and normal subgroups}. We combine these descriptions with the verification of conditions outlined in Subsections \ref{subsec: bruhat tits theory} and \ref{subsec: parahoric group schemes} to compute the levels at the boundary. The results are summarized in Proposition \ref{prop: main prop for section BT theory}.

\begin{remark}\label{rmk: choose Q standard}
 All maximal $F$-split tori are conjugated by $G(F)$, and there is a one-to-one correspondence between the reduced apartments $A_{\red}$ in $B_{\red}$ and these tori. When considering the rational boundary components of the given Shimura datum, we choose a specific representative element $Q$ from each conjugacy class of admissible parabolic subgroups. These representatives are selected such that each $Q$ corresponds to the standard parabolic subgroup associated with a fixed choice of $(S, J \subset \Delta)$ for some $J$.
 
 Let $K_{x, (r)} \subset G(F)$ be a type $(\ast)$ group with respect to $x \in A_{\red}(G, S)$ and $r\in\R_{\geq 0}$ and let $g \in G(F)$. Due to Lemma \ref{lemma: enhanced Iwasawa decomposition}, $g$ decomposes as $ulnk$, $l \in L(F)$, $u\in R_uQ(F)$, $n \in N(F)$, $k \in P_{x}^{0}$. The level group at boundary is of the form:
\begin{equation}\label{eq: level group at boundary}
    \pi(Q(F) \cap gK_{x, (r)}g^{-1}) = \pi(Q(F) \cap ulK_{\nu_N(n)(x), (r)}(ul)^{-1} ) = l\pi(Q(F) \cap K_{\nu_N(n)(x), (r)}) l^{-1},
\end{equation}
here $\nu_N(n)(x) \in A_{\red}(G, S)$.
\end{remark}

\subsubsection{Special points}
Let us recall the following definitions introduced in \cite{tits1979reductive}: 

Let $\theta \in \Phi_{\aff}$, $r_{\theta}$ denotes the reflection along the vanishing affine hyperplane $H_{\theta}:=\theta^{-1}(0)$. Let $W_{\aff}$ be the affine Weyl group of the affine root system $\Phi_{\aff}$, generated by all $r_{\theta}, \theta \in \Phi_{\aff}$. $W_{\aff}$ acts transitively on all alcoves. Let $\alc$ be a fixed alcove whose simplicial closure contains a fixed special point. Let $H_i$ be the walls bounding $\mathfrak{a}$, and $\theta_i$ be the unique affine root with $H_{\theta_i} = H_i$, $\frac{1}{2}\theta_i \not\in \Phi_{\aff}$, then the set $\lrbracket{\theta_i}_i$ is called a \emph{basis} of $\Phi_{\aff}$ associated with $\mathfrak{a}$. Given $(A_{\red}, \Phi_{\aff})$, one could construct a \emph{local dynkin diagram} $\Delta_{\aff}$ as in \cite[\S 1.8]{tits1979reductive}, the vertices are those affine roots $\theta_i$ in a fixed basis. The local dynkin diagram is independent of the choice of $\mathfrak{a}$.

Let $x \in A_{\red}$, since $W_{\aff}$ acts transitively on alcoves, there exists $w\in W_{\aff}$ such that $wx \in \Bar{\mathfrak{a}}$. Let $I_x:=\lrbracket{\theta_i|\ wx\not\in H_i}$, $\Phi_{\aff, x}\subset \Phi_{\aff}$ be the set of elements which vanish at $x$, $\Phi_i \subset \Phi$ be the vector part of $\Phi_{\aff, x}$, $W_x \subset W_{\aff}$ be the set of reflections $r_{\theta}$ for those $\theta\in\Phi_{\aff, x}$. The weyl group of $\Phi_x$ is the vector part of $W_x$, and its associated ordinary dynkin diagram is isomorphic to the sub-diagram of $\Delta_{\aff}$ by deleting all vertexes belonging to $I_x$ and the edges containing those vertexes. Note that $\Phi_x$ does not depend on the choice of the affine apartment containing $x$.

\begin{definition-proposition}\label{def: special points}
     A vertex $x \in B_{\red}$ is called a \emph{special} point if $\Phi$ and $\Phi_x$ have the same Weyl group. There always exists a special point.
\end{definition-proposition}
Note that if $x \in B_{\red}(G, F)$ is special in $B_{\red}(G, F_1)$ (under the canonical embedding $B_{\red}(G, F) \to B_{\red}(G, F_1)$) for some unramified extension $F_1|F$, then $x$ is special in $B_{\red}(G, F)$.
\begin{definition}{{\cite[\S 1.10.2]{tits1979reductive}, \cite[\S G]{kaletha2021supercuspidal}}}\label{def: various special points}
    Given a vertex $x \in B_{\red}(G, F)$,
    \begin{enumerate}
        \item $x$ is \emph{hyperspecial} if and only if $G$ splits over an unramified extension $F_1$ of $F$ and $x$ is special in $B_{\red}(G, F_1)$.
        \item $x$ is \emph{absolutely special} if and only if $x$ is special in $B_{\red}(G, F_1)$ for any finite Galois extension $F_1|F$.
        \item $x$ is \emph{superspecial} if and only if $x$ is special in $B_{\red}(G, F_1)$ for any finite unramified extension $F_1|F$.
    \end{enumerate}
\end{definition}
\begin{remark}\leavevmode
\begin{enumerate}
    \item When $G$ is unramified, then all these three notions agree.
    \item When $G$ is quasi-split, then $G$ is unramified if and only if there exists a hyperspecial point in $B_{\red}(G, F)$. 
    \item When $G$ is ramified, there is no hyperspecial point, the latter two notions are replacements of the hyperspecial points.
    \item If there exists a superspecial point in $B_{\red}(G, F)$, then $G$ is quasi-split. See \cite[\S G]{kaletha2021supercuspidal} for detailed discussion.
\end{enumerate}
\end{remark}
\begin{definition}{{\cite[\S 1.10.2]{tits1979reductive}}}\leavevmode
    \begin{enumerate}
        \item $G$ is \emph{residually quasi-split} over $F$ if there exists a $\Sigma_0$-invariant alcove $\alc' \subset B_{\red}(G, \bF)$.
        \item $G$ is \emph{residually split} over $F$ if there exists a $\Sigma_0$-fixed alcove $\alc' \subset B_{\red}(G, \bF)$.
    \end{enumerate}
\end{definition}
\begin{remark}
    In \cite[\S 1.10.3, 1.10.4]{tits1979reductive}, it is asserted that when the residue field of $F$ is finite, the group $G$ is residually quasi-split, and when the residue field is algebraically closed, $G$ is residually split. The condition of being residually quasi-split implies that $G$ is quasi-split. As explained in \cite[\S 3.5.2]{tits1979reductive}, when $G$ is residually quasi-split (resp. residually split), the maximal reductive quotient of the special fiber of $\Ggc_x$, denoted as $\ovl{\Ggc_{x}}/R_u\ovl{\Ggc_{x}}$, is quasi-split (resp. split), and this elucidates the use of the adverb \emph{residually}.
\end{remark}

\subsubsection{Levi subgroups}\label{subsubsec: levi subgroups}
Let $x \in A_{\red}(G, S)$, $(Q, L)$ be in the standard position with respect to $(S, J \subset \Delta)$. 

Consider the canonical $N_L(F)$-equivariant map $\pi: A_{\red}(G, S) \to A_{\red}(L, S)$, it extends canonically to $L(F)$-equivariant maps $\Tilde{\pi}: L(F)A_{\red}(G, S) \to B_{\red}(L, F)$ and $\Tilde{\pi}: L(F)A_{\ext}(G, S) \to B_{\ext}(L, F)$, see \cite[Lemma 2.1.4]{landvogt2000some}. On the other hand, there is an $L(F)$-equivariant embedding $\iota_L: B_{\ext}(L, F) \to B_{\ext}(G, F)$. The buildings $B_{\ext}(L, F)$ and $B_{\ext}(G, F)$ share an apartment $A_{\ext}(L, S) = A_{\ext}(G, S)$ with different affine hyperplanes. Given two different $L(F)$-equivariant embeddings $\iota_{L, 1}, \iota_{L, 2}: B_{\ext}(L, F) \to B_{\ext}(G, F)$, there is a unique $v_0 \in V_0(L, S)$ such that $\iota_{L, 1}(x, v) = \iota_{L, 2}(x, v + v_0)$ for any $x \in B_{\red}(L, F)$ and $v \in V_0(L, S)$, see \cite[Proposition 2.1.5]{landvogt2000some}. Note that the construction there is based on \cite[\S 7.6.5, 7.6.6]{bruhat1972groupes} and \cite[\S 4.2.17]{bruhat1984groupes}. 

Let $x_0 \in V_0(G, S)$, $\Tilde{x} =  (x, x_{0}) \in A_{\ext}(G, S)$, then $\pr(\Tilde{x}) = x$. Consider the composition:
\begin{equation*}
  s_L:  L(F)A_{\ext}(G, S) \stackrel{\Tilde{\pi}}{\to} B_{\ext}(L, F) \stackrel{\iota_L}{\to} B_{\ext}(G, F),
\end{equation*}
by adjusting $\iota_L$ with some $v \in V_0(L, S)$, we assume $\Tilde{x} \in A_{\ext}(G, S)$ maps to itself under $s_L$. Let $x_L = \pi(x)$, $\Tilde{x}_L := \Tilde{\pi}(\Tilde{x}) = (x_L, x_{L, 0}) \in B_{\ext}(L, F)$ for some $x_{L, 0} \in V_0(L, S)$, $\Tilde{x} = \iota_L(\Tilde{x}_L) \in B_{\ext}(G, F)$. Since $\hat{P}_{x_L}^1 = \hat{P}_{\Tilde{x}_L}$, $\hat{P}_{x}^1 = \hat{P}_{\Tilde{x}}$, and $\iota_L: B_{\ext}(L, F) \to B_{\ext}(G, F)$ is $L(F)$-equivariant,

\begin{corollary}\label{corollary: P^1_x cap L}
    $\hat{P}_{x_L}^1 = \hat{P}_{x}^1 \cap L(F)$.
\end{corollary}

Let $f: \wdh{\Phi} \to \wdt{\R}$ be a concave function, $f_L = f|_{\wdh{\Phi}}: \wdh{\Phi}_L \to \wdt{\R}$ is also a concave function.

\begin{proposition}\label{proposition: Levi of parahoric is parahoric}
    Let $K_{x, (r)}$ and $K_{x_L, (r)}$ be subgroups of $G(F)$ and of $L(F)$ of the same type $(\ast)$ with respect to $x \in  A_{\red}(G, S)$ and $x_L \in  A_{\red}(L, S)$ respectively, and $r \geq 0$, then $K_{x, (r)} \cap L(F) = \pi( K_{x, (r)} \cap Q(F)) = K_{x_L, (r)}$.
\end{proposition}
\begin{proof}
    Due to Proposition \ref{proposition: MP-group satisfies two conditions} and \ref{lem: parahoric subgroup satisfies the condition}, $K_{x, (r)}$ and $K_{x_L, (r)}$ satisfy the condition \ref{eq: U^-U^+N}. Proposition \ref{proposition: henniart factors} implies that $K_{x, (r)} \cap L(F) = \pi( K_{x, (r)} \cap Q(F))$. 
    
    When $K_{x, (r)} = \hat{P}_x^1$ and $K_{x_L, (r)} = \hat{P}_{x_L}^1$, then $K_{x, (r)} \cap L(F) = K_{x_L, (r)}$ follows from Corollary \ref{corollary: P^1_x cap L}. In all the other cases, the condition \ref{eq: U^+U^-U^-Z} is also valid. Due to Corollary \ref{cor: when K cap L is K_L}, it suffices to show:
    \begin{equation}
        \forall\alpha \in \Phi_{J, \red},\ K_{x_L, (r)} \cap U_{\alpha}(F) = K_{x, (r)} \cap  U_{\alpha}(F), \quad Z(F) \cap K_{x_L, (r)} = Z(F) \cap K_{x, (r)}.
    \end{equation}
    
    Firstly, we consider the unipotent groups $U_{\alpha}$. Let $x_1, x_2 \in A_{\red}(G, S)$ such that $\pi(x_1) = \pi(x_2)$ under $\pi: A_{\red}(G, S) \to A_{\red}(L, S)$, then $x_1 -  x_2 \in V_0':= X_*(S_J/Z_{c, \spl}) \otimes \R$, where $S_J = (\cap_{\alpha \in J}\Ker \alpha)^{\circ}$. Thus, for any $\alpha \in \Phi_J$, $\alpha(x_1 - x_2) = 0$. Therefore, given any $\alpha \in \Phi_J$, $\alpha(x_1 - x_2) = \alpha(\pi(x_1) - \pi(x_2))$. In particular, $f_x(\alpha) = -\alpha(x - x_o) = -\alpha(x_L - \pi(x_o)) = f_{x_L}(\alpha)$, $K_{x_L, (r)} \cap U_{\alpha}(F) = K_{x, (r)} \cap  U_{\alpha}(F)$.

    Secondly, we consider the minimal levi subgroup $Z$, the centralizer of maximal $F$-split torus $S$. When $K_{x, (r)} = P^1_x$, $K_{x_L, (r)} = P^1_{x_L}$, since the closure of $Z$ in $\Gg_x$ (resp. $\Gg_{x_L}$, here $x_L \in A_{\red}(L, S)$) is $\ZZ$, due to \cite[\S 5.2.4]{bruhat1984groupes},
    \begin{equation}\label{eq: Z(F)^1, intersection}
        Z(F)^1 = P^1_x \cap Z(F) = P^1_{x_L} \cap Z(F).
    \end{equation} 
    When $K_{x, (r)} = P^0_x$, $K_{x_L, (r)} = P^0_{x_L}$, it suffices to show
    \begin{equation}\label{eq: Z(F)_0, intersection}
        Z(F)_0 = P^0_x \cap Z(F) = P^0_{x_L} \cap Z(F).
    \end{equation}
    Recall that $Z'$ is the centralizer of a maximal $\bF$-torus $T \supset S$ defined over $F$, due to \cite[Lemma 5]{haines2008parahoric}, it satisfies $Z'(\bF)_0 = P^0_{\iota(x)} \cap Z'(\bF)$, where $\iota: \Bui_{\ext}(G, F) \to \Bui_{\ext}(G, \bF)$ is the canonical embedding. In particular, we have $Z(\bF)_0 = P^0_{\iota(x)} \cap Z(\bF)$, by applying the previous paragraphs to the case when $L$ is replaced with $Z$ and $F$ is replaced with with $\bF$. Finally, replace $G$ by $L$, and take the $\Sigma_0$-invariants, we get equations \ref{eq: Z(F)_0, intersection}.

    Finally, when $K_{x, (r)} = G(F)_{x, r}$ is a Moy-Prasad subgroup, it suffices to show $G(\bF)_{x, f} \cap Z(\bF) = L(\bF)_{x_L, f} \cap Z(\bF)$. We verify this over $\bF$ Utilizing the factorizations \ref{eq: U^+U^-U^-Z} over $\bF$. The identity between the unipotent parts follows from the above discussion, and the identity between the maximal torus part $U_0 = Z'$ follows from the definition $f(0) = f_L(0)$.
\end{proof}
\begin{corollary}\label{cor: good intersection with L}
    Let $x \in \Bui_{\ext}(G, F)$, then $K_{x, (r)} \cap L(F) = \pi(K_{x, (r)} \cap Q(F))$ if $x \in \Bui_{\ext}(L, F)$.
\end{corollary}
\begin{proof}
    If $x \in \Bui_{\ext}(L, F)$, then we can find a maximal $F$-split torus $S' \subset L$ such that $x \in A_{\ext}(L, S')$, we also regard $x \in A_{\ext}(G, S')$. Since all maximal $F$-split torus are conjugated, pick $g \in L(F)$ such that $gS'g^{-1} = S \subset L$ is at standard position. We rewrite $x_L := x \in \Bui_{\ext}(L, F)$, and $x \in \Bui_{\ext}(G, F)$. Apply Proposition \ref{proposition: Levi of parahoric is parahoric}, we have 
    \[ K_{gx, (r)} \cap L(F) = \pi(K_{gx, (r)} \cap Q(F)) = K_{gx_L, (r)}.\] 
    Since $gK_{x_L, (r)}g^{-1} = K_{gx_L, (r)}$, $gK_{x, (r)}g^{-1} = K_{gx, (r)}$, and $g \in L(F)$, we have 
    \[ K_{x, (r)} \cap L(F) = \pi(K_{x, (r)} \cap Q(F)) = K_{x_L, (r)}.\]
    
\end{proof}
\begin{remark}
    Let $\FF \subset A_{\red}(G, S)$ be a facet, there exists a unique facet $\FF_L \subset A_{\red}(L, S)$ containing $\FF$. \cite[Lemma 4.1.1]{haines2010satake} proved that $\Ggc_{\FF}(\OO_F) \cap L(F) = \Ggc_{\FF_L}(\OO_F)$. The proof is based on some detailed calculations involving Iwahori-Weyl groups and Kottwitz maps and is carried out over $\bF$. Proposition \ref{proposition: Levi of parahoric is parahoric} and Corollary \ref{cor: good intersection with L} offer an alternative perspective that is also applicable at deeper levels.
\end{remark}

\begin{lemma}\label{lemma: levi of special points}
    Let $x \in A_{\red}(G, S)$. When $x$ is a special (resp. absolutely special, hyperspecial, superspecial) point, then $x_L \in A_{\red}(L, S)$ is also special (resp. absolutely special, hyperspecial, superspecial) point.
\end{lemma}
\begin{proof}
    It suffices to prove the case when $x$ is special over $F$, since all the other cases are on various field extensions $F_1|F$, we could apply same arguments over $F_1$. In the hyperspecial case, $L$ is unramified since $G$ is unramified, see Lemma \ref{lemma: splitting fields are keepting the same}. Let $L = \lrbracket{Z, U_{\Phi_J}}$ be as before, it follows from the definition that the affine roots $\Phi_{\aff, L}$ are those elements in $\Phi_{\aff, G}$ whose linear parts are in $\Phi_J$, thus $x \in A_{\red}(G, S)$ is special implies $x \in A_{\red}(L, S)$ follows from the definition.
\end{proof}
\begin{lemma}\label{lemma: residually split, boundary}
    When $G$ is residually split, then $L$ is residuelly split.
\end{lemma}
\begin{proof}

    Let $\alc' \subset A_{\red}(G, T)$ be a $\Sigma_0$-invariant alcove containing $\alc$, and $\alc_L' \subset A_{\red}(L, T)$ be the unique alcove containing $\alc'$, then $\alc_L \subset \alc_L'$, where $\alc_L \subset A_{\red}(L, S)$ is the unique alcove containing $\alc$. Fix a basis $\Delta_{\aff}'$ of $\Phi_{\aff}'$, since the walls bounding $\alc_L'$ form a subset of the walls bounding $\alc'$ (with those walls whose associated affine roots not in $\Phi_{\aff, L}'$ removed), we get a basis $\Delta_{\aff, L}'$ of $\Phi_{\aff, L}'$ which is a subset of $\Delta_{\aff}'$. Note that the affine roots $\Phi_{\aff}$ are those nonconstant restrictions $\alpha'|_{A_{\red}(G, S)}$ with $\alpha'\in\Phi_{\aff}'$. Let $\Delta_{\aff, \mathrm{dist}}' := \Delta_{\aff}' \cap \Phi_{\aff}$, it is $\Sigma_0$-invariant and is the \emph{distinguished set} introduced in \cite[\S 1.11]{tits1979reductive}, $G$ is residually split over $F$ if and only if $\Sigma_0$ acts trivially on $\Delta_{\aff, \mathrm{dist}}'$, see \emph{loc. cit.} The key point is that, $\Delta_{\aff, L, \mathrm{dist}}' := \Delta_{\aff, L}' \cap \Phi_{\aff, L}$ is a subset of $\Delta_{\aff, \mathrm{dist}}'$, thus $G$ is residually split implies that $L$ is residually split.

\end{proof}

\subsubsection{Normal reductive subgroups}\label{subsub: normal subgroups}

In this subsection, we consider a normal subgroup $G_h \subset L$, and using the notations introduced in subsection \ref{subsec: general facts}.

\begin{lemma}{{\cite[\S 2]{landvogt2000some}}}\leavevmode
\begin{enumerate}
        \item There is a canonical $G_h(\bF) \times G_l'(\bF)$-$\Sigma_0$-equivariant bijection 
        \[ B_{\ext}(G_h, \bF) \times B_{\ext}(G_l', \bF)  \to B_{\ext}(L, \bF). \]
    \item There is a canonical $L(\bF)$-$\Sigma_0$-equivariant map $\iota_{\pi}: B_{\ext}(L, \bF) \to B_{\ext}(G_l, \bF)$. Moreover, there is a $G_l'(\bF)$-$\Sigma_0$-equivariant bijection $B_{\ext}(G_l', \bF) \to B_{\ext}(G_l, \bF)$ which provides a section to $\iota_{\pi}$.
\end{enumerate}
\end{lemma}
Therefore, there is a canonical $G_h(\bF) \times G_l'(\bF)$-$\Sigma_0$-equivariant bijection 
\begin{equation}\label{eq: building decomposition}
    B_{\ext}(G_h, \bF) \times B_{\ext}(G_l, \bF)  \to B_{\ext}(L, \bF). 
\end{equation}
The same statements hold with $B_{\ext}$ replaced with $B_{\red}$ for precisely the same reasons. Lemma \ref{lemma: good sub and quotient for centralizers and normalizers} together with the facts $\Phi_{\aff} \cong \Phi_{h, \aff} \cup \Phi_{l', \aff}$, $\Phi_{l', \aff} \cong \Phi_{l, \aff}$ show that the bijection \ref{eq: building decomposition} keeps poly-simplicial structures: Given a facet $\FF \subset B_{\red}(L, \bF)$, then there exists a unique pair of facets $(\FF_h, \FF_l) \subset B_{\red}(G_h, \bF) \times B_{\red}(G_l, \bF)$ such that $\FF_h \times \FF_l = \FF$. $\FF$ is $\Sigma_0$-equivariant if and only if both $\FF_h$ and $\FF_l$ are $\Sigma_0$-equivariant.

Let $\lrbracket{x_L} \in B_{\ext}(L, F)$ be any point, under the bijection \ref{eq: building decomposition}, $\lrbracket{x_L} = \lrbracket{x_h} \times \lrbracket{x_l} \in B_{\ext}(G_h, F) \times B_{\ext}(G_l, F)$. It follows from the above discussion that
\begin{corollary}\label{corollary: sub and quotient of special or residually split is of the same type}\leavevmode
\begin{enumerate}
    \item $x_L$ is a special (resp. absolutely special, superspecial, hyperspecial) point if and only if both $x_h$ and $x_l$ are special (resp. absolutely special, superspecial, hyperspecial) point.
    \item $L$ is residually split if and only if both $G_h$ and $G_l'$ (or equivalently, $G_l$) are residually split.
\end{enumerate}
\end{corollary}

The exact sequence $1 \to G_h \to L \to G_l \to 1$ induces an exact sequence $1 \to \pi_1(G_h) \to \pi_1(L) \to \pi_1(G_l) \to 1$ of algebraic fundamental groups (which we recall in subsection \ref{subsec: general facts}). Taking $I$-coinvariant spaces of these Galois modules, we have an exact sequence $\pi_1(G_h)_I \to \pi_1(L)_I \to \pi_1(G_l)_I \to 1$. In general, $\pi_1(G_h)_I \to \pi_1(L)_I$ might not be an injection, thus the following inclustion might be strict:
\begin{equation}
    G_h(\bF)_0 \subset L(\bF)_0 \cap G_h(\bF)
\end{equation}
On the other hand, recall that $L(\bF)^1$ (resp. $G_h(\bF)^1$) is the kernel of $\nu_L: L(\bF) \to X_*(Z_{c, \spl})\otimes \Q$ (resp. $\nu_{G_h}: G_h(\bF) \to X_*(Z_{c, h, \spl})\otimes \Q$). Since $X_*(Z_{c, h, \spl})\otimes \Q \to X_*(Z_{c, \spl})\otimes \Q$ is an injection, we have:
\begin{equation}
  G_h(\bF)^1 =  L(\bF)^1 \cap G_h(\bF).
\end{equation}
These implies that
\begin{equation}\label{equation: intersects with L1 is Gh1}
    G_h(F)_0 \subset L(F)_0 \cap G_h(F),\quad G_h(F)^1 = L(F)^1 \cap G_h(F).
\end{equation}

Apply these arguments to $1 \to Z_h \to Z_L \to Z_l \to 1$, where $Z_h$, $Z_L$, $Z_l$ are centralizers of maximal $F$-split tori $S_h$, $S_L$, $S_l$ of $G_h$, $L$, $G_l$ respectively, here we pick $S_L$ such that $x_L \in A_{\red}(L, S_L)$, thus $x_h \in A_{\red}(G_h, S_h)$, $x_l \in A_{\red}(G_l, S_l)$. We have

\begin{equation}\label{equation: intersects with Z_L1 is Z_h1}
    Z_L(F)_0 \cap Z_h(F) \supset Z_h(F)_0, \quad Z_L(F)^1 \cap Z_h(F) = Z_h(F)^1.
\end{equation}

On the other side, consider the surjection $L \to G_l$ as well as $Z_L \to Z_l$, they are surjective since $H^1(L, G_h) = H^1(L, Z_h) = \lrbracket{0}$. Since both $\pi_1(L)_I \to \pi_1(G_l)_I$ and $X_*(Z_{c, \spl})\otimes \Q \to X_*(Z_{c, l, \spl}) \otimes \Q$ are surjective, and $\nu_{G_h}$, $\kappa_{G_h}$ are surjective thus we have
\begin{equation}
    \pi_l(L(\bF)^1) = G_l(\bF)^1,\quad \pi_l(L(\bF)_0) = G_l(\bF)_0.
\end{equation}
Similarly,
\begin{equation}\label{equation: quotient of Z_L1}
    \pi_l(Z_L(\bF)^1) = Z_l(\bF)^1,\quad \pi_l(Z_L(\bF)_0) = Z_l(\bF)_0.
\end{equation}

\begin{proposition}\label{proposition: sub and quotient of BT group is BT group}\leavevmode
    \begin{enumerate}
        \item When $K_{x_L, (r)}$ is of type $(\ast)$ as defined in \ref{def: type ast groups, in an appartment}, then $K_{x_L, (r)} \cap G_h(F) = K_{x_h, (r)}$ (except for the case when $K_{x_L, (r)}$ is parahoric), and $\pi_l(K_{x_L, (r)}) \subset K_{x_l, (r)}$.
        \item $\Ggc_{x_L}(\OO_F) \cap G_h(F) = \Ggc_{x_h}(\OO_F)$ if and only if $Z_h(F)_0 = Z_L(F)_0 \cap Z_h(F)$. In particular, $\Ggc_{x_L}(\OO_F) \cap G_h(F) = \Ggc_{x_h}(\OO_F)$ when $G$ is unramified or when $L \to G_l$ has a section.
        \item Let $\KK_{x_L, (r)}$ (resp. $\KK_{x_l, (r)}$) be the $\OO_{\bF}$-points of the associated smooth models $\Gg_{x_L}$ (resp. $\Gg_{x_l}$) or $\Ggc_{x_L}$ (resp. $\Ggc_{x_l}$), then $\pi_l(\KK_{x_L, (r)}) = \KK_{x_l, (r)}$.
    \end{enumerate}
\end{proposition}
\begin{proof}
    We show these properties by examining these four types \ref{def: type ast groups, in an appartment} separately.
    \begin{enumerate}
        \item When $K_{x_L, (r)} = \Ggh_{x_L}(\OO_F) = \hat{P}_{\pr(x_L)}^1 = \hat{P}_{x_L}$, $K_{x_h, (r)} = \Ggh_{x_h}(\OO_F) = \hat{P}_{\pr(x_h)}^1 = \hat{P}_{x_h}$, then $\hat{P}_{x_L} \cap G_h(F) = \hat{P}_{x_h}$ since there exists a $G_h(\bF) \times G_l'(\bF)$-equivariant bijection \ref{eq: building decomposition} and since $G_h$-action on $B_{\ext}(G_l, F)$ is trivial. Similarly, $\pi_l(\hat{P}_{x_L}) \subset \hat{P}_{x_l}$ since $B_{\ext}(L, F) \to B_{\ext}(G_l, F)$ is $L(F)$-equivariant.
        \item When $K_{x_L, (r)} = \Gg_{x_L}(\OO_F) = P_{x_L}^1$, $K_{x_h, (r)} = \Gg_{x_h}(\OO_F) = P_{x_h}^1$ (resp. $K_{x_l, (r)} = \Gg_{x_l}(\OO_F) = P_{x_l}^1$), and $\KK_{x_L, (r)} = \Gg_{x_L}(\OO_{\bF}) = P_{\iota(x_L)}^1$, $\KK_{x_l, (r)} = \Gg_{x_l}(\OO_{\bF}) = P_{\iota(x_l)}^1$, they satisfy conditions \ref{eq: U^-U^+N} and \ref{eq: U^+U^-U^-Z} due to Lemma \ref{lem: parahoric subgroup satisfies the condition}, we could apply Corollary \ref{cor: when K_L cap G_h is K_h}: On one hand, let $\alpha \in \Phi_h$ (resp. $\alpha \in \Phi_l' \cong \Phi_l$), since $\Phi_L = \Phi_h \sqcup \Phi_l'$ is an orthogonal decomposition, then $\alpha(\lrbracket{x_{h, 1}, x_{l, 1}} - \lrbracket{x_{h, 2}, x_{l, 2}}) = \alpha(x_{h, 1} - x_{h, 2})$ (resp. $= \alpha(x_{l, 1} - x_{l, 2})$, thus $f_{x_L}(\alpha) = f_{x_h}(\alpha)$ (resp. $= f_{x_l}(\alpha)$), $K_{x_L, (r)} \cap U_{\alpha}(F) = K_{x_h, (r)} \cap U_{\alpha}(F)$ (resp. $\pi_l(K_{x_L, (r)}\cap U_{\alpha}(F)) = K_{x_h, (r)} \cap U_{\alpha}(F)$). On the other hand, due to identifications \ref{eq: Z(F)^1, intersection}, $P_{x_L}^1 \cap Z_L(F) = Z_L(F)^1$, $P_{x_h}^1 \cap Z_h(F) = Z_h(F)^1$ (resp. $P_{x_l}^1 \cap Z_l(F) = Z_l(F)^1$), and $P_{\iota(x_L)}^1 \cap Z_L(\bF) = Z_L(\bF)^1$, $P_{\iota(x_l)}^1 \cap Z_l(\bF) = Z_l(\bF)^1$, thus $P_{x_L}^1 \cap Z_h(F) = P_{x_h}^1 \cap Z_h(F)$ (resp. $\pi_l(P_{x_L}^1 \cap Z_L(F)) \subset P_{x_l}^1 \cap Z_l(F)$) and $\pi_l(P_{\iota(x_L)}^1 \cap Z_L(\bF)) = P_{\iota(x_l)}^1 \cap Z_l(\bF)$ due to \ref{equation: intersects with Z_L1 is Z_h1} and \ref{equation: quotient of Z_L1}.
        \item When $K_{x_L, (r)} = \Ggc_{x_L}(\OO_F) = P_{x_L}^0$, $K_{x_h, (r)} = \Ggc_{x_h}(\OO_F) = P_{x_h}^0$ (resp. $K_{x_l, (r)} = \Ggc_{x_l}(\OO_F) = P_{x_l}^0$), and $\KK_{x_L, (r)} = \Ggc_{x_L}(\OO_{\bF}) = P_{\iota(x_L)}^0$, $\KK_{x_L, (r)} = \Ggc_{x_L}(\OO_{\bF}) = P_{\iota(x_l)}^0$. The calculation is similar to the previous case, and use the identifications \ref{eq: Z(F)_0, intersection},\. Let us explain the part $(2)$: due to Corollary \ref{cor: when K_L cap G_h is K_h}, $P_{x_L}^0 \cap G_h(F) = P_{x_h}^0$ if and only if $P_{x_L}^0 \cap Z_h(F) = P_{x_h}^0 \cap Z_h(F)$ if and only if $Z_L(F)_0 \cap Z_h(F) = Z_h(F)_0$. If $G$ is unramified, then $Z_L, Z_h, Z_l$ are all unramified due to Lemma \ref{lemma: splitting fields are keepting the same}, then taking $I$-coinvariant spaces of the exact sequence $1 \to \pi_1(Z_h) \to \pi_1(Z_L) \to \pi_1(Z_l) \to 1$ does not change anything, $\pi_1(Z_h)_I \to \pi_1(Z_L)_I$ is an injection. If $L \to G_l$ has a section, then $\pi_1(Z_h)_I \to \pi_1(Z_L)_I$ is also an injection. In both cases, $Z_L(F)_0 \cap Z_h(F) = Z_h(F)_0$.
        \item When $K_{x_L, (r)} = L(F)_{x_L, r}$, $K_{x_h, (r)} = G_h(F)_{x_h, r}$ (resp. $K_{x_l, (r)} = G_l(F)_{x_h, r}$) with $r>0$. The calculation is similar to the previous case, but over $\bF$. For the $G_h$-part, it suffices to show, for any torus $Z'_h \subset Z'_L$, $Z'_h(\bF)^{\MPF}_r = Z'_L(\bF)^{\MPF}_r \cap Z'_h(\bF)$ for all $r > 0$. For the $G_l$-part, it suffices to show, for any torus $Z'_L \to Z_l'$, $\pi_l(Z'_L(\bF)^{\MPF}_r) \subset Z'_l(\bF)^{\MPF}_r$ for all $r > 0$. These follow from the definition in Remark \ref{remark: MP independent of embedding}.
    \end{enumerate}
\end{proof}
\begin{remark}
    We can replace $G(F)_{x, r}$ with any $G(F)_{x, f}$ with $f$ being a concave function. Let $f_L: \wdh{\Phi}_L \to \wdt{\R}$ be a concave function, then $f_h := f_L|_{\wdh{\Phi}_h}: \wdh{\Phi}_h \to \wdt{\R}$ and $f_l := f_L|_{\wdh{\Phi}_l'}: \wdh{\Phi}_l'\cong \wdh{\Phi}_l \to \wdt{\R}$ are also concave functions. When we replace $G(F)_{x, r}$ with $G(F)_{x, f}$, Proposition \ref{proposition: Levi of parahoric is parahoric} holds as well, and Proposition \ref{proposition: sub and quotient of BT group is BT group} holds if and only if $Z'_h(\bF)_{f_L(0)} = Z'_L(\bF)_{f_L(0)} \cap Z'_h(\bF)$.
\end{remark}
\begin{example}\label{example: norm torus}
    Let $F|F_0$ be a quadratic ramified extension, $I = \Gal(F|F_0)$, let $1 \to G_h \to L \to G_l \to 1$ be
    \[ 1 \to \Res_{F|F_0}^{(1)}\Gm \to \Res_{F|F_0}\Gm \to \Gm \to 1, \]
    then $\pi_1(G_h)_I \to \pi_1(L)_I \to \pi_1(G_l)_I \to 1$ is 
    \[ \Z/2\Z \stackrel{0}{\to} \Z \to \Z \to 0. \]
    In particular, $L(F)_0 \cap G_h(F) = G_h(F)^1$, and $G_h(F)^1$ strictly contains $G_h(F)_0$.
\end{example}

\begin{corollary}\label{corollary: if L splits and K special type then K splits}
    If $K:= K_{x_L, (r)} \subset L(F)$ is of type $(\ast)$ as defined in \ref{def: type ast groups, in an appartment}, and if $1 \to G_h \to L \to G_l \to 1$ has a section inducing $L \cong G_h \times G_l$, then $K = (K \cap G_h(F))(K \cap G_l(F))$.
\end{corollary}
\begin{proof}
   Note that $L \cong G_h \times G_l$ implies $Z_L \cong Z_h \times Z_l$, $Z_h(F)_0 = Z_L(F)_0 \cap Z_h(F)$. Consider the natural embedding $K_{x_L, (r)} \cap G_l(F) \hookrightarrow \pi_l(K_{x_L, (r)})$. The left hand side is $K_{x_l, (r)}$, the right hand side is contained in $K_{x_l, (r)}$, thus $K \cap G_l(F) = \pi_l(K)$. Given $k \in K$, let $k_1 \in K \cap G_l(F)$ be a lifting of $\pi_l(k) \in \pi_l(K)$, then $kk_1^{-1} \in K \cap G_h(F)$, thus $K = (K \cap G_h(F))(K \cap G_l(F))$.
\end{proof}

\subsubsection{Isogenies}

Let us say a few words about isogenies. This would not be used in the rest of the paper.

\begin{proposition}\label{prop: isogeny of tori, parahoric}
    Let $\varphi: T_1 \to T_2$ be an isogeny of tori, then $T_1(\bF)_1 = \varphi^{-1}(T_2(\bF)_1)$.
\end{proposition}
\begin{proof}
    By functoriality, $\varphi(T_1(\bF)_1) \subset  T_2(\bF)_1$. $X^*(T_2) \to X^*(T_1)$ is injective with finite cokernel, then $X^*(T_2)^I \to X^*(T_1)^I$ is also injective with finite cokernel, then 
    \[ \Hom(X^*(T_1)^I, \Z) \to \Hom(X^*(T_2)^I, \Z) \]
    is injective. Therefore, $T_1(\bF)_1 = \varphi^{-1}(T_2(\bF)_1)$.
\end{proof}

\begin{remark}\label{example: isogeny, non-parahoric}
    We keep notations from last example \ref{example: norm torus}. Let $T_1 = \Res^{(1)}_{F|F_0} \Gm \times (\Gm)_{F_0}$, $T_2 = \Res_{F|F_0}\Gm$, then the natural morphism $\varphi: T_1 \to T_2$ is an isogeny, $\pi_1(T_1) \to \pi_1(T_2)$ is injective with finite cokernel, but $\pi_1(T_1)_I \to \pi_1(T_2)_I$ is no longer injective. In other words, $\varphi^{-1}(T_2(\bF)_0) \subsetneq T_1(\bF)_0$.
\end{remark}

\begin{proposition}
    Let $\varphi: G_1 \to G_2$ be an isogeny between reductive groups, then $\varphi^{-1}(G_2(\bF)_1) = G_1(\bF)_1$.
\end{proposition}
\begin{proof}
    One can use another definition of $\nu_G$ (\ref{equation: nuG another definition}),
    \[ \nu_{G_i}:\quad G_i(\bF) \to \Hom_{\Z}(X^*(G), \Z) \otimes \R = X_*(Z_{i, c, \spl}) \otimes \R, \]
    it suffices to show  $X_*(Z_{1, c, \spl}) \otimes \R \to X_*(Z_{2, c, \spl}) \otimes \R$ is injective. This follows from that $Z_{1, c, \spl} \to Z_{2, c, \spl}$ is an isogeny.
\end{proof}

\begin{remark}
    Let $G_1 \to G_2 \to \cdots \to G_n$ be a chain of isogenies among semisimple groups, $G_1$ is simply-connected, $G_n$ is adjoint. It is worth noting that $Z_1$ and $Z_n$ (the centralizer of a maximal $\bF$-split $F$-torus) in $G_1$ and $G_n$ respectively are induced (\cite[\S 4.4.16]{bruhat1984groupes}), thus $X_*(Z_1)_I$ is torsion-free and $X_*(Z_1)_I \to X_*(Z_n)_I$ is injective. Then $G_i(\bF)^1 = G_i(\bF)$ for all $i$, and the index of $G_i(\bF)^1/G_i(\bF)_0$ is non-decreasing when $i$ increase (starting from $G_1(\bF)^1/G_1(\bF)_0 = \lrbracket{0}$).
\end{remark}

\subsubsection{Summary}

\begin{proposition}\label{prop: main prop for section BT theory}
    Let $x \in A_{\red}(G, S)$, $(Q, L)$ be a standard pair with respect to $(S, J \subset \Delta)$, $G_h \subset L$ be a normal subgroup, $G_l = L/G_h$ be its quotient, $\pi: Q \to L$, $\pi_l: L \to G_l$ be projections to their Levi quotients, $x_L$ be the canonical image of $x$ under the projection $A_{\red}(G, S) \to A_{\red}(L, S)$, and $x' = \lrbracket{x'_h, x'_l} \in B_{\red}(G_h, F) \times B_{\red}(G_l, F)$ be the image of a given point $x' \in B_{\red}(L, F)$ under the bijection \ref{eq: building decomposition}. Let $r \geq 0$, $K_{x, (r)} \subset G(F)$ be a group of type $(\ast)$ (see \ref{definition: type ast} and \ref{def: type ast groups, in a building}), let $g \in G(F)$, decompose $g = ulnk$, $l \in L(F)$, $u\in R_uQ(F)$, $n \in N(F)$, $k \in P_{x}^{0}$ (see Remark \ref{rmk: choose Q standard}), then

    \begin{enumerate}
        \item $\pi(Q(F) \cap gK_{x, (r)}g^{-1}) = K_{x_{g, L}, (r)}$, where $x_{g, L} = l\cdot (\nu_N(n)(x))_L$ (see \ref{def: type ast groups, in an appartment}).
        \item $K_{x_{g, L}, (r)} \cap G_h(F) \supset K_{(x_{g, L})_h, (r)}$ and the equality holds at least in following cases:
        \begin{enumerate}
            \item if $K_{x, (r)}$ is not a parahoric subgroup.
            \item if $K_{x, (r)}$ is a parahoric subgroup, and if $G$ is unramified.
            \item if $K_{x, (r)}$ is a parahoric subgroup, and if there exists a section $G_l \to L$ of $\pi_l$.
        \end{enumerate}
        \item $\pi_l(K_{x_{g, L}, (r)}) \subset K_{(x_{g, L})_l, (r)}$.
        \item When $x$ is a special (resp. absolutely special, superspecial, hyperspecial) point, then $x_{g, L}$, $(x_{g, L})_h$, $(x_{g, L})_l$ are special (resp. absolutely special, superspecial, hyperspecial) points.
        \end{enumerate}
\end{proposition}
\subsubsection{Adelic version}
\begin{lemma}{{\cite[\S 3.9.1]{tits1979reductive}}}\label{lemma: global BT groups}
   Let $G$ be a $\Q$-reductive group, $K \subset G(\A_f)$ be an open compact subgroup, then for all but finitely many primes $p$, $K_p \subset G(\Q_p)$ is a hyperspecial subgroup. That is to say, $K_p = \GGc_{x_p}(\Z_p) = \GG_{x_p}(\Z_p)$ for some hyperspecial point $x_p \in B_{\ext}(G, \Q_p)$. Moreover, let $G \hookrightarrow \GL_n := \GL(V)$ be a closed embedding, $V_{\Z} \subset V$ be a full lattice, $\mathcal{GL}_n:=\GL(V_{\Z})$ be the integral model of $\GL_n$ over $\Spec \Z$. Let $\GG$ be the schematic closure of $G$ in $\mathcal{GL}_n$, then for all but finitely many primes $p$, $\GG \otimes_{\Z} \Z_p = \GG_{x_p}$, given $x_p$ as above.
\end{lemma}

\begin{definition-proposition}\label{def: ast-type group, global}
    A subgroup $K \subset G(\A_f)$ is called a \emph{($\ast$)-type subgroup} if for every prime $p$ its $p$-part $K_p \subset G(\Q_p)$ is a group type of $(\ast)$ as in Definition \ref{def: type ast groups, in a building}, and for all but finitely many primes, $K_p$ is a hyperspecial subgroup. Due to Lemma \ref{lemma: global BT groups}, 
    $(\ast)$-type subgroups form a basis of open compact neighbourhood at identity.
\end{definition-proposition}

\subsubsection{Group schemes at boundary}\label{subsec: group schemes at boundary}

Due to Proposition \ref{proposition: sub and quotient of BT group is BT group}, $\Ggc_{x_L}(\OO_F) \cap G_h(F) \supset \Ggc_{x_h}(\OO_F)$ might be strict. Nevertheless, $\Gg_{x_L}(\OO_F) \cap G_h(F) = \Gg_{x_h}(\OO_F)$, thus
\begin{equation}
   \Ggc_{x_h}(\OO_F) \subset \Ggc_{x_L}(\OO_F) \cap G_h(F) \subset \Gg_{x_h}(\OO_F).
\end{equation}

Let $\Ggf_{x_h, 1}$ be the schematic closure of $G_h$ in $\Ggc_{x_L}$, it is the unique model of $G_h$ over $\OO_F$ which represents the following functor on the category of flat $\OO_F$-algebras (see \cite[\S 2.2]{yu2015smooth}):
\[  R \mapsto G_h(R[1/p]) \cap \Ggc_{x_L}(R). \]
In particular, $\Ggf_{x_h, 1}(\OO_{\bF}) = \Ggc_{x_L}(\OO_{\bF}) \cap G_h(\bF)$, $\Ggf_{x_h, 1}(\OO_{F}) = \Ggc_{x_L}(\OO_{F}) \cap G_h(F)$. Since $\Ggf_{x_h, 1}$ is affine by construction, let $\Ggf_{x_h}$ be the smoothening of $\Ggf_{x_h, 1}$, $\Ggf_{x_h}$ is the unique smooth model of $G_h$ over $\OO_F$ such that $\Ggf_{x_h}(\OO_{\bF}) = \Ggf_{x_h, 1}(\OO_{\bF}) $, $\Ggf_{x_h}(\OO_{F}) = \Ggf_{x_h, 1}(\OO_{F})$, see \cite[\S 2.5]{yu2015smooth}.

Therefore we have a series of smooth models of $G_h$ over $\OO_F$:
\begin{equation}\label{eq: G x_h different levels}
   \Ggf_{x_h} \to \Gg_{x_h}, \quad \Ggc_{x_h} = (\Ggf_{x_h})^{\circ} = (\Gg_{x_h})^{\circ}.
\end{equation}

When we work with various strata in the Kisin-Pappas integral models, we usually need to assume $\Ggh_x = \Ggc_x$. The assumption $\Ggh_x = \Ggc_x$ is equivalent to $\Ggh_x = \Gg_x$ and $\Gg_x = \Ggc_x$.

By checking the $\OO_{\bF}$-points of smooth models, it follows from Proposition \ref{prop: main prop for section BT theory} that $\Ggh_{x_{g, L}} = \Gg_{x_{g, L}}$ and $\Gg_{x_{g, L}} = \Ggc_{x_{g, L}}$, and we furthur have $\Ggh_{x_{g, h}} = \Gg_{x_{g, h}}$ as well as $\Gg_{x_{g, l}} = \Ggc_{x_{g, l}}$. However, it might not be true that $\Gg_{x_{g, h}} = \Ggc_{x_{g, h}}$ and $\Ggh_{x_{g, l}} = \Gg_{x_{g, l}}$. Nevertheless, due to Proposition \ref{prop: main prop for section BT theory}, 
\begin{corollary}
   Assume $\Ggh_x = \Ggc_x$, we have $\Ggh_{x_{g, h}} = \Ggc_{x_{g, h}}$ at least in one of the following cases:
\begin{enumerate}
    \item $G$ is unramified.
    \item $1 \to G_h \to L \to G_l \to 1$ has a splitting.
\end{enumerate}
\end{corollary}

\subsection{Extending propositions to the boundary}\label{subsec: general facts}

In this subsection, we show that some properties about $G$ are inherited by its Levi subgroups and normal reductive subgroups.

In the setting of Kisin-Pappas models \cite{kisin2018integral}, the assumption is the following:
\begin{definition}\label{general condition}
    We say $(p, G)$ is under the Kisin-Pappas setting if $G$ splits over a tamely ramified extension of $\Q_p$, $p>2$, and $p\nmid |\pi_1(G^{\der})|$.
\end{definition}
In the setting of Kisin-Pappas-Zhou models \cite{kisin2024integral}, the assumption is the following:
\begin{definition}\label{general condition, KPZ}
    We say $(p, G)$ is under the Kisin-Pappas-Zhou setting if $G$ is $R$-smooth (see \cite[Definition 2.4.3]{kisin2024independence}), $p>2$, and $p\nmid |\pi_1(G^{\der})|$.
\end{definition}

\begin{remark}\label{rmk: R smooth}
    In \cite{kisin2024independence} and \cite{kisin2024integral}, the condition that $G$ splits over a tamely ramified extension of $\Q_p$ can be replaced by the weaker condition that $G$ is $R$-smooth (\cite[\S 2.1.4]{kisin2024integral}). In this subsection, we see that the condition \ref{general condition} is inherited by its Levi $L$ and normal subgroup $G_h$ of $L$. However, $G$ being $R$-smooth does not imply $G_h$ is $R$-smooth (although it does implies $L$ is $R$-smooth), since a torus $T$ being $R$-smooth does not imply its subtori are $R$-smooth. Nevertheless, when we are under the setting \ref{general condition}, $G_h$ splits over a tamely ramified extension of $\Q_p$, then $G_h$ is $R$-smooth.
\end{remark}

     \begin{proposition}\label{lemma: splitting fields are keepting the same}
          Let $G$ be a reductive group over $F$ which satisfies the property $(\ast)$, $Q$ be a parabolic $F$-subgroup of $G$. Then any normal reductive $F$-subgroup $G_h \subset L:=Q/R_uQ$ (including $L$ itself) satisfies the property $(\ast)$. Here $(\ast)$ could be one of the followings:
          \begin{enumerate}
              \item Is quasi-split over $F$.
              \item Splits over a fixed field extension $F_1$ of $F$.
              \item Is residually split.
          \end{enumerate}
     \end{proposition}
     \begin{proof}
     \begin{enumerate}
         \item Since $G$ is quasi-split over $F$, let $B$ be a Borel $F$-subgroup of $G$. Since all the Borel $F$-subgroups are $G(F)$-conjugated, we assume $B \subset Q$ such that unipotent radical $R_uQ$ of $Q$ is contained in the unipotent radical $R_uB$ of $B$, we define $B_L = B/R_uQ \subset L$, it is connected and solvable. Since solvability satisfies the \emph{three for the price of two} property: $H$ is solvable if and only if $N$ and $H/N$ are solvable for some (thus for all) normal subgroups $N$ of $H$, $B_L$ is maximal among those connected solvable subgroups of $L:=Q/R_uQ$, itself is a Borel subgroup of $L$, thus $L$ is quasi-split over $F$.
         
         Let $B_h:= (B_L \cap G_h)^{\circ}$, the intersection of closed subschemes commutes with base change, we do everything over $\ovl{F}$. Due to \cite[Proposition 11.14(2)]{borel2012linear}, $B_h$ is a Borel subgroup of $G_h$, and is defined over $F$, $G_h$ is quasi-split over $F$.
         
         \item We replace $F_1$ with $F$, it suffices to prove when $G$ splits over $F$, then $G_h$ splits over $F$. Let $T$ be a maximal $F$-split torus of $G$ ($T$ is also a maximal torus of $G$), all such tori are $G(F)$-conjugated, thus we could assume $T \subset L$, $T \subset L$ is a $F$-split maximal torus of $L$. Consider $T_h:= T \cap G_h$, it is a maximal torus of $G_h$ (see also Lemma \ref{lemma: good sub and quotient for centralizers and normalizers}), and it is split since $X_*(T_h) \hookrightarrow X_*(T)$ is $\Gal(\ovl{F}|F)$-invariant.

         \item See Lemma \ref{lemma: residually split, boundary} and Corollary \ref{corollary: sub and quotient of special or residually split is of the same type}.
           \end{enumerate}
     \end{proof}

     Next, let us recall the construction of the algebraic fundamental group of a reductive group $G$ following \cite{borovoi1998abelian}.
     
     Let $G$ be a reductive group over $F$, $T \subset G_{\ovl{F}}$ be a maximal torus. Let $\rho: \wdt{G} \to G^{\der} \to G$ be the simply connected cover, and $\wdt{T}$ be the preimage of $T$, it is a maximal torus of $\wdt{G}$. The algebraic fundamental group with respect to the pair $(G, T)$ is defined as
     \begin{equation}
         \pi_1(G, T) = X_*(T) / \rho_*(\wdt{T}).
     \end{equation}
     
     Given $\gamma \in \Sigma = \Gal(\ovl{F}/F)$, there exists $g_{\gamma}\in G(\ovl{F})$ such that $g_{\gamma}(\gamma(T))g_{\gamma}^{-1} = T$, $\Sigma$ acts on $X_*(T)$ as
     \begin{equation}\label{eq: twist action Galois}
         \gamma \cdot \mu = \Ad(g_{\gamma}) \circ \gamma(\mu),\quad \gamma \in \Sigma.
     \end{equation}

     If $T' = gTg^{-1}$ for some $g \in G(\ovl{F})$, then $\ad(g)$ induces a canonical $\Sigma$-invariant isomorphism $\pi_1(G, T) \rightiso \pi_1(G, T')$, and \cite[Lemma 1.2]{borovoi1998abelian} showed that this is independent of the choice of $g$. We abbreviate $\pi_1(G, T)$ as $\pi_1(G)$ once we fix a maximal torus. If the torus $T$ is defined over $F$, then the Galois action on $\pi_1(T)=X_*(T)$ is compatible with the Galois action on $\pi_1(G)$, this is also due to \cite[Lemma 1.2]{borovoi1998abelian}.

     \begin{proposition}\label{lemma: p does not devide fund group of Gh}
        Keep the notations as established in Lemma \ref{lemma: splitting fields are keepting the same}. If $p \nmid |\pi_1(G^{\der})|$, then $p\nmid |\pi_1(G_h^{\der})|$.
     \end{proposition}
     \begin{proof}
         This is a problem over algebraic closed fields, we assume $G$ is split. We first work with $L$. $\wdt{L} \to G^{\der}$ factors through $\wdt{G}$ canonically. It suffices to show $\wdt{L} \to \wdt{G}$ is a group embedding. When $\wdt{L} \to \wdt{G}$ is a group embedding, then $\mu_L \to \mu_G$ as well as $\pi_1(L^{\der}) \to \pi_1(G^{\der})$ are embeddings.

         There is a one-to-one correspondence between the set of parabolic subgroups (resp. maximal tori) of $G$ and the set of parabolic subgroups (resp. maximal tori) of $G^{\der}$ given by $Q \mapsto Q\cap G^{\der}$ (resp. $T \mapsto T \cap G^{\der}$). Let $L \subset Q$, one can find a Levi subgroup $\hat{L}$ of $Q \cap G^{\der}$ such that $\hat{L}^{\der} = L^{\der}$. Without lose of generality, we assume $G = G^{\der}$. The central isogeny $\wdt{G} \to G$ maps the set of parabolic subgroups of $\wdt{G}$ to the set of parabolic subgroups of $G$. Let $Q'$ (resp. $L'$) be the pullback of $Q$ (resp. $L$). $L'$ is a Levi subgroup of $\wdt{G}$. It is well-known that the derived group of a Levi subgroup of a simply-connected semisimple group is again simply-connected, for example, see \cite[Proposition 12.14]{malle2011linear}. Thus $L^{\prime\der}$ is simply-connected, $L^{\prime\der} \to L^{\der}$ is a simply-connected covering, we have $\wdt{L} = L^{\prime\der} \hookrightarrow \wdt{G}$.
         
         Now, consider $G_h$. $G_h^{\der}$ is a normal subgroup of $L^{\der}$, let $\ovl{G}_l$ be the quotient $L^{\der}/G_h^{\der}$, one has an exact sequence:
         \[ 1 \to \pi_1(G_h^{\der}) \to \pi_1(L^{\der}) \to \pi_1(\ovl{G}_l) \to 1,  \]
         see \cite[Lemma 1.5]{borovoi1998abelian}. If $p \nmid |\pi_1(L^{\der})|$, then $p\nmid |\pi_1(G_h^{\der})|$.
     \end{proof}
  
  \begin{lemma}\label{lemma: KZ are preserved}
      When $G$ is acceptable \cite[Definition 3.1.2]{kisin2024independence} (resp. essentially tame \cite[Definition 3.1.4]{kisin2024integral}, classical), its Levi subgroups $L$ and normal subgroup $G_h$ are acceptable (resp. essentially tame, classical).
  \end{lemma}
  \begin{proof}
      This follows from Lemma \ref{lemma: splitting fields are keepting the same} and the facts that, parabolic subgroups of Weil-restrictions are Weil-restrictions of parabolic subgroups, $L$ and $L^{\ad}$ have the same splitting field, $G_h^{\ad}$ is a direct factor of $L^{\ad}$.
  \end{proof}

  \begin{lemma}\label{lemma: group scheme embeddings, levi}\leavevmode
  \begin{enumerate}
      \item Let $L \subset G$ be a Levi subgroup, $S$ be a maximal $F$-split torus of $G$ contained in $L$, and $x \in A_{\red}(G, S)$, $x_L \in A_{\red}(L, S)$ be its image under the canonical projection $A_{\red}(G, S) \to A_{\red}(L, S)$ (see subsection \ref{subsubsec: levi subgroups} for detailed discussions), in particular, $\Ggh_x(\OO_{\bF})\cap L(\bF) = \Ggh_{x_L}(\OO_{\bF})$, and we have a natural morphism between group schemes $\Ggh_{x_L} \to \Ggh_x$. Then the closure of $L$ in $\Ggh_x$ is $\Ggh_{x_L}$. In other words, $\Ggh_{x_L} \to \Ggh_x$ is a closed embedding.
      \item Let $G_h \subset L$ be a normal subgroup, $x_h$ be the projection of $x_L$ (see subsection \ref{subsub: normal subgroups} for detailed discussions), in particular, $\Ggh_{x_L}(\OO_{\bF})\cap G_h(\bF) = \Ggh_{x_h}(\OO_{\bF})$, and we have a natural morphism between group schemes $\Ggh_{x_h} \to \Ggh_{x_L}$. If there exists a maximal $\bF$-split torus of $G_h$ is $R$-smooth, then the closure of $G_h$ in $\Ggh_{x_L}$ is $\Ggh_{x_h}$. In other words, $\Ggh_{x_h} \to \Ggh_{x_L}$ is a closed embedding.
  \end{enumerate}
  \end{lemma}
  \begin{proof}
      \begin{enumerate}
          \item Let $\mathcal{L}_x$ be the closure of $L$ in $\Ggh_x$, then $\mathcal{L}_x(\OO_{\bF}) = L(\bF) \cap \Ggh_{x}(\OO_{\bF})$, $\Ggh_{x_L} \to \mathcal{L}_x$ is the unique smooth affine group scheme with generic fiber $L$ and $\Ggh_{x_L}(\OO_{\bF}) = \mathcal{L}_{x}(\OO_{\bF})$. Note that $\mathcal{L}_x$ contains the open big cell $\UU_{x_L}^+ \times \ZZ \times \UU_{x_L}^-$ (cf. \ref{eq: open big cell}) of $\Ggh_{x_L}$, since $\UU_{x_L}^{\pm}$, $\ZZ$ are the schematic closures of $U_{\Phi_L^{\pm}}$ and $Z = Z_G(S) = Z_L(S)$ in $\Ggh_x$ respectively, see Remark \ref{rmk: combine to positively closed set} and \ref{rmk: combine to positively closed set, over F}. Since $\UU_{x_L}^{\pm}$, $\ZZ$ are smooth due to \cite[Corollary 2.2.5]{bruhat1984groupes}, then $\mathcal{L}_x$ is smooth, and  $\Ggh_{x_L} \rightiso \mathcal{L}_x$.
          \item Let $\mathcal{G}_{h, x_L}$ be the closure of $G_h$ in $\Ggh_{x_L}$, then $\mathcal{G}_{h, x_L}(\OO_{\bF}) = G_h(\bF) \cap \Ggh_{x_L}(\OO_{\bF})$, $\Ggh_{x_h} \to\mathcal{G}_{h, x_L}$ is the unique smooth affine group scheme with generic fiber $G_h$ and $\Ggh_{x_h}(\OO_{\bF}) = \mathcal{G}_{h, x_L}(\OO_{\bF})$. If there exists a maximal $\bF$-split torus of $G_h$ is $R$-smooth, then any maximal $\bF$-split torus of $G_h$ is $R$-smooth, then $Z_h = Z_{G_h}(T_h)$ is $R$-smooth, here $T_h$ is a maximal $\bF$-split $F$-torus, the closure $\ZZ_h$ of $Z_h$ in $\ZZ$ is the finite-type N\'eron model of $Z_h$ due to \cite[Lemma 2.4.4]{kisin2024independence}. Therefore, $\mathcal{G}_{h, x_L}$ contains the big open cell $\UU_{x_h}^+ \times \ZZ_h \times \UU_{x_h}^-$ of $\Ggh_{x_h}$ (the open big cell defined using $T_h$ instead of $S_h$), thus $\Ggh_{x_h} \rightiso \mathcal{G}_{h, x_L}$.
      \end{enumerate}
  \end{proof}

\subsection{Functoriality of level groups}\label{subsec: functoriality of levis}

Let $G_1 \to G_2$ be an embedding of $F$-reductive group. Let $Q_2$ be a parabolic subgroup of $G_2$ with unipotent radical $W_2$, such that $Q_1 = G_1 \cap Q_2$ is a parabolic subgroup of $G_1$ with unipotent radical $W_1 = G_1 \cap W_1$, let $\pi_1: Q_1 \to L_1$, $\pi_2: Q_2 \to L_2$ be Levi quotients, then there is a natural embedding of Levis $L_1 \to L_2$, with $L_1 = G_1 \cap L_2$ when we fix sections $L_1 \to Q_1$ and $L_2 \to Q_2$ properly.

By the main result of \cite{landvogt2000some}, there exist $G_1(\bF)$-$\Gal(\bF|F)$-equivariant embedding $\iota_G: \Bui_{\ext}(G_1, \bF) \to \Bui_{\ext}(G_2, \bF)$ and $L_1(\bF)$-$\Gal(\bF|F)$-equivariant embedding $\iota_L: \Bui_{\ext}(L_1, \bF) \to \Bui_{\ext}(L_2, \bF)$. An elementary result \cite[Proposition 2.1.5]{landvogt2000some} shows that for $i = 1, 2$, there exists $L_i(\bF)$-$\Gal(\bF|F)$-equivariant embedding $\Bui_{\ext}(L_i, \bF) \to \Bui_{\ext}(G_i, \bF)$ and moreover, the images of $\Bui_{\ext}(L_i, \bF)$ in $\Bui_{\ext}(G_i, \bF)$ are well-defined.

Let us briefly recall the construction of $\iota_G$. Let $T_1$ be a maximal $\bF$-split torus of $G_1$, $x_1 \in A_{\ext}(G_1, T_1)$ be a special point, let $x_2 \in \Bui_{\ext}(G_2, \bF)$ be a point such that $\hat{P}^1_{x_1} \subset \hat{P}^1_{x_2}$ ($\hat{P}^1_{x_i} \subset G_i(\bF)$ are the stablizers of $x_i \in \Bui_{\ext}(G_i, \bF)$). We say $x_2$ has property $(\mathrm{TOR})$ if there exists a maximal $\bF$-split torus $T_2 \subset G_2$ which contains $T_1$ and $x_2 \in A_{\ext}(G_2, T_2)$. Since $\bF$ is strictly henselian, $x_2$ always has property $(\mathrm{TOR})$ due to \cite[Proposition 2.4.1]{landvogt2000some}. When $x_2$ has property $(\mathrm{TOR})$, there exists a unique injective morphism $\iota_{x_2, T_2}: A_{\ext}(G_1, T_1) \to A_{\ext}(G_2, T_2)$ which maps $x_1$ to $x_2$ and is compatible with $X_*(T_1) \otimes \R \to X_*(T_2) \otimes \R$. Due to \cite[Proposition 2.2.4]{landvogt2000some}, the composition $\iota_{x_2}: A_{\ext}(G_1, T_1) \to \Bui_{\ext}(G_2, \bF)$ is independent of the choice of $T_2$. 
    
    We say $x_2$ has property $(\mathrm{STAB})$ if $\hat{P}^1_{y} \subset \hat{P}^1_{\iota_{x_2}(y)}$ for all $y \in A_{\ext}(G_1, T_1)$, and $(\mathrm{CENT})$ if $\iota_{x_2}$ is $Z_1(\bF)$-equivariant, $Z_1 = Z_{G_1}(T_1)$ is the centralizer of $T_1$. When $x_2$ has property $(\mathrm{TOR})$, $(\mathrm{STAB})$ and $(\mathrm{CENT})$, then $\iota_{x_2}$ is $N_1(\bF)$-equivariant, $N_1 = N_{G_1}(T_1)$ is the normalizaer of $T_1$, and moreover it extends to a $G_1(\bF)$-equivariant toral embedding $\iota_{x_2}: \Bui_{\ext}(G_1, \bF) \to \Bui_{\ext}(G_2, \bF)$, see \cite[Proposition 2.2.8, 2.2.9]{landvogt2000some}.

    In \cite{landvogt2000some}, the author proved that we can always modify $x_2$ such that $x_2$ has property $(\mathrm{TOR})$, $(\mathrm{STAB})$ and $(\mathrm{CENT})$ and that $\iota_{x_2}$ becomes $\Gal(\bF|F)$-equivariant.

\begin{proposition}\label{prop: good levi BT building embedding}
    Let $x_1 \in \Bui_{\ext}(G_1, \bF)$ be a special point, $x_2 = \iota_G(x_1) \in \Bui_{\ext}(G_2, \bF)$, assume $x_i \in \Bui_{\ext}(L_i, \bF)$. Then there exists a unique $L_1(\bF)$-$\Gal(\bF|F)$-equivariant embedding $\iota_L: \Bui_{\ext}(L_1, \bF) \to \Bui_{\ext}(L_2, \bF)$ commutes with $\iota_G$ when we fix embeddings $\Bui_{\ext}(L_i, \bF) \to \Bui_{\ext}(L_i, \bF)$ that map $x_i$ to $x_i$ for $i = 1, 2$. In particular, $\iota_G$ maps $\Bui_{\ext}(L_1, \bF) $ into $ \Bui_{\ext}(L_2, \bF)$.
\end{proposition}
\begin{proof}
    This uniqueness comes from \cite[Proposition 2.2.10]{landvogt2000some}. Now let us show the existence. By definition, $x_2$ satisfies $(\mathrm{TOR})$, $(\mathrm{STAB})$ and $(\mathrm{CENT})$. We regard $x_1, x_2$ as elements in $\Bui_{\ext}(L_1, \bF)$, $\Bui_{\ext}(L_2, \bF)$ respectively, and denote them by $x_{L, 1}$, $x_{L, 2}$. Then $\hat{P}^1_{x_{L, i}} = \hat{P}^1_{x_i} \cap L_i(\bF)$, in particular, $\hat{P}^1_{x_{L, 1}} \subset \hat{P}^1_{x_{L, 2}}$. Moreover, $x_{L, 1}$ is a special point due to Lemma \ref{lemma: levi of special points}. Note that $x_{L, 2}$ satisfies $(\mathrm{TOR})$, $(\mathrm{STAB})$ and $(\mathrm{CENT})$ with respect to $x_{L, 1}$: $(\mathrm{TOR})$ comes for free since $\bF$ is strictly henselian. $(\mathrm{STAB})$ comes from the fact that $\iota_{x_2}: A_{\ext}(G_1, T_1) \to \Bui_{\ext}(G_2, \bF)$ (we fix $T_1 \subset L_1$ and regard $A_{\ext}(G_1, T_1) = A_{\ext}(L_1, T_1)$ as affine spaces) factors through $\Bui_{\ext}(L_2, \bF)$. $(\mathrm{CENT})$ comes from the fact that $Z_{L_1}(T_1) \subset Z_{G_1}(T_1)$.
    
    Then we have a $L_1(\bF)$-equivariant embedding $\iota_L: \Bui_{\ext}(L_1, \bF) \to \Bui_{\ext}(L_2, \bF)$ which maps $x_{L_1}$ to $x_{L, 2}$. Consider the diagram
\[\begin{tikzcd}
	{\Bui_{\ext}(L_1, \bF)} & {\Bui_{\ext}(L_2, \bF)} \\
	{\Bui_{\ext}(G_1, \bF)} & {\Bui_{\ext}(G_2, \bF)}
	\arrow["{\iota_L}", from=1-1, to=1-2]
	\arrow["{x_{L, 1} \mapsto x_1}", from=1-1, to=2-1]
	\arrow["{x_{L, 2} \mapsto x_2}"', from=1-2, to=2-2]
	\arrow["{\iota_G}", from=2-1, to=2-2]
\end{tikzcd}\]
     it is commutative due to \cite[Proposition 2.2.10]{landvogt2000some}. Since all horizontal and vertical morphisms are injective, $\iota_L$ has to be $\Gal(\bF|F)$-equivariant.
\end{proof}

\begin{lemma}\label{lemma: good levi embedding, lemma 1}
    Fix an embedding $\iota_G$ which maps $\Bui_{\ext}(L_1, \bF)$ into $\Bui_{\ext}(L_2, \bF)$. Let $x_1 \in  \Bui_{\ext}(L_1, \bF)$, $x_2 = \iota_G(x_1) \in  \Bui_{\ext}(L_2, \bF)$, let $\KK_{x_i, (r)} \subset G_i(\bF)$ be of type $(\ast)$, assume $\KK_{x_1, (r)} = G_1(\bF) \cap \KK_{x_2, (r)}$ (this is automatically true when $\KK_i = \Ggh_{x_i}(\OO_{\bF})$), then for any $g \in G_1(\bF)$,
    \[  \pi_1(g\KK_{x_1, (r)}g^{-1} \cap Q_1(\bF)) = \pi_2(g\KK_{x_2, (r)}g^{-1} \cap Q_2(\bF)) \cap L_1(\bF). \]
    The same statement is true when we replace $\bF$ with $F$.
\end{lemma}
\begin{proof}
    Assume $x_1 \in A_{\ext}(G_1, T_1)$ for some maximal $\bF$-split torus $T_1$, without lose of generality, we assume $L_1$ is at standard position with respect to $T_1$. Factor $g = ulnk$ as in Remark \ref{eq: level group at boundary}, where $u \in W_1(\bF) \subset W_2(\bF)$, $l \in L_1(\bF) \subset L_2(\bF)$, $n \in N_1(\bF)$ ($N_1 = N_{G_1}(T_1)$ is the normalizer), and $k \in \Ggh_{x_1}(\OO_{\bF}) \subset \Ggh_{x_2}(\OO_{\bF})$. Due to Corollary \ref{cor: good intersection with L}, we have
  \[ \pi_1(g\KK_{x_1, (r)}g^{-1} \cap Q_1(\bF)) = l \pi_1(\KK_{n \cdot x_1, (r)} \cap Q_1(\bF))l^{-1} = l (\KK_{n \cdot x_1, (r)} \cap L_1(\bF))l^{-1},  \] 
     \[ \pi_2(g\KK_{x_2, (r)}g^{-1} \cap Q_2(\bF)) = l \pi_2(\KK_{n \cdot x_2, (r)} \cap Q_2(\bF))l^{-1} = l (\KK_{n \cdot x_2, (r)} \cap L_2(\bF))l^{-1}, \]
     where $n \cdot x_1 = \nu_{N_1}(n)(x_1) \in A_{\ext}(G_1, T_1)$ by definition, and $n \cdot x_2 \in \Bui_{\ext}(L_2, \bF)$ since $\iota_G$ is $G_1(\bF)$-equivariant and $\iota_G$ maps $\Bui_{\ext}(L_1, \bF)$ into $\Bui_{\ext}(L_2, \bF)$. Since $\KK_{x_1, (r)} = G_1(\bF) \cap \KK_{x_2, (r)}$, then $\KK_{n \cdot x_1, (r)} = G_1(\bF) \cap \KK_{n \cdot x_2, (r)}$, $\KK_{n \cdot x_1, (r)} \cap L_1(\bF) = \KK_{n \cdot x_2, (r)} \cap L_1(\bF)$, we are done.
\end{proof}

\begin{lemma}\label{lemma: good levi embedding, lemma 2}
     Fix an embedding $\iota_G$ which maps $\Bui_{\ext}(L_1, \bF)$ into $\Bui_{\ext}(L_2, \bF)$. We furthur assume that there exist maximal $\bF$-split tori $T_1, T_2$ of $L_1$, $L_2$ respectively such that $T_1 \subset T_2$, $\iota_G$ maps $A_{\ext}(G_1, T_1)$ into $A_{\ext}(G_2, T_2)$, $N_1 \subset N_2$, $Z_1 \subset Z_2$, where $Z_i = Z_{G_i}(T_i)$, $N_i = N_{G_i}(T_i)$ are centralizer and normalizer of $T_i$ respectively.

    Let $x_i \in A_{\ext}(G_i, T_i)$, $\KK_{x_i, (r)} \subset G_i(\bF)$ be of type $(\ast)$. Assume $\KK_{x_1, (r)} = G_1(\bF) \cap \KK_{x_2, (r)}$ (but $x_2$ does not need to be the image of $x_1$ under $\iota_G$). Then for any $g \in G_1(\bF)$,
    \[  \pi_1(g\KK_{x_1, (r)}g^{-1} \cap Q_1(\bF)) = \pi_2(g\KK_{x_2, (r)}g^{-1} \cap Q_2(\bF)) \cap L_1(\bF).  \]
    The same statement is true when we replace $\bF$ with $F$.
\end{lemma}
\begin{proof}
    Without lose of generality, we assume $L_1$ is at standard position with respect to $T_1$. Factor $g = ulnk$ as in Remark \ref{eq: level group at boundary}, where $u \in W_1(\bF) \subset W_2(\bF)$, $l \in L_1(\bF) \subset L_2(\bF)$, $n \in N_1(\bF) \subset N_2(\bF)$ (by assumption), and $k \in \Ggh_{x_1}(\OO_{\bF}) \subset \Ggh_{x_2}(\OO_{\bF})$. Then, due to Corollary \ref{cor: good intersection with L},
    \[ \pi_1(g\KK_{x_1, (r)}g^{-1} \cap Q_1(\bF)) = l \pi_1(\KK_{\nu_{N_1}(n)(x_1), (r)} \cap Q_1(\bF))l^{-1} = l (\KK_{\nu_{N_1}(n)(x_1), (r)} \cap L_1(\bF))l^{-1},  \]
    \[ \pi_2(g\KK_{x_2, (r)}g^{-1} \cap Q_2(\bF)) = l \pi_2(\KK_{\nu_{N_2}(n)(x_2), (r)} \cap Q_2(\bF))l^{-1} = l (\KK_{\nu_{N_2}(n)(x_2), (r)} \cap L_2(\bF))l^{-1}. \]
    Since $\KK_{x_1, (r)} = G_1(\bF) \cap \KK_{x_2, (r)}$, and $\iota_G$ is $G_1(\bF)$-equivariant, then $\KK_{\nu_{N_1}(n)(x_1), (r)} = G_1(\bF) \cap \KK_{\nu_{N_2}(n)(x_2), (r)}$, $\KK_{\nu_{N_1}(n)(x_1), (r)} \cap L_1(\bF) = \KK_{\nu_{N_2}(n)(x_2), (r)} \cap L_1(\bF)$, we are done.
\end{proof}
\begin{remark}\label{rmk: good levi embeddings, only need normal subgroups}
    In both Lemma \ref{lemma: good levi embedding, lemma 1} and \ref{lemma: good levi embedding, lemma 2}, if $\KK_{x_2, (r)} \subset \Ggh_{x_2}(\OO_{\bF})$ is a normal subgroup and $\Ggh_{x_1}(\OO_{\bF}) = G_1(\bF) \cap \Ggh_{x_2}(\OO_{\bF})$, then the statements hold when we replace $\KK_{x_1, (r)}$ with $\KK_1 := G_1(\bF) \cap \KK_{x_2, (r)}$. Let us denote $\KK_2:= \KK_{x_2, (r)}$. $\KK_1 \subset \Ggh_{x_1}(\OO_{\bF})$ is a normal subgroup. Factor $g = ulnk$ as usual, $k \in \Ggh_{x_1}(\OO_{\bF})$. $\KK_1 = G_1(\bF) \cap \KK_2$ implies that $n\KK_1n^{-1} = G_1(\bF) \cap n\KK_2n^{-1}$, $n\KK_1n^{-1} \cap L_1(\bF) = L_1(\bF) \cap (n\KK_2n^{-1} \cap L_2(\bF))$. Consider the inclusions:

\[\begin{tikzcd}
	{\pi_1(g\KK_1g^{-1} \cap Q_1(\bF))} & {\pi_2(g\KK_2g^{-1} \cap Q_2(\bF)) \cap L_1(\bF)} \\
	{l (n\KK_1n^{-1} \cap L_1(\bF))l^{-1}} & {l (n\KK_2n^{-1} \cap L_2(\bF))l^{-1} \cap L_1(\bF)}
	\arrow[hook, from=1-1, to=1-2]
	\arrow[equals, from=1-2, to=2-2]
	\arrow[hook, from=2-1, to=1-1]
	\arrow[equals, from=2-1, to=2-2]
\end{tikzcd}\]
where the right vertical map is an identity as in the proof of Lemma \ref{lemma: good levi embedding, lemma 1} and \ref{lemma: good levi embedding, lemma 2}. This implies the left and the upper inclusions are identity.
\end{remark}

\subsection{Remarks on quasi-parahoric groups}\label{subsec: quasi-parahoric subgroups}

  This subsection would not be used in the remaining of the paper. We introduce some interesting observations which might be useful when we study general quasi-parahoric level subgroups. 
  
  The group schemes used in \cite{kisin2018integral} are $\Ggc_x \subset \Ggh_x$. The group schemes used in \cite{haines2009base} and \cite{haines2010satake} are $\Ggc_x \subset \Gg_x$.
\begin{remark}\label{remark: when is Gg equals Ggh}\leavevmode
\begin{enumerate}
    \item If the simplicial closure of $x \in A_{\red}(G, S)$ contains a special vertex, then $\Ggh_{x}(\OO_{\bF}) = \Gg_{x}(\OO_{\bF})$, $\Ggh_{x}(\OO_{F}) = \Gg_{x}(\OO_{F})$. In other words, $\Ggh_x = \Gg_x$, see \cite[\S 1.7]{bruhat1984groupes}. This identification was proved in \cite[Lemma 7.0.2]{haines2009corrigendum} when $G$ is split, and later generalized in \cite[Lemma 2.11.1]{boumasmoud2021tale} for general $G$. In particular, when $x$ is in the interior of an alcove, $\Ggc_x$ is called an Iwahori group scheme, and $\Gg_x = \Ggh_x$.
    \item When $G$ is simply connected, then $\Ggc_x = \Gg_x = \Ggh_x$, see \cite[Proposition 4.6.32]{bruhat1984groupes}.
    \item When $G$ is unramified, then $\Ggc_x = \Gg_x$. This is because $G$ splits over $\bF$, then $Z'(\bF)_0 = Z'(\bF)^1$ (see \cite[Corollary 2.3.2]{haines2009base}). In this case, $Z(F)_0 = Z(F)^1$ (see the following Proposition \ref{lemma: bijection, Z' and Z}).
    \item If $\Ggc_x = \Gg_x$ for some $x \in A_{\red}(G, S)$, then $Z(\bF)_0 = Z(\bF)^1$, $\pi_1(Z)_I$ is torsion-free, and $\Ggc_y = \Gg_y$ for any $y \in A_{\red}(G, S)$.
\end{enumerate}
\end{remark}
\subsubsection{}

    Let $\KKf \subset G(\bF)$ be a quasi-parahoric subgroup controlled by $\Ggc_x(\OO_{\bF}) \subset \KKf \subset \Ggh_x(\OO_{\bF})$ that is stable under $\sigma$-action, by the main result of \cite{bruhat1984groupes}, there exists a unique smooth affine group scheme $\Ggf_{x}$ over $\OO_{F}$ with generic fiber $G$ such that $\Ggf_{x}(\OO_{\bF}) = \KKf$. Also, by \cite{bruhat1984groupes}, $\Ggc_x(\OO_{\bF}) \subset \Ggf_x(\OO_{\bF}) \subset \Ggh_x(\OO_{\bF})$ induces morphisms $\Ggc_x \to \Ggf_x \to \Ggh_x$, and $\Ggc_x$ is the neutral connected component of $\Ggf_x$. 

   In case $\KKc:= \Ggc_x(\OO_{\bF}) \subset \KKf \subset \KKn:= \Gg_x(\OO_{\bF})$, let $Z(\bF)^{\flat} = \KKf \cap Z(\bF)$, then $Z(\bF)_0 \subset Z(\bF)^{\flat} \subset Z(\bF)^1$. Since $\KKc$ and $\KKn$ satisfy conditions \ref{eq: U^-U^+N} and \ref{eq: U^+U^-U^-Z} due to Lemma \ref{lem: parahoric subgroup satisfies the condition}, then $\KKf$ also satisfies conditions \ref{eq: U^-U^+N} and \ref{eq: U^+U^-U^-Z} due to Corollary \ref{cor: U, X, f, good factorization}. Moreover:
\begin{equation}\label{eq: nice decomposition for Kc Kf Kn}
    \KKc = U^{+}_xU^{-}_xU^{+}_x Z(\bF)_0, \quad \KKf = U^{+}_xU^{-}_xU^{+}_x Z(\bF)^{\flat}, \quad  \KKn = U^{+}_xU^{-}_xU^{+}_x Z(\bF)^1.
\end{equation}

\subsubsection{}
Recall $Z = Z_G(S)$ is a minimal Levi subgroup, $Z' = Z_G(T)$ is a maximal torus inside $Z$.
\begin{proposition}\label{lemma: bijection, Z' and Z}\leavevmode
   \begin{enumerate}
       \item $Z'(\bF)^1/Z'(\bF)_0 \to Z(\bF)^1/Z(\bF)_0$ is an isomorphism.
       \item $Z'(F)^1/Z'(F)_0 \to Z(F)^1/Z(F)_0$ is an isomorphism.
   \end{enumerate}
\end{proposition}
\begin{proof}
   \begin{enumerate}
       \item Due to \cite[Lemma 4.59]{haines2014stable}, $Z'(\bF)^1 \cap G(\bF)_0 = Z'(\bF)_0$, thus
       \begin{equation*}
           Z'(\bF)^1/Z'(\bF)_0 \to Z(\bF)^1/Z(\bF)_0
       \end{equation*}
       is injective. It suffices to show the surjectivity.
   
    Let $x \in A_{\red}(G, S)$, consider $\Ggc_x$ and $\Gg_x$. Base change the big cell \ref{eq: big cell for parahoric} to $\OO_{\bF}$, it is again an open embedding. We have open embeddings
    \begin{equation*}
      \UU_{\Phi^{\prime, +}, x} \times \ZZ' \times \UU_{\Phi^{\prime, -}, x} \hookrightarrow  (\UU_{\Phi^+, x})_{\OO_{\bF}} \times (\ZZ)_{\OO_{\bF}} \times (\UU_{\Phi^-, x})_{\OO_{\bF}} \hookrightarrow (\Gg_x)_{\OO_{\bF}}.
    \end{equation*}
    Therefore, due to Lemma \ref{lem: special fiber and integral points}, we have
    \begin{equation*}
        \Gg_x(\OO_{\bF}) = U_{\Phi^{\prime, +}, x}U_{\Phi^{\prime, -}, x}U_{\Phi^{\prime, +}, x}Z'(\bF)^1 = U_{\Phi^+, x}U_{\Phi^-, x}U_{\Phi^+, x}Z(\bF)^1.
    \end{equation*}
    Similarly, we have 
    \begin{equation*}
        \Ggc_x(\OO_{\bF}) = U_{\Phi^{\prime, +}, x}U_{\Phi^{\prime, -}, x}U_{\Phi^{\prime, +}, x}Z'(\bF)_0 = U_{\Phi^+, x}U_{\Phi^-, x}U_{\Phi^+, x}Z(\bF)_0.
    \end{equation*}
    Therefore $\Gg_x(\OO_{\bF})/\Ggc_x(\OO_{\bF}) = Z'(\bF)^1/Z'(\bF)_0 = Z(\bF)^1/Z(\bF)_0$.
     \item Consider the exact sequence 
     \begin{equation*}
         1 \to Z'(\bF)_0 \to Z'(\bF)^1 \to Z'(\bF)^1/Z'(\bF)_0 \to 1,
     \end{equation*}
     \begin{equation*}
         1 \to Z(\bF)_0 \to Z(\bF)^1 \to Z(\bF)^1/Z(\bF)_0 \to 1,
     \end{equation*}
     take $\lrangle{\sigma}$-invariants of these exact sequences, due to Lang's theorem,
     \begin{equation*}
         H^1(\lrangle{\sigma}, Z'(\bF)_0) = H^1_{\et}(\Spec \Z_p, (\ZZ')^{\circ}) = 0,
     \end{equation*}
     \begin{equation*}
         H^1(\lrangle{\sigma}, Z(\bF)_0) = H^1_{\et}(\Spec \Z_p, \ZZ^{\circ}) = 0,
     \end{equation*}
     therefore 
     \begin{equation*}
         Z'(F)^1/Z'(F)_0=(Z'(\bF)^1/Z'(\bF)_0)^{\lrangle{\sigma}} = (Z(\bF)^1/Z(\bF)_0)^{\lrangle{\sigma}}=Z(F)^1/Z(F)_0.
     \end{equation*}
   \end{enumerate}
\end{proof}

\subsubsection{}

 Recall that $T$ is a maximal $\bF$-split torus of $G$ containing $S$ and $T$ is defined over $F$, $Z = Z_G(T)$ is the centralizer of $T$ (we use $Z$ instead of $Z'$!), $Z$ is a maximal torus defined over $F$. Due to \cite[Lemma 4.59]{haines2014stable}, 
\begin{equation}\label{eq: Z to G, inj}
    Z(\bF)^1/Z(\bF)_0 \hookrightarrow G(\bF)^1/G(\bF)_0
\end{equation}
is an injection, i.e., $G(\bF)_0 \cap Z(\bF)^1 = Z(\bF)_0$. 

In general, it is not true that $G(\bF)^1 \cap Z(\bF) = Z(\bF)^1$, therefore it might not be true that $G(\bF)_0 \cap Z(\bF) = Z(\bF)_0$. For example, when $G$ is a simply connected semisimple group, then $G(\bF)^1 = G(\bF)_0 = G(\bF)$, but for any maximal torus $T \subset G$, $T(\bF)_0$ is the unique parahoric subgroup of $T(\bF)$ which is indeed a proper subgroup. Nevertheless, we have the following:

\begin{proposition}\label{prop: Z_G stablizer group}
    Let $Z_G$ be the neutral connected component of the center of $G$, then $G(\bF)^1 \cap Z_G(\bF) = Z_G(\bF)^1$, and $G(\bF)_0 \cap Z_G(\bF)$ is included in any parahoric subgroup of $G(\bF)$.
\end{proposition}
\begin{proof}
    Consider $G^{\tor}:= G/G^{\der}$, we have an isogeny of tori $\varphi: Z_G \to G^{\tor}$. By functoriality, $G(\bF)^1 \cap Z_G(\bF)$ maps to $G^{\tor}(\bF)^1$, thus $G(\bF)^1 \cap Z_G(\bF) \subset \varphi^{-1}(G^{\tor}(\bF)^1) = Z_G(\bF)^1$, see Proposition \ref{prop: isogeny of tori, parahoric}. Another way to prove this is to use Proposition \ref{proposition: sub and quotient of BT group is BT group}, note that $Z_G$ is normal in $G$.
    
    Let $\KKc$ be any parahoric subgroup with respect to a point $x \in B_{\red}(G, \bF)$, then it falls in some $A_{\red}(G, T)$ for some maximal $\bF$-split torus $T$ (not necessarily defined over $F$), let $Z$ be its centralizer, which is a maximal torus, we have $Z_G \subset Z$, and
    \[ G(\bF)_0 \cap Z_G(\bF) = G(\bF)_0 \cap Z_G(\bF)^1 = G(\bF)_0 \cap Z(\bF)^1 \cap Z_G(\bF)^1 = Z(\bF)_0 \cap Z_G(\bF)^1. \]
    In particular, $G(\bF)_0 \cap Z_G(\bF) \subset Z(\bF)_0 \subset \KKc$, we are done.
\end{proof}
\begin{remark}
    It might be possible that $Z_G(\bF)_0 \subset G(\bF)_0 \cap Z_G(\bF)$ is a strict inclusion, as we mentioned in Remark \ref{example: isogeny, non-parahoric}.
\end{remark}

\subsubsection{}
Let $\Lambda_G:= \pi_1(G)_I$, $\pi_0(\Ggf_x) := \KKf/\KKc \stackrel{\kappa_G}{\hookrightarrow} \Lambda_G$, $\Pi_{\Ggf_x} = \Ker (\pi_0(\Ggf_x)_{\sigma} \to (\Lambda_G)_{\sigma})$. $\Pi_{\Ggf_x}$ is an important index set which appears in \cite{pappas2024integral} and \cite{daniels2024conjecture} when one considers relations between objects defined under parahoric subgroups $\KKc = \Ggc_x(\OO_{\bF})$ and under general quasi-parahoric subgroups $\KKf = \Ggf_x(\OO_{\bF})$. We give a creterion when it is trivial.

   Let $W_F \subset \Sigma$ be the Weil group whose restriction on $\bF$ is an integer power of $\sigma$. We have an exact sequence $1 \to I \to W_F \to \lrangle{\sigma} \to 1$. Note that $\pi_1(Z)_{W_F} = \pi_1(Z)_{\Sigma}$, $\pi_1(G)_{W_F} = \pi_1(G)_{\Sigma}$.
    \begin{lemma}\label{lem: kernel of Z to G, torsion free}
       The kernel of $\pi_1(Z)_{W_F} \to \pi_1(G)_{W_F}$ is torsion-free.
    \end{lemma}
    \begin{proof}
       Let $\wdt{Z}$ be the preimage of the maximal torus $Z \cap G^{\der}$ along the simply connected cover $\rho: \wdt{G} \to G^{\der}$, recall that $\pi_1(G)  =  X_*(Z) / \rho_*(X_*(\wdt{Z}))$.

        Taking $W_F$-coinvariant spaces, we have an exact sequence:
        \begin{equation}\label{eq: pi_1(G)_Sigma}
            X_*(\wdt{Z})_{W_F} \to X_*(Z)_{W_F} \to \pi_1(G)_{W_F} \to 0.
        \end{equation}

        For a $\Z$-module $M$, the torsion part $M_{\tor}$ is $\Tor^1_{\Z}(M, \Q/\Z)$, thus the exact sequence \ref{eq: pi_1(G)_Sigma} induces an exact sequence:
       \begin{equation*}
           X_*(\wdt{Z})_{W_F, \tor} \to X_*(Z)_{W_F, \tor} \to \pi_1(G)_{W_F, \tor}.
       \end{equation*}

       Since $X_*(\wdt{Z})_{W_F}$ is torsion-free ($\wdt{Z}$ is an induced torus, see \cite[\S 4.4.16]{bruhat1984groupes}), $X_*(\wdt{Z})_{W_F, \tor} = 0$, thus $X_*(Z)_{W_F, \tor} \to \pi_1(G)_{W_F, \tor}$ is injective, the kernel of $\pi_1(Z)_{W_F} \to \pi_1(G)_{W_F}$ is torsion-free.
       \end{proof}
    
    \begin{proposition}
         $\Pi_{\Gg_x}$ is trivial.
    \end{proposition}
    \begin{proof}
        Recall that $Z(\bF)^1/Z(\bF)_0 \rightiso \KKn/\KKc$ (\ref{eq: nice decomposition for Kc Kf Kn}), we have 
        \begin{equation}
            \pi_0(\ZZ):=Z(\bF)^1/Z(\bF)_0 \rightiso \KKn/\KKc = \pi_0(\Gg_x).
        \end{equation}
        Consider the commutative diagram:
\[\begin{tikzcd}
	{\pi_0(\Gg_x)_{\sigma}} & {\pi_1(G)_{W_F, \tor}} \\
	{\pi_0(\ZZ)_{\sigma}} & {\pi_1(Z)_{W_F, \tor}}
	\arrow["\cong", from=2-1, to=1-1]
	\arrow[from=1-1, to=1-2]
	\arrow[from=2-1, to=2-2]
	\arrow[from=2-2, to=1-2]
\end{tikzcd}\]
      Note that $\pi_0(\ZZ) \rightiso (\Lambda_{Z})_{\tor}$, where $\Lambda_{Z} = \pi_1(Z)_I$ is a $\Sigma_0$-module, finitely generated over $\Z$, thus $\Lambda_{Z} \cong \Lambda_{Z, \tor} \oplus \Lambda_{Z, \mathrm{tor-free}}$, and both the torsion part $\Lambda_{Z, \tor}$ and torsion-free part $\Lambda_{Z, \mathrm{tor-free}}$ are $\Sigma_0$-submodules of $\Lambda_{Z}$, thus $(\Lambda_{Z})_{\sigma} \cong (\Lambda_{Z, \tor})_{\sigma} \oplus (\Lambda_{Z, \mathrm{tor-free}})_{\sigma}$. In particular, $\pi_0(\ZZ)_{\sigma} \rightiso (\Lambda_{Z, \tor})_{\sigma} \hookrightarrow (\Lambda_{Z})_{\sigma, \tor} = \pi_1(Z)_{W_F, \tor}$, the lower horizontal map in the diagram is an injection. Then the Proposition follows from Lemma \ref{lem: kernel of Z to G, torsion free}.
    \end{proof}

        \begin{remark}\label{remark: H^1 et GGc}
       The proof of \cite[Lemma 3.1.1]{pappas2024integral} shows that $\pi_0(\GG)_{\sigma} = H^1_{\et}(\Spec \Z_p, \GG)$ and $\pi_1(G)_{W_F, \tor} = H^1_{\et}(\Spec \Q_p, G)$, then the above proposition implies that \cite[Assumption 25.3.1]{scholze2020berkeley} ($H^1_{\et}(\Spec \Z_p, \GG) \to H^1_{\et}(\Spec \Q_p, G)$ is injective) holds when $\GGc(\Z_p)$ is contained in some special parahoric subgroups, see Remark \ref{remark: when is Gg equals Ggh}.
    \end{remark}

\section{Compactifications of Hodge-type Shimura varieties}\label{sec: compactifications of HT Shimura variety}

In this section, we apply the calculations from the previous section to establish properties of the boundary structures in the compactifications of integral models of Hodge-type Shimura varieties with quasi-parahoric level structures.

In subsection \ref{subsec: Compactification theory over generic fiber}, we recall the notations and the strata under consideration, following \cite{pink1989arithmetical}. 

In subsection \ref{subsec: level groups at boundary}, we compute the level groups associated with pure and mixed Shimura data at the boundary, using the full machinery developed in the previous section.

In subsection \ref{subsec: good embeddings}, we prove that starting with a good Siegel embedding as in \cite{kisin2018integral} ensures that the corresponding boundary Siegel embeddings are also good. This result is crucial for analyzing the boundaries of central leaves, Newton strata, and other strata.

In subsection \ref{subsec: integral models}, we recall the various integral models under consideration, including the Kisin-Pappas-Zhou integral models and Pappas-Rapoport integral models. 

In subsection \ref{subsection: Compactification theory over integral base}, we recall the construction of compactifications of integral models of Hodge-type Shimura varieties, following \cite{pera2019toroidal}, and establish some structural properties required for our purposes.

Finally, in subsection \ref{subsec: summary}, we summarize the key properties of the boundary structures in the compactifications of integral models of Hodge-type Shimura varieties with quasi-parahoric level structures.

\subsection{Compactification theory over generic fiber}\label{subsec: Compactification theory over generic fiber}

\subsubsection{Shimura varieties}

Let $G$ be a connected reductive group over $\Q$ and $X$ be a conjugacy class of homomorphisms of algebraic groups $h: \DS \to G_{\R}$ over $\R$ such that $(G, X)$ is a Shimura datum in the sense of \cite{deligne1979varietes}. Here $\DS$ is the Deligne torus. Let $K \subset G(\A_f)$ be an open compact subgroup (also called a \emph{level structure}), consider the double coset space
\[ \shu{K}(G, X)(\CC) = G(\Q)\backslash X \times (G(\A_f) / K) \]
     in which $G(\Q)$ acts on $X \times (G(\A_f)/K)$ diagonally. 
  
     For each $h \in X$, we have a \emph{Hodge cocharacter} $\mu_h$ of $G_{\CC}$: $\mu_h(z) := h_{\CC}(z, 1)$. Let $E = E(G, X)$ be the reflex field. Thanks to the theory of \emph{canonical model}, $\shu{K}(G, X)$ has a model defined over $E(G, X)$. For all neat $K' \subset K$, the transition morphisms $\shu{K}(G, X) \to \shu{K'}(G, X)$ are finite regular morphisms of algebraic varieties defined over $E(G, X)$. Fix $K_p$, let $\shu{K_p}(G, X):= \prolim_{K^p} \shu{K}(G, X)$. 

     In \cite[Definition 2.1]{pink1989arithmetical}, the definition of pure Shimura varieties is slightly more general than \cite{deligne1979varietes}, it allows $X$ being replaced with a $G(\R)$-equivalent map $X \to \Hom(\DS, G_{\R})$ with finite fibers. Nevertheless, there is no difference in the case of Hodge-type Shimura varieties, see \cite[Corollary 1.10]{wu2025arith}.

\subsubsection{Set up}

     In the remaining part of this section, we mainly follow the notions from \cite{pink1989arithmetical} and \cite{pera2019toroidal}. 
     \begin{enumerate}
         \item  Since we start with a reductive group $G$, not a general linear algebraic group $P$, we use $(P, W, U, V)$ instead of $(P_1, W_1, U_1, V_1)$ (the latter is used in \cite{pink1989arithmetical}) when we describe the rational boundary strata.
         \item The groups $(Q, P, W, U)$ in \cite{pera2019toroidal} are the groups $(P, Q, U, W)$ in \cite{pink1989arithmetical}. We follow the notations in \cite{pink1989arithmetical}.
     \end{enumerate}

    Fix a $\Q$-admissible parabolic subgroup $Q \subset G$, for any $x\in X$, there exists a unique morphism $h_{\infty, x}: \DS_{\CC} \to Q_{\CC}$ ($h_{\infty, x} := \omega_x\circ h_{\infty}$ in Pink's notations) which attaches a weight filtration associated with $Q$ to the Hodge filtration induced by $h_x: \DS \to G_{\R}$, see \cite[Proposition 4.6, Proposition 4.12]{pink1989arithmetical}. Let $P$ be the smallest normal $\Q$-subgroup of $Q$ such that $h_{\infty, x}$ factors through $P_{\CC}$, due to \cite[Claim 4.7]{pink1989arithmetical}, $P$ is independent of the choice of $x \in X$, and there is a $Q(\R)$-equivariant morphism $\varphi: X \to \Hom(\DS_{\CC}, P_{\CC})$ maps $x$ to $h_{\infty, x}$. Consider the conjugation action $\Ad_Q \circ h_{\infty, x}$ on $\Lie P$ and $\Lie Q$, let $W_{-n}$ be the weight space of weight $\leq -n$ under this action. Let $W$ (resp. $W_Q$) be the unipotent radical of $P$ (resp. of $Q$), the proof of \cite[Lemma 4.8]{pink1989arithmetical} showed $\Lie W = W_{-1}(\Lie W_Q)$. Since $\Lie W_Q$ is of weight $\leq -1$ when $G$ is a reductive group, thus $W = W_Q$. Let $U = \exp(W_{-2}(\Lie W_1))$, it's a normal $\Q$-subgroup of $Q$, only depends on $(G, X)$ and $Q$. \cite[Lemma 4.8]{pink1989arithmetical} also showed  $W_{\leq -3}(\Lie P) = 0$, thus $U$ is the center of $W$. Let $V = W/U$, $\Lie U$ and $\Lie V$ are the $-2$- and $-1$-graded pieces of the weight filtration on $\Lie P$ under $\Ad_P\circ h_{\infty, x}$.

    One has a $Q(\R)$-equivariant morphism
    \begin{equation}\label{eq: from interior to boundary}
        X \stackrel{(\pi_0, \varphi)}{\longrightarrow} \pi_0(X) \times \Hom(\DS_{\CC}, P_{\CC}),\quad x\mapsto ([x], h_{\infty, x}).
    \end{equation}
    Fix a connected component $X^+$ of $X$, its image under $\varphi$ lies in the same $P(\R)U(\CC)$-orbit, one gets a $P(\R)U(\CC)$-orbit $X_{P, +}$ which only depends on $(G, X), Q$ and $X^+$. Let $X' \subset X$ be the preimage of $X_{P, +}$, then $X' \to X_{P, +}$ is injective and holomorphic (\cite[Corollary 4.13]{pink1989arithmetical}). We call $(P, X_{P, +})$ a \emph{rational boundary component} of $(G, X)$, it is a mixed Shimura datum in the sense of \cite[Definition 2.1]{pink1989arithmetical}. When modulo $W$, it provides the Shimura datum of a boundary stratum of the minimal compactification of the Shimura variety associated with $(G, X)$. We explain this in details later.

    A \emph{cusp label representative} $\Phi = ((P, X_{P, +}), g)$ consists of a rational boundary component $(P, X_{P, +})$ and an element $g \in G(\A_f)$. Fix a cusp label representative $\Phi$, we define a list of notations:
    \[ \ovl{P} = P / U,\quad \ovl{Q} = Q / U,\quad G_h = P/ W,\quad L = Q/W, \quad G_l = L/G_h, \]
    \[ K_{\Phi} = gKg^{-1},\quad K_{\Phi, P}=K_{\Phi}\cap P(\A_f),\quad  K_{\Phi, Q}=K_{\Phi}\cap Q(\A_f)   \]
    \[ K_{\Phi, \ovl{P}}=K_{\Phi, P}/(U(\A_f)\cap K_{\Phi, P}), \quad K_{\Phi, \ovl{Q}}=K_{\Phi, Q}/(U(\A_f)\cap K_{\Phi, Q}), \]
    \[ K_{\Phi, h}=  K_{\Phi, P}/(W(\A_f)\cap K_{\Phi, P}), \]
    \[ X_{\Phi} = X_{P, +},\quad X_{\Phi, \ovl{P}}=U(\CC)\backslash X_{P, +}, \quad X_{\Phi, h}=W(\R)U(\CC)\backslash X_{P, +},  \]
    \[  \shu{K_{\Phi, P}}(P, X_{\Phi})(\CC) = P(\Q)\backslash X_{\Phi}\times (P(\A_f)/K_{\Phi, P}),  \]
    \[  \shu{K_{\Phi, \ovl{P}}}(\ovl{P}, X_{\Phi, \ovl{P}})(\CC) = \ovl{P}(\Q)\backslash X_{\Phi, \ovl{P}}\times (\ovl{P}(\A_f)/K_{\Phi, \ovl{P}}), \]
    \[  \shu{K_{\Phi, h}}(G_h, X_{\Phi, h})(\CC) = G_h(\Q)\backslash X_{\Phi, h}\times (G_h(\A_f)/K_{\Phi, h}). \]
    Note that $X_{\Phi}$, $X_{\Phi, \ovl{P}}$, $X_{\Phi, h}$ only depend on $(P, X_{P, +})$, not on $g$. $g$ affects the level groups.

    Let $\gamma \in G(\Q)$, $p_2 \in P_2(\A_f)$, we say $\Phi_1\stackrel{(\gamma, p_2)_K}{\longrightarrow}{\Phi_2}$ if $\gamma P_1 \gamma^{-1} \subset P_2$ (or equivalently $\gamma U_1 \gamma^{-1} \supset U_2$), $\gamma\cdot X'_1 \subset X'_2$, and $g_2K = p_2^{-1}\gamma g_1 K$. And we say $\Phi_1 \sim \Phi_2$ ($\Phi_1$ is equivalent to $\Phi_2$) if there exists $\gamma$ and $p_2$ such that $\Phi_1\stackrel{(\gamma, p_2)_K}{\longrightarrow}{\Phi_2}$, and $\gamma Q_1 \gamma^{-1} = Q_2$, $\gamma\cdot X_1' = X_2'$. Let $\Cusp_K:=\Cusp_K(G, X)$ be the set of equaivalent classes of $\Phi$.  
    
    We have a tower of spaces:
    \begin{equation}\label{eq: tower over C}
        \shu{K_{\Phi, P}}(P, X_{\Phi})(\CC) \to \shu{K_{\Phi, \ovl{P}}}(\ovl{P}, X_{\Phi, \ovl{P}})(\CC) \to \shu{K_{\Phi, h}}(G_h, X_{\Phi, h})(\CC).
    \end{equation}
    These are all mixed Shimura varieties. In particular, $\shu{K_{\Phi, h}}(G_h, X_{\Phi, h})(\CC)$ is a pure Shimura variety in the sense of \cite[Definition 2.1]{pink1989arithmetical}, and is also a pure Shimura variety in the sense of \cite{deligne1979varietes} when $(G, X)$ is of Hodge-type.

\subsubsection{A detailed look into this Tower}\label{sec: a detailed look into this tower}

    Let $K_{\Phi, W} = K_{\Phi, P} \cap W(\A_f)$, $K_{\Phi, V}$ be its image in $V(\A_f)$, let $\wdt{K}_{\Phi, W}$ be the projection of $(Z(P)(\Q)^0 \times W(\A_f)) \cap K_{\Phi, P}$ under $Z(P) \times W \to W$. Similarily, one defines $\wdt{K}_{\Phi, U}$ and $\wdt{K}_{\Phi, V}$ (regard $V$ as a normal subgroup of $\ovl{P}$ and replace $K_{\Phi, P}$ with $K_{\Phi, \ovl{P}}$). Also, let $\wdt{\Gamma}_{\Phi, (\ast)} = \wdt{K}_{\Phi, (\ast)} \cap (\ast)(\Q)$.

    Since $Z(P) = Z(\ovl{P}) = Z(G_h)$, and $K$ is neat (thus every involved open compact subgroups are neat), \cite[Corollary 3.12]{pink1989arithmetical} implies that:
     \begin{proposition}\label{prop: first properties about the tower}
   $\shu{K_{\Phi, P}}(P, X_{\Phi})(\CC) \to \shu{K_{\Phi, \ovl{P}}}(\ovl{P}, X_{\Phi, \ovl{P}})(\CC)$ (resp. $\shu{K_{\Phi, \ovl{P}}}(\ovl{P}, X_{\Phi, \ovl{P}})(\CC) \to \shu{K_{\Phi, h}}(G_h, X_{\Phi, h})(\CC)$) is in a canonical way a holomorphic torsor under the holomorphic family of abelian complex Lie groups $\shu{\wdt{K}_{\Phi, U}\rtimes K_{\Phi, \ovl{P}}}(P, X_{\Phi})(\CC) \to \shu{K_{\Phi, \ovl{P}}}(\ovl{P}, X_{\Phi, \ovl{P}})(\CC)$ (resp. $\shu{\wdt{K}_{\Phi, V}\rtimes K_{\Phi, h}}(\ovl{P}, X_{\Phi, \ovl{P}})(\CC) \to \shu{K_{\Phi, h}}(G_h, X_{\Phi, h})(\CC).$)

     \end{proposition}

    \begin{remark}
        Consider the following condition: 
    \begin{quote}
        $Z(P)^{\circ}$ is an almost product of a $\Q$-split torus and a compact-type torus.
    \end{quote}
     This is true when $(G, X)$ is a Hodge-type Shimura datum. Since any arithmetic subgroup in an almost product of $\Q$-split torus or in a compact-type torus is finite, if it is neat (thus torsion-free), then it is trivial. Therefore, $\wdt{K}_{\Phi, W} = K_{\Phi, W}$, $\wdt{K}_{\Phi, U} = K_{\Phi, U}$ and $\wdt{K}_{\Phi, V} = \pi_{W, V}(\wdt{K}_{\Phi, W}) = K_{\Phi, V}$, similariy for $\wdt{\Gamma}_{\Phi, (\ast)}$.
    \end{remark}

    Let $x \in X_{\Phi}$, the Hodge-strucure on $\Lie V$ induced by $h_{\infty, x}$ is of type $\lrbracket{(-1, 0), (0, -1)}$ and is polarizable due to \cite[Proposition 2.14(b), Proposition 3.20]{pink1989arithmetical}. Recall that we have an exact sequence of unipotent groups $1 \to U \to W \to V \to 1$ and $U$, $V$ are abelian. $\Lie V$ as a constant vector bundle $\shu{K_{\Phi, h}}(G_h, X_{\Phi, h})(\CC)$ has a variation of Hodge structure, we denote it by $\mathbf{V}_{\MH}(\Phi)$. Note that $K_{\Phi, V} \subset V(\A_f)$ is a $K_{\Phi, P}$-stable $\wdh{\Z}$-lattice, it refines $\mathbf{V}_{\MH}(\Phi)$ to a $\Z$-Hodge structure $\mathbf{V}_{\MH}(\Phi)_{\Z}$. Thus one has a canonical family of abelian varieties $A_K(\Phi)(\CC) \to \shu{K_{\Phi, h}}(G_h, X_{\Phi, h})(\CC)$ whose relative integral homology is identified with $\mathbf{V}_{\MH}(\Phi)_{\Z}$. Acturally, $A_K(\Phi)(\CC) = \shu{\wdt{K}_{\Phi, V}\rtimes K_{\Phi, h}}(\ovl{P}, X_{\Phi, \ovl{P}})(\CC)$, see \cite[\S 3.14, 3.22]{pink1989arithmetical}. $\shu{K_{\Phi, \ovl{P}}}(\ovl{P}, X_{\Phi, \ovl{P}})(\CC) \to \shu{K_{\Phi, h}}(G_h, X_{\Phi, h})(\CC)$ is a torsor under $A_K(\Phi)(\CC)$.

    On the other hand, $\shu{K_{\Phi, P}}(P, X_{\Phi})(\CC) \to \shu{K_{\Phi, \ovl{P}}}(\ovl{P}, X_{\Phi, \ovl{P}})(\CC)$ is a torsor whose fiber is identifed with $\wdt{\Gamma}_{\Phi, U}\backslash U(\CC)$. $\wdt{\Gamma}_{\Phi, U}\backslash U(\CC)$ is an algebraic torus over $\CC$. Moreover, let $\mathbf{B}_K(\Phi) = \Hom(\Z(1), \wdt{\Gamma}_{\Phi, U}) = \wdt{\Gamma}_{\Phi, U}(-1) \subset U(\Q)(-1)$, $\mathbf{E}_K(\Phi)$ be the torus defined over $\Z$ with cocharacter group $\mathbf{B}_K(\Phi)$, $P$ acts on $U$ as a character and the action factors through $\ovl{P}$, then
    \begin{equation}\label{eq: Xi to C, factorization}
        \shu{\wdt{K}_{\Phi, V}\rtimes K_{\Phi, h}}(\ovl{P}, X_{\Phi, \ovl{P}})(\CC) \cong \mathbf{E}_K(\Phi)(\CC) \times \shu{K_{\Phi, \ovl{P}}}(\ovl{P}, X_{\Phi, \ovl{P}})(\CC),
    \end{equation}
    this isomorphism is given in \cite[Example 3.16]{pink1989arithmetical}. We do not iterate here.
    
    There are canonical models $\shu{K_{\Phi, P}}(P, X_{\Phi})$, $\shu{K_{\Phi, \ovl{P}}}(\ovl{P}, X_{\Phi, \ovl{P}})$, $\shu{K_{\Phi, h}}(G_h, X_{\Phi, h})$ and $A_K(\Phi)$ defined over the same reflex field $E=E(G, X)$ due to \cite[\S 11]{pink1989arithmetical}, and the above towers descend to a tower over $E$:
    \begin{equation}\label{eq: tower, over E}
         \shu{K_{\Phi, P}}(P, X_{\Phi}) \to \shu{K_{\Phi, \ovl{P}}}(\ovl{P}, X_{\Phi, \ovl{P}}) \to \shu{K_{\Phi, h}}(G_h, X_{\Phi, h}). 
    \end{equation}

\subsubsection{Minimal compactification}\label{subsubsec: min compactifications}
    Now we describe the boundary of the minimal compactification $\shum{K}(G, X)$: we abbreviate $ W(\R)U(\CC)\backslash X_{P, +}$ as $X_{P, +}/W =: X_{\Phi, h}$ (see \cite[Proposition 2.9]{pink1989arithmetical} and the Remark below it), one first defines
    \begin{equation}\label{eq: X* decomposition}
        X^* := \bigsqcup X_{P, +}/W,
    \end{equation}
    the disjoint union is over all rational boundary components $(P, X_{P, +})$ of $(G, X)$, and then attaches the so-called Satake topology on $X^*$, see \cite[\S 6.2]{pink1989arithmetical}. Finally, let
    \[ \shum{K}(G, X)(\CC):= G(\Q)\backslash X^* \times (G(\A_f)/K) \]
    and endow it with the quotient topology. It turns out to be a projective normal complex space.
    
    Fix a representative $(P, X_{P, +})$ in a $G(\Q)$-conjugacy class of rational boundary components of $(G, X)$, let $Q$ be the admissable parabolic $\Q$-subgroup of $G$ associated with $P$. This conjugacy class contributes to $\shum{K}(G, X)(\CC)$ the following subspace: 
    \begin{equation}\label{eq: I_K, def}
        \Stab_{Q(\Q)}(X_{P, +})\backslash (X_{P, +}/W) \times (G(\A_f)/K) 
    \end{equation}
 
     Fix $I_K$ a set of representative elements in the double cosets 
     \[ \Stab_{Q(\Q)}(X_{P, +})P(\A_f)\backslash G(\A_f)/K, \] 
     let $Q(\Q)^{\star}:=\Stab_{Q(\Q)}(X_{P, +})$, it includes $P(\Q)$ by definition. Then
     \begin{align*}
          Q(\Q)^{\star}\backslash (X_{P, +}/W) &\times (G(\A_f)/K) =\\  &\bigsqcup_{g\in I} (Q(\Q)^{\star}\cap (P(\A_f) gKg^{-1}))\backslash (X_{P, +}/W) \times (P(\A_f) gK/K)   
     \end{align*}
     Note that $P$ is normal in $Q$, any element $q = p gkg^{-1} \in Q(\Q)^{\star}\cap (P(\A_f) gKg^{-1})$ acts on the second factor $p_1 gK/K \in P(\A_f) gK/K$ as follows:
     \begin{equation}\label{equation: hecke action, exact formula}
         qp_1(gK/K) = qp_1q^{-1}qgK/K= qp_1q^{-1}p gkg^{-1}g K/K = (qp_1q^{-1}) p (gK/K).
     \end{equation}

     Let $g\in I_K$, $\Phi:= ((P, X_{P, +}), g)$ forms a cusp label representative. Identify $P(\A_f)/K_{\Phi, P}$ with $P(\A_f)\cdot gK/K$ via the bijection $[x]\to [x\cdot g]$. Combining Lemma \ref{lemma: shimura variety reduces to levi level} and Equation \ref{equation: Lambda from parabolic to levi}, we have furthur identification
     \begin{align*}
         \Lambda_{\Phi, K}\backslash (P(\Q)\backslash X_{\Phi, h} \times (P(\A_f)\cdot gK/K)) &\cong \Lambda_{\Phi, K}\backslash (P(\Q)\backslash X_{\Phi, h} \times (P(\A_f)/K_{\Phi, P})) \\
         &\cong \Lambda_{\Phi, K}\backslash \shu{K_{\Phi, h}}(G_h, X_{\Phi, h})(\CC).
     \end{align*}
     
 By calculating the index set, for example, see \cite[Proposition 1.4]{wu2025arith}, we have:

     \[\shum{K}(G, X)(\CC) = \bigsqcup_{[\Phi] \in \Cusp_K(G, X)} \Lambda_{\Phi, K}\backslash \shu{K_{\Phi, h}}(G_h, X_{\Phi, h})(\CC).\]

     In order to apply results from Section \ref{sec: Bruhat-Tits buildings}, we remove the unipotent part $W$ of $P$: Since $W(\Q)$ acts trivially on $X_{\Phi, h}$, $W(\Q)$ is dense in $W(\A_f)$ (due to the strong approximation theorem for unipotent groups), $P(\A_f)/K_{\Phi, P}$ is discrete, we have
          \begin{lemma}\label{lemma: shimura variety reduces to levi level}
         \[ P(\Q)\backslash X_{\Phi, h} \times (P(\A_f)/K_{\Phi, P}) \rightiso G_h(\Q)\backslash X_{\Phi, h} \times (G_h(\A_f) / K_{\Phi, h}) \]
     \end{lemma}
     
     Let $K_{\Phi, L}:= \pi_{Q, L}(K_{\Phi, Q})$, $L(\Q)^{\star} = \pi(Q(\Q)^*)$. Since $W(\R)$ is connected, then $W(\Q) \subset Q(\Q)^{\star}$, $L(\Q)^{\star} = Q(\Q)^* \cap L(\Q)$. There is a natural identification
     \begin{equation}\label{equation: Lambda from parabolic to levi}
         \Lambda_{\Phi, K} = P(\Q) \backslash (Q(\Q)^{\star} \cap P(\A_f)K_{\Phi, Q}) \rightiso G_h(\Q)\backslash(L(\Q)^{\star}\cap G_h(\A_f) K_{\Phi, L}).
     \end{equation}
     Since there are only finitely many conjugacy classes of rational boundary components of $(G, X)$, and $I_K$ is finite due to approximation theorems, the boundary of $\shum{K}(G, X)(\CC)$ is a finite union of quotients of pure Shimura varieties. Moreover, these quotient are finite quotients. In particular, when $\shu{K}(G, X)$ is a Hodge-type Shimura variety, then its boundary data in $\shum{K}(G, X)$ are finite quotients of Hodge-type Shimura varieties.

     It follows from direct calculation that
     \begin{equation}\label{lemma: subquotient and quotient sub for unipotent radicals}
          K_{\Phi, h}:=\pi_{P, W}(P(\A_f) \cap gKg^{-1}) = \pi_{Q, W}(Q(\A_f) \cap gKg^{-1}) \cap G_h(\A_f) = K_{\Phi, L} \cap G_h(\A_f).
     \end{equation}

\subsubsection{Lambda-action}\label{subsec: triviality of lambda-action}

  It is natural to ask when does $\Lambda_{\Phi, K}$ acts trivially on $\shu{K_{\Phi, h}}(G_h, X_{\Phi, h})$, if so, this would simplify Pink's formula \cite{pink1992}. One can check this over $\CC$-points. 
  
  Let $G_l'$ be a normal subgroup of $L$ such that $G_h \times G_l' \to L$ is an almost product, $G_l' \to G_l$ is an isogeny. We start with $K_{\Phi, L} = \pi_{Q, L}(gKg^{-1} \cap Q(\A_f))$, $\pi_l:=\pi_{L, G_l}: L \to G_l$, $G_l'(\Q)^{\star} = G_l'(\Q) \cap L(\Q)^{\star}$, $G_l(\Q)^{\star} = \pi_l(L(\Q)^{\star})$. 

  \begin{proposition}\label{proposition: Lambda acts trivially a necessary condition}
         If the open compact subgroup $K_{\Phi, L} \subset L(\A_f)$ satisfies the following two conditions: 
         \begin{equation}\label{equation: first equation}
             L(\Q)^{\star}\cap G_h(\A_f)K_{\Phi, L} = G_h(\Q)(L(\Q)^{\star} \cap K_{\Phi, L}),
         \end{equation}
         \begin{equation}\label{equation: second equation}
             \pi_l(G_l'(\Q)^{\star} \cap K_{\Phi, L}) = \pi_l(L(\Q)^{\star} \cap K_{\Phi, L}),
         \end{equation}
          then $\Lambda_{\Phi, K}$ acts trivially on $\shu{K_{\Phi, h}}(G_h, X_{\Phi, h})(\CC)$.
      \end{proposition}
      \begin{proof}
        By conditions,
         \[ \Lambda_{\Phi, K} = G_h(\Q)\backslash (L(\Q)^{\star}\cap G_h(\A_f)K_{\Phi, L}) = =\pi_l(G_l'(\Q)^{\star} \cap K_{\Phi, L}). \]
         Let $q \in G_l'(\Q)^{\star} \cap K_{\Phi, L} \subset K_{\Phi, L}$ be a lifting of an element in $\Lambda_{\Phi, K}$. The action of $q$ on $(G_h(\Q)\backslash X_{\Phi, h} \times (G_h(\A_f)/K_{\Phi, h})$ maps $[x, p]$ to $q[x, p]=[qx, qp]=[qx, qpq^{-1}q]=[qx, qpq^{-1}]$. Note that $G_l'$ and $G_h$ commutes, one has $q[x, p]=[qx, qpq^{-1}]=[x, p]$, $\Lambda_{\Phi, K}$ acts trivially on $\shu{K_{\Phi, h}}(G_h, X_{\Phi, h})(\CC)$.
      \end{proof}

    \begin{corollary}\label{cor: Lambda trivial action}
       Assume $1 \to G_h \to L \to G_l \to 1$ splits as $L \cong G_h \times G_l$ (i.e. $G_l' \rightiso G_l$), and let $K$ be an $(\ast)$-type group (cf. Definition \ref{def: ast-type group, global}), then $\Lambda_{\Phi, K}$ acts trivially on $\shu{K_{\Phi, h}}(G_h, X_{\Phi, h})$ for each $[\Phi] \in \Cusp_K(G, X)$.
   \end{corollary}
   \begin{proof}
        Due to Proposition \ref{prop: main prop for section BT theory}, when $K \subset G(\A_f)$ is a $(\ast)$-type subgroup, then every $K_{\Phi, L}$ is a $(\ast)$-type subgroup. Then the Corollary follows from Proposition \ref{proposition: Lambda acts trivially a necessary condition} and Corollary \ref{corollary: if L splits and K special type then K splits}.
   \end{proof}
   \begin{remark}
       This splitting condition holds for PEL type Shimura varieties, this gives a group theoretical criterion supplement to the moduli theoretical criterion as in \cite{lan2013arithmetic}. Note that for general PEL-type Shimura data, the group $G_{h}$ used in \cite{pink1989arithmetical} is slightly different from the one used in \cite{lan2013arithmetic}. See \cite[\S 9.6]{mao2024} for details.
   \end{remark}

    Moreover, we show that two conditions in Proposition \ref{proposition: Lambda acts trivially a necessary condition} are also necessary conditions.
     \begin{proposition}{\cite[Proposition 9.4.6]{mao2024}}
          Let $(G, X)$ be a Hodge type Shimura datum, $K$ be a neat level group, then $\Lambda_{\Phi, K}$ acts trivially on $\shu{K_{\Phi, h}}(G_h, X_{\Phi, h})(\CC)$ if and only if equations \ref{equation: first equation} and \ref{equation: second equation} hold.
    \end{proposition}
    
    \begin{corollary}{\cite[Corollary 9.5.5]{mao2024}}
   Assume $G_l^{\ad}$ has no compact-type $\Q$-simple factor. If for any $K_0 \subset G(\A_f)$ open-compact subgroup, there always exists $K \subset K_0$ such that the $\Lambda_{\Phi, K}$-action on $\shu{K_{\Phi, h}}(G_h, X_{\Phi, h})$ is trivial, then $G_l^{\prime,\der} \rightiso G_l^{\der}$.
   \end{corollary}

\subsubsection{Toroidal compactifications}\label{subsubsec: toroidal compactifications}

     We reveal how the tower \ref{eq: tower, over E} is related with toroidal compactification. Since $\shu{K_{\Phi, P}}(P, X_{\Phi})(\CC) \to \shu{K_{\Phi, \ovl{P}}}(\ovl{P}, X_{\Phi, \ovl{P}})(\CC)$ is a $\mathbf{E}_K(\Phi)(\CC)$-torsor, and the cocharacter group $\mathbf{B}_K(\Phi) = \wdt{\Gamma}_{\Phi, U}(-1)$ spans $\mathbf{B}_K(\Phi)_{\R} = U(\R)(-1)$. Let $\sigma \subset U(\R)(-1)$ be a rational polyhedral cone, then we could form a twisted torus embedding (\cite[\S 2.1.17]{pera2019toroidal}):
    \[ \shu{K_{\Phi, P}}(P, X_{\Phi})(\CC) \to \shu{K_{\Phi, P}}(P, X_{\Phi}, \sigma)(\CC)\] 
   over $\shu{K_{\Phi, \ovl{P}}}(\ovl{P}, X_{\Phi, \ovl{P}})(\CC)$, we denote by $Z_{K_{\Phi, P}}(P, X_{\Phi}, \sigma)(\CC)$ the closed stratum. The twisted torus embedding as well as the closed stratum all descend over $E = E(G, X)$. 

    There is a notion $(\Phi_1, \sigma_1) \stackrel{(\gamma, p_2)_K}{\longrightarrow} (\Phi_2, \sigma_2)$ (resp. $(\Phi_1, \sigma_1) \sim (\Phi_2, \sigma_2)$) extending $\Phi_1 \stackrel{(\gamma, p_2)_K}{\longrightarrow} \Phi_2$ (resp. $\Phi_1 \sim \Phi_2$). Let $\Lambda_{K, \Phi}(\sigma) \subset \Lambda_{K, \Phi}$ be the subgroup preserving $\sigma$, then it induces strata-preserving automorphism of $\shu{K_{\Phi, P}}(P, X_{\Phi})$. $\Lambda_{\Phi, K}(\sigma)$ is automatically trivial in the Hodge-type case, see \cite[Lemma 2.1.20]{pera2019toroidal}. 
    
    Fix a collection of rational polyhedral cone decompositions $\Sigma = \lrbracket{\Sigma(\Phi)}$. The compactification $\shuc{K, \Sigma}$ is stratified by locally closed subschemes $Z_K(\Upsilon)$ with respect to $\Upsilon \in \Cusp_K(\Sigma):=\Cusp_K(G, X, \Sigma)$ (the equivalent classes of $(\Phi, \sigma)$). Fix a cusp label representative $\Phi$ in each equivalence class $[\Phi]$, then those $\Upsilon = [\Phi, \sigma] \in \Cusp_K(\Sigma)$ over $[\Phi]$ are indexed by $[\sigma] \in \Lambda_{K, \Phi}\backslash\Sigma(\Phi)^+$.
    
    We have $\Lambda_{\Phi, K}(\sigma)\backslash Z_{K_{\Phi, P}}(P, X_{\Phi}, \sigma) \hookrightarrow \shuc{K, \Sigma}(G, X)$ which is isomorphic to $Z_K(\Upsilon) \hookrightarrow \shuc{K, \Sigma}(G, X)$, $\Upsilon = [(\Phi, \sigma)] \in \Cusp_K(\Sigma)$, it extends to an isomorphism between the formal completion of $ \Lambda_{\Phi, K}(\sigma)\backslash \shu{K_{\Phi, P}}(P, X_{\Phi, P})$ at $\Lambda_{\Phi, K}(\sigma)\backslash Z_{K_{\Phi, P}}(P, X_{\Phi}, \sigma)$ and the formal completion of $\shuc{K, \Sigma}(G, X)$ at $Z_K(\Upsilon)$.

 \subsubsection{Functoriality}\label{subsubsec: functorialities}

   Let $\iota=(\psi_{\iota}, \phi_{\iota}): (G_1, X_1) \to (G_2, X_2)$ be an embedding of Shimura data. Let $\Phi_1 = ((P_1, X_{1, P_1, +}), g_1)$ be a cusp label representative of $(G_1, X_1)$, $\Phi_2:= ((P_2, X_{P_2, +}), g_2)$ be the induced cusp label representative of $(G_2, X_2)$, where $g_2 = g_1$ and $(P_2, X_{2, P_2, +}) = \iota_*((P_1, X_{1, P_1, +}))$. Let $K_1 \subset G_1(\A_f)$, $K_2 \subset G_2(\A_f)$ such that $K_1 \subset  K_2$. The natural inclusion $K_1 \subset K_2$ induces various morphisms $K_{1, \Phi_1, (\ast)_1} \to K_{2, \Phi_2, (\ast)_2}$, where $(\ast)$ could be $P, Q, \ovl{P}, \ovl{Q}, h, W, U, V$. There are functorial maps between various mixed Shimura varieties and strata on compactifications associated with $(G_1, X_1)$ and those associated with $(G_2, X_2)$, these data include $\shu{K}(G, X)$, $\shu{K_{\Phi, P}}(P, X_{\Phi})$, $\shu{K_{\Phi, \ovl{P}}}(\ovl{P}, X_{\Phi, \ovl{P}})$, $\shu{K_{\Phi, h}}(G_h, X_{\Phi, h})$, $\shu{K_{\Phi, P}}(P, X_{\Phi}, \sigma)$, $Z_{K_{\Phi, P}}(P, X_{\Phi}, \sigma)$, $Z_{K}(\Upsilon)$, $\shuc{K, \Sigma}(G, X)$, $\shum{K}(G, X)$, etc. All of these could be found in \cite{pink1989arithmetical}.

\subsection{Level groups at boundary}\label{subsec: level groups at boundary}

\subsubsection{Functoriality}\label{subsubsec: initial step}

Let $(G_1, X_1) \to (G_2, X_2)$ be an embedding of Shimura data. Due to \cite{landvogt2000some}, there is a $G_1(\bQ)$- $\Gal(\bQ|\Q_p)$-equivariant toral embedding between the extended Bruhat-Tits buildings
\begin{equation}\label{eq: Bruhat Tits building embedding}
    \iota_*: B_{\ext}(G_1, \bQ)\hookrightarrow B_{\ext}(G_2, \bQ)
\end{equation}
induced by the embedding $\iota: G_1 \hookrightarrow G_2$. $\iota_*$ maps a fixed point $x\in B_{\ext}(G_1, \bQ)$ to $y=\iota_*(x)\in B_{\ext}(G_2, \bQ)$. Then $\Ggh_{1, x}(\bZ_p) = \Ggh_{2, y}(\bZ_p) \cap G_1(\bQ)$, there is a morphism of group schemes $\Ggh_{1, x} \to \Ggh_{2, y}$ over $\Z_p$, which might not be a closed embedding. Note that $G_1(\bQ)$-$\Gal(\bQ|\Q_p)$-equivariant toral embedding $B_{\ext}(G_1, \bQ)\hookrightarrow B_{\ext}(G_2, \bQ)$ is uniquely determined by the image of a fixed special point due to \cite[Proposition 2.2.10]{landvogt2000some}. 

Let $T_1$ be a maximal $\bQ$-split $\Q_p$-torus of $G_1$ such that $x \in A_{\ext}(G_1, T_1)$. In each $[\Phi_1] \in \Cusp_{K_1}(G_1, X_1)$, as in Remark \ref{rmk: choose Q standard}, we could fix a cusp label representative $\Phi_1 = ((P_1, X_{P_1, +}), g)$ such that $(Q_1, L_1) = (Q_{1, \Phi_1}, L_{1, \Phi_1})$ is in the standard position respect to $(T_1, J_1 \subset \Delta_1)$ for some $J_1$. In particular, $x \in \Bui_{\ext}(L_{1, \Phi_1}, \bQ)$ for all such choice. 

\begin{lemma}\label{lemma: y is in all levi}
    For all induced cusp labels $\Phi_2:= \iota_*\Phi_1$, $y \in \Bui_{\ext}(L_{2, \Phi_2}, \bQ)$.
\end{lemma}
\begin{proof}
    Let $Z_1$ be the centralizer of $T_1$, it is a maximal torus over $\Q_p$. Let $\lambda_{\Phi_1}: \Gm \to Z_1$ be the cocharacter which induce $Q_{1, \Phi_1}$. On the other hand, let $T_2$ be a maximal $\bQ$-split $\Q_p$-torus of $G_2$ which contains $T_1$ and $y \in A_{\ext}(G_2, T_2)$, let $Z_2$ be the centralizer of $T_2$. Let $\Phi_2:= \iota_*\Phi_1$ be the induced cusp label, then $Q_{2, \Phi_2}:= \iota_*Q_{1, \Phi_1}$ is induced by $\lambda_{\Phi_2}$, where $\lambda_{\Phi_2}: \Gm \stackrel{\lambda_{\Phi_1}}{\to} Z_1 \to Z_2$. In particular, $L_{2, \Phi_2}$ contains $Z_2$ thus contains $T_2$ for all induced cusp labels $\Phi_2:= \iota_*\Phi_1$. This implies that $y \in \Bui_{\ext}(L_{2, \Phi_2}, \bQ)$.
\end{proof}

   \begin{corollary}\label{corollary: Bruhat Tits group are still intersection from GSp side}
        Let $K_{1, p} = \Ggh_{1, x}(\Z_p)$, $K_{2, p} = \Ggh_{2, y}(\Z_p)$. Then $K_{1, \Phi_1, L_1, p} :=  \pi(g_pK_{1, p}g_p^{-1}\cap Q_1(\Q_p))$ equals $K_{2, \Phi_2, L_2, p} \cap L_1(\Q_p)$. In particular, with the help of Lemma \ref{lemma: subquotient and quotient sub for unipotent radicals}, 
        \[K_{1, \Phi_1, h, p} = K_{1, \Phi_1, L_1, p} \cap G_{1, h}(\Q_p) = K_{2, \Phi_2, L_2, p} \cap L_1(\Q_p) \cap G_{1, h}(\Q_p) = K_{2, \Phi_2, h, p} \cap G_{1, h}(\Q_p).\]
        The same statement holds true when we replace $(K_{1, p}, K_{2, p})$ with $(K_{1, p}(n), K_{2, p}(n))$, where $K_{2, p}(n):= \ker (\Ggh_{2, y}(\Z_p) \to \Ggh_{2, y}(\Z_p/p^n\Z_p))$, and $K_{1, p}(n):= G_1(\Q_p) \cap K_{2, p}(n)$.
     \end{corollary}
     \begin{proof}
         This follows from Proposition \ref{prop: good levi BT building embedding}, Lemma \ref{lemma: good levi embedding, lemma 1} and \ref{lemma: y is in all levi}, Remark \ref{rmk: good levi embeddings, only need normal subgroups}.
     \end{proof}

     We can formulate $x_{g, h} \in \Bui_{\red}(G_{1, \Phi_1, h}, \Q_p)$, $x_{g, L} \in \Bui_{\red}(L_{1, \Phi_1}, \Q_p)$, $y_{g, h} \in \Bui_{\red}(G_{2, \Phi_2, h}, \Q_p)$, $y_{g, L} \in \Bui_{\red}(L_{2, \Phi_2}, \Q_p)$ induced by $\pr(x) \in  A_{\red}(G_1, T_1)$ and $\pr(y) \in A_{\red}(G_2, T_2)$ respectively as in Proposition \ref{prop: main prop for section BT theory}. Then $K_{1, \Phi_1, h, p} = \Ggh_{x_{g, h}}(\Z_p)$, $K_{1, \Phi_1, L_1, p} = \Ggh_{x_{g, L}}(\Z_p)$, $K_{2, \Phi_2, h, p} = \Ggh_{y_{g, h}}(\Z_p)$, $K_{2, \Phi_2, L_2, p} = \Ggh_{y_{g, L}}(\Z_p)$. Moreover, by choosing the embedding $\Bui_{\ext}(L_{1, \Phi_1}, \bQ) \to \Bui_{\ext}(G_1, \bQ)$ properly, we lift $x_{g, L} \in \Bui_{\ext}(L_{1, \Phi_1}, \Q_p)$ which has image $x_g \in \Bui_{\ext}(G_1, \Q_p)$, where $x_g = lnk(x) = l\cdot \nu_N(n)(x)$, we also regard $x_{g, h} \in \Bui_{\ext}(G_{1, \Phi_1, h}, \Q_p)$ by canonical projection from $L_{\Phi}$. We can do this to $y_{g, L}$, $y_{g, h}$ as well in a compatible way: let $\iota(x_g) = y_g := l \cdot (\nu_N(n)(y))$. Note that we do not use the decomposition $g = u_2l_2n_2k_2$ in $G_2(\Q_p)$, but use $g = ulnk$ in $G_1(\Q_p)$, thus $n$ might not be in $N$, $\nu_N(n)(y)$ might not be in $A_{\ext}(G_2, T_2)$. Nevertheless, we still have $\nu_N(n)(y) \in \Bui_{\ext}(L_2, \bQ)$ due to Proposition \ref{prop: good levi BT building embedding}, in particular, we formulate $y_{g, L}$ and $y_{g, h}$ in Proposition \ref{prop: main prop for section BT theory} using a differenta affine apartment. Choosing a different affine apartment which $y_g$ lies on does not change the group scheme $\Ggh_{y_{g, L}}$ and $\Ggh_{y_{g, h}}$.

     \begin{definition}\label{def: nice embedding}
   Let $G_1 \to G_2$ be an embedding of $\Q_p$-reductive groups, and fix $K_1 \subset G_1(\Q_p)$, $K_2 \subset G_2(\Q_p)$, $K_1 \subset K_2$. 
   \begin{enumerate}
       \item We say the embedding $(G_1, K_1) \hookrightarrow (G_2, K_2)$ is \emph{nice} if $K_i = \Ggf_{x_i}(\Z_p)$ are quasi-parahoric for some $x_i \in \Bui_{\ext}(G_i, \Q_p)$, such that $\Ggh_{x_1}(\bZ_p) = G(\bQ) \cap \Ggh_{x_2}(\bZ_p)$.
       \item We say a nice embedding $(G_1, K_1) \hookrightarrow (G_2, K_2)$ is \emph{very nice} if moreover $x_2$ is the image of $x_1$ under some $G_1(\bQ)$-$\Gal(\bQ|\Q_p)$-equivariant morphism $\iota: \Bui_{\ext}(G_1, \bQ) \to \Bui_{\ext}(G_2, \bQ)$.
   \end{enumerate}
   In case $G_1 \to G_2$ comes from an embedding of Shimura data $(G_1, X_1) \to (G_2, X_2)$, we say $(G_1, X_1, K_{1, p}) \hookrightarrow (G_2, X_2, K_{2, p})$ is \emph{nice} (resp. \emph{very nice}) if $(G_{1, \Q_p}, K_{1, p}) \hookrightarrow (G_{2, \Q_p}, K_{2, p})$ is \emph{nice} (resp. \emph{very nice}).
\end{definition}

\begin{corollary}\label{cor: very nice emb implies very nice emb}
    Let $\iota: (G_1, X_1, K_{1, p}) \hookrightarrow (G_2, X_2, K_{2, p})$ be a very nice embedding with stablizer quasi-parahoric groups $K_{1, p}$, $K_{2, p}$, where $K_{1, p} = \Ggh_{x_i}(\Z_p)$. Fix an apartment $x_1 \in A_{\ext}(G_1, S_1)$ for a maximal $\Q_p$-split torus $S_1$, take a representative $\Phi_1 = ((P_{\Phi_1}, X_{\Phi_1, +}), g)$ in each class $[\Phi_1] \in \Cusp_{K_1}(G_1, X_1)$ such that $(Q_{\Phi_1}, L_{\Phi_1})$ is at standard position with $(S_1, J_1)$ for some $J_1 \subset \Delta_1$, and $\Phi_2 = \iota_*\Phi_1$. Then all these boundary embeddings 
    \begin{equation}\label{eq: boundary embeddings}
        (G_{1, \Phi_1, h}, X_{1, \Phi_1, h}, K_{1, \Phi_1, h, p}) \hookrightarrow (G_{2, \Phi_2, h}, X_{2, \Phi_2, h}, K_{2, \Phi_2, h, p})
    \end{equation}
    are very nice embeddings, with $K_{i, \Phi_i, h, p} = \Ggh_{x_{i, g, h}}(\Z_p)$. Moreover, as in Proposition \ref{prop: main prop for section BT theory}, when we replace $(K_{1, p}, K_{2, p})$ with type $(\ast)$-groups $(K_{1, x_1, (r)}, K_{2, x_2, (r)})$ such that $K_{1, x_1, (r)} = G_1(\Q_p) \cap K_{2, x_2, (r)}$ (which is automatically true when they are stablizer quasi-parahoric), then $(K_{1, \Phi_1, h, p}, K_{2, \Phi_2, h, p}) = (K_{1, x_{1, g, h}, (r)}, K_{2, x_{2, g, h}, (r)})$, $K_{1, x_{1, g, h}, (r)} = G_{1, \Phi_1, h}(\Q_p) \cap K_{2, x_{2, g, h}, (r)}$.
\end{corollary}

\subsubsection{Modifications}\label{subsubsec: modifications}

In practice, one needs to replace $G_1 = \GSp(V, \psi)$ with a larger group $G_2 = \GSp(V', \psi')$ and replace $y_1:= y \in \Bui_{\ext}(G_1, \bQ)$ by a hyperspecial point $y_2 \in \Bui_{\ext}(G_2, \bQ)$. It is not true that this replacement induces a $G_1(\bQ)$-$\Gal(\bQ|\Q_p)$-equivariant morphism $\Bui_{\ext}(G_1, \bQ) \hookrightarrow \Bui_{\ext}(G_2, \bQ)$ which maps $y_1$ to $y_2$, but we still want Corollary \ref{corollary: Bruhat Tits group are still intersection from GSp side} holds true.

The replacement was done in two steps, first step was introduced in \cite{kisin2018integral}, second step use Zarhin's trick. Let us briefly recall the construction as in \cite[\S 1.1.11]{kisin2018integral}.

Given $y_1 \in \Bui_{\ext}(G_1, \bQ) \subset \Bui_{\ext}(\GL(V), \bQ)$, it associates with a graded almost self-dual lattice chain $(\lrbracket{\Lambda^{\bullet}}, c)$ with $\Lambda^{mr + i} = p^m\Lambda^i$ for $m \in \Z, 0 \leq i < r$,
\begin{equation}\label{eq: self-dual lattice chain}
    \Lambda^{r-1} \subset \cdots \subset \Lambda^0 \subset (\Lambda^0)^{\vee} \subset \cdots \subset (\Lambda^{r-1})^{\vee} \subset p^{-1}\Lambda^{r-1}, 
\end{equation}
$c(\Lambda^{\vee}) = -c(\Lambda) + m$ for some $m \in \Z$, $c(p^n\Lambda) = c(\Lambda) + n$, and for $a = 0$ or $1$, we have $(\Lambda^i)^{\vee} = \Lambda^{-i-a}$ for each $i$. Let $V' = \oplus^{r-1}_{i = -(r-1)-a} V$, $\psi' = \oplus^{r-1}_{i = -(r-1)-a} \psi$, then we have the natural diagonal embedding $\GSp(V, \psi) \hookrightarrow \GSp(V', \psi') \subset \GL(V')$ which factors through $\prod^{\prime, r-1}_{i = -(r-1)-a} \GSp(V, \psi)$, here $\prod^{\prime, r-1}_{i = -(r-1)-a} \GSp(V, \psi)$ is the subgroup of $\prod^{r-1}_{i = -(r-1)-a} \GSp(V, \psi)$ consists of elements
\[ (g_{-(r-1)-a}, \dots, g_{r-1}) \in \prod^{r-1}_{i = -(r-1)-a} \GSp(V, \psi), \quad \nu(g_{-(r-1)-a}) = \cdots = \nu(g_{r-1}).  \]
Consider the lattice $\Lambda' = \oplus^{r-1}_{i = -(r-1)-a} \Lambda^i \subset V'$, by replace it with a scalar multiple, we have $\Lambda' \subset \Lambda^{\prime\vee}$. Let $y_2 \in \Bui_{\ext}(G_2, \bF)$ ($G_2:= \GSp(V', \psi')$) be the point associated with $(\lrbracket{\Lambda'}, c')$, where $c'(\Lambda') = \Sigma_{i = -(r-1) - a}^{r-1} c(\Lambda^i)$. We have
\begin{align*}
    \GSP_{y_1}(\Z_p) &= \bigcap_i \GL(\Lambda_i) \cap \GSp(V, \psi)(\Q_p) \\ &=  \GL(\Lambda') \cap \GSp(V, \psi)(\Q_p) = \GSP_{y_2}(\Z_p) \cap \GSp(V, \psi)(\Q_p),
\end{align*}
and $\GG_1:=\GSP_{y_1} \hookrightarrow \GG_2:=\GSP_{y_2}$ is a closed embedding. We can furthur extend $\GSp(V, \psi) \hookrightarrow \GSp(V', \psi')$ to a Siegel embedding of Shimura data $\iota': (G_1, X_1) \hookrightarrow (G_2, X_2)$. We apply the functoriality between the boundary strata to $\shu{K_1}(G_1, X_2) \to \shu{K_2}(G_2, X_2)$, where $K_{1, p} = \GG_1(\Z_p)$, $K_{2, p} = \GG_2(\Z_p)$, $K_1^p \subset G_1(\A_f^p)$ and $K_2^p \subset G_2(\A_f^p)$ are properly chosen.

From Lemma \ref{lemma: y is in all levi}, for those $\Phi_1 = \iota_* \Phi$ coming from $\Cusp_K(G, X)$ (we write $G^{\dd}$ as $G_1$), $y_1 \in A_{\ext}(G_1, T_1)$ and $y_1$ is contained in all of these $\Bui_{\ext}(L_{\Phi_1}, \bQ)$.

\begin{lemma}\label{lem: modification, G_h}
        Let $K_{1, p} = \GG_1(\Z_p)$, $K_{2, p} = \GG_2(\Z_p)$. Let $\Phi_1 = ((P_1, X_{P_1, +}), g)$ be a cusp label representative in $\Cusp_{K_1}(G_1, X_1)$ coming from $\Cusp_K(G, X)$, $\Phi_2 = \iota'_*\Phi_1 = ((P_2, X_{P_2, +}), g) \in \Cusp_{K_2}(G_2, X_2)$. Then 
        
        \[ K_{\Phi_1, L_1, p} :=  \pi(g_pK_{1, p}g_p^{-1}\cap Q_1(\Q_p)) = K_{\Phi_2, L_2, p}\cap L_2(\Q_p).\]
        
        In particular, 
        \[K_{\Phi_1, h, p} = K_{\Phi_1, L_1, p} \cap G_{1, h}(\Q_p) = K_{\Phi_2, L_2, p} \cap L_1(\Q_p) \cap G_{1, h}(\Q_p) = K_{\Phi_2, h, p} \cap G_{1, h}(\Q_p).\]
        The same statement holds true when we replace $(K_{1, p}, K_{2, p})$ with $(K_{1, p}(n), K_{2, p}(n))$, where $K_{i, p}(n):= \ker (\GSP_{y_i}(\Z_p) \to \GSP(\Z_p/p^n\Z_p))$.
\end{lemma}
\begin{proof}
    The group embedding $G_1 \hookrightarrow G_2$ is an almost diagonal embedding and factors through $G_1' := \prod^{\prime}_{i \in I} G_1$, $I = \lrbracket{-(r-1)-a, \cdots, r-1}$, in particular, it maps $T_1$ diagonally to $T_1':= \prod^{\prime}_{i \in I} T_1 \subset \prod^{\prime}_{i \in I} G_1 \subset G_2$. Here $G_1' = \prod^{\prime}_{i \in I} G_1 \subset \prod G_1$ only differs by a center, the reduced Bruhat-Tits buildings are canonically identified. As in \cite[Proposition 2.1.6, 2.1.7]{landvogt2000some}, there is a canonical map $\iota^1:  \Bui_{\ext}(G_1, \bQ) \to \Bui_{\ext}(G_1', \bQ)$ which is $G_1(\bF)$-$\Gal(\bQ|\Q_p)$-equivariant. Since $\prod^{\prime}_{i \in I} G_1 \to G_2$ is a Levi subgroup, as in \cite[Proposition 2.1.5]{landvogt2000some}, there is a $G_1'(\bQ)$-$\Gal(\bQ|\Q_p)$-equivariant embedding $\iota^2: \Bui_{\ext}(G_1', \bQ) \to \Bui_{\ext}(G_2, \bQ)$ whose image is well-defined. By construction, one can choose a maximal $\bQ$-split torus $T_2$ contains $T_1'$, such that $Z_2 = Z_{G_2}(T_2)$ is contained in all $Q_{\Phi_2}$ for those $\Phi_2 = \iota_*\Phi_1$, see the proof of Lemma \ref{lemma: y is in all levi}. In particular $\iota^2\circ\iota^1$ maps $\Bui_{\ext}(L_{\Phi_1}, \bQ)$ to $\Bui_{\ext}(L_{\Phi_2}, \bQ)$ for all $\Phi_1 = \iota_*\Phi$.

     Under $\iota^1$, the image of $y_1$ is $y_1' = (y_1, y_1, \cdots, y_1) \in A_{\ext}(G_1', T_1')$ (to be more precise, it is only a diagonal embedding in the reduced affine apartment, we ignore the difference here and later). Since $G_1$ and $G_2$ are split, $T_1$ and $T_1'$ are indeed maximal tori of $G_1$ and $G_1'$ respectively, thus $T_2 = T_1'$ is a maximal torus of $G_2$, we can regard $y_1'$ as a point in $y_1'' \in A_{\ext}(G_2, T_2)$, and we can choose the embedding $\iota^2$ that identifies $A_{\ext}(G_1', T_1')$ with $A_{\ext}(G_2, T_2)$ (but with different affine hyperplane structures).

     In $A_{\ext}(G_1', T_1')$, we can find another point $z = (z_i)_{i \in I}$ where each $z_i \in A_{\ext}(G_1, T_1)$ is associated with the period lattice $(\lrbracket{\Lambda^i}, c)$ generated by a single $\Lambda^i$. Note that $z_i$ are indeed in the same appartment $A_{\ext}(G_1, T_1)$ with $y_1$, and the collection $\lrbracket{z_i}_{i\in I}$ are the vertices of the facet defined by $y_1$, this becomes more clear in the case of $\GL(V)$. Then $y_2 = \iota^2_*(z) \in A_{\ext}(G_2, T_2)$ is the point we need.

     Now we can apply Proposition \ref{prop: good levi BT building embedding}, Lemma \ref{lemma: good levi embedding, lemma 2} and Remark \ref{rmk: good levi embeddings, only need normal subgroups}.
\end{proof}

\begin{remark}\label{rmk: intermediate steps, G_h}
    In the proof of Lemma \ref{lem: modification, G_h}, as intermediate steps, we see that the same statement holds true when we replace $(y_2, G_2)$ with $(\prod_{i \in I} z_i, \prod_{i \in I} G_1)$, which more or less follow directly from Lemma \ref{lemma: good levi embedding, lemma 2}.
\end{remark}

One can furthur use Zarhin's trick, take $V'' = \oplus V^{\prime \oplus 8}$, $\Lambda'' = \Lambda^{\prime\oplus 4} \oplus \Lambda^{\prime\vee, \oplus 4} \subset V''$, and then there is an alternating form on $V''$ which makes $\Lambda''$ self-dual, see \cite[\S 4.5.9]{pera2012toroidal} for details. Therefore, when we consider the level groups and Bruhat-Tits buildings, we analyze exactly in the way above (use twisted diagonal embeddings). 

\begin{remark}\label{rmk: points correspondence}
   As the lines above Definition \ref{def: nice embedding}, in the proof of Lemma \ref{lem: modification, G_h} as well as the step on Zarhin's trick, we can similarly formulate $x_{g, h}$, $x_{g, L}$, $y_{g, h}$, $y_{g, L}$ from $x$ and $y$ respectively. Although $x$ is not mapped to $y$, we still find a torus $T_1 \subset T_2$ such that $x \in A_{\red}(G_1, T_1)$, $y \in A_{\red}(G_2, T_2)$. We factor $g = ulnk \in G_1(\Q_p)$, in this situation, $W_1 \subset W_2$, $L_1 \subset L_2$, $N_1 \subset N_2$, $K_1 \subset K_2$, we can define $x_g = l\nu_N(n)\cdot x$, $y_g = l\nu_N(n)\cdot y$, which coincide with $x_{g, L}$ and $y_{g, L}$ respectively. 
\end{remark}

\begin{definition}\label{def: adjusted embedding}
   Let $(G, X) \to (G^{\dd}, X^{\dd})$ be a Siegel embedding, $K_p = \Ggf_x(\Z_p)$ and $K_p^{\dd} = \GSP_y(\Z_p)$ are quasi-parahoric with respect to $x \in \Bui_{\ext}(G, \Q_p)$ and $y \in \Bui_{\ext}(G^{\dd}, \Q_p)$. We say the embedding of the triples $(G, X, K_p) \hookrightarrow (G^{\dd}, X^{\dd}, K^{\dd}_p)$ is \emph{adjusted} if $y$ is hyperspecial and there exists an intermediate Siegel embedding $(G, X) \hookrightarrow (G^{\prime\dd}, X^{\prime\dd}) \hookrightarrow (G^{\dd}, X^{\dd})$ such that 
   \begin{enumerate}
       \item there exists an embedding $\iota_*:  \Bui_{\ext}(G, \Q_p) \to \Bui_{\ext}(G^{\prime\dd}, \Q_p)$ which is $G(\bQ)$-$\Gal(\bQ|\Q_p)$-equivariant and maps $x$ to a point $y' \in \Bui_{\ext}(G^{\prime\dd}, \Q_p)$,
       \item the point $y'$ associates with a chain of almost self-dual lattice, and $y$ is related with $y'$ in the sense of step \ref{subsubsec: modifications} (in short, applying \cite{kisin2018integral} to combine that chain of lattices to a single lattice, and then use Zarhin's trick to make the lattice self-dual).
   \end{enumerate}
\end{definition}
\begin{remark}
   Given a Siegel embedding $\iota: (G, X, K_p) \hookrightarrow (G^{\dd}, X^{\dd}, K^{\dd}_p)$:
    \begin{enumerate}
        \item If $\iota$ is an adjusted Siegel embedding, then we are in the setting \ref{def: compactification, setting}, where we have good theory for compactifications $\Shumc{K, \Sigma}(G, X)$ and $\Shumm{K}(G, X)$.
        \item If $\iota$ is an adjusted Siegel embedding, then $\iota$ is also a nice Siegel embedding in the sense of definition \ref{def: nice embedding} where we have a good theory of integral model $\Shum{K}(G, X)$, see subsection \ref{subsec: pappas-rapoport integral models} in the sense of \cite{pappas2024p}.
        \item Under the condition \ref{general condition}, the Siegel embedding $\iota$ used to define Kisin-Pappas integral models is an adjusted Siegel embedding, see subsection \ref{subsec: KP group embeddings}.
        \item  Under the condition \ref{general condition, KPZ}, the Siegel embedding $\iota$ used to define Kisin-Pappas-Zhou integral models is an adjusted Siegel embedding, see \ref{subsec: KPZ group embeddings}
    \end{enumerate}
\end{remark}

Combine Corollary \ref{corollary: Bruhat Tits group are still intersection from GSp side} and Lemma \ref{lem: modification, G_h}, we have
 \begin{corollary}\label{corollary: Bruhat Tits group are still intersection from GSp side, 2}
        Given an adjusted Siegel embedding $(G, X, K_p) \hookrightarrow (G^{\dd}, X^{\dd}, K^{\dd}_p)$, assume $K_p$ is moreover a stablizer quasi-parahoric subgroup, i.e., $K_p = \Ggh_x(\Z_p)$. Then for all $\Phi^{\dd} = \iota_*\Phi$ and $n \geq 0$,
        \[ K_{\Phi, L, p} :=  \pi(g_pK_pg_p^{-1}\cap Q(\Q_p)) = K^{\dd}_{\Phi^{\dd}, L^{\dd}, p} \cap L(\Q_p). \]
        In particular,
        \[K_{\Phi, h, p} = K_{\Phi, L, p} \cap G_h(\Q_p) = K^{\dd}_{\Phi^{\dd}, L^{\dd}, p} \cap L(\Q_p) \cap G_h(\Q_p) = K^{\dd}_{\Phi^{\dd}, h, p} \cap G_h(\Q_p),\]
        thus $(G_{\Phi, h}, X_{\Phi, h}, K_{\Phi, h}) \to (G_{\Phi^{\dd}, h}^{\dd}, X_{\Phi^{\dd}, h}^{\dd}, K^{\dd}_{\Phi^{\dd}, h})$ are nice embeddings for all $\Phi^{\dd} = \iota_*\Phi$ in the sense of Definition \ref{def: nice embedding}. The same statement holds true when we replace $(K_p, K^{\dd}_p)$ with $(K_p(n), K^{\dd}_p(n))$.
     \end{corollary}

\begin{lemma}\label{lem: boundary adjusted embedding}
    Given an adjusted Siegel embedding $(G, X, K_p) \hookrightarrow (G^{\dd}, X^{\dd}, K^{\dd}_p)$, then for all $\Phi^{\dd} = \iota_*\Phi$, $(G_{\Phi, h}, X_{\Phi, h}, K_{\Phi, h, p}) \to (G_{\Phi^{\dd}, h}^{\dd}, X_{\Phi^{\dd}, h}^{\dd}, K_{\Phi^{\dd}, h, p}^{\dd})$ are adjusted Siegel embeddings.
\end{lemma}
\begin{proof}
    By definition, $(G, X, K_p) \hookrightarrow (G^{\dd}, X^{\dd}, K^{\dd}_p)$ factors through $(G^{\prime\dd}, X^{\prime\dd}, K^{\prime\dd}_p)$, where $(G, X, K_p) \to (G^{\prime\dd}, X^{\prime\dd}, K^{\prime\dd}_p)$ is a very nice Siegel embedding and we denote the image of $x$ by $y'$, and $(G^{\prime\dd}, X^{\prime\dd}, K^{\prime\dd}_p) \to (G^{\dd}, X^{\dd}, K^{\dd}_p)$ is the Siegel embedding introduced in Step \ref{subsubsec: modifications}.
    
    It follows from Corollary \ref{cor: very nice emb implies very nice emb} that $(G_{\Phi, h}, X_{\Phi, h}, K_{\Phi, h, p}) \to (G_{\Phi^{\dd}, h}^{\prime\dd}, X_{\Phi^{\dd}, h}^{\prime\dd}, K_{\Phi^{\dd}, h, p}^{\prime\dd})$ is a very nice embedding, and we denote the image of $x_{g, h}$ as $y_{g, h}'$. We need to show $(G_{\Phi^{\dd}, h}^{\prime\dd}, X_{\Phi^{\dd}, h}^{\prime\dd}, K_{\Phi^{\dd}, h, p}^{\prime\dd}) \to (G_{\Phi^{\dd}, h}^{\dd}, X_{\Phi^{\dd}, h}^{\dd}, K_{\Phi^{\dd}, h, p}^{\dd})$ is the Siegel embedding introduced in Step \ref{subsubsec: modifications}.

    To save notations, we omit $(\ast)^{\dd}$ in the remaining of the proof and omit the grading $c$ in the pair $(\Lambda^{\bullet}, c)$ since $c$ changes accordingly with the change of $\Lambda^{\bullet}$ ad has no impact on the stablizer group schemes. By definition, the Siegel embedding $(G', X', K'_p) \to (G, X, K_p)$ can be furthur factored as $(G', X', K_p') \to (G'', X'', K_p'') \to (G, X, K_p)$: As in \ref{eq: self-dual lattice chain}, $y'$ associates with a periodic chain of almost self-dual lattice $\Lambda^{\bullet}$ in $(V', \psi')$, and its corresponding $y'' \in \Bui(G'', \Q_p)$ associates with $\Lambda'' = \oplus_{i= -(r-1)-a}^{r-1} \Lambda^i \subset V''$. By replacing it with a scalar multiple, we can assume $\Lambda'' \subset \Lambda^{\prime\prime\vee}$. By applying Zarhin's trick, $y$ associated with a self-dual lattice $\Lambda := \Lambda^{\prime\prime \oplus 4} \bigoplus \Lambda^{\prime\prime\vee \oplus 4} \subset V:= V^{\prime\prime \oplus 8}$ by adjusting the alternating form on $V$. 

    Consider $y_{g, h}'$, $y_{g, h}''$, $y_{g, h}$. From the proofs of Lemma \ref{lemma: good levi embedding, lemma 2}, \ref{lem: modification, G_h} and \ref{lemma: standard position of Q and L}, we have seen that, $(y', y_{g, h}')$ corrsponds to $(\Lambda^{\bullet}, \Lambda^{\bullet}_{g, h})$, where $\Lambda^{i}_{g, h}:= \gr_{-1}(g' \Lambda^i)$, $g' = lnk$ (we write $g = ulnk$ in $G'(\Q_p)$ as usual), $g'y' \in \Bui_{\ext}(L', \Q_p)$. Similarly $(y'', y_{g, h}'')$ corrsponds to (periodic chains generated by) $(\Lambda^{\prime\prime}, \Lambda^{\prime\prime}_{g, h})$. In particular, the corresponding $y_{g, h}' \mapsto y_{g, h}''$ is presented by combining lattices as in $y' \mapsto y''$. Similarly, the corresponding $y_{g, h}'' \mapsto y_{g, h}$ is also presented by combining lattices as in $y'' \mapsto y$, and $y_{g, h}''$ is already a hyperspecial point.
\end{proof}

\subsubsection{About unipotent groups}\label{subsubsec: about unipotent groups}

In this subsection, we fix an adjusted Siegel embedding $(G, X, K_p) \hookrightarrow (G^{\dd}, X^{\dd}, K^{\dd}_p)$. 

  Given $(Q^{\dd}, L^{\dd})$, there exists an isotropic subspace $V_1 \subset V$ such that the chain $V^{\bullet}: 0=V_0 \subset V_1 \subset V_1^{\bot}:= V_2\subset V_3=V$ defines $Q^{\dd}:= \Stab_{\GSp(V)}(V^{\bullet})$ and $L^{\dd}:= \Stab_{\GSp(V)}(\gr(V^{\bullet}))$ ($\gr^i(V) := V_{i+1}/V_{i}$), and $W^{\dd}:=\{ g\in Q|\ g|_{\gr(V^{\bullet})} = \identity \}$ is the unipotent radical of $Q^{\dd}$, $G^{\dd}_h:= \{ g\in L |\ g|_{V_1 \oplus V_3/V_2 }= \identity \}$, $G^{\dd}_l:= \{ g\in L |\ g|_{V_2/V_1} = \identity \}$.

    Given any point $z \in B_{\ext}(G^{\dd}, \Q_p)$, it associates with an almost self-dual lattice $\Lambda^{\bullet} \subset V$ with a grading $c$ as in \ref{eq: self-dual lattice chain}.
    
    \begin{definition}
        We say $(Q^{\dd}, L^{\dd})$ is at the standard position with respect to a point $z\in B_{\ext}(G^{\dd}, \Q_p)$ if the lattices $\Lambda^i$ associated with $z$ satisfies $\Lambda^i = \oplus \gr^{\bullet}(\Lambda^i)$ for all $i$, here the graded filtration is the weight filtration induced by $V^{\bullet}$.
    \end{definition}

      \begin{lemma}\label{lemma: standard position of Q and L}
        $(Q^{\dd}, L^{\dd})$ is at the standard position with respect to $z$ if $z$ is in the canonical image of $B_{\ext}(L^{\dd}, \Q_p)$.
    \end{lemma}
    \begin{proof}
        Recall that, $z$ is in the facet enclosed by vertices $y_0, \dots, y_{r-1}$, where $y_i$ associates with the sub-periodic chain generated by $\Lambda^i \subset (\Lambda^i)^{\vee}$, and $z \in  B_{\ext}(L^{\dd}, \Q_p)$ implies that $y_i \in B_{\ext}(L^{\dd}, \Q_p)$ for all $i \in \lrbracket{0, \dots, r-1}$. It suffices to work with a single lattice $\Lambda \subset \Lambda^{\vee}$. 
        
        Since $L^{\dd} \cong \GSp(V_2/V_1) \times \GL(V_1)$, and its action on $V/V_2$ is determined by the action on $V_1$, use the canonical embedding $\Bui_{\ext}(\GSp, \Q_p) \hookrightarrow \Bui_{\ext}(\GL, \Q_p)$ under $\GSp \hookrightarrow \GL$, it suffice to work with a single lattice $\Lambda$ together with the weight filtration $0 \subset V_1 \subset V_2 \subset V$. Then the result is classical.
    \end{proof}
    
      Let $\pi:= \pi_{W, V}: W \to V$, since $U$, $V$ are abelian, the commutator operator on $W$ induces an alternating bilinear form $\Psi: V \times V \to U$, $\Psi$ determines the extension $1 \to U \to W \to V \to 1$ up to an isomorphism. Given $U, V, \Psi$, one can recover $W=U\times^{\Psi} V$ as follows: let $W = U\times V$ as a $\Q$-variety, and define the product $(u,v )(u', v') = (u + u' + \frac{1}{2}\Psi(v, v'), v+v')$. Although there is no section (as groups) $V \to W$ which splits the extension $1 \to U \to W \to V \to 1$, we do have a section (as sets) $V \to W$. Let $w_0=(u_0, v_0) \in W(\Q_p)$, $w=(u, v) \in W(\Q_p)$, then $(u_0, v_0)^{-1} = -(u_0, v_0)$, $(u_0, v_0)(u, v)(u_0, v_0)^{-1} = (u + \Psi(v_0, v), v)$. The conjugation actions of $Q$ on both $U$ and $V$ factor through $L$, and the only possible action of $l \in L(\Q_p)$ on $(u, v) \in W = U \times V$ is $g((u, v)) = (g(u), g(v))$. We have the following commutative diagram:
\[\begin{tikzcd}
	V & V & U \\
	{V^{\ddagger}} & {V^{\ddagger}} & {U^{\ddagger}}
	\arrow["\times"{description}, draw=none, from=1-1, to=1-2]
 \arrow["\times"{description}, draw=none, from=2-1, to=2-2]
	\arrow["\Psi", shorten <=2pt, shorten >=2pt, from=1-2, to=1-3]
	\arrow["{\Psi^{\ddagger}}", shorten <=2pt, shorten >=3pt, from=2-2, to=2-3]
	\arrow[from=1-3, to=2-3]
	\arrow[shorten <=1pt, shorten >=1pt, from=1-1, to=2-1]
	\arrow[shorten <=1pt, shorten >=1pt, from=1-2, to=2-2]
\end{tikzcd}\]

   Let $(\Lambda, c_{\Lambda})$ be the lattice associated with a vertex point $y$. Then $K_p = \GSP_y(\Z_p)$ is the subgroup of $\GSp(\Q_p)$ that stablizes $\Lambda$, and $K^{\dd}_{W^{\dd}, p} := K^{\dd}_p \cap W^{\dd}(\Q_p)$ is the subgroup of $W^{\dd}(\Q_p)$ stablizes $\Lambda$. 

   If $(Q^{\dd}, L^{\dd})$ is at the standard position with respect to $y$, we use the concrete matrices description in \cite[\S 4.1.3]{lan2017example} (the notations $U$ and $W$ are used reversely there), the set of integral points $K^{\dd}_{W^{\dd}, p}$ (resp. $K^{\dd}_{V^{\dd}, p}$) is the group $W^{\dd}(\Z_p)$ (resp. $V^{\dd}(\Z_p)$) defined there. In particular, the following natural map is a bijection:
   \begin{equation}\label{eq: natural bijection, standard V}
       K_{W^{\dd}, p}^{\dd}(n) \cap V^{\dd}(\Q_p) = K_p^{\ddagger}(n) \cap V^{\dd}(\Q_p) \to \pi^{\dd}(K_p^{\ddagger}(n) \cap W^{\dd}(\Q_p))=:K_{V^{\dd}, p}^{\dd}(n)
   \end{equation}
   Given a general point $z$, we pick the vertices $y_0, y_1, \dots, y_{r-1}$ of its closure facet as in the proof of Lemma \ref{lemma: standard position of Q and L}, note that $\GSP_{z}(\Z_p) = \bigcap_{i = 0}^{r-1} \GSP_{y_i}(\Z_p)$, since image of intersections is contained in the intersection of images, the above equation \ref{eq: natural bijection, standard V} is true when $K_p = \GSP_z(\Z_p)$.

   \begin{lemma}\label{lemma: standard pos implies sub is qup}
       If $(Q^{\dd}, L^{\dd})$ is at the standard position with respect to $z$, then $K_{V, p}(n) = K^{\dd}_{V^{\dd}, p}(n) \cap V(\Q_p)$ for $n \geq 0$.
   \end{lemma}
   \begin{proof}
       We identify $W$ with $U \times^{\Psi} V$ and $W^{\ddagger}$ with $U^{\ddagger} \times^{\Psi^{\ddagger}} V^{\ddagger}$. Since $K_p(n) = K_p^{\dd}(n) \cap G(\Q_p)$, we have $K_{W, p}(n) = K_{W^{\dd}, p}^{\dd}(n) \cap W(\Q_p)$, the lemma follows from the bijection \ref{eq: natural bijection, standard V}:
       \begin{align*}
           K_{W, p}(n) \cap V(\Q_p) \subset K_{V, p}(n) \subset K^{\dd}_{V^{\dd}, p}(n) \cap V(\Q_p) &= K_{W^{\dd}, p}^{\dd}(n) \cap V^{\dd}(\Q_p) \cap V(\Q_p) \\ &= K_{W, p}(n) \cap V(\Q_p)
       \end{align*}
   \end{proof}

   \begin{corollary}\label{corollary: equality of cpt groups on V}
       $K^{\ddagger}_{\Phi^{\dd}, V^{\ddagger}, p}(n) \cap V(\Q_p) = K_{\Phi, V, p}(n)$ for each $[\Phi] \in \Cusp_K(G, X)$.
   \end{corollary}
    \begin{proof}
       We need to separate the calculation into two steps as in steps \ref{subsubsec: initial step} and \ref{subsubsec: modifications}, thus we write $G^{\prime\dd} = \GSp(V')$ and $G^{\dd} = \GSp(V)$, and we add $\prime$ to everything related with $\GSp(V')$. We need to show $K^{\prime\ddagger}_{\Phi^{\prime\dd}, V^{\prime\ddagger}, p}(n) \cap V(\Q_p) = K_{\Phi, V, p}(n)$ for each $\Phi^{\prime\dd} = (\iota'\circ\iota)_*\Phi$ with $[\Phi] \in \Cusp_K(G, X)$, where $\iota: (G, X) \to (\GSp(V, \psi), S^{\pm})$, $\iota': (\GSp(V, \psi), S^{\pm}) \to (\GSp(V', \psi'), S^{\prime\pm})$.
       
       First, let $\pi^{\prime\ddagger}: W^{\prime\ddagger} \to V^{\prime\ddagger}$ and $\pi^{\ddagger}: W^{\ddagger} \to V^{\ddagger}$. Let $g_p = q_1m_1k_1$, where $q_1 \in Q^{\dd}(\Q_p)$, $m_1 \in N^{\dd}(\Q_p)$, $k_1 \in K^{\dd}_p$. Then 
\[K^{\prime\ddagger}_{\Phi^{\prime\dd}, V^{\prime\ddagger}, p}(n) := \pi^{\prime\dd}(g_pK^{\prime\dd}_p(n)g_p^{-1} \cap W^{\prime\dd}(\Q_p)) = q_1\pi^{\prime\dd}(m_1K_p^{\prime\dd}(n)m_1^{-1} \cap W^{\prime\dd}(\Q_p))q_1^{-1}.\] 
 Since $m_1K_p^{\prime\dd}m_1^{-1} = m_1\GSP_{y'}(\Z_p)m_1^{-1} = \GSP_{y''}(\Z_p)$, $y'' = \nu_{N^{\prime\dd}}(m_1)(y') \in A_{\ext}(L^{\prime\dd}, \bQ)$. Since $y''$ is again a hyperspecial point due to Proposition \ref{prop: main prop for section BT theory}, then we apply Lemma \ref{lemma: standard position of Q and L} and \ref{lemma: standard pos implies sub is qup}, thus
 \begin{align*}
    K_{\Phi^{\dd}, V^{\dd}, p}^{\dd}(n) &\subset K^{\prime\ddagger}_{\Phi^{\prime\dd}, V^{\prime\ddagger}, p}(n) \cap V^{\dd}(\Q_p) \\ 
    &= q_1(\pi^{\prime\dd}(m_1K_p^{\prime\dd}(n)m_1^{-1} \cap W^{\prime\dd}(\Q_p)) \cap V^{\dd}(\Q_p))q_1^{-1} \\ &= q_1(\pi^{\dd}(m_1K^{\dd}_p(n)m_1^{-1} \cap W^{\dd}(\Q_p)))q_1^{-1} \\ &= \pi^{\dd}(q_1m_1K^{\dd}_p(n)m_1^{-1}q_1^{-1} \cap W^{\dd}(\Q_p)) = K_{\Phi^{\dd}, V^{\dd}, p}^{\dd}(n).
 \end{align*}
 In particular, $K_{\Phi^{\dd}, V^{\dd}, p}^{\dd}(n) = K^{\prime\ddagger}_{\Phi^{\prime\dd}, V^{\prime\ddagger}, p}(n) \cap V^{\dd}(\Q_p)$.

  Second, let $\pi: W \to V$ and rewrite $g_p = qmk$, where $q \in Q(\Q_p)$, $m \in N(\Q_p)$, $k \in K_p$, 
\[K^{\ddagger}_{\Phi^{\dd}, V^{\ddagger}, p}(n) := \pi^{\dd}(g_pK^{\dd}_p(n)g_p^{-1} \cap W^{\dd}(\Q_p)) = q\pi^{\dd}(mK_p^{\dd}(n)m^{-1} \cap W^{\dd}(\Q_p))q^{-1}.\]   
Since $mK_p^{\dd}(n)m^{-1} = m\GSP_{y}(\Z_p)m^{-1} = \GSP_{m\cdot y}(\Z_p)$, where $m \cdot y \in \Bui_{\ext}(G^{\dd}, \bQ)$ is the image of $m \cdot x$ under the $G(\bQ)$-$\Gal(\bQ|\Q_p)$-equivariant embedding $\iota_*: \Bui_{\ext}(G, \bQ) \to \Bui_{\ext}(G^{\dd}, \bQ)$. Recall that, our choice of $x, y=\iota_*(x)$ automatically implies that $y$ is included in the image of $\Bui_{\ext}(L^{\dd}, \bQ)$. Apply Lemma \ref{lemma: standard position of Q and L} and \ref{lemma: standard pos implies sub is qup} again (where we replace a single lattice at standard position with a chain of lattices at standard position), we have $K_{V, p}(n) = K^{\dd}_{V^{\dd}, p}(n) \cap V(\Q_p)$ for $n \geq 0$. Therefore,  
 \begin{align*}
    K_{\Phi, V, p}(n) &\subset K^{\ddagger}_{\Phi^{\dd}, V^{\ddagger}, p}(n) \cap V(\Q_p) = q(\pi^{\dd}(mK_p^{\dd}(n)m^{-1} \cap W^{\dd}(\Q_p)) \cap V(\Q_p))q^{-1} \\ &= q(\pi(mK_p(n)m^{-1} \cap W(\Q_p)))q^{-1} = \pi(qmK_p(n)m^{-1}q^{-1} \cap W(\Q_p)) = K_{\Phi, V, p}(n).
 \end{align*}
  In particular, $ K_{\Phi, V, p}(n) = K^{\ddagger}_{\Phi^{\dd}, V^{\ddagger}, p}(n) \cap V(\Q_p)$. Combine these two steps together, we are done.
    \end{proof}     
\begin{lemma}\label{lemma: Kc and K have same W part}
    $K_{\Phi, W, p}(n) = K^{\circ}_{\Phi, W, p}(n)$. In particular,  $K_{\Phi, V, p}(n) = K^{\circ}_{\Phi, V, p}(n)$.
\end{lemma}
\begin{proof}
   It suffices $K_{\Phi, W, p}(n) \subset \ker\kappa_G$. Since $\ker\kappa_G$ is a normal subgroup, we let assume $g_p = 1$. Given any unipotent subgroup $W \subset G$, the composition $W(\Q_p) \to G(\Q_p) \stackrel{\kappa_G}{\to} \pi_1(G)_I^{\Sigma_0}$ factors through $\pi_1(W)_I^{\Sigma_0}$. It follows from definition that $\pi_1(W)$ is trivial, thus $W(\Q_p) \subset \Ker \kappa_G$. In particular,
    \[ K_p^{\circ}(n) \cap W(\Q_p) = K_p(n) \cap \Ker\kappa_G \cap W(\Q_p) = K_p(n) \cap W(\Q_p),  \]
    thus $K_{W, p}(n) = K^{\circ}_{W, p}(n)$.
 \end{proof}

\subsection{Good embeddings and very good embeddings}\label{subsec: good embeddings}
\subsubsection{Set up}
Let $\lrbracket{\mu}$ be the $G_{\ovl{\Q}}$-conjugacy class of a Hodge cocharacter, we choose embedding an $\ovl{\Q} \hookrightarrow \ovl{\Q}_p$ and regard $\lrbracket{\mu}$ as a $G_{\ovl{\Q}_p}$-conjugacy class. Consider the embedding:
\begin{equation}\label{eq: local embedding to GL}
    (G_{\Q_p}, \lrbracket{\mu}) \hookrightarrow (\GSp(V_{\Q_p}, \psi),  \lrbracket{\mu^{\dd}})) \hookrightarrow (\GL(V_{\Q_p}), \lrbracket{\mu^{\dd}})),
\end{equation}
the (local) Hodge type triple $(G_{\Q_p}, \lrbracket{\mu})$ (in the sense of \cite[Definition 3.1.2]{kisin2024integral}) comes from global Hodge-type Shimura datum $(G, \lrbracket{\mu})$, thus satisfies the standard assumption in \cite[\S 3.1]{kisin2024integral} and are acceptable in \cite[Definition 3.1.2]{kisin2024independence}, see \cite[Remark 3.1.5]{kisin2024integral}.

In the following part of this subsection, we only consider the local case (that comes from global Hodge-type Shimura data), thus will abbreviate $(G_{\Q_p}, \lrbracket{\mu})$ as $(G, \lrbracket{\mu})$, and we denote $E$ the completion of the reflex field $E(G, X)$ at a prime $v | p$.

Let $\GG:=\Ggh_x$ be a stablizer quasi-parahoric group scheme of $G$ associated with $x \in \Bui_{\ext}(G, \Q_p)$. One has a local model $\mathbb{M}^{\loc}_{\GG, \mu}$ associated with the triple $(G, \lrbracket{\mu}, \GG)$, it is a proper flat $\OO_E$-scheme with $\GG$-action, and has generic fiber $G/P_{\mu}$. It follows from definition that $\mathbb{M}^{\loc}_{\GG, \mu} = \mathbb{M}^{\loc}_{\GGc, \mu}$. $\mathbb{M}^{\loc}_{\GG, \mu}$ is unique and realizes the $v$-sheaf theoretic local model $\mathbf{M}^{v}_{\GG, \mu}$ over $\Spd \OO_E$ defined in \cite{scholze2020berkeley}. The existence of such $\mathbb{M}^{\loc}_{\GG, \mu}$ was proved in \cite{anschutz2022p} and \cite{gleason2024tubular}. In \cite{kisin2018integral}, the local models used there were constructed in \cite{pappas2013local}, which realizes the $v$-sheaf theoretic local model $\mathbf{M}^{v}_{\GG, \mu}$. In the situation when $(G, \lrbracket{\mu})$ only satisfies the standard assumption, \cite{kisin2024integral} gave a specific construction of $\mathbb{M}^{\loc}_{\GG, \mu}$ that is more useful in their setting.

Recall that if $(G_1, K_1) \hookrightarrow (G_2, K_2)$ is a nice embedding, then by Bruhat-Tits theory there exists a group scheme morphism $\GG_1 \to \GG_2$, where $\GG_i = \Ggh_{x_i}$ for $i = 1, 2$. It might not be a closed embedding.

\begin{definition}{\cite[Definition 3.1.6]{kisin2024independence}}\label{def: good embeddings}
    Let $\mu_1$ be a miniscule cocharacter of $G_1$, $\mu_2$ be its composition under $G_1 \to G_2$. Given a morphism between triples $(G_1, \mu_1, K_1) \hookrightarrow (G_2, \mu_2, K_2)$, it is a \emph{good} embedding if
    \begin{enumerate}
        \item $(G_1, K_1) \hookrightarrow (G_2, K_2)$ is a nice embedding,
        \item $\GG_1 \to \GG_2$ is a closed embedding of group schemes.
        \item $ \mathbb{M}^{\loc}_{\GG_1, \mu_1} \to \mathbb{M}^{\loc}_{\GG_2, \mu_2} \otimes \OO_E$ is a closed embedding.
    \end{enumerate}
    In case $G_1 \to G_2$ comes from an embedding of Shimura data $(G_1, X_1) \to (G_2, X_2)$, we say $(G_1, X_1, K_{1, p}) \hookrightarrow (G_2, X_2, K_{2, p})$ is \emph{good} if $(G_{1, \Q_p}, \mu_{1, \ovl{\Q}_p}, K_{1, p}) \hookrightarrow (G_{2, \Q_p}, \mu_{2, \ovl{\Q}_p}, K_{2, p})$ is \emph{good}, here $\mu_1$ is a Hodge-cocharacter of $(G_1, X_1)$ and $\mu_2$ is the composition of $\mu_1$ and $G_1 \to G_2$.
\end{definition}

\subsubsection{First Properties}\label{subsubsec: first properties good embedding}

Fix an adjusted Siegel embedding. In steps \ref{subsubsec: initial step} and \ref{subsubsec: modifications}, we mentioned that, under the embedding \ref{eq: local embedding to GL}, the image $y \in \Bui_{\ext}(\GSp(V, \psi), \Q_p)$ of $x$ is associated with a periodic chain of lattices with a grading, replace $((V, \psi), (\Lambda^{\bullet}, c))$ with $((V', \psi'), (\Lambda', c'))$, we get a point $y' \in \Bui_{\ext}(\GSp(V', \psi'), \Q_p)$ which is associated with a lattice $\Lambda' \subset \Lambda^{\prime\vee}$. Use Zarhin's trick, replace  $((V', \psi'), (\Lambda', c'))$ with $((V'', \psi''), (\Lambda'', c''))$, we get a point $y'' \in \Bui_{\ext}(\GSp(V'', \psi''), \Q_p)$ which is associated with a self-dual lattice $\Lambda''$. We identify $y, y', y''$ with their canonical images in the Bruhat-tits building of $\GL(V)$, $\GL(V')$ and $\GL(V'')$ respectively.

\begin{lemma}\label{lemma: change setting, good embeddings}
    In this process, we have successive closed embeddings of local models:
\begin{equation}\label{eq: successsive embeddings, local models}
    \mathbb{M}^{\loc}_{\GSP_{y}, \mu^{\dd}} \hookrightarrow \mathbb{M}^{\loc}_{\GSP_{y'}, \mu^{\dd}} \hookrightarrow \mathbb{M}^{\loc}_{\GSP_{y''}, \mu^{\dd}}\quad (\textit{induced by}\ \mathbb{M}^{\loc}_{\GLL_{y}, \mu^{\dd}} \hookrightarrow \mathbb{M}^{\loc}_{\GLL_{y'}, \mu^{\dd}} \hookrightarrow \mathbb{M}^{\loc}_{\GLL_{y''}, \mu^{\dd}})
\end{equation}
and closed embeddings of group schemes:
\begin{equation}\label{eq: successsive embeddings, group schemes}
   \GSP_y \hookrightarrow \GSP_{y'} \hookrightarrow \GSP_{y''}\quad (\textit{induced by}\ \GLL_y \hookrightarrow \GLL_{y'} \hookrightarrow \GLL_{y''})
\end{equation}
\end{lemma}
\begin{proof}
Note that both $ \mathbb{M}^{\loc}_{\GSP_{y}, \mu^{\dd}} \to \mathbb{M}^{\loc}_{\GLL_{y}, \mu^{\dd}}$ and $\GSP_y \to \GLL_y$ are closed embeddings (also for $y'$, $y''$), then we only need to consider the $\GLL$ case. For local models, we use the Grassmannian interpretations as in \cite[\S 2.3.4, 2.3.15]{kisin2018integral}. For group schemes, we use \cite[Lemma 2.3.1]{kisin2024integral} (\cite[\S 3.8, 3.9]{bruhat1984schemas}).
\end{proof}

\begin{lemma}
  Let $(G_1, \mu_1, K_1) \stackrel{f}{\hookrightarrow} (G_2, \mu_2, K_2) \stackrel{g}{\hookrightarrow} (G_3, \mu_3, K_3)$, if $g$ and the composition $g\circ f$ is a good embedding, then $f$ is a good embedding.
\end{lemma}
\begin{proof}
    This follows from definition. $f$ satisfies conditions $(2)$ and $(3)$ in Definition \ref{def: good embeddings} since $g\circ f$ does. $f$ satisfies condition $(1)$ in Definition \ref{def: good embeddings} since $g\circ f$ and $g$ do. 
\end{proof}

\subsubsection{Siegel embedding in Kisin-Pappas case}\label{subsec: KP group embeddings}

In this subsection, we assume the setting \ref{general condition}. Given a Siegel embedding $(G, X) \hookrightarrow (G^{\dd}, X^{\dd})$, under the setting \ref{general condition}, in \cite{kisin2018integral}, the authors constructed a $G(\bQ)$-$\Gal(\bQ|\Q_p)$-equivariant toral embedding $\iota_*: \Bui_{\ext}(G, \bQ) \to \Bui_{\ext}(G^{\dd}, \bQ)$ such that for all $x \in \Bui_{\ext}(G, \Q_p)$, the group morphism $\GG_{x} \to \GG^{\dd}_{\iota(x)}$ is a closed embedding, see \cite[Proposition 1.3.3, Lemma 2.3.3]{kisin2018integral}. Since $\iota_*$ is uniquely determined by the image of a prescribed special point $x \in \Bui_{\ext}(G, \bQ)$ due to \cite[Proposition 2.2.10]{landvogt2000some}, we regard $\iota_*$ as a special realization of the functorial morphism constructed abstractly in \cite{landvogt2000some}. On the other hand, let $\mu^{\dd}$ be the composition of $\mu$ under $G \to \GSp$, \cite[Proposition 2.3.7]{kisin2018integral} shows that under the setting \ref{general condition}, $\mathbb{M}^{\loc}_{\GG_x, \mu} \to \mathbb{M}^{\loc}_{\GG_{\iota(x)}^{\dd}, \mu^{\dd}} \otimes \OO_E$ is a closed embedding. Therefore, let $K_p = \GG_x(\Z_p)$, $K_p^{\dd} = \GG^{\dd}_{\iota(x)}(\Z_p)$, the very nice embedding $(G, K_p) \hookrightarrow (G^{\dd}, K^{\dd}_p)$ gives a good embedding $(G, X, K_p) \hookrightarrow (G^{\dd}, X^{\dd}, K^{\dd}_p)$.

After this, they apply the process we recalled in subsection \ref{subsubsec: modifications} to replace the Siegel Shimura datum $(G^{\dd}, X^{\dd})$ by another Siegel Shimura datum and replace the chain of almost self-dual lattices associated to $\iota(x)$ by a single lattice that is contained in its dual. We can furthur apply Zarhin's trick which we recalled in subsection \ref{subsubsec: modifications}. Then we apply Lemma \ref{lemma: change setting, good embeddings}. In summary, 
\begin{corollary}
    The good Siegel embedding $(G, X, K_p) \to (G^{\dd}, X^{\dd}, K_p^{\dd})$ constructed in \cite{kisin2018integral} is an adjusted Siegel embedding after doing Zarhin's trick.
\end{corollary}

\subsubsection{Siegel embedding in Kisin-Pappas-Zhou case}\label{subsec: KPZ group embeddings}

In this subsection, we work under the setting \ref{general condition, KPZ}. We show the following:
\begin{proposition}\label{prop: KPZ embeddings are adjusted Siegel embeddings}
    The good Siegel embedding $(G, X, K_p) \hookrightarrow (G^{\dd}, X^{\dd}, K^{\dd}_p)$ constructed in \cite[Lemma 4.1.4]{kisin2024independence} is an adjusted Siegel embedding after doing Zarhin's trick, up to replace $x \in \Bui_{\ext}(G, \Q_p)$ by a nearby point $x' \in \Bui_{\ext}(G, \Q_p)$ such that $\Ggc_x = \Ggc_{x'}$ (see \cite[Theorem 3.3.25]{kisin2024integral}) and replace $K_p$ by $K'_p = \Ggh_{x'}(\Z_p)$.
\end{proposition}

By chasing the constructions in \cite[Lemma 4.1.4]{kisin2024independence}, which is based on \cite[\S 2.4.4, Theorem 3.3.25]{kisin2024integral}, \cite[Proposition 3.1.9, 3.1.12]{kisin2024independence}, it suffices to show, in each step, the level groups are chosen in a way that the associated points in the Bruhat-Tits buildings can be induced via some maps among the buildings. This Proposition follows from applying following Corollary \ref{cor: lattice over split groups} as well as  \cite[Theorem 2.1.1, 2.2.1]{landvogt2000some} (taking fields extensions), \cite[Lemma 4.5, Proposition 4.6]{haines2020test} (taking Weil restrictions), \cite[Proposition 2.1.5]{landvogt2000some} (taking embedding of Levi subgroups), \cite[Proposition 2.1.6]{landvogt2000some} (taking products):

We work with $F$, a finite extension with $\Q_p$. 
\begin{lemma}
    Let $G_1 \hookrightarrow G_2$ be split reductive groups over $F$, let $x_1 \in \Bui_{\ext}(G_1, F)$ be a special point, $x_2 \in \Bui_{\ext}(G_2, F)$ be a point such that the stablizer group $P_{x_1} \subset G_1(F)$ of $x_1$ is contained in the stablizer group $P_{x_2} \subset G_2(F)$ of $x_2$, then there exists a unique $G_1(\bF)$-$\Gal(\bF|F)$-equivariant toral embedding $\Bui_{\ext}(G_1, \bF) \to \Bui_{\ext}(G_2, \bF)$ which maps $x_1$ to $x_2$. 
\end{lemma}
\begin{proof}
    The existence follows from the proof of \cite[Theorem 2.2.1]{landvogt2000some}, note that $x_2 \in \Bui_{\ext}(G_2, \bF)$ is $\Gal(\bF|F)$-equivariant and automatically satisfies $(\mathrm{TOR})$, $(\mathrm{STAB})$, $(\mathrm{CENT})$ since $G_1$ and $G_2$ are split over $F$. The uniqueness follows from \cite[Proposition 2.2.10]{landvogt2000some}.
\end{proof}
\begin{corollary}\label{cor: lattice over split groups}
    Let $G$ be a split reductive group over $F$, $\rho: G \to \GL(V)$ be a faithful representation. Let $x \in \Bui_{\ext}(G, F)$ be a special point, $\Lambda \subset V$ be a $\OO_F$-lattice such that $P_x \subset \GL(\Lambda)$, then one can attach a grading $c$ to the period chain generated by $\Lambda$, $y = (\Lambda^{\bullet}, c) \in \Bui_{\ext}(\GL(V), F)$, and there exists a unique $G(\bF)$-$\Gal(\bF|F)$-equivariant toral embedding $\Bui_{\ext}(G, \bF) \to \Bui_{\ext}(\GL(V), \bF)$ which maps $x$ to $y$. The same statement is true when we replace $\GL(V)$ by $\GSp(V, \psi)$ and replace $\Lambda$ by $\Lambda \subset \Lambda^{\vee}$, an $\OO_F$-lattice contained in its dual under $\psi$.
\end{corollary}

\subsubsection{Boundary embeddings are good}

Fix an adjusted Siegel embedding $(G, X, K_p) \to (G^{\dd}, X^{\dd}, K_p^{\dd})$ with $K_p$ being stablizer quasi-parahoric group $K_p:=\Ggh_x(\Z_p)$ and $K^{\dd}_p = \Ggh_y(\Z_p)$. We furthur assume $(G, \mu, K_p) \to (G^{\dd}, \mu^{\dd}, K_p^{\dd})$ is a good embedding.

Our goal is to show the morphism on the triples between boundary strata are all good embeddings under some week conditions. Let $\mu_h$ be a Hodge cocharacter of $(G_{\Phi, h}, X_{\Phi, h})$, due to \cite[Proposition 12.1]{pink1989arithmetical}, the Hodge cocharacters $\mu_h$ and $\mu$ (associated with $X$, the interior Shimura datum) are conjugated under $Q_{\Phi}$ and their conjugacy classes define the same reflex field. In particular, we can a priori pick a Hodge cocharacter $\mu$ in its $G(\ovl{\Q}_p)$-conjugacy class and assume $\mu = \mu_h$. Since we only consider the conjugacy classes of $\mu$ and $\mu_h$, we always do this manipulation.

 For each cusp label representative $\Phi$, we have a Siegel embedding $(G_{\Phi, h}, X_{\Phi, h}) \hookrightarrow (G^{\dd}_{\Phi^{\dd}, h}, X_{\Phi^{\dd}, h}^{\dd})$, we denote $\mu_h^{\dd}$ by the composition of $\mu_h$ under $G_{\Phi, h} \to G^{\dd}_{\Phi^{\dd}, h}$, then $\mu_h^{\dd}$ and $\mu^{\dd}$ are in the same  $\GSp(V, \psi)(\ovl
 {\Q}_p)$-conjugacy class. Due to Lemma \ref{lem: boundary adjusted embedding}, the embedding $(G_{\Phi, h}, K_{\Phi, h}) \hookrightarrow (G^{\dd}_{\Phi^{\dd}, h}, K^{\dd}_{\Phi^{\dd}, h})$ is an adjusted Siegel embedding. Let $x_{g, h} \in \Bui_{\ext}(G_{\Phi, h}, \Q_p)$ and $y_{g, h} \in \Bui_{\ext}(G^{\dd}_{\Phi^{\dd}, h}, \Q_p)$ be induced by $x \in \Bui_{\ext}(G, \Q_p)$ and $y \in \Bui_{\ext}(G^{\dd}, \Q_p)$ respectively in Remark \ref{rmk: points correspondence}.

\begin{proposition}\label{prop: good embedding}
   Under the setting \ref{general condition} (resp. setting \ref{general condition, KPZ} with extra assumption that $G_{\Phi, h}$ are $R$-smooth for all $[\Phi] \in \Cusp_K(G, X)$), given an adjusted Siegel embedding $(G, X, K_p) \to (G^{\dd}, X^{\dd}, K_p^{\dd})$ which is a good embedding, then it induces adjusted Siegel embeddings $(G_{\Phi, h}, X_{\Phi, h}, K_{\Phi, h}) \to (G_{\Phi^{\dd}, h}^{\dd}, X_{\Phi^{\dd}, h}^{\dd}, K^{\dd}_{\Phi^{\dd}, h})$ which are good embeddings.
\end{proposition}
\begin{proof}
   Boundary morphisms are again adjusted Siegel embeddings due to Lemma \ref{lem: boundary adjusted embedding}. We need to show boundary morphisms are good. The first, second, third conditions of being a good morphism (Definition \ref{def: good embeddings}) are given in Corollary \ref{corollary: Bruhat Tits group are still intersection from GSp side, 2}, Lemma \ref{lemma: closed group embedding} and \ref{lemma: closed local model embedding} respectively.
\end{proof}

We fix a $\Phi$ from now, and omit it in the index to save notations.
\begin{lemma}\label{lemma: closed group embedding}
    We follow the assumptions in Proposition \ref{prop: good embedding}. The group scheme homomorphism $\Ggh_{x_{g, h}} \to \Ggh_{y_{g, h}}$ is a closed embedding.
\end{lemma}
\begin{proof}
    In the setting \ref{def: adjusted embedding}, the construction of $\Ggh_{x_{g, h}} \to \Ggh_{y_{g, h}}:= \Ggh_{y_{g, h}'}$ is separated into two steps. In the second step \ref{subsubsec: modifications}, the group scheme homomorphism $\Ggh_{y_{g, h}} \to \Ggh_{y_{g, h}'}$ is a closed embedding due to Lemma \ref{lemma: change setting, good embeddings}. Therefore, we focus on the fisrt step \ref{subsubsec: initial step}, where there exists a $G(\bQ)$-$\Gal(\bQ|\Q_p)$- equivariant morphism $\iota: \Bui_{\ext}(G, \bQ) \to \Bui_{\ext}(G^{\dd}, \bQ)$ such that $y = \iota(x)$, $y$ might not be a hyperspecial point.

   Recall that in Corollary \ref{cor: very nice emb implies very nice emb}, we have commutative diagram of group schemes homomorphisms
\[\begin{tikzcd}
	{\Ggh_{y_{g, h}}} & {\Ggh_{y_{g, L}}} & {\Ggh_{y_g}} \\
	{\Ggh_{x_{g, h}}} & {\Ggh_{x_{g, L}}} & {\Ggh_{x_g}}
	\arrow[from=1-1, to=1-2]
	\arrow[from=1-2, to=1-3]
	\arrow[from=2-1, to=1-1]
	\arrow[from=2-1, to=2-2]
	\arrow[from=2-2, to=1-2]
	\arrow[from=2-2, to=2-3]
	\arrow[from=2-3, to=1-3]
\end{tikzcd}\]
   where the right vertical morphism is a closed embedding since it is $g'$-conjugation of the closed embedding $\Ggh_{x} \to \Ggh_{y}$, where $g' = lnk$ as in the lines above Definition \ref{def: nice embedding}. The lower morphisms are again closed embeddings due to Lemma \ref{lemma: group scheme embeddings, levi}. Note that under the setting \ref{general condition}, $G$ splits over a tamely ramified extension of $\Q_p$, then $G_h$ splits over a tamely ramified extension of $\Q_p$ due to Lemma \ref{lemma: splitting fields are keepting the same}, thus is automatically $R$-smooth.
\end{proof}

\begin{lemma}\label{lemma: closed local model embedding}
     We follow the assumptions in Proposition \ref{prop: good embedding}. The morphism between local models $\mathbb{M}^{\loc}_{\Ggh_{x_{g, h}}, \mu_h} \to \mathbb{M}^{\loc}_{\Ggh_{y_{g, h}}, \mu_h^{\dd}}$ is a closed embeddings.
\end{lemma}
\begin{proof}
    As in the proof of Lemma \ref{lemma: closed group embedding}, we only need to consider the fisrt step \ref{subsubsec: initial step}, use functoriality of local models, we have following commutative diagram:
\[\begin{tikzcd}
	{\mathbb{M}^{\loc}_{\Ggh_{y_{g, h}}, \mu_h^{\dd}}} & {\mathbb{M}^{\loc}_{\Ggh_{y_{g, L}}, \mu_h^{\dd}}} & {\mathbb{M}^{\loc}_{\Ggh_{y_{g}}, \mu_h^{\dd}}} \\
	{\mathbb{M}^{\loc}_{\Ggh_{x_{g, h}}, \mu_h}} & {\mathbb{M}^{\loc}_{\Ggh_{x_{g, L}}, \mu_h}} & {\mathbb{M}^{\loc}_{\Ggh_{x_{g}}, \mu_h}}
	\arrow[from=1-1, to=1-2]
	\arrow[from=1-2, to=1-3]
	\arrow[from=2-1, to=1-1]
	\arrow[from=2-1, to=2-2]
	\arrow[from=2-2, to=1-2]
	\arrow[from=2-2, to=2-3]
	\arrow[from=2-3, to=1-3]
\end{tikzcd}\]
 where the right vertical map is a closed embedding sine it is the $g'$-conjugation of the closed embedding $\mathbb{M}^{\loc}_{\Ggh_{x}, \mu_h} \to \mathbb{M}^{\loc}_{\Ggh_{y}, \mu_h^{\dd}}$. We need to show lower morphisms are closed embeddings. The first morphism $\mathbb{M}^{\loc}_{\Ggh_{x_{g, h}}, \mu_h} \to \mathbb{M}^{\loc}_{\Ggh_{x_{g, L}}, \mu_h}$ is an isomorphism), since $L^{\ad} = G_h^{\ad} \times G_l^{\ad}$ and $\mu_h$ factors through $G_h$, see \cite[Proposition 21.5.1]{scholze2020berkeley}, \cite[\S 3.4.2]{kisin2024integral}. The second morphism $i_{L, G}: \mathbb{M}^{\loc}_{\Ggh_{x_{g, L}}, \mu_h} \to \mathbb{M}^{\loc}_{\Ggh_{x_{g}}, \mu_h}$ is a closed embedding can be checked via Zariski main theorem: we take $\mathbb{M}^{\loc, \prime}$ as the normalization of the closure of the generic fiber of $\mathbb{M}^{\loc}:=\mathbb{M}^{\loc}_{\Ggh_{x_{g, L}}, \mu_h}$ in $\mathbb{M}^{\loc}_{\Ggh_{x_{g}}, \mu_h}$, and to show the proper birational morphism between the noetherian normal schemes $\mathbb{M}^{\loc} \to \mathbb{M}^{\loc, \prime}$ is an isomorphism, it suffices to show $i_{L, G}$ is an injection on special fibers. Note that on the perfection of the special fibers, one use the descriptions of Witt affine Grassmannian and use stratifications of affine Schubert varieties, see \cite[Theorem 1.1]{anschutz2022p}.
\end{proof}

\begin{remark}
As explained in \cite{kisin2024integral}, the construction in \cite{kisin2018integral} actually relies on a notion called \emph{very good embedding} (\cite[Definition 5.2.5]{kisin2024integral}): choosing a good embedding is not sufficient for the construction in \cite{kisin2018integral}, one really needs to choose a very good embedding. We do not recall the precise definition, just note that good embeddings are very good embeddings in most of the cases, especially when $(G^{\ad}, \mu^{\ad})$ coming from a global Hodge type Shimura datum has no factor of type $D^{\mathbb{H}}$, see \cite[Theorem 4.4.3, Proposition 5.3.10]{kisin2024integral}. In this paper, we simply ignore this issue since we are not working with local model diagrams.
\end{remark}

\subsection{Integral models with quasi-parahoric level structures}\label{subsec: integral models}

\subsubsection{Relative normalizations}

  In this subsection we only consider locally noetherian schemes. We list some very useful propositions of relative normalizations.

         \begin{lemma}\label{lemma: rela norm of rela norm is rela norm}
         Let $Y \to Z \to X$ be qcqs morphisms, let $Z_X \to X$ (resp. $Y_Z$, $Y_X$) be the relative normalization of $X$ in $Z$ (resp. $Z_X$ in $Y$, $X$ in $Y$), then there is a canonical isomorphism $Y_X = Y_Z$. 

     \end{lemma}
     \begin{proof}
         Due to \cite[\href{https://stacks.math.columbia.edu/tag/035I}{Tag 035I}]{stacks-project}, there is a unique morphism $Y_X \to Y_Z$ such that $Y \to Y_Z$ factors through $Y_X$, and $Y_X \to Y_Z$ is the relative normalization of $Y_Z$ in $Y$. On the other hand, due to \cite[\href{https://stacks.math.columbia.edu/tag/0BXA}{Tag 0BXA}]{stacks-project}, the relative normalization of $Y_Z$ in $Y$ is $Y_Z$ itself, thus $Y_X = Y_Z$.
     \end{proof}

      \begin{lemma}{\cite[\href{https://stacks.math.columbia.edu/tag/03GV}{Tag 03GV}]{stacks-project}}\label{lemma: relative normalization and smooth base change}
        Relative normalizations of qcqs morphisms commute with smooth base change. In particular, it commutes with any open embedding.
    \end{lemma}

     In practice, when we apply relative normalizations to those schemes related with integral models of Shimura varieties, we will apply the following lemma:
  
  \begin{lemma}{\cite[\href{https://stacks.math.columbia.edu/tag/03GR}{Tag 03GR}]{stacks-project}}\label{lemma: nagata and finite morphism}
    Let $f: X \to S$ be a morphism, $S$ is Nagata (all locally of finite type schemes over $\Spec \OO_{E(v)}$ or $E$ are Nagata, see \cite[\href{https://stacks.math.columbia.edu/tag/035B}{Tag 035B}]{stacks-project}), $f$ is finite-type, and $X$ is reduced, then the relative normalization $S'\to S$ of $f$ is finite, in particular, separated.
  \end{lemma}

          \begin{lemma}\label{lemma: rela norm and normalization}
         Let $X \to Y$ be a qcqs morphism, and assume $X$ is normal. Let $Z$ be the relative normalization of $Y$ in $X$, then the usual normalization of $Z$ (\cite[\href{https://stacks.math.columbia.edu/tag/035E}{Tag 035E}]{stacks-project}) is $Z$ itself. In particular, $Z$ is normal.
     \end{lemma}
     \begin{proof}
         Let $Z'$ be the usual normalization of $Z$. Due to \cite[\href{https://stacks.math.columbia.edu/tag/0AXP}{Tag 0AXP}]{stacks-project}, every irreducible component of $X$ dominants an irreducible component of $Z$, then there exists a unique factorization $X \to Z' \to Z$ due to \cite[\href{https://stacks.math.columbia.edu/tag/035Q}{Tag 035Q}]{stacks-project}. On the other hand, $Z' \to Z$ and $Z \to Y$ are integral, thus there exists a unique factorization $X \to Z \to Z'$, which provides an inverse morphism to $Z' \to Z$, thus $Z'$ is canonically isomorphic to $Z$.
     \end{proof}

    \begin{lemma}\label{lemma: rela correct generic fiber}
         Let $U_1 \to U_2$ be an integral morphism (e.g. a finite morphism, a closed embedding), $U_2 \to X_2$ be an open embedding, $X_1$ be the relative normalization of $X_2$ in $U_1$, then $U_1 = X_1 \times_{X_2} U_2$.
     \end{lemma}
     \begin{proof}
         Consider the fiber products:
\[\begin{tikzcd}
	{U_1\times_{X_2}U_2 = U_1} & {X_1\times_{X_2}U_2} & {U_2} \\
	{U_1} & {X_1} & {X_2}
	\arrow[from=1-1, to=1-2]
	\arrow[from=1-1, to=2-1]
	\arrow[from=1-2, to=1-3]
	\arrow[from=1-2, to=2-2]
	\arrow[from=1-3, to=2-3]
	\arrow[from=2-1, to=2-2]
	\arrow[from=2-2, to=2-3]
\end{tikzcd}\]
          Since relative normalizations commute with open embeddings, then $X_1 \times_{X_2} U_2$ is the relative normalization of $U_2$ in $U_1$, which is $U_2$ itself due to \cite[\href{https://stacks.math.columbia.edu/tag/03GP}{Tag 03GP}]{stacks-project}.
     \end{proof}
     
     \begin{corollary}\label{cor: correct double rela}
          Let $U_1 \to U_2$ be an integral morphism , $U_2 \to Y_2$, $Y_2 \to X_2$ be open embeddings, $X_1$ (resp. $Y_1$) be the relative normalization of $X_2$ (resp. $Y_2$) in $U_1$, then $X_1$ is the relative normalization of $X_2$ in $Y_1$, and following diagrams are all cartesian.
\[\begin{tikzcd}
	{U_1} & {Y_1} & {X_1} \\
	{U_2} & {Y_2} & {X_2}
	\arrow[from=1-1, to=1-2]
	\arrow[from=1-1, to=2-1]
	\arrow[from=1-2, to=1-3]
	\arrow[from=1-2, to=2-2]
	\arrow[from=1-3, to=2-3]
	\arrow[from=2-1, to=2-2]
	\arrow[from=2-2, to=2-3]
\end{tikzcd}\]
     \end{corollary}
     \begin{proof}
         This follows from Lemma \ref{lemma: rela norm of rela norm is rela norm} and \ref{lemma: rela correct generic fiber}.
     \end{proof}
     \begin{remark}\label{rmk: can take rela norm interiorly}
         In practice, we will apply this corollary to diagrams:
\[\begin{tikzcd}
	{\shu{K}(G, X)} & {\shuc{K, \Sigma}(G, X) (\textit{resp.}\ \Shum{K}(G, X))} & {\Shumc{K, \Sigma}(G, X)} \\
	{\shu{K^{\dd}}(G^{\dd}, X^{\dd})} & {\shuc{K^{\dd}, \Sigma^{\dd}}(G^{\dd}, X^{\dd}) (\textit{resp.}\ \Shum{K^{\dd}}(G^{\dd}, X^{\dd}))} & {\Shumc{K^{\dd}, \Sigma^{\dd}}(G^{\dd}, X^{\dd})}
	\arrow[from=1-1, to=1-2]
	\arrow[from=1-1, to=2-1]
	\arrow[from=1-2, to=1-3]
	\arrow[from=1-2, to=2-2]
	\arrow[from=1-3, to=2-3]
	\arrow[from=2-1, to=2-2]
	\arrow[from=2-2, to=2-3]
\end{tikzcd}\]
     \end{remark}

     \begin{lemma}\label{lemma: finite morphism does not impact rela norm}
         Let $Z_1, Z_2$ be normal schemes locally of finite type and flat over $\Spec \OO_{E(v)}$, $X_1, X_2$ be their generic fibers over $E$ respectively. Assume we have a finite morphism $Z_1 \to Z_2$. Let $X \to X_1$ be a finite morphism over $E$, $X$ be normal, and $Z$ be the relative normalization of $Z_1$ in $X$, then $Z$ is the relative normalization of $Z_2$ in $X$, $Z \to Z_1$ is finite, $Z$ is normal, and has generic fiber $X$. We have the following diagram, where each square is Cartesian.
\[\begin{tikzcd}
	X & {X_1} & {X_2} \\
	Z & {Z_1} & {Z_2}
	\arrow[from=1-1, to=1-2]
	\arrow[from=1-1, to=2-1]
	\arrow[from=1-2, to=1-3]
	\arrow[from=1-2, to=2-2]
	\arrow[from=1-3, to=2-3]
	\arrow[from=2-1, to=2-2]
	\arrow[from=2-2, to=2-3]
\end{tikzcd}\]
     \end{lemma}
     \begin{proof}
         The right square is a Cartesian diagram by definition. Let $Z_1'$ be the relative normalization of $Z_2$ in $X_1$. By functoriality, there exists a morphism $\pi: Z_1' \to Z_1$. $\pi$ is a birational morphism since $Z_1'$ has generic fiber $X_1$ due to Lemma \ref{lemma: rela correct generic fiber}, $Z_1'$ is normal due to Lemma \ref{lemma: rela norm and normalization}, $\pi$ is a finite morphism since $Z_1' \to Z_2$ is a finite morphism due to Lemma \ref{lemma: nagata and finite morphism}. Therefore, $\pi$ is a birational finite morphism between normal noetherian schemes. In particular, $\pi$ is an isomorphism due to Zariski main theorem. By Lemma \ref{lemma: rela norm of rela norm is rela norm}, $Z$ is the relative normalization of $Z_2$ in $X$. For the same reasons, $Z \to Z_1$ is finite, $Z$ is normal, and has generic fiber $X$.
     \end{proof}

    \begin{lemma}\label{lem: rela norm finite surjection}
        Let $f_i: X_i \to S$ satisfy the conditions in Lemma \ref{lemma: nagata and finite morphism} and moreover $f_i$ is finite, for $i = 1, 2$, assume there exists a morphism $f: X_1 \to X_2$ over $S$ such that $f$ is an \'etale torsor, then the induced morphism $S_1' \to S_2'$ is finite and surjective.
    \end{lemma}
    \begin{proof}
        The morphism $S_1' \to S_2'$ comes from \cite[\href{https://stacks.math.columbia.edu/tag/035J}{Tag 035J}]{stacks-project}. Since both $S_1' \to S$ and $S_2' \to S$ are finite due to Lemma \ref{lemma: nagata and finite morphism}, then $S_1' \to S_2'$ is finite. \'Etale locally over $X_2$, there exist sections $X_2 \to X_1$ since $X_1 \to X_2$ is an \'etale torsor. We can refine such \'etale covering by an \'etale covering of $S$, see \cite[\href{https://stacks.math.columbia.edu/tag/04DQ}{Tag 04DQ}]{stacks-project}. Since relative normalizations commute with \'etale base change (see Lemma \ref{lemma: relative normalization and smooth base change}), then \'etale locally over $S$ we have sections $S_2' \to S_1'$, thus $S_1' \to S_2'$ is surjective.
    \end{proof}

\subsubsection{Kisin-Pappas-Zhou integral models}\label{subsec: kisin-pappas integral models}

\begin{enumerate}
    \item Under the setting \ref{general condition}, the good embedding $(G, X, K_p) \to (G^{\dd}, X^{\dd}, K^{\dd}_p)$ constructed in \cite{kisin2018integral} is an adjusted Siegel embedding, see subsection \ref{subsec: KP group embeddings}.
    \item Under the weaker setting \ref{general condition, KPZ}, the good embedding $(G, X, K_p) \to (G^{\dd}, X^{\dd}, K^{\dd}_p)$ constructed in \cite{kisin2024independence} and \cite{kisin2024integral} is an adjusted Siegel embedding, see subsection \ref{subsec: KPZ group embeddings}.
\end{enumerate}
   
    In either case, fix such an adjusted Siegel embedding $(G, X, K_p) \to (G^{\dd}, X^{\dd}, K^{\dd}_p)$. Let $K=K_pK^p\subset G(\A_f)$ (resp. $K^{\dd}=K^{\dd}_pK^{\dd, p} \subset \GSp(\A_f)$) be a neat open compact subgroup such that $K_p=\Ggh_x(\Z_p)$ and $K^p \subset G(\A_f^p)$ (resp. $K^{\dd}_p=\GSP_y(\Z_p)$ and $K^{\dd, p}\subset \GSp(\A_f^p)$). Let $\shu{K}:=\shu{K}(G, X)$ and $\shu{K^{\dd}}:=\shu{K^{\dd}}(\GSp, S^{\pm})$. Given $K^p\subset G(\A_f^p)$, By \cite[Lemma 2.1.2]{kisin2010integral}, we could choose a suitable neat $K^{\dd, p}\subset \GSp(\A_f^p)$ such that $\iota$ induces a closed embedding of Shimura varieties $\shu{K}\hookrightarrow\shu{K^{\dd}}\otimes E$ over the reflex field $E=E(G, X)$. 
      
      Recall that $\shu{K^{\dd}}$ has a smooth integral model $\Shum{K^{\dd}}$ over $\Spec \Z_{(p)}$ which parametrizes abelian schemes with principal polarizations and away-from-$p$ level structures. Let $v|p$ be a prime of $\OO_E$, by taking the relative normalization of $\Shum{K^{\dd}} \otimes_{\Z_{(p)}} \OO_{E(v)}$ in $\shu{K}$, we get an integral model $\Shum{K}$ defined over $\Spec \OO_{E(v)}$ which is normal due to Lemma \ref{lemma: rela norm and normalization}, together with a finite morphism $\Shum{K}\to\Shum{K^{\dd}}$ due to Lemma \ref{lemma: nagata and finite morphism}. Also, let us denote by $\Kf=\Kf_pK^p\subset G(\A_f)$ where $\Kf_p=\Ggf_x(\Z_p)$, and denote by $\Shum{\Kf}$ the relative normalization of $\Shum{K}$ in $\shu{\Kf}$. 

\begin{remark}\label{remark: y hyperspecial setting}
    In constructions in \cite{kisin2018integral}, \cite{kisin2024independence}, \cite{kisin2024integral}, and \cite{pappas2024p}, $K_p^{\dd}$ does not need to be hyperspecial, we can use any lattice $\Lambda \subset \Lambda^{\vee}$. We can even use any parahoric subgroup $K_p^{\dd}$ as in \cite[\S 4.5.2]{pappas2024p}. We choose $y$ to be hyperspecial in order to apply the setting \ref{def: compactification, setting} when we consider the toroidal compactifications. This simplifies a lot the modular interpretations on the boundary strata of Siegel Shimura varieties. This does not make difference for $\Shum{K}$, since the morphism between the integral models before and after applying Zarhin's trick is finite, see \cite[\S 4.5.9]{pera2012toroidal}, and then we apply Lemma \ref{lemma: finite morphism does not impact rela norm}.
\end{remark}

     \begin{remark}\label{remark: two normalizations equal}\leavevmode
      \begin{enumerate}
          \item In \cite{kisin2018integral}, \cite{kisin2024integral}, the authors first took $\Shum{K}^{-}$ as the integral closure of $\shu{K}$ in $\Shum{K^{\dd}} \otimes \OO_{E(v)}$ and then took $\Shum{K}$ as the usual normalization of $\Shum{K}^-$. This coincides with the relative normalization construction by the proof of Lemma \ref{lemma: rela norm and normalization}.
        \item When $\shu{K} \to \shu{K^{\dd}}$ is not a closed embedding ($K^{\dd, p}$ can be freely chosen, but need to contain $K^p$), we can still define $\Shum{K}$ as the relative normalization of $\Shum{K^{\dd}}$ in $\shu{K}$: Since $\shu{K} \to \shu{K^{\dd}}$ factors through $\shu{K} \hookrightarrow \shu{K^{\prime\dd}} \to \shu{K^{\dd}}$ for some $K^{\prime\dd}$, where $K^{\prime\dd}_p = K^{\dd}_p$, $K^{\prime\dd, p} \subset K^{\dd, p}$, $\shu{K} \hookrightarrow \shu{K^{\prime\dd}}$ is a closed embedding, and $\shu{K^{\prime\dd}} \to \shu{K^{\dd}}$ is a finite surjection. Apply Lemma \ref{lemma: finite morphism does not impact rela norm}, the relative normalization of $\Shum{K^{\dd}}$ in $\shu{K}$, is the relative normalization of $\Shum{K^{\prime\dd}}$ in $\shu{K}$.
        \item Similarly, when $K'_p = K_p \subset G(\Q_p)$, $K^{\prime, p} \subset K^p \subset G(\A_f^p)$, $\Shum{K'}$ can be either defined using via $\shu{K'} \hookrightarrow \Shum{K^{\prime\dd}}$, or via $\shu{K'} \to \shu{K} \to \Shum{K}$.
                 \item Due to Lemma \ref{lemma: rela norm of rela norm is rela norm}, $\Shum{\Kf}$ coincides with the relative normalization of $\Shum{K^{\dd}} \otimes \OO_{E(v)}$ in $\shu{\Kf}$.
      \end{enumerate}         
     \end{remark}

\subsubsection{Pappas-Rapoport integral model}\label{subsec: pappas-rapoport integral models}

Let $\iota: (G, X, K_p) \hookrightarrow (G^{\dd}, X^{\dd}, K^{\dd}_p)$ be a nice Siegel embedding. Due to \cite[Theorem 4.3.1, 4.5.2]{pappas2024p} (for $K_p = \Ggh_x(\Z_p)$ being stablizer quasi-parahoric subgroup) and \cite[Theorem 4.1.8, 4.1.12]{daniels2024conjecture} (for $K_p = \Ggf_x(\Z_p)$ being a general quasi-parahoric subgroup), the relative normalization $\Shum{K}(G, X)$ of $\Shum{K^{\dd}}(G^{\dd}, X^{\dd})$ in $\shu{K}(G, X)$ is the integral model in the sense of \cite[Conjecture 4.2.2]{pappas2024p}. In particular, $\Shum{K}(G, X)$ is independent of choice of Siegel embeddings due to \cite[Theorem 4.2.4]{pappas2024p} and \cite[Corollary 4.1.9]{daniels2024conjecture}.

\begin{remark}\label{rmk: not need to be Siegel embedding}
   Here $\iota$ does not need to be very nice. Moreover, from the quoted theorems, given a nice embedding $(G_1, X_1, K_{1, p}) \hookrightarrow (G_2, X_2, K_{2, p})$, if the Shimura varieties $\lrbracket{\shu{K_2}(G_2, X_2)}_{K_2^p}$ associated with $(G_2, X_2, K_{2, p})$ have canonical integral models $\lrbracket{\Shum{K_2}(G_2, X_2)}_{K_2^p}$, then the relative normalizations $\lrbracket{\Shum{K_1}(G_1, X_1)}_{K_1^p}$ of $\lrbracket{\Shum{K_2}(G_2, X_2)}_{K_2^p}$ in $\lrbracket{\shu{K_1}(G_1, X_1)}_{K_1^p}$ are also canonical integral models. In this case, we say $\lrbracket{\Shum{K_1}(G_1, X_1)}_{K_1^p}$ is defined via the nice embedding $(G_1, X_1, K_{1, p}) \hookrightarrow (G_2, X_2, K_{2, p})$.
\end{remark}
 
 Since an adjusted Siegel embedding is a nice embedding, Kisin-pappas integral models (resp. Kisin-Pappas-Zhou integral models) (if exists) are also Pappas-Rapoport integral models, they are independent of choice of Siegel embeddings. Kisin-pappas integral models (resp. Kisin-Pappas-Zhou integral models) support local model diagrams, which implies that \'etale locally such integral models can be determined by associated local models, using the deformation of $p$-divisible groups. This property can not be seen in Pappas-Rapoport construction.

\subsection{Compactification theory over integral base}\label{subsection: Compactification theory over integral base}

   We recall the constructions and results from \cite{pera2019toroidal}.
   \begin{definition}\label{def: compactification, setting}
         Recall the setting in \cite{pera2019toroidal}: $(G, X) \hookrightarrow (G^{\dd}, X^{\dd}) = (\GSp(V, \psi), S^{\pm})$ is a Siegel embedding, $K \subset K^{\dd}$, $\Lambda \subset V_{\Q}$ be a self-dual lattice, $K^{\dd}$ stablizes $\Lambda \otimes \hat{\Z}$ and $K^{\dd}_p$ is the stablizer subgroup of $\Lambda \otimes \Z_p$, in which case $K_p^{\dd}$ is a hyperspecial group. 
     \end{definition}

\subsubsection{Boundary tower}\label{subsubsec: lambda action}
   When we talk about mixed Shimura varieties at boundary, once we fix the level groups $K$, $K^{\ddagger}$ and the cusp label representatives $\Phi$ and $\Phi^{\dd}$, to save notation, we usually omit the indexes $K$, $K^{\ddagger}$, $g$, when they are clear from the context.  
   
   Under the setting \ref{def: compactification, setting}, one can formulate moduli interpretations of integral models
   \begin{equation}\label{equation: tower of GSp over integral base}
       \Xi^{\dd}:=\Shum{K^{\dd}_{\Phi^{\dd}, P^{\dd}}}(P^{\dd}, X_{\Phi^{\dd}})\to C^{\dd}:=\Shum{K^{\dd}_{\Phi^{\dd}, \ovl{P}^{\dd}}}(P^{\dd}, X_{\Phi^{\dd}, \ovl{P}})\to \Zb^{\bigsur, \dd}:=\Shum{K^{\dd}_{\Phi^{\dd}, h}}(G_h^{\dd}, X_{\Phi^{\dd}, h}^{\dd}),
   \end{equation}
   with generic fibers 
   \begin{equation}\label{equation: tower of GSp over generic fiber}
      \Xi^{\dd}_{\eta}:= \shu{K^{\dd}_{\Phi^{\dd}, P^{\dd}}}(P^{\dd}, X_{\Phi^{\dd}})\to C^{\dd}_{\eta}:=\shu{K^{\dd}_{\Phi^{\dd}, \ovl{P}^{\dd}}}(P^{\dd}, X_{\Phi^{\dd}, \ovl{P}})\to \Zb^{\bigsur, \dd}_{\eta}:=\shu{K^{\dd}_{\Phi^{\dd}, h}}(G_h^{\dd}, X_{\Phi^{\dd}, h}^{\dd}),
   \end{equation}
   respectively, see \cite[\S 2.2]{pera2019toroidal}.

   \begin{remark}\label{rmk: cusp labels coincide}
       In Siegel case, the cusp labels used in \cite{lan2016compactifications} coincide with the cusp labels used in generic fiber (\cite{lan2012comparison}), and coincide with the cusp labels used in \cite{pera2019toroidal}. This is a nontrivial fact, and is not true for general PEL type Shimura varieties. This correspondence was hiden behind the statements in \cite{pera2019toroidal}, especially in \cite[\S 2.2]{pera2019toroidal}. Please refer to \cite[Lemma A.4]{wu2025arith} for a more thorough explanation.
   \end{remark}

    The action of $\Lambda^{\ddagger}:=\Lambda_{\Phi^{\dd}, K^{\dd}}$ on the tower \ref{equation: tower of GSp over generic fiber} extends to the tower \ref{equation: tower of GSp over integral base}. The action of $\Lambda^{\ddagger}$ is compactible with the one in \cite[Errata \S 6.2.4]{lan2013arithmetic}, where $\mathrm{M}^{\Phi^{\ddagger}}=\Zb^{\bigsur, \ddagger}$, $\mathrm{M}^{\Zb^{\ddagger}}=\Zb^{\ddagger}$. If $K^{\ddagger}$ is chosen as some principle level subgroup $K^{\ddagger}(n)$, then $\Lambda^{\ddagger}$ acts on $\Zb^{\bigsur, \ddagger}$ trivially, for example, see Corollary \ref{cor: Lambda trivial action}.

    Recall we have the towers of mixed Shimura varieties over the reflex field $E=E(G, X)$ defined in \ref{eq: tower, over E} and \ref{equation: tower of GSp over generic fiber}. By functoriality, we have a finite morphism which respects the natural tower structures on each side:
    \[ \iota: \Xi_{\eta}:= \shu{K_{\Phi, P}}(P, X_{\Phi}) \to E\otimes_{\Q}\Xi^{\dd}_{\eta} = E\otimes_{\Q}\shu{K_{\Phi^{\dd}, P^{\dd}}}(P^{\dd}, X_{\Phi^{\dd}}).  \]     

     Take the relative normalization of
     \[ \OO_{E(v)}\otimes_{\Z_{(p)}}\Xi^{\ddagger}\to \OO_{E(v)}\otimes_{\Z_{(p)}}C^{\ddagger}\to \OO_{E(v)}\otimes_{\Z_{(p)}}\Zb^{\bigsur, \ddagger}\]
     in \ref{eq: tower, over E}, we get a tower of integral models over $\OO_{E(v)}$:
     \begin{equation}\label{equation: tower of G over integral base}
        \Xi:= \Shum{K_{\Phi, P}}(P, X_{\Phi}) \to C:=\Shum{K_{\Phi, \ovl{P}}}(\ovl{P}, X_{\Phi, \ovl{P}}) \to \Zb^{\bigsur}:=\Shum{K_{\Phi, h}}(G_h, X_{\Phi, h}).  
     \end{equation}
    
    Due to the functoriality of relative normalizations (\cite[\href{https://stacks.math.columbia.edu/tag/035J}{Tag 035J}]{stacks-project}), the $\Lambda := \Lambda_{\Phi, K}$-action on the tower \ref{eq: tower, over E} extends to the tower \ref{equation: tower of G over integral base}. Let us denote by $\Zb = \Zb^{\bigsur}/\Lambda$.

    Since $\Zb^{\bigsur}$ is quasi-projective, such a quotient exists (either as a geometric quotient, or a fppf quotient, or a categorical quotient, they are identical with each other in this case). The quotient map $\Zb^{\bigsur} \to \Zb$ is finite and surjective.

   \subsubsection{Twisted torus embedding}
   
     Fix a $\Phi \in [\Phi] \in \Cusp_K(G, X)$. Let $\sigma^{\dd} \subset U^{\dd}(\R)(-1) = \mathbf{B}_{K^{\dd}}(\Phi^{\dd}) \otimes \R$ be a rational polyhedral cone, we have a twisted torus embedding
     \begin{equation}
       \Xi^{\dd}:=  \Shum{K^{\dd}_{\Phi^{\dd}, P^{\dd}}}(P^{\dd}, X_{\Phi^{\dd}}) \to \Xi^{\dd}(\sigma^{\dd}):=\Shum{K^{\dd}_{\Phi^{\dd}, P^{\dd}}}(P^{\dd}, X_{\Phi^{\dd}}, \sigma^{\dd})
     \end{equation}
      Let $\Xi^{\dd}_{\sigma^{\dd}}:=\mathcal{Z}_{K^{\dd}_{\Phi^{\dd}, P^{\dd}}}(P^{\dd}, X_{\Phi^{\dd}}, \sigma^{\dd})$ be the associated closed subscheme of $\Xi^{\dd}(\sigma^{\dd})$. The morphisms and schemes are defined over $\Spec \Z_{(p)}$.

     Let $\sigma = \sigma^{\dd} \cap U(\R)(-1)$, and let $\Xi(\sigma):=\Shum{K_{\Phi, P}}(P, X_{\Phi}, \sigma)$ be the relative normalization of $ \Xi^{\dd}(\sigma^{\dd}) \otimes \OO_{E(v)}$ in $\Xi_{\eta}:=\shu{K_{\Phi, P}}(P, X_{\Phi})$, it has generic fiber $\Xi(\sigma)_{\eta}:=\shu{K_{\Phi, P}}(P, X_{\Phi}, \sigma)$. Let $\Xi_{\sigma}:=\mathcal{Z}_{K_{\Phi, P}}(P, X_{\Phi}, \sigma)$ be the relative normalization of $\Xi^{\dd}_{\sigma^{\dd}} \otimes \OO_{E(v)}$ in $\Xi_{\sigma, \eta}:= Z_{K_{\Phi, P}}(P, X_{\Phi}, \sigma)$, then we obtain a finite morphism $\Xi_{\sigma} \to \Xi(\sigma)$ which extends the morphism $\Xi_{\sigma, \eta} \to \Xi(\sigma)_{\eta}$ defined over $E$. 
     
     Fix rational polyhedral cone decompositions $\Sigma$, $\Sigma^{\dd}$ of $(G, X, K)$, $(G^{\dd}, X^{\dd}, K^{\dd})$ respectively. We say $\Sigma$ \emph{is induced by} $\Sigma^{\ddagger}$ if for each $\Phi$ and for each $\sigma \in \Sigma(\Phi)^{+}$, $\sigma$ is exactly the preimage of some $\sigma^{\ddagger} \in \Sigma^{\dd}(\Phi^{\dd})^{+}$ (and we denoted by $\sigma^{\dd} = \iota_*\sigma$). one can always refine the given $\Sigma$ and $\Sigma^{\ddagger}$ such that $\Sigma$ is induced by $\Sigma^{\ddagger}$. Even more, one can even request $(\Sigma, \Sigma^{\dd})$ to be strictly compatible in the sense of \cite[Definition 4.6]{lan2022closed} due to \cite[Proposition 4.10]{lan2022closed}.

     Fix a good enough $\Sigma^{\dd}$ (i.e. $\Sigma^{\dd}$ is smooth, complete, finite, admissable, and has no self-intersection). We do not study in details the properties of toroidal compactifications depending on the cone decompositions, thus we always assume the cone decomposition to be good enough. Let $\Sigma$ be induced by $\Sigma^{\dd}$. Let $\Shumc{K^{\dd}, \Sigma^{\dd}}$ be the Chai-Faltings compactification of $\Shum{K^{\dd}}$, and let $\Shumc{K, \Sigma}$ be the relative normalization of $\Shumc{K^{\dd}, \Sigma^{\dd}} \otimes \OO_{E(v)}$ in $\shuc{K, \Sigma}$. Let $\Upsilon = [(\Phi, \sigma)] \in \Cusp_K(\Sigma)$. Let $\Upsilon^{\dd} = [(\Phi^{\dd}, \sigma^{\dd})] = [(\iota_*\Phi, \iota_*\sigma)] \in \Cusp_{K^{\dd}}(\Sigma^{\dd})$, and $\mathcal{Z}_K(\Upsilon)$ be the relative normalization of $\mathcal{Z}_{K^{\dd}}(\Upsilon^{\dd}) \otimes \OO_{E(v)}$ in $Z_K(\Upsilon)$. Then there exists a canonical isomorphism $\Xi_{\sigma} = \mathcal{Z}_{K_{\Phi, P}}(P, X_{\Phi}, \sigma) \rightiso \mathcal{Z}_K(\Upsilon)$ extending the canonical isomorphism $\Xi_{\sigma, \eta} = Z_{K_{\Phi, P}}(P, X_{\Phi}, \sigma) \rightiso Z_K(\Upsilon)$, and $\mathcal{Z}_K(\Upsilon) \to \mathcal{Z}_{K^{\dd}}(\Upsilon^{\dd}) \otimes \OO_{E(v)} \to \Shumc{K^{\dd}, \Sigma^{\dd}} \otimes \OO_{E(v)}$ factors through $\mathcal{Z}_K(\Upsilon) \to \Shumc{K, \Sigma}$ which extends the locally closed embedding $Z_K(\Upsilon) \to \shuc{K, \Sigma}$.

  \subsubsection{Minimal compactifications}
     We recall the integral model $\Shumm{K}$ of the  minimal compactification $\shum{K}$. See \cite[\S 5]{pera2019toroidal} for details. Let $\omega(\Sigma)$ be the pullback of the Hodge bundle $\omega^{\dd}(\Sigma^{\dd})$ on $\Shumc{K^{\dd}, \Sigma^{\dd}}$, it is an ample bundle on $\Shumc{K, \Sigma}$, and is independent to the choice of $\Sigma$ and $\Sigma^{\dd}$. 
     
     $\bigoplus_{n\geq 0} H^0(\Shumc{K, \Sigma}, \omega(\Sigma)^n)$ is a finitely generated $\OO_{E(v)}$-algebra. Let $N$ be an integer large enough such that $\omega(\Sigma)^N$ is generated by its global sections over $\Shumc{K, \Sigma}$. Consider the projection:
     \begin{equation}
         \oint_{K, \Sigma}: \Shumc{K, \Sigma} \to \Shumm{K}:=
     \Proj(\bigoplus_{n\geq 0} H^0(\Shumc{K, \Sigma}, \omega(\Sigma)^n),
     \end{equation}
   it is the Stein factorization of the proper map:
   \begin{equation}\label{equation: stein factorization, minimal compactifications}
       \Shumc{K, \Sigma} \to H^0(\Shumc{K, \Sigma}, \omega(\Sigma)^N).
   \end{equation}
   The composition $\Xi_{\sigma} = \mathcal{Z}_{K_{\Phi, P}}(P, X_{\Phi}, \sigma) \to \Shumc{K, \Sigma} \stackrel{\oint_{K, \Sigma}}{\to} \Shumm{K}$
   factors through $\Lambda_{\Phi, K} \backslash \Zb^{\bigsur}(\Phi) \to \Shumm{K}$. Due to the Zariski main theorem, $\Lambda_{\Phi, K} \backslash\Zb^{\bigsur}(\Phi) \to \Shumm{K}$ is a locally closed embedding, and only depends on the class $[\Phi] \in \Cusp_K(G, X)$, see \cite[Proposition 5.2.7, Corollary 5.2.9]{pera2019toroidal}, thus we denote it by $\mathcal{Z}([\Phi]) \to \Shumm{K}$. Such defined $\Shumc{K, \Sigma}$, $\Shumm{K}$ are normal and projective over $S = \Spec \OO_{E(v)}$, and the complement of $\Shum{K}$ in $\Shumm{K}$ is also flat over $S$.

   \subsubsection{Summary}
   The main theorems in \cite{pera2019toroidal} could be summarized as in \cite[Proposition 2.1.2]{lan2018compactifications}, the following theorem is not the full list:

\begin{theorem}{{\cite[Theorem 4.1.5, Theorem 5.2.11]{pera2019toroidal}, \cite[Propersition 2.1.3]{lan2018compactifications}}}\label{theorem: axiomatic descriptions of compactifications}

    Fix $(\Sigma, \Sigma^{\dd})$ as above. There are proper surjective morphisms between noetherian schemes $\oint_{K, \Sigma}: \Shumc{K, \Sigma}\to\Shumm{K}$ over $S$ that are functorial with respect to $\Sigma$. Denote by $\imin{K}: \Shum{K} \to \Shumm{K}$, $\itor{K, \Sigma}: \Shum{K} \to \Shumc{K, \Sigma}$, then $\oint_{K, \Sigma}\circ\imin{K} = \itor{K, \Sigma}$.

    The compactification $\Shumc{K, \Sigma}$ (resp. $\Shumm{K}$) is stratified by locally closed subschemes $\mathcal{Z}(\Upsilon)$ (resp. $\mathcal{Z}([\Phi])$) with respect to $\Upsilon \in \Cusp_K(G, X, \Sigma)$ (resp. $[\Phi] \in \Cusp_K(G, X)$), and has an open dense subscheme $\Shum{K}$ when $\mathcal{Z}([\Phi]) = \Shum{K}$ and $\sigma = \lrbracket{0}$. Fix a cusp label representative $\Phi$ in each equivalence class $[\Phi]$, then those $\Upsilon = [\Phi, \sigma] \in \Cusp_K(\Sigma)$ over $[\Phi]$ are indexed by $[\sigma] \in \Lambda_{K, \Phi}\backslash\Sigma(\Phi)^+$.
    
    The preimage $\oint_{K, \Sigma}^{-1}(\mathcal{Z}([\Phi]))=\bigsqcup_{[\sigma] \in \Lambda_{K, \Phi}\backslash\Sigma(\Phi)^+} \mathcal{Z}([\Phi, \sigma])$ (set theoretically) and the closure of $\mathcal{Z}([\Phi, \sigma])$ (resp. $\mathcal{Z}([\Phi])$) is the union of $\mathcal{Z}([\Phi', \sigma'])$ (resp. $\mathcal{Z}([\Phi'])$) such that $[\Phi', \sigma'] \preceq [\Phi, \sigma]$ (resp. $[\Phi'] \preceq [\Phi]$).
   
   For each $\Zb^{\bigsur}(\Phi)$, there is a projective surjective morphism $C(\Phi)\to\Zb^{\bigsur}(\Phi)$. Also, there is a $\mathbf{E}_K(\Phi)$-torsor $\Xi(\Phi)\to C(\Phi)$ where $\mathbf{E}_K(\Phi)$ is the pullback of a split torus over $\Spec \Z$. For each $\sigma\in\Sigma(\Phi)$, there is a toroidal embedding $\Xi(\Phi)\hookrightarrow\Xi(\Phi, \sigma)$ and a specified closed subscheme $\Xi_{\sigma}(\Phi)$ of $\Xi(\Phi,\sigma)$. $\Xi(\Phi, \sigma)$ has a stratification by $\{\Xi_{\tau}(\Phi)\}_{\tau}$ where $\tau$ runs over the set of all faces of $\sigma$. For each $\sigma\in\Sigma(\Phi)^+$, there is a canonical isomorphism $\Xi_{\sigma}(\Phi)\cong\mathcal{Z}([\Phi, \sigma])$. The following diagram is commutative:
\begin{equation}\label{diag: strata in two sides}
    \begin{tikzcd}
	{\mathcal{Z}([\Phi, \sigma])} &&& {\mathcal{Z}([\Phi])} \\
	{\Xi_{\sigma}(\Phi)} & {C(\Phi)} & {\Zb^{\bigsur}(\Phi)} & {\Zb(\Phi)}
	\arrow["{\oint_{K, \Sigma}}", from=1-1, to=1-4]
	\arrow[from=2-1, to=2-2]
	\arrow[from=2-2, to=2-3]
	\arrow[from=2-3, to=2-4]
	\arrow["\cong", from=2-1, to=1-1]
	\arrow["\cong"', from=2-4, to=1-4]
\end{tikzcd}
\end{equation}
   
        The local structures of toroidal compactifications along boundary strata $\mathcal{Z}([\Phi, \sigma])$ could be described as follows:

     	Fix a $\Phi$ from now on. We omit $\Phi$ in the subscript to save notions. Consider the full toroidal embedding and its formal completion:
        \begin{equation}\label{eq: general compactification, Xi}
         \overline{\Xi}=\underset{\sigma\in\Sigma(\Phi)^+}{\bigcup} \Xi( \sigma), \quad  \mathfrak{X}=(\overline{\Xi})^{\wedge}_{\underset{\tau\in\Sigma(\Phi)^{+}}{\cup}\Xi_{\tau}} 
        \end{equation}
     	Note that $\mathfrak{X}$ has open coverings $\mathfrak{X}_{\sigma}^{\circ}$ for $\sigma$ running over all $\sigma\in\Sigma(\Phi)^{+}$, where
        \begin{equation}\label{eq: mathfrak X, definition}
             \mathfrak{X}_{\sigma}^{\circ}=(\Xi(\sigma))_{\underset{\tau\in\Sigma(\Phi)^{+},\ \bar{\tau}\subset\bar{\sigma}}{\cup}\Xi_{\tau}}^{\wedge}
        \end{equation}
     	The composition $\mathfrak{X}_{\sigma}^{\circ}\to\mathfrak{X}\to\Shumc{K, \Sigma}$ induces an isomorphism
       \begin{equation}\label{eq: mathfrak X, identification}
           \mathfrak{X}_{\sigma}^{\circ}\cong(\Shumc{K, \Sigma})_{\underset{\tau\in \Sigma(\Phi)^{+},\ \bar{\tau}\subset\bar{\sigma}}{\cup} \mathcal{Z}([\Phi, \tau])}^{\wedge} 
       \end{equation}
     	which extends the canonical boundary isomorphism $\Xi_{\sigma}\cong\mathcal{Z}([\Phi, \sigma])$.
     \end{theorem}

\begin{remark}\label{rmk: morphism between toroidal, bad}
    It is worth noting that in the construction, $K_p$ does not need to be (quasi-)parahoric, and $K_p \subset K^{\dd}_p \cap G(\Q_p)$ is the only restriction for the level $K_p$ at prime $p$. In particular, we do not even require $(G, X, K_p) \to (G^{\dd}, X^{\dd}, K_p^{\dd})$ to be a nice morphism. For different $K, K' \subset G(\A_f)$ and $g \in G(\A_f)$ such that $gK'g^{-1} \subset K$, let $\Sigma'$ refines $g^*\Sigma$, then the Hecke action $\shuc{K', \Sigma'}(G, X) \to \shuc{K, \Sigma}(G, X)$ extends to $\Shumc{K', \Sigma'}(G, X) \to \Shumc{K, \Sigma}(G, X)$ due to \cite[Proposition 4.1.13]{pera2019toroidal}. Here both $\Shumc{K', \Sigma'}(G, X)$ and $\Shumc{K, \Sigma}(G, X)$ are constructed as relative normalizations of same $\Shumc{K^{\dd}, \Sigma^{\dd}}(G^{\dd}, X^{\dd})$. Since both $\Shumc{K, \Sigma}(G, X)$ and $\Shumc{K', \Sigma'}(G, X)$ are proper over $\OO_{E(v)}$, $\Shumc{K', \Sigma'}(G, X) \to \Shumc{K, \Sigma}(G, X)$ is also proper. Indeed, this is a finite morphism by using Lemma \ref{lemma: nagata and finite morphism}, since both $\Shumc{K', \Sigma'}(G, X)$ and $\Shumc{K, \Sigma}(G, X)$ are finite over $\Shumc{K^{\dd}, \Sigma^{\dd}}(G^{\dd}, X^{\dd})$. 
\end{remark}

   \begin{remark}\label{rmk: morphism between minimal, bad}
          In particular, given $gK'g^{-1} \subset K \subset G(\A_f)$ as usual, the Hecke action $\shum{K'}(G, X) \to \shum{K}(G, X)$ extends to $\Shumm{K'}(G, X) \to \Shumm{K}(G, X)$ due to \cite[Proposition 5.2.13]{pera2019toroidal}. Since both $\Shumm{K}(G, X)$ and $\Shumm{K'}(G, X)$ are proper over $\OO_{E(v)}$, $\Shumm{K'}(G, X) \to \Shumm{K}(G, X)$ is also proper. Indeed, this is a finite morphism, see Remark \ref{rmk: morphism between toroidal, bad}.
   \end{remark}

    Of course, in many cases, we do not want to use same $\Shumc{K^{\dd}, \Sigma^{\dd}}(G^{\dd}, X^{\dd})$ in the construction, then even the existence of such morphism $\Shumc{K', \Sigma'}(G, X) \to \Shumc{K, \Sigma}(G, X)$ becomes a question, see subsection \ref{subsec: change of parahoric}.

   \begin{proposition}\label{prop: minimal compactification, rela normalization}
       $\Shumm{K}(G, X)$ has generic fiber $\shum{K}(G, X)$, and is the relative normalization of $\Shumm{K^{\dd}}(G^{\dd}, X^{\dd}) \otimes \OO_{E(v)}$ in $\shum{K}(G, X)$.
   \end{proposition}
   \begin{proof}
       Due to \cite[Theorem 5.1.1]{lan2012comparison}, $\Shumm{K^{\dd}}(G^{\dd}, X^{\dd})$ has generic fiber $\shum{K^{\dd}}(G^{\dd}, X^{\dd})$. Using the same comparison on the pullback of Hodge bundles (see the proof of \cite[Theorem 5.1.1]{lan2012comparison}), $\Shumm{K}(G, X)$ has generic fiber $\shum{K}(G, X)$. This can also be seen from \cite[Theorem 5.2.11(3)]{pera2019toroidal}. 
       
       Apply \cite[Proposition 3.4 case (4)]{lan2022closed}, $\Shumm{K}(G, X) \to \Shumm{K^{\dd}}(G^{\dd}, X^{\dd})$ is a finite morphism.
       
       Let $S$ be the relative normalization of $\Shumm{K^{\dd}}(G^{\dd}, X^{\dd}) \otimes \OO_{E(v)}$ in $\shum{K}(G, X)$, due to \cite[\href{https://stacks.math.columbia.edu/tag/035I}{Tag 035I}]{stacks-project} and \cite[\href{https://stacks.math.columbia.edu/tag/03GR}{Tag 03GR}]{stacks-project}, there exists a finite morphism $S \to \Shumm{K}(G, X)$. Since it is a birational, finite morphism between normal noetherian schemes, it is an isomorphism by Zariski main theorem.
   \end{proof}

    \subsubsection{Lambda action at quasi-parahoric levels}

    Recall that in subsection \ref{subsubsec: lambda action}, we have a finite surjection $\Zb^{\bigsur} \to \Zb$ given by the quotient of $\Lambda$.
    \begin{proposition}{\cite[Corollary 4.27]{wu2025arith}}\label{prop: Zb finite etale}
        The action of $\Lambda$ on $\Zb^{\bigsur}$ is \'etale. In particular, $\Zb^{\bigsur} \to \Zb$ is finite \'etale.
    \end{proposition}
    Note that in \cite[Corollary 4.27]{wu2025arith}, he proved the quotient by $\Lambda$ is \'etale in abelian-type case and for each level subgroup $K_p$, using the theory of twisting of $1$-motif. 

\subsubsection{Abelian schemes torsor}

    In subsection \ref{sec: a detailed look into this tower}, we recall that $\shu{K_{\Phi, \ovl{P}}}(\ovl{P}, X_{\Phi, \ovl{P}}) \to \shu{K_{\Phi, h}}(G_h, X_{\Phi, h})$ is a torsor under the abelian scheme $A_K(\Phi) = \shu{K_{\Phi, V}\rtimes K_{\Phi, h}}(\ovl{P}, X_{\Phi, \ovl{P}})$ ($K_{\Phi, W} = \wdt{K}_{\Phi, W}$ in Hodge-type case). Consider the weight filtrations on $\Lie(P)$ and $\Lie(P^{\ddagger})$, we see that $U \to W \to W^{\ddagger}$ factors through $U^{\ddagger}$, thus $W \cap U^{\ddagger} = U$, $V = W/U \to V^{\ddagger}=W^{\ddagger}/U^{\ddagger}$ is an embedding. The natural morphism $K_{\Phi, W} \to K_{\Phi^{\dd}, W^{\dd}} \to K_{\Phi^{\dd}, V^{\dd}}$ factors through $K_{\Phi, V} \to K_{\Phi^{\dd}, V^{\dd}}$, which induces a morphism of abelian schemes $A_K(\Phi) \to A_{K^{\dd}}(\Phi^{\dd})|_{\shu{K_{\Phi, h}}(G_h, X_{\Phi, h})}$. Since $A_{K^{\dd}}(\Phi^{\dd})$ extends to an abelian scheme $\ab_{K^{\dd}}(\Phi^{\dd})$  over the integral base $\Spec \Z_{(p)}$, we take the relative normalization of $\ab_{K^{\dd}}(\Phi^{\dd})|_{\Shum{K_{\Phi, h}}(G_h, X_{\Phi, h})}$ in $A_K(\Phi)$, and denote it by $\ab_K(\Phi)$.
    \begin{proposition}{{\cite[Proposition 4.2.2]{pera2019toroidal}}}\label{prop: Pera 4.2.2}
       When $K_{\Phi, V, p} = K^{\dd}_{\Phi^{\dd}, V^{\dd}, p} \cap V(\Q_p)$, then $\ab_K(\Phi) \to \Shum{K_{\Phi, h}}(G_h, X_{\Phi, h})$ is an abelian scheme.
    \end{proposition}

\begin{corollary}\label{cor: pera, 4.1.5}
    Assume $(G, X, K_p) \to (G^{\dd}, X^{\dd}, K^{\dd}_p)$ is an adjusted Siegel embedding, $K_p$ being a general quasi-parahoric subgroup, then for each $[\Phi] \in \Cusp_K(G, X)$, then $K^{\ddagger}_{\Phi^{\dd}, V^{\ddagger}, p} \cap V(\Q_p) = K_{\Phi, V, p}$, and $\Shum{K_{\Phi, \ovl{P}}}(\ovl{P}, X_{\Phi, \ovl{P}}) \to \Shum{K_{\Phi, h}}(G_h, X_{\Phi, h})$ is a torsor under an abelian scheme $\ab_K(\Phi) \to \Shum{K_{\Phi, h}}(G_h, X_{\Phi, h})$. The same statement holds true when we replace $K_p$ by $K_p(n)$.
\end{corollary}
\begin{proof}
      This is \cite[Theorem 4.1.5(1)]{pera2019toroidal}. We verify the condition $K^{\ddagger}_{\Phi^{\dd}, V^{\ddagger}, p} \cap V(\Q_p) = K_{\Phi, V, p}$ for each $[\Phi] \in \Cusp_K(G, X)$ using the calculation in subsection \ref{subsubsec: about unipotent groups}, especially Corollary \ref{corollary: equality of cpt groups on V} and Lemma \ref{lemma: Kc and K have same W part}. 
\end{proof}

\subsubsection{Fiber products at boundary}

  Fix a stabalizer quasi-parahoric subgroup $K_p$. For each boundary component $\Zb \subset \Shumm{K}$, we have tower $\Xi_K\to C_K\to \Zb^{\bigsur}_K \to \Zb_K$ (\ref{equation: tower of G over integral base}) over $\Spec \OO_{E(v)}$. When we replace $K$ with a general quasi-parahoric subgroup $\Kf$ (and use same Siegel embedding and $K^{\dd}\subset G(\A_f)$), we similarly get a tower $\Xi_{\Kf}\to C_{\Kf}\to \Zb_{\Kf}^{\bigsur} \to \Zb_{\Kf}$. We can say more about the relations between $C_K \to \Zb^{\bigsur}_K$ and $C_{\Kf}\to \Zb^{\bigsur}_{\Kf}$. Consider the following commutative graph:
     \begin{equation}\label{diagram: commutative diagram for abelian scheme torsor}
\begin{tikzcd}[sep=small]
	{\shu{\Kf_{\Phi, \ovl{P}}}(\ovl{P}, X_{\Phi, \ovl{P}}) } && {C_{\Kf}} \\
	& {\shu{K_{\Phi, \ovl{P}}}(\ovl{P}, X_{\Phi, \ovl{P}}) } &&& {C_K} \\
	{\shu{\Kf_{\Phi, h}}(G_h, X_{\Phi, h})} && {\Zb_{\Kf}^{\bigsur}} \\
	& {\shu{K_{\Phi, h}}(G_h, X_{\Phi, h})} &&& {\Zb_K^{\bigsur}}
	\arrow[from=1-1, to=1-3]
	\arrow[from=1-1, to=2-2]
	\arrow[from=1-1, to=3-1]
	\arrow[from=1-3, to=2-5]
	\arrow[from=1-3, to=3-3]
	\arrow[from=2-2, to=2-5]
	\arrow[from=2-2, to=4-2]
	\arrow[from=2-5, to=4-5]
	\arrow[from=3-1, to=3-3]
	\arrow[from=3-1, to=4-2]
	\arrow[from=3-3, to=4-5]
	\arrow[from=4-2, to=4-5]
\end{tikzcd}
 \end{equation}
     \begin{proposition}\label{prop: cartesian products}
     	The leftmost and the rightmost squares of the commutative graph \eqref{diagram: commutative diagram for abelian scheme torsor} are cartesian products. The same is true when we replace $K_p$ and $\Kf_p$ by $K_p(n)$ and $\Kf_p(n):= K_p(n) \cap \Kf_p$ respectively.
     \end{proposition}
     \begin{proof}
     Let $C$ be the cartesian product of $C_{K}\rightarrow\Zb_{K}^{\bigsur}\leftarrow \Zb_{\Kf}^{\bigsur}$. Since the relative normalizations commute with base change along open embeddings \cite[\href{https://stacks.math.columbia.edu/tag/035K}{Tag 035K}]{stacks-project}, and the relative normalization of $X$ in $Y$ along a finite morphism $Y\to X$ is $Y$ \cite[\href{https://stacks.math.columbia.edu/tag/03GP}{Tag 03GP}]{stacks-project}, the leftmost square of \eqref{diagram: commutative diagram for abelian scheme torsor} is the fiber of the rightmost square over $E=E(G,X)$ and the generic fiber $C_{E}$ of $C$ is the cartesian product of 
      \[\shu{\Kf_{\Phi, h}}(G_h, X_{\Phi, h}) \rightarrow \shu{K_{\Phi, h}}(G_h, X_{\Phi, h}) \leftarrow \shu{K_{\Phi, \ovl{P}}}(\ovl{P}, X_{\Phi, \ovl{P}}).\]
     	
     First, the leftmost square of \eqref{diagram: commutative diagram for abelian scheme torsor} is a cartesian product, i.e. $\shu{\Kf_{\Phi, \ovl{P}}}(\ovl{P}, X_{\Phi, \ovl{P}}) \to C_{E}$ is an isomorphism. It follows from \cite[3.10, 3.11]{pink1989arithmetical} (first base change over $\CC$), the conditions there could be verified in our case:
     	\begin{enumerate}
     		\item Lemma \ref{lemma: Kc and K have same W part} says $\Kf_{\Phi, V, p} = K_{\Phi, V, p}$, thus 
       \[ \Kf_{\Phi, h, p}\times_{K_{\Phi, h, p}}K_{\Phi, \ovl{P}, p} = (\pi_{\ovl{P}, G_h}|_{K_{\Phi, \ovl{P}, p}})^{-1}(\Kf_{\Phi, h, p}) = \Kf_{\Phi, \ovl{P}, p}\]
     		\item The center of $G_h$ is an almost product of a $\Q$-split torus with a torus of compact type defined over $\Q$.
    	\end{enumerate}

      Second, consider the morphism $C_{\Kf}\to C_K$. Due to Corollary \ref{cor: pera, 4.1.5}, $C_{K}\to\Zb_{K}^{\bigsur}$ is an abelian scheme torsor, then $C\to\Zb^{\bigsur}_{\Kf}$ is also an abelian scheme torsor over a normal base $\Zb^{\bigsur}_{\Kf}$, then $C$ is normal. Due to Corollary \ref{lem: rela norm finite surjection}, $C_{\Kf} \to C_K$ is finite and surjective. Similarly, $\Zb_{\Kf}^{\bigsur} \to \Zb_{K}^{\bigsur}$ is finite and surjective, thus $C \to C_K$ is finite and surjective. Then $C_{\Kf} \to C_K$ is finite. Use Zariski main theorem, $C_{\Kf} \to C_K$ is finite and birational between noetherian normal schemes, then it is an isomorphism.
     \end{proof}

\subsection{Summary}\label{subsec: summary}

We have seen that the boundary strata of an integral model $\Shumm{K}(G, X)$ with quasi-parahoric level are again finite quotients of integral models with quasi-parahoric levels. We summarize the results from this section as follows.

\begin{theorem}\leavevmode\label{thm: boundary of kp is again kp}

\begin{enumerate}
    \item (over generic fiber, classical results) Let $(G, X)$ be a Hodge-type Shimura datum, $E = E(G, X)$, $\iota: (G, X) \to (G^{\dd}, X^{\dd})$ be a Siegel embedding, $K \subset G(\A_f)$, $K^{\dd} \subset \GSp(\A_f)$ be neat open compact subgroups such that $K \subset K^{\dd}$. Then $\iota$ induces a finite morphism $\shu{K} \to \shu{K^{\dd}} \otimes E$ and extends to a finite morphism between their minimal compactifications $\iota^{\min}: \shum{K} \to \shum{K^{\dd}} \otimes E$. The boundaries of $\shum{K}$ and of $\shum{K^{\dd}}$ are stratified by 
    \begin{equation*}
        \bigsqcup_{[\Phi]  \in \Cusp_K(G, X)} \Lambda_{\Phi, K}\backslash\shu{K_{\Phi, h}}(G_h, X_{\Phi, h}),\quad \bigsqcup_{[\Phi^{\dd}]  \in \Cusp_{K^{\dd}}(G^{\dd}, X^{\dd})} \Lambda_{\Phi^{\dd}, K^{\dd}}\backslash\shu{K^{\dd}_{\Phi^{\dd}, h}}(G^{\dd}_h, X_{\Phi^{\dd}, h}^{\dd}),
    \end{equation*}
    respectively. The action of $\Lambda_{\Phi, K}$ on $\shu{K_{\Phi, h}}$ and the action of $\Lambda_{\Phi^{\dd}, K^{\dd}}$ on $\shu{K^{\dd}_{\Phi^{\dd}, h}}$ factors through finite quotients.
    
    For each cusp label representative $\Phi$ in $[\Phi] \in \Cusp_K(G, X)$, let $\Phi^{\dd} = \iota_*\Phi$, then $\iota$ induces a Siegel embedding $\iota_{\Phi}: (G_h, X_{\Phi, h}) \to (G_h^{\dd}, X_{\Phi^{\dd}, h}^{\dd})$, $(G_h, X_{\Phi, h})$ is a Hodge-type Shimura datum. Note that $K_{\Phi, h} \subset K^{\dd}_{\Phi^{\dd}, h}$, $\iota_{\Phi}$ induces a finite morphism
    \begin{equation*}
        \shu{K_{\Phi, h}}:= \shu{K_{\Phi, h}}(G_h, X_{\Phi, h}) \to \shu{K^{\dd}_{\Phi^{\dd}, h}} \otimes E :=\shu{K^{\dd}_{\Phi^{\dd}, h}}(G^{\dd}_h, X_{\Phi^{\dd}, h}^{\dd}) \otimes E,
    \end{equation*}
    modulo the $\Lambda$-actions, it is compatible with the morphism induced by $\iota^{\min}$.

    \item (over integral base, \cite{pera2019toroidal}) We moreover assume $(G, X, K_p) \to (G^{\dd}, X^{\dd}, K_p^{\dd})$ is an adjusted Siegel embedding (Definition \ref{def: adjusted embedding}). Let $\Shum{K^{\dd}}$ be the smooth integral model, $\Shum{K}$ be the relative normalization of $\Shum{K^{\dd}} \otimes \OO_{E(v)}$ (we drop the $\otimes \OO_{E(v)}$-part in the following) in $\shu{K}$, here we do not need $\iota$ to be a closed embedding (Remark \ref{remark: two normalizations equal}). Then $\iota^{\min}$ extends to a finite morphism $\Shumm{K} \to \Shumm{K^{\dd}}$. The boundaries of $\Shumm{K}$ and of $\Shumm{K^{\dd}}$ are stratified by 
    \begin{equation*}
        \bigsqcup_{[\Phi]  \in \Cusp_K(G, X)} \Lambda_{\Phi, K}\backslash\Shum{K_{\Phi, h}}(G_h, X_{\Phi, h}),\quad \bigsqcup_{[\Phi^{\dd}]  \in \Cusp_{K^{\dd}}(G^{\dd}, X^{\dd})} \Lambda_{\Phi^{\dd}, K^{\dd}}\backslash\Shum{K^{\dd}_{\Phi^{\dd}, h}}(G^{\dd}_h, X_{\Phi^{\dd}, h}^{\dd}),
    \end{equation*}
    respectively. Here $K^{\dd}_{\Phi^{\dd}, h, p}$ is again hyperspecial (Proposition \ref{prop: main prop for section BT theory}), $\Shum{K^{\dd}_{\Phi^{\dd}, h}}$ is the smooth integral model similarly defined as $\Shum{K^{\dd}}$ (\cite[\S 4.1.1]{pera2019toroidal}), $\Shum{K_{\Phi, h}}$ is the relative normalization of $\Shum{K^{\dd}_{\Phi^{\dd}, h}}$ in $\shu{K_{\Phi, h}}$ induced by $\iota_{\Phi}$. The quotients $\Shum{K_{\Phi, h}}(G_h, X_{\Phi, h}) \to \Lambda_{\Phi, K}\backslash\Shum{K_{\Phi, h}}(G_h, X_{\Phi, h})$ are finite \'etale (Proposition \ref{prop: Zb finite etale}), and $\Shum{K_{\Phi, \ovl{P}}}(\ovl{P}, X_{\Phi, \ovl{P}}) \to \Shum{K_{\Phi, h}}(G_h, X_{\Phi, h})$ are torsors under abelian schemes (Corollary \ref{cor: pera, 4.1.5}).
    \item (level groups) If $(G, X, K_p) \to (G^{\dd}, X^{\dd}, K_p^{\dd})$ is an adjusted Siegel embedding (Definition \ref{def: adjusted embedding}) (resp. very nice embedding (Definition \ref{def: nice embedding})), then $(G_{\Phi, h}, X_{\Phi, h}, K_{\Phi, h}) \to (G_{\Phi^{\dd}, h}^{\dd}, X_{\Phi^{\dd}, h}^{\dd}, K^{\dd}_{\Phi^{\dd}, h})$ is an adjusted Siegel embedding (Lemma \ref{lem: boundary adjusted embedding}) (resp. very nice embedding (Corollary \ref{cor: very nice emb implies very nice emb})). In group theory aspect, let $x \in \Bui_{\ext}(G, \Q_p)$, $y \in \Bui_{\ext}(G^{\dd}, \Q_p)$, $K_p := \Ggh_x(\Z_p)$, $K_p^{\dd} = \Ggh_{y}(\Z_p)$, $K_p = G(\Q_p) \cap K^{\dd}_p$ (resp. $y$ is the image of $x$ under some $\iota_{G}: G(\bQ)$-$\Gal(\bQ|\Q_p)$-equivariant embedding $\Bui_{\ext}(G, \bQ) \to \Bui_{\ext}(G^{\dd}, \bQ)$). Then $K_{\Phi, h, p} = \Ggh_{x_{g, h}}(\Z_p)$ for some $x_{g, h} \in B_{\ext}(G_h, \Q_p)$ determined by $x$ in Proposition \ref{prop: main prop for section BT theory}, and similarly $K_{\Phi^{\dd}, h, p}^{\dd} = \Ggh_{y_{g, h}}(\Z_p)$ for some $y_{g, h} \in B_{\ext}(G_h^{\dd}, \Q_p)$ determined by $y$, $y_{g, h}$ is a hyperspecial point, and $K_{\Phi, h, p} = K_{\Phi^{\dd}, h, p}^{\dd} \cap G_h(\Q_p)$ (resp. $y_{g, h}$ is the image of $x_{g, h}$ under some $\iota_{G_h}: G_{h}(\bQ)$-$\Gal(\bQ|\Q_p)$-equivariant embedding $\Bui_{\ext}(G_{h}, \bQ) \to \Bui_{\ext}(G_{h}^{\dd}, \bQ)$ compatible with $\iota_G$). Moreover, Proposition \ref{prop: main prop for section BT theory} shows that if $\Kn_p = \Gg_x(\Z_p)$, $\Kc_p = \Ggc_x(\Z_p)$, then $(\Kn_p)_{\Phi, h} = \Gg_{x_{g, h}}(\Z_p)$,  $\Ggc_{x_{g, h}}(\Z_p) \subset (\Kc_p)_{\Phi, h} \subset \Gg_{x_{g, h}}(\Z_p)$. In particular, if $K_p$ is a quasi-parahoric subgroup, then $K_{\Phi, h, p}$ is a quasi-parahoric subgroup.

    \item (properties inherited by boundary, subsection \ref{subsec: general facts}) If $x$ is a special (resp. hyperspecial, absolutely special, superspecial) point, then $x_{g, h}$ is also a special (resp. hyperspecial, absolutely special, superspecial) point. If $G$ satisfies the property $(\ast)$, then $G_h$ also satisfies the property $(\ast)$, here $(\ast)$ can be any of the followings:
    \begin{itemize}
        \item Is quasi-split over $\Q_p$.
        \item Splits over a fixed field extension $F_1$ of $\Q_p$.
        \item Is residually split over $\Q_p$.
        \item Satisfies the condition \ref{general condition} posed on $(p, G)$.
    \end{itemize}

    \item (Pappas-Rapoport integral models, subsection \ref{subsec: pappas-rapoport integral models}) Let $(G, X, K_p) \to (G^{\dd}, X^{\dd}, K_p^{\dd})$ be an adjusted Siegel embedding, $\Shum{K}(G, X)$ is the Pappas-Rapoport integral model \cite[\S 4.5]{pappas2024p} defined via the nice embedding $(G, X, K_p) \to (G^{\dd}, X^{\dd}, K_p^{\dd})$ (Remark \ref{rmk: not need to be Siegel embedding}), then $\Shum{K_{\Phi, h}}(G_h, X_{\Phi, h})$ are Pappas-Rapoport integral models for all $[\Phi] \in \Cusp_K(G, X)$ defined via the adjusted Siegel embeddings $(G_{\Phi, h}, X_{\Phi, h}, K_{\Phi, h, p}) \to (G_{\Phi^{\dd}, h}^{\dd}, X_{\Phi^{\dd}, h}^{\dd}, K^{\dd}_{\Phi^{\dd}, h, p})$ (Lemma \ref{lem: boundary adjusted embedding}).

    \item (Kisin-Pappas(-Zhou) integral models, subsection \ref{subsec: kisin-pappas integral models}) Assume condition \ref{general condition} (resp. condition \ref{general condition, KPZ}). Let $(G, X, K_p) \to (G^{\dd}, X^{\dd}, K_p^{\dd})$ be an adjusted Siegel embedding that is a good embedding, which exists, see subsection \ref{subsec: KP group embeddings} (resp. \ref{subsec: KPZ group embeddings}). Let $\Shum{K}(G, X)$ be the Kisin-Pappas (resp. Kisin-Pappas-Zhou) integral model defined via the good embedding $(G, X, K_p) \to (G^{\dd}, X^{\dd}, K_p^{\dd})$. Then $\Shum{K_{\Phi, h}}(G_h, X_{\Phi, h})$ are Kisin-Pappas integral (resp. Kisin-Pappas-Zhou) models defined via the good embeddings $(G_{\Phi, h}, X_{\Phi, h}, K_{\Phi, h, p}) \to (G_{\Phi^{\dd}, h}^{\dd}, X_{\Phi^{\dd}, h}^{\dd}, K^{\dd}_{\Phi^{\dd}, h, p})$ for all $[\Phi] \in \Cusp_K(G, X)$ (resp. under the extra assumption that $G_{\Phi, h}$ are $R$-smooth for all $[\Phi]$). (Proposition \ref{prop: good embedding}).
\end{enumerate}
\end{theorem}

\begin{remark}
    We expect if the adjusted Siegel embedding $(G, X, K_p) \to (G^{\dd}, X^{\dd}, K_p^{\dd})$ is a very good embedding, then the adjusted Sigel embeddings $(G_{\Phi, h}, X_{\Phi, h}, K_{\Phi, h, p}) \to (G_{\Phi^{\dd}, h}^{\dd}, X_{\Phi^{\dd}, h}^{\dd}, K^{\dd}_{\Phi^{\dd}, h, p})$ are all very good embeddings for all $[\Phi] \in \Cusp_K(G, X)$.
\end{remark}

\section{Uniqueness of compactifications and change-of-parahorics}\label{sec: uniqueness of compactifications}

\subsection{Uniqueness of toroidal and minimal compactifications}

In this section, we extend the uniqueness result of integral models to toroidal compactifications. Let $\Shumc{K, \Sigma}(G, X)_1$ and $\Shumc{K, \Sigma}(G, X)_2$ be two integral models of $\shuc{K, \Sigma}(G, X)$ defined via relative normalizations as in \cite{pera2019toroidal} with possibly different adjusted Siegel embeddings, and we add script $\mathcal{Y}_i$ to the following strata $\mathcal{Y}$ related with $\Shumc{K, \Sigma}(G, X)_i$ in the construction: $\mathcal{Y}:=\Shum{K}(G, X)$, $\Shumc{K, \Sigma}(G, X)$, $\Shumm{K}(G, X)$ or $\mathcal{Y} = \mathcal{Z}(\Upsilon)$ for all $\Upsilon = [\Phi, \sigma] \in \Cusp_K(G, X, \Sigma)$, or $\mathcal{Y} = \mathcal{Z}([\Phi])$, $\Xi(\Phi)$, $\Xi(\Phi)(\sigma)$, $\Xi(\Phi)_{\sigma}$, $\ovl{\Xi}(\Phi)$, $\mathfrak{X}^{\circ}_{\sigma}(\Phi)$, $\mathfrak{X}(\Phi)$  $C(\Phi)$, $\Zb^{\bigsur}(\Phi)$, $\Zb(\Phi)$ for all $[\Phi] \in \Cusp_K(G, X)$, $\sigma \in \Sigma(\Phi)^+$, see Theorem \ref{theorem: axiomatic descriptions of compactifications} for notations. We abbreviate $\Shumc{i, K, \Sigma}:= \Shumc{K, \Sigma}(G, X)_i$, $\Shumm{i, K} = \Shumm{K}(G, X)_i$.

Consider the composition of diagonal embedding and open embeddings
\[  Y \stackrel{\Delta}{\longrightarrow} Y \times_S Y \longrightarrow \mathcal{Y}_1 \times_{\mathcal{S}} \mathcal{Y}_2, \]
let $\mathcal{Y}$ be the relative normalization of $\mathcal{Y}_1 \times \mathcal{Y}_2$ in $Y$, and $\pi_{Y, i}: \mathcal{Y} \to \mathcal{Y}_i$ be the natural projections, here the base $\mathcal{S}$ (with generic fiber $S$) is taking as the following:
\begin{compactitem}
    \item $\mathcal{S} = \OO_{E(v)}$ when $\mathcal{Y} = \Shum{K}, \Shumc{K, \Sigma}, \Shumm{K}, \mathcal{Z}(\Upsilon), \mathcal{Z}([\Phi]), \Zb^{\bigsur}(\Phi), \Zb(\Phi)$,
    \item $\mathcal{S} = C(\Phi)$ when $\mathcal{Y} = \Xi(\Phi), \Xi_{\sigma}(\Phi), \Xi(\Phi)(\sigma)$,
    \item $\mathcal{S} = \Zb^{\bigsur}(\Phi)$ when $\mathcal{Y} = C(\Phi)$, $\ab_K(\Phi)$. 
\end{compactitem}

\begin{definition}\label{def: natural isomorphism}
 In this subsection, we say there exists a \emph{natural isomorphism} $\mathcal{Y}_1 \cong \mathcal{Y}_2$ if $\pi_{Y, i}$ are isomorphisms for both $i = 1, 2$.
\end{definition}
\begin{remark}
    In particular, natural isomorphism extends the identity morphism on the generic fiber $Y = Y_1 = Y_2$. On the other hand, since all strata $\mathcal{Y}$ listed above are reduced and separated over the base scheme $S$, then there exists at most one isomorphism $\mathcal{Y}_1 \cong \mathcal{Y}_2$ which extends the identity morphism on the generic fiber $Y_1 = Y_2$.
\end{remark}

Once we fix a $[\Phi]$, we omit $\Phi$ in the writing if it is clear. The main result is the following:
\begin{proposition}\label{prop: uniqueness of toroidal compactification}
   Assume $(G, X, K_p) \hookrightarrow (G_i^{\dd}, X_i^{\dd}, K_{i, p}^{\dd})$ are two \emph{adjusted} Siegel embeddings, and $K_p$ is quasi-parahoric (not necessarily stablizer quasi-parahoric).
   \begin{enumerate}
       \item Fix any $[\Phi] \in \Cusp_K(G, X)$, $\sigma \in \Sigma(\Phi)^+$. There exist natural isomorphisms fit into the commutative diagram:
\[\begin{tikzcd}[row sep=tiny]
	{\Xi_{1, \sigma}} & {\Xi_1(\sigma)} & {\Xi_1} \\
	&&& {C_1} & {\Zb_1^{\bigsur}} & {\Zb_1} \\
	{\Xi_{2, \sigma}} & {\Xi_2(\sigma)} & {\Xi_2} \\
	&&& {C_2} & {\Zb_2^{\bigsur}} & {\Zb_2}
	\arrow[hook, from=1-1, to=1-2]
	\arrow[from=1-1, to=2-4]
	\arrow["\cong", from=1-1, to=3-1]
	\arrow[from=1-2, to=2-4]
	\arrow["\cong", from=1-2, to=3-2]
	\arrow[hook', from=1-3, to=1-2]
	\arrow[from=1-3, to=2-4]
	\arrow["\cong", from=1-3, to=3-3]
	\arrow[from=2-4, to=2-5]
	\arrow["\cong", from=2-4, to=4-4]
	\arrow[from=2-5, to=2-6]
	\arrow["\cong", from=2-5, to=4-5]
	\arrow["\cong", from=2-6, to=4-6]
	\arrow[hook, from=3-1, to=3-2]
	\arrow[from=3-1, to=4-4]
	\arrow[from=3-2, to=4-4]
	\arrow[hook', from=3-3, to=3-2]
	\arrow[from=3-3, to=4-4]
	\arrow[from=4-4, to=4-5]
	\arrow[from=4-5, to=4-6]
\end{tikzcd}\]
        These isomorphisms are functorial on $\sigma$ and glue to $\ovl{\Xi}_1 \cong \ovl{\Xi}_2$, and extend to strata-preserving isomorphisms $\mathfrak{X}_{1, \sigma}^{\circ} \cong \mathfrak{X}_{2, \sigma}^{\circ}$, $\mathfrak{X}_1 \cong \mathfrak{X}_2$.
      \item There exist natural isomorphisms $\mathcal{Z}_1(\Upsilon) \cong \mathcal{Z}_2(\Upsilon)$ for all $\Upsilon \in \Cusp_K(G, X, \Sigma)$, and exists natural strata-preserving isomorphism $\Shumc{1, K, \Sigma} \cong \Shumc{2, K, \Sigma}$.
      \item The above isomorphisms fit into the commutative diagram
\[\begin{tikzcd}
	{\mathfrak{X}_{1, \sigma}^{\circ}} & {(\Shumc{1, K, \Sigma})_{\underset{\tau\in \Sigma(\Phi)^{+},\ \bar{\tau}\subset\bar{\sigma}}{\cup} \mathcal{Z}_1([\Phi, \tau])}^{\wedge} } \\
	{\mathfrak{X}_{2, \sigma}^{\circ}} & {(\Shumc{2, K, \Sigma})_{\underset{\tau\in \Sigma(\Phi)^{+},\ \bar{\tau}\subset\bar{\sigma}}{\cup} \mathcal{Z}_2([\Phi, \tau])}^{\wedge} }
	\arrow["\cong", from=1-1, to=1-2]
	\arrow["\cong", from=1-1, to=2-1]
	\arrow["\cong", from=1-2, to=2-2]
	\arrow["\cong", from=2-1, to=2-2]
\end{tikzcd}\]
   \item Let $\Shumm{i, K}$ be the integral model of $\shum{i, K}:=\shum{K}(G, X)_i$ defined as the Stein factorization $\Shumc{i, K, \Sigma} \to H^0(\Shumc{i, K, \Sigma}, \omega_i(\Sigma)^{N_i})$ (or defined using relative normalizations, see Proposition \ref{prop: minimal compactification, rela normalization}), then there exist isomorphisms $\mathcal{Z}_1([\Phi]) \cong \mathcal{Z}_2([\Phi])$ compatible with $\Zb_1(\Phi) \cong \Zb_2(\Phi)$ for all $[\Phi] \in \Cusp_K(G, X)$, and there exists a natural strata-preserving isomorphism $\Shumm{1, K} \cong \Shumm{2, K}$ which fits into the commutative diagram:
   \begin{equation}\label{eq: comm diagram tor min}
\begin{tikzcd}
	{\Shumc{1, K, \Sigma}} & {\Shumc{2, K, \Sigma}} \\
	{\Shumm{1, K}} & {\Shumm{2, K}}
	\arrow["\cong", no head, from=1-1, to=1-2]
	\arrow["{\oint_{1, K, \Sigma}}", from=1-1, to=2-1]
	\arrow["{\oint_{2, K, \Sigma}}", from=1-2, to=2-2]
	\arrow["\cong", no head, from=2-1, to=2-2]
\end{tikzcd}
   \end{equation}
\item All the above isomorphisms extend the identity morphism on the their generic fibers.
   \end{enumerate}
   
\end{proposition}

We prove this theorem piece by piece in this subsection. 

\begin{remark}\label{rmk: real conditions we need}
    In the proof we will see that, besides the axiomatic descriptions in Theorem \ref{theorem: axiomatic descriptions of compactifications}, what we really need in the proof are following two conditions, they are guaranteed when $(G, X, K_p) \hookrightarrow (G_i^{\dd}, X_i^{\dd}, K_{i, p}^{\dd})$ are adjusted Siegel embeddings for $i = 1, 2$.
    \begin{enumerate}
        \item There exist natural isomorphisms $\Zb_1^{\bigsur}(\Phi) \cong \Zb_2^{\bigsur}(\Phi)$ for all $[\Phi] \in \Cusp_K(G, X)$, including $\Shum{1, K} \cong \Shum{2, K}$. 
        \item $C_i(\Phi) \to \Zb_i^{\bigsur}(\Phi)$ are abelian-scheme torsors for all $[\Phi] \in \Cusp_K(G, X)$, $i = 1, 2$.
    \end{enumerate}
\end{remark}

\subsubsection{}

Fix a $[\Phi]$. Let $\mathcal{Y}_i = \Zb^{\bigsur}_i(\Phi) = \Shum{i, K_{\Phi, h}}$. Due to part $(5)$ of Theorem \ref{thm: boundary of kp is again kp}, for each $i = 1, 2$, $\Shum{i, K_{\Phi, h}}$ is an integral model of the Hodge-type Shimura variety $\shu{K_{\Phi, h}}$ with quasi-parahoric level structure defined via the nice embeddings $(G_{\Phi, h}, X_{\Phi, h}, K_{\Phi, h}) \to (G_{\Phi^{\dd}, h}^{\dd}, X_{\Phi^{\dd}, h}, K^{\dd}_{\Phi^{\dd}, h})$, then due to the main result of \cite{pappas2024p} and \cite{daniels2024conjecture}, $\Shum{K_{\Phi, h}} \to \Shum{K_{1, \Phi, h}}$ and  $\Shum{K_{\Phi, h}} \to \Shum{K_{2, \Phi, h}}$ are isomorphisms due to \cite[Theorem 4.2.4]{pappas2024p} and \cite[Corollary 4.1.9]{daniels2024conjecture}. In other words, there is a natural isomorphism $\Shum{K_{1, \Phi, h}} \cong \Shum{K_{2, \Phi, h}}$.

\subsubsection{}\label{subsec: A Y}

Let $\mathcal{Y}_i = \ab_{i, K}(\Phi)$. Note that we are using two different adjusted Siegel embeddings $\iota_i: (G, X) \to (G^{\dd}_i, X^{\dd}_i)$, for each $i$, we have fixed $\Phi_i^{\dd} = \iota_{i, *}(\Phi)$. Then $\ab_{i, K}(\Phi)$ is the relative normalization of $\ab_{i, K^{\dd}}(\Phi_i^{\dd})|_{\Zb_i^{\bigsur}}$ in $A_{K}(\Phi)$, by construction.

Recall that over $\shu{K_{\Phi, h}}(\CC)$, there is a variation of Hodge structure on the constant vector bundle $\Lie V_{\Phi}$ which we denote by $\mathbf{V}_{\MH}(\Phi)$, and $K_{\Phi, V} \subset V_{\Phi}(\A_f)$ refines $\mathbf{V}_{\MH}(\Phi)$ to a $\Z$-Hodge structure $\mathbf{V}_{\MH}(\Phi)_{\Z}$, which is the relative integral homology of $A_K(\Phi) =  \shu{K_{\Phi, V}\rtimes K_{\Phi, h}}(\ovl{P}, X_{\Phi, \ovl{P}}) \to \shu{K_{\Phi, h}}$. The homomorphism $A_K(\Phi) \to A_{i, K}(\Phi_i^{\dd})|_{\shu{K_{\Phi, h}}}$ is associated with $V_{\Phi} \hookrightarrow V_{i, \Phi^{\dd}}^{\dd}$, and $K_{\Phi, V, p} = K_{\Phi^{\dd}_i, V_i^{\dd}, p} \cap V_{\Phi}(\Q_p)$ due to Corollary \ref{corollary: equality of cpt groups on V}.

We pull back $\ab_{i, K}(\Phi) \to \Zb_i^{\bigsur}(\Phi)$ along the isomorphism $\pi_{\Zb^{\bigsur}, i}: \Zb^{\bigsur}(\Phi) \rightiso \Zb_i^{\bigsur}(\Phi)$ and still denote it by $\ab_{i, K}(\Phi)$. We apply the proofs of \cite[Proposition 4.2.2, Lemma A.3.8]{pera2019toroidal} to
\[ A_K(\Phi) \stackrel{\Delta}{\to} \ab_{1, K}(\Phi) \times \ab_{2, K}(\Phi) |_{\Zb^{\bigsur}(\Phi)} \to \ab_{1, K}(\Phi_1^{\dd}) \times \ab_{2, K}(\Phi_2^{\dd}) |_{\Zb^{\bigsur}(\Phi)}, \]
and note that along $V_{\Phi} \stackrel{\Delta}{\hookrightarrow} V_{\Phi} \times V_{\Phi} \hookrightarrow  V_{1, \Phi^{\dd}}^{\dd} \times V_{2, \Phi^{\dd}}^{\dd}$,
\[ K_{\Phi, V, p} = (K_{\Phi^{\dd}_1, V_1^{\dd}, p} \times K_{\Phi^{\dd}_2, V_2^{\dd}, p} \cap (V_{\Phi} \times V_{\Phi}) (\Q_p)) \cap V_{\Phi}(\Q_p) \]
then $\ab_K(\Phi) \to \Zb^{\bigsur}(\Phi)$ is an abelian scheme.

Now we omit the index $\Phi$. We have projections $\pi_{\ab_K, i}: \ab_K \to \ab_{i, K}$ of abelian schemes over $\Zb^{\bigsur}$. By the functoriality of relative normalizations (\cite[\href{https://stacks.math.columbia.edu/tag/035J}{Tag 035J}]{stacks-project}), $\pi_{\ab_K, i}$ is a homomorphism between group schemes. Over each connected component of $\Zb^{\bigsur}$, $\ab_K$, $\ab_{i, K}$ have same relative dimensions, thus $\pi_{\ab_K, i}$ is an isogeny, in particular, it is finite. As $\Zb^{\bigsur}$ is normal, then $\ab_K$ and $\ab_{i, K}$ are normal since they are smooth over $\Zb^{\bigsur}$. By Zariski main theorem, $\pi_{\ab_K, i}$ is an isomorphism. In other words, there is a natural isomorphism $\ab_{1, K} \cong \ab_{2, K}$ compatible with $\Zb_1^{\bigsur} \cong \Zb_2^{\bigsur}$.

\subsubsection{}\label{subsec: C Y}

Let $\mathcal{Y}_i = C_i = \Shum{i, K_{\Phi, \ovl{P}}}$. Due to Corollary \ref{cor: pera, 4.1.5}, $C_i \to \Zb_i^{\bigsur}$ is a torsor under abelian scheme $\ab_{i, K}$. We pull back $C_i$ along $\pi_{\Zb_i^{\bigsur}}: \Zb^{\bigsur} \rightiso \Zb_i^{\bigsur}$ and still denote it by $C_i$.

Due to the result of previous subsection, and since relative normalization commutes with smooth base change, then $\pi_{C, i}: C \to C_i$ is $\pi_{\ab_K, i}: \ab_K \rightiso \ab_{i, K}$ \'etale locally over $\Zb^{\bigsur}$, thus $\pi_{C, i}$ is an isomorphism. In other words, there is a natural isomorphism $C_1 \cong C_2$ compatible with $\Zb_1^{\bigsur} \cong \Zb_1^{\bigsur}$ and $\ab_{1, K} \cong \ab_{2, K}$.

\subsubsection{}\label{subsec: Xi Y}

Let $\mathcal{Y}_i = \Xi_i = \Shum{i, K_{\Phi, P}}$, $\Xi_i \to C_i$ is an $\mathbf{E}_K(\Phi)$-torus torsor due to \cite[Theorem 4.1.5(4)]{pera2019toroidal}, here $\mathbf{E}_K(\Phi)$ is the constant torus over $\Z$ with cocharacter group $\mathbf{B}_K(\Phi):= (U_{\Phi}(\Q) \cap K_{\Phi, U})(-1) \subset U_{\Phi}(\Q)(-1)$. Apply the arguments from previous subsection \ref{subsec: C Y}, by checking \'etale locally, there is a natural isomorphism $\Xi_1 \cong \Xi_2$ compatible with $C_1 \cong C_2$. Similarly, we have natural ismorphisms $\Xi_{1, \sigma} \cong \Xi_{2, \sigma}$ and $\Xi_1(\sigma) \cong \Xi_2(\sigma)$ (use definitions of twisted torus embedding, \cite[\S 2.1.17]{pera2019toroidal}), then they glue to natural isomorphism $\ovl{\Xi}_1 \cong \ovl{\Xi}_2$ (by checking over open covering $\bigcup_{\sigma\in \Sigma(\Phi)^+} \Xi_i(\sigma)$), thus extend to $\mathfrak{X}_{1, \sigma}^{\circ} \cong \mathfrak{X}_{2, \sigma}^{\circ}$, $\mathfrak{X}_1 \cong \mathfrak{X}_2$. We finished the proof of part (1) of Proposition \ref{prop: uniqueness of toroidal compactification}.

\subsubsection{}\label{subsec: Z_i Upsilon, Y}

Let $\mathcal{Y}_i = \mathcal{Z}_i(\Upsilon)$ for some $\Upsilon = [\Phi, \sigma] \in \Cusp_K(G, X, \Sigma)$. Use functoriality and the canonical isomorphisms $\mathcal{Z}_i(\Upsilon) \cong \Xi_{i, \sigma}(\Phi)$ which extend the canonical isomorphism $\mathcal{Z}(\Upsilon)_{\eta} \cong \Xi_{\sigma}(\Phi)_{\eta}$ on the generic fiber, we have
\[\begin{tikzcd}
	{\mathcal{Z}(\Upsilon)_{\eta}} & {\mathcal{Z}(\Upsilon)} & {\mathcal{Z}_1(\Upsilon) \times \mathcal{Z}_2(\Upsilon)} & {\mathcal{Z}_i(\Upsilon)} \\
	{\Xi_{\sigma}(\Phi)_{\eta}} & {\Xi_{\sigma}(\Phi)} & {\Xi_{1, \sigma}(\Phi) \times \Xi_{2, \sigma}(\Phi)} & {\Xi_{i, \sigma}(\Phi)}
	\arrow[from=1-1, to=1-2]
	\arrow["\cong", from=1-1, to=2-1]
	\arrow[from=1-2, to=1-3]
	\arrow["\cong", from=1-2, to=2-2]
	\arrow["{\pr_i}", from=1-3, to=1-4]
	\arrow["\cong", from=1-3, to=2-3]
	\arrow["\cong", from=1-4, to=2-4]
	\arrow[from=2-1, to=2-2]
	\arrow[from=2-2, to=2-3]
	\arrow["{\pr_i}", from=2-3, to=2-4]
\end{tikzcd}\]
$\pi_{\Xi_{\sigma}, i}: \Xi_{\sigma}(\Phi) \to \Xi_{i, \sigma}(\Phi)$ is an isomorphism implies that $\pi_{\mathcal{Z}(\Upsilon), i}: \mathcal{Z}(\Upsilon) \to \mathcal{Z}_i(\Upsilon)$ is an isomorphism. In other words, there is a natural isomorphism $\mathcal{Z}_1(\Upsilon) \cong \mathcal{Z}_2(\Upsilon)$ compatible with $\Xi_{1, \sigma}(\Phi) \to \Xi_{2, \sigma}(\Phi)$.

\subsubsection{}\label{ref: subsec toroidal Y}

Let $\mathcal{Y}_i = \Shumc{i, K, \Sigma}$. Use functoriality, consider
\[\begin{tikzcd}
	{\mathcal{Z}(\Upsilon)} \\
	& {\mathcal{Z}(\Upsilon)'} & {\mathcal{Z}_1(\Upsilon)\times\mathcal{Z}_2(\Upsilon)} \\
	& {\Shumc{K, \Sigma}} & {\Shumc{1, K, \Sigma} \times \Shumc{2, K, \Sigma}}
	\arrow[from=1-1, to=2-2]
	\arrow[from=1-1, to=2-3]
	\arrow[from=1-1, to=3-2]
	\arrow[from=2-2, to=2-3]
	\arrow[from=2-2, to=3-2]
	\arrow["\lrcorner"{anchor=center, pos=0.125, rotate=45}, draw=none, from=2-2, to=3-3]
	\arrow[from=2-3, to=3-3]
	\arrow[from=3-2, to=3-3]
\end{tikzcd}\]
where $\mathcal{Z}(\Upsilon)'$ is the fiber product. Let $\wdt{\mathcal{Z}(\Upsilon)}$ be the normalization of $\mathcal{Z}(\Upsilon)'$, then we have a morphism $\pi: \mathcal{Z}(\Upsilon) \to \wdt{\mathcal{Z}(\Upsilon)}$. Since the projections $\pi_{\mathcal{Z}(\Upsilon), i}: \mathcal{Z}(\Upsilon) \to \mathcal{Z}_i(\Upsilon)$ are isomorphisms, and $\mathcal{Z}(\Upsilon)$ is separated over $\OO_{E(v)}$ (both $\mathcal{Z}_i(\Upsilon)$ are separated and $\mathcal{Z}(\Upsilon) \to \mathcal{Z}_1(\Upsilon)\times\mathcal{Z}_2(\Upsilon)$ is finite), then $\mathcal{Z}(\Upsilon) \to \mathcal{Z}_1(\Upsilon)\times\mathcal{Z}_2(\Upsilon)$ is a closed embedding. Therefore, both $\mathcal{Z}(\Upsilon) \to \mathcal{Z}(\Upsilon)'$ and $\pi: \mathcal{Z}(\Upsilon) \to \wdt{\mathcal{Z}(\Upsilon)}$ are closed embeddings. In particular, $\pi$ is an isomorphism due to Zariski main theorem, thus the normalization $\wdt{\mathcal{Z}(\Upsilon)} \to \mathcal{Z}(\Upsilon)'$ is a closed embedding, then $\wdt{\mathcal{Z}(\Upsilon)} = \mathcal{Z}(\Upsilon)'_{\red}$. 

Therefore, the projections $\pi_{\Shumc{K, \Sigma}, i}: \Shumc{K, \Sigma} \to \Shumc{i, K, \Sigma}$ map $\Shum{K}$ isomorphically onto $\Shum{i, K}$ (relative normalizations commute with open embeddings, and then use the main result of \cite{pappas2024p} and \cite{daniels2024conjecture} with the proof of \cite[Theorem 4.2.4]{pappas2024p} and \cite[Corollary 4.1.9]{daniels2024conjecture}), map $\mathcal{Z}(\Upsilon)$ homeomorphically onto $\mathcal{Z}_i(\Upsilon)$. 

Next, we need to show $\Shumc{K, \Sigma}$ is covered by $\bigcup_{\Upsilon \in \Cusp_K(G, X, \Sigma)} \mathcal{Z}(\Upsilon)$. This shows that $\pi_{\Shumc{K, \Sigma}, i}$ are bijections.

Since $\mathcal{Z}_1(\Upsilon) \cong \mathcal{Z}_2(\Upsilon)$ for all $\Upsilon \in \Cusp_K(G, X, \Sigma)$. Given a point $x_1 \in \mathcal{Z}_1(\Upsilon)$, we can identify it with a point $x_2 \in \mathcal{Z}_2(\Upsilon)$. Use the canonical isomorphisms $\mathcal{Z}_i(\Upsilon) \cong \Xi_i(\Phi)_{\sigma}$, we identify $x_i$ with $y_i \in \Xi_i(\Phi)_{\sigma}$. There exist \'etale neighourhoods $(\ovl{U}_i, z_i)$ of $(\Shumc{i, K, \Sigma}, x_i)$ together with \'etale morphisms $(\ovl{U}_i, z_i) \to (\Xi_i(\Phi)(\sigma), y_i)$. These \'etale morphisms match the stratifications on the source and the target. We can assume $\ovl{U}_1$ and $\ovl{U}_2$ have same generic fiber. Note that the image of the diagonal embedding $\shuc{K, \Sigma} \to \shuc{K, \Sigma} \times \shuc{K, \Sigma}$ is covered by such \'etale neighourhoods $\ovl{U}_1 \times \ovl{U}_2$ when varing $(x_1, x_2)$ together with the image of $\shu{K} \to \shu{K} \times \shu{K}$, and since taking relative normalizations commutes with \'etale base change, we can check the covering $\Shumc{K, \Sigma} = \bigcup_{\Upsilon \in \Cusp_K(G, X, \Sigma)} \mathcal{Z}(\Upsilon)$ over \'etale neighourhoods $\ovl{U}_1 \times \ovl{U}_2$. Then we can apply the results in the step \ref{subsec: Xi Y} and \ref{subsec: Z_i Upsilon, Y}.

Since $\Shumc{i, K, \Sigma} \to \Spec \OO_{E(v)}$ are proper, $\Shumc{K, \Sigma} \to \Shumc{1, K, \Sigma} \times \Shumc{2, K, \Sigma}$ is finite, then $\pi_{\Shumc{K, \Sigma}, i}$ are proper. By Zariski main theorem, $\pi_{\Shumc{K, \Sigma}, i}$ are isomorphisms, and they induce isomorphisms
\[  (\Shumc{K, \Sigma})_{\underset{\tau\in \Sigma(\Phi)^{+},\ \bar{\tau}\subset\bar{\sigma}}{\cup} \mathcal{Z}([\Phi, \tau])}^{\wedge}  \rightiso (\Shumc{i, K, \Sigma})_{\underset{\tau\in \Sigma(\Phi)^{+},\ \bar{\tau}\subset\bar{\sigma}}{\cup} \mathcal{Z}_i([\Phi, \tau])}^{\wedge} . \]
We finish the proof of part $(2)$ of Proposition \ref{prop: uniqueness of toroidal compactification}.

\subsubsection{}

We need to show that $f_1 = f_2$ for the unique isomorphism $f_i$ that fits into the commutative diagram:
\[\begin{tikzcd}
	{\mathfrak{X}_{\sigma}^{\circ}} & {(\Shumc{K, \Sigma})_{\underset{\tau\in \Sigma(\Phi)^{+},\ \bar{\tau}\subset\bar{\sigma}}{\cup} \mathcal{Z}([\Phi, \tau])}^{\wedge} } \\
	{\mathfrak{X}_{i, \sigma}^{\circ}} & {(\Shumc{i, K, \Sigma})_{\underset{\tau\in \Sigma(\Phi)^{+},\ \bar{\tau}\subset\bar{\sigma}}{\cup} \mathcal{Z}_i([\Phi, \tau])}^{\wedge} }
	\arrow["{f_i}", from=1-1, to=1-2]
	\arrow["\cong", from=1-1, to=2-1]
	\arrow["\cong", from=1-2, to=2-2]
	\arrow["\cong", from=2-1, to=2-2]
\end{tikzcd}\]

Since $f_1 = f_2$ over the generic fiber (the vertical isomorphisms are simply identity morphisms), then they are identical over an open dense formal subscheme of a reduced subcheme. Since the target is separated over $\Spec \OO_{E(v)}$, then $f_1 = f_2$. Note that the notion of a morphism being separate works well for formal schemes, and \cite[\href{https://stacks.math.columbia.edu/tag/01KM}{Tag 01KM}]{stacks-project} works as well. We finish the proof of part $(3)$ of Proposition \ref{prop: uniqueness of toroidal compactification}.

\subsubsection{}

Let $\mathcal{Y}_i = \Shumm{i, K}$. Use functoriality, consider
\[\begin{tikzcd}
	{\Zb^{\bigsur}} \\
	& {\Zb^{\bigsur\prime}} & {\Zb^{\bigsur}_1\times\Zb^{\bigsur}_2} \\
	& {\Shumm{K}} & {\Shumm{1, K} \times \Shumm{2, K}}
	\arrow[from=1-1, to=2-2]
	\arrow[from=1-1, to=2-3]
	\arrow[from=1-1, to=3-2]
	\arrow[from=2-2, to=2-3]
	\arrow[from=2-2, to=3-2]
	\arrow["\lrcorner"{anchor=center, pos=0.125, rotate=45}, draw=none, from=2-2, to=3-3]
	\arrow[from=2-3, to=3-3]
	\arrow[from=3-2, to=3-3]
\end{tikzcd}\]
where $\Zb^{\bigsur\prime}$ is the fiber product. Let $\wdt{\Zb^{\bigsur}}$ be the normalization of $\Zb^{\bigsur\prime}$. Similar arguments as in step \ref{ref: subsec toroidal Y} show that $\Zb^{\bigsur} \cong \wdt{\Zb^{\bigsur}}$ is finite over $\Zb^{\bigsur\prime}$. 

Given different $[\Phi], [\Phi'] \in \Cusp_K(G, X)$, the images of $\Zb^{\bigsur}(\Phi)$, $\Zb^{\bigsur}(\Phi')$ in $\Shumm{K}$ are disjoint by definition. We need to show $\Shumm{K}$ is covered by the images of $\Zb^{\bigsur}(\Phi)$ with varing $[\Phi] \in \Cusp_K(G, X)$, this will show that $\pi_{\Shumm{K}, i}: \Shumm{K} \to \Shumm{K, i}$ are quasi-finite. Since $\pi_{\Shumm{K}, i}$ are proper, then are isomorphisms.

By Proposition \ref{prop: minimal compactification, rela normalization} and functoriality of relative normalization, we have a morphism $\oint_{K, \Sigma}: \Shumc{K, \Sigma} \to \Shumm{K}$ which fits into the commutative diagram:
\[\begin{tikzcd}
	{\Shumc{K, \Sigma}} && {\Shumc{1, K, \Sigma} \times \Shumc{2, K, \Sigma}} \\
	{\Shumm{K}} && {\Shumm{1, K}\times \Shumm{2, K}}
	\arrow[from=1-1, to=1-3]
	\arrow["{\oint_{K, \Sigma}}", from=1-1, to=2-1]
	\arrow["{\oint_{1, K, \Sigma} \times \oint_{2, K, \Sigma}}", from=1-3, to=2-3]
	\arrow[from=2-1, to=2-3]
\end{tikzcd}\]
Since $\oint_{K, \Sigma}$ is proper and birational, it is surjective. In particular, any point in the boundary of $\Shumm{K}$ is in the image of some $\mathcal{Z}(\Upsilon) \cong \Xi_{\sigma}(\Phi)$, which is in the image of some $\Zb^{\bigsur}(\Phi)$. This also proves the commutative diagram \ref{eq: comm diagram tor min}.

\subsubsection{}

Let $\mathcal{Y}_i = \Zb_i(\Phi) = \Lambda_{\Phi, K}\backslash \Shum{i, K_{\Phi, h}}$. By functoriality, we have commutative diagram
\[\begin{tikzcd}
	{\Zb^{\bigsur}} & {\Zb^{\bigsur}_1 \times \Zb^{\bigsur}_2} \\
	\Zb & {\Zb_1 \times \Zb_2}
	\arrow[from=1-1, to=1-2]
	\arrow[from=1-1, to=2-1]
	\arrow[from=1-2, to=2-2]
	\arrow[from=2-1, to=2-2]
\end{tikzcd}\]
Since $\Zb^{\bigsur} \to \Zb^{\bigsur}_1 \times \Zb^{\bigsur}_2$ is a closed embedding and $\Zb^{\bigsur}_i \to \Zb_i$ are finite morphism, then $\Zb^{\bigsur} \to \Zb$ is a finite morphism. On the other hand, by functoriality, $\Lambda_{\Phi, K}$ acts on $\Zb^{\bigsur}$. We denote $\ovl{\Zb} = \Lambda_{\Phi, K}\backslash \Zb^{\bigsur}$, then we have a finite surjection $\Zb^{\bigsur} \to \ovl{\Zb}$, and $\Zb^{\bigsur} \to \Zb$ factors through $\ovl{\Zb} \to \Zb$, then $\ovl{\Zb} \to \Zb$ is finite surjective. By Zariski main theorem, $\ovl{\Zb} \cong \Zb$.

On the other hand, consider the diagram 
\[\begin{tikzcd}
	\Zb \\
	& {\Zb'} & {\Zb_1 \times \Zb_2} \\
	& {\Shumm{K}} & {\Shumm{1, K} \times \Shumm{2, K}}
	\arrow[from=1-1, to=2-2]
	\arrow[from=1-1, to=2-3]
	\arrow[from=1-1, to=3-2]
	\arrow[from=2-2, to=2-3]
	\arrow[from=2-2, to=3-2]
	\arrow["\lrcorner"{anchor=center, pos=0.125, rotate=45}, draw=none, from=2-2, to=3-3]
	\arrow[from=2-3, to=3-3]
	\arrow[from=3-2, to=3-3]
\end{tikzcd}\]
where $\Zb'$ is the fiber product. Let $\wdt{\Zb}$ be the normalization of $\Zb'$. Since $\Zb \to \wdt{\Zb}$ is a finite morphism, by Zariski main theorem, $\Zb \cong \wdt{\Zb}$, and $\Zb \to \Zb'$ is a finite morphism, $\Zb' \to \Shumm{K}$ is a locally closed embedding. Consider $\pi_{\Zb, i}: \Zb \to \Zb_i$. Since both $\Zb^{\bigsur}_i \to \Zb_i$ and $\Zb^{\bigsur} \to \Zb$ are quasi-finite and surjective, and $\pi_{\Zb^{\bigsur}, i}: \Zb^{\bigsur} \cong \Zb^{\bigsur}_i$ is an isomorphism, then $\pi_{\Zb, i}$ is quasi-finite and surjective. To show $\pi_{\Zb, i}$ is an isomorphism, it suffices to show it is proper. Let $\ovl{\Zb}_i$ be the closure of $\Zb_i$ in $\Shumm{i, K}$, by removing $\ovl{\Zb}_i\backslash \Zb_i$ from $\Shumm{i, K}$, we can assume $\Zb_i$ is closed in $\Shumm{i, K}$, then $\Zb \to \Shumm{K} \cong \Shumm{i, K}$ is finite, $\pi_{\Zb, i}: \Zb \to \Zb_i$ is finite, we are done. This shows that $\Zb \to \Shumm{K}$ is a locally closed embedding, and we denote it by $\mathcal{Z}([\Phi])$. The natural projections $\pi_{\Shumm{K}, i}: \Shumm{K} \to \Shumm{i, K}$ are isomorphisms, they map $\Shum{K}$ isomorphically onto $\Shum{i, K}$, map $\mathcal{Z}([\Phi])$ isomorphically onto $\mathcal{Z}_i([\Phi])$. We finish the proof of Proposition \ref{prop: uniqueness of toroidal compactification}.

\subsection{Change-of-Parahoric: Constructions}\label{subsec: change of parahoric}

    Let $(G, X)$ be a Hodge-type Shimura datum. Let $K_{1, p} \subset K_{2, p}$ be two quasi-parahoric subgroups of $G(\Q_p)$, choose $K_1^p \subset K_2^p$, let $K_1 = K_{1, p}K_1^p$, $K_2 = K_{2, p}K_2^p$, then we have a morphism $\Shum{K_1} \to \Shum{K_2}$ which extends the finite morphism $\shu{K_1} \to \shu{K_2}$, due to \cite[Corollary 4.3.2]{pappas2024p} and \cite[Corollary 4.1.10]{daniels2024conjecture}. 

    In this subsection, we prove the following:
\begin{theorem}\label{Theorem: change of parahoric}
    Let $K_{1, p} \subset K_{2, p}$ be two quasi-parahoric subgroups of $G(\Q_p)$.
    \begin{enumerate}
        \item The morphism $\Shum{K_1} \to \Shum{K_2}$ extends (uniquely) to a proper morphism $\pi^{\min}: \Shumm{K_1} \to \Shumm{K_2}$ and $\pi^{\min, -1}(\Shum{K_2}) = \Shum{K_1}$.
        \item Let $\Sigma_1$ refine $\Sigma_2$ such that we have a proper morphism $\shuc{K_1, \Sigma_1} \to \shuc{K_2, \Sigma_2}$ extending the projection $\shu{K_1} \to \shu{K_2}$, then there is a (unique) proper morphism $\pi^{\tor}: \Shumc{K_1, \Sigma_1} \to \Shumc{K_2, \Sigma_2}$ extending both $\shuc{K_1, \Sigma_1} \to \shuc{K_2, \Sigma_2}$ and $\Shum{K_1} \to \Shum{K_2}$, and makes the following diagram commutes:
\[\begin{tikzcd}
	{\Shumc{K_1, \Sigma_1}} & {\Shumc{K_2, \Sigma_2}} \\
	{\Shumm{K_1}} & {\Shumm{K_2}}
	\arrow["{\pi^{\tor}}", from=1-1, to=1-2]
	\arrow["{\oint_{K_1, \Sigma_1}}", from=1-1, to=2-1]
	\arrow["{\oint_{K_2, \Sigma_2}}", from=1-2, to=2-2]
	\arrow["{\pi^{\min}}", from=2-1, to=2-2]
\end{tikzcd}\]
        
    \end{enumerate}
\end{theorem}

\begin{remark}
     From Proposition \ref{prop: uniqueness of toroidal compactification}, we see that the construction of $\Shumc{K, \Sigma}$ and $\Shumm{K}$ are independent of the choice of an adjusted Siegel embedding just as $\Shum{K}$, thus it makes sense to speak the Proposition without mentioning the adjusted Siegel embedding we are using. In constructions of $\Shum{K_i}$, $\Shumc{K_i, \Sigma_i}$ and $\Shumm{K_i}$, we use adjusted Siegel embeddings instead of random ones as in \cite{pera2019toroidal}. For $i = 1, 2$, we need to use different adjusted Siegel embeddings. Otherwise the inteior integral models $\Shum{K_i}$ are not the correct models, and the constructions of $\pi^{\tor}$ and $\pi^{\min}$ are easy consequence from functoriality, see Remark \ref{rmk: morphism between toroidal, bad} and \ref{rmk: morphism between minimal, bad}.
\end{remark}
\begin{remark}\label{rmk: remaining of the proof, change of parahoric}\leavevmode
    \begin{enumerate}
        \item The uniquenss of the extensions as well as the commutativity of the diagram come from the fact that given $f_1, f_2: X \to Y$, $X, Y$ flat, separated over $\Spec \OO_{E(v)}$, if $f_1 = f_2$ over the generic fiber $X_{\eta}$ (which is open and dense) of the reduced scheme $X$, then $f_1 = f_2$. 
        \item The properness of both $\pi^{\min}$ and $\pi^{\tor}$ comes trivially from the fact that $\Shumm{K_1}$, $\Shumm{K_2}$, $\Shumc{K_1, \Sigma_1}$, $\Shumc{K_2, \Sigma_2}$ are proper over $\Spec \OO_{E(v)}$. This in turns implies that $\pi: \Shum{K_1} \to \Shum{K_2}$ is proper.
    \end{enumerate}
\end{remark}

\subsubsection{}
Let $(G, X) \hookrightarrow (G^{\dd}, X^{\dd}) = (\GSp(V, \psi), S^{\pm})$ be a Siegel embedding, $x_1$, $x_2 \in \Bui_{\ext}(G, \Q_p)$ be two points associated with the quasi-parahoric subgroups $K_{1, p}$ and $K_{2, p}$ respectively. Fix a $G(\bQ)$-$\Gal(\bQ|\Q_p)$-equivariant embedding $\iota: \Bui_{\ext}(G, \Q_p) \to \Bui_{\ext}(G^{\dd}, \Q_p)$, and let $y_1 = \iota(x_1)$, $y_2 = \iota(x_2)$. Let $\mathfrak{g}_1$, $\mathfrak{g}_2$, $\mathfrak{h}_1$, $\mathfrak{h}_2$ be the facets that contain $x_1$, $x_2$, $y_1$, $y_2$ respectively. Since $K_1 \subset K_2$, then the closure of $\mathfrak{g}_1$ contains $\mathfrak{g}_2$, and  the closure of $\mathfrak{h}_1$ contains $\mathfrak{h}_2$.

As in \cite[\S 1.1.11]{kisin2018integral}, $\mathfrak{h}_1$ corresponds to a periodic chain of $\Z_p$-lattices in $V_{\Q_p}$
\[ \cdots \subset p\Lambda_0 \subset \Lambda_r \subset \Lambda_{r-1} \subset \cdots \subset \Lambda_0 \subset \Lambda_0^{\vee} \subset \cdots, \]
and $\mathfrak{h}_2$ corresponds to a periodic sub-chain of $\Z_p$-lattices:
\[ \cdots \subset p\Lambda_{s_1} \subset \Lambda_{s_t} \subset \Lambda_{s_{t-1}} \subset \cdots \subset \Lambda_{s_1} \subset \Lambda_{s_1}^{\vee} \subset \cdots, \]
for some $\lrbracket{s_1, \dots, s_t} \subset \lrbracket{1, 2, \dots, r}$. Let $J^1 = \lrbracket{1, 2, \dots, r}$, $J^2 = \lrbracket{s_1, s_2, \dots, s_t}$, $J^1 \subset J^2$

Let $K_{i, p}^{\dd} \subset \GSp(\Q_p)$ be the stablizer group of such chain of lattices. In other words, $K_{i, p}^{\dd} = \GSP_{y_i}(\Z_p)$. Choose $K_i^{\dd, p}$ properly, and let $K_i^{\dd} = K_{i, p}^{\dd}K_i^{\dd, p}$, we have a commutative diagram with horizontal morphisms being closed embeddings:
\begin{equation}\label{eq: change of parahoric, generic}
\begin{tikzcd}
	{\shu{K_1}(G, X)} & {\shu{K_1^{\dd}}(G^{\dd}, X^{\dd})} \\
	{\shu{K_2}(G, X)} & {\shu{K_2^{\dd}}(G^{\dd}, X^{\dd})}
	\arrow[from=1-1, to=1-2]
	\arrow[from=1-1, to=2-1]
	\arrow[from=1-2, to=2-2]
	\arrow[from=2-1, to=2-2]
\end{tikzcd}
\end{equation}

The choice of $\lrbracket{\Lambda_j}_{j \in J^i}$ (or more precisely the choice of $\Z$-lattices $\lrbracket{V_{j, \Z}}_{j \in J^i}$, $V_{j, \Z} \subset V_{\Q}$, $V_{j, \Z} \otimes \Z_{(p)} = \Lambda_j \cap V_{\Q}$) gives a moduli interpretation of $\shu{K_i^{\dd}}(G^{\dd}, X^{\dd})$, which says that $\shu{K_i^{\dd}}(G^{\dd}, X^{\dd})$ classifies chains of $p$-isogenies between polarized abelian schemes $A_j$ with $K_i^{\dd, p}$-level structures. Use this moduli interpretention, one can construct integral models $\Shum{K_i^{\dd}}(G^{\dd}, X^{\dd})$ of $\shu{K_i^{\dd}}(G^{\dd}, X^{\dd})$.

It follows from the moduli interpretation that $\shu{K_1^{\dd}}(G^{\dd}, X^{\dd}) \to \shu{K_2^{\dd}}(G^{\dd}, X^{\dd})$ extends to $\Shum{K_1^{\dd}}(G^{\dd}, X^{\dd}) \to \Shum{K_2^{\dd}}(G^{\dd}, X^{\dd})$. We take $\Shum{K_i}(G, X)$ to be the relative normalization of $\Shum{K_i^{\dd}}(G^{\dd}, X^{\dd})$ in $\shu{K_i}(G, X)$, since $(G, X, K_i) \to (G^{\dd}, X^{\dd}, K_i^{\dd})$ are very nice embeddings, $\Shum{K_i}(G, X)$ satisfy the conditions in \cite[Conjecture 4.2.2]{pappas2024p} due to \cite[Theorem 4.5.2]{pappas2024p} and \cite[Theorem 4.1.12]{daniels2024conjecture}. Due to \cite[Corollary 4.3.2]{pappas2024p} and \cite[Corollary 4.1.10]{daniels2024conjecture}, the morphism $\shu{K_1}(G, X) \to \shu{K_2}(G, X)$ extends to $\Shum{K_1}(G, X) \to \Shum{K_2}(G, X)$.

As in the Remark \ref{rmk: remaining of the proof, change of parahoric}, the diagram \ref{eq: change of parahoric, generic} extends to a commutative diagram:
    \begin{equation}\label{eq: change of parahoric, interior}
\begin{tikzcd}
	{\Shum{K_1}(G, X)} & {\Shum{K_1^{\dd}}(G^{\dd}, X^{\dd})} \\
	{\Shum{K_2}(G, X)} & {\Shum{K_2^{\dd}}(G^{\dd}, X^{\dd})}
	\arrow[from=1-1, to=1-2]
	\arrow[from=1-1, to=2-1]
	\arrow[from=1-2, to=2-2]
	\arrow[from=2-1, to=2-2]
\end{tikzcd}
    \end{equation}

\subsubsection{}

There are two ways to define $\Shumc{K_i^{\dd}, \Sigma_i^{\dd}}(G^{\dd}, X^{\dd})$ and $\Shumm{K_i^{\dd}}(G^{\dd}, X^{\dd})$ of $\Shum{K_i^{\dd}}(G^{\dd}, X^{\dd})$. One comes from choosing an adjusted Siegel embedding, this give correct definition of $\Shumc{K_i, \Sigma_i}(G, X)$ and $\Shumm{K_i}(G, X)$ via furthur relative normalizations. One comes from the construction in \cite{lan2016compactifications}, which gives functorial morphisms \ref{eq: toroidal, intermediate}, \ref{eq: intermediate, minimal}, and \ref{eq: change of parahoric, compactification} by furthur relative normalizations. We show that these two constructions coincide each other.

Let us recall the construction in \cite{lan2016compactifications}. In our case (see \cite[Example 13.12]{lan2016compactifications} for details), for each lattice $\Lambda_j$ ($j \in J^1$), one can associated it with a Shimura variety $\Shum{H_j}(G^{\dd}_j, X^{\dd}_j)$ over $\Spec \Z_{(p)}$, where $(G^{\dd}, X^{\dd}) \cong (G^{\dd}_j, X^{\dd}_j)$, $H_{j, p} \subset G^{\dd}_j(\Q_p)$ is the stablizer of the lattice $\Lambda_j$, and $H_j^p \subset G^{\dd}_j(\A_f^p)$ is properly chosen. $(G^{\dd}_j, X^{\dd}_j)$ is denoted by $(G_{j, \mathrm{aux}}, X_{j, \mathrm{aux}})$ in that paper. $\Shum{H_j}(G^{\dd}_j, X^{\dd}_j)$ classifies a single $A_j$ with away from $p$-level structures instead of classifying a chain of $p$-isogenies between $\lrbracket{A_j}_{j \in J^i}$. Use moduli interpretations, we have a commutative diagram:
\begin{equation}
\begin{tikzcd}
	{\Shum{K_1^{\dd}}(G^{\dd}, X^{\dd})} & {\prod_{j \in J^1} \Shum{H_j^{\dd}}(G^{\dd}_j, X^{\dd}_j)} \\
	{\Shum{K_2^{\dd}}(G^{\dd}, X^{\dd})} & {\prod_{j \in J^2} \Shum{H_j^{\dd}}(G^{\dd}_j, X^{\dd}_j)}
	\arrow[from=1-1, to=1-2]
	\arrow[from=1-1, to=2-1]
	\arrow[from=1-2, to=2-2]
	\arrow[from=2-1, to=2-2]
\end{tikzcd}    
\end{equation}

\begin{lemma}\label{lemma: interior, rela norm}
    The horizontal morphisms are finite morphisms, and $\Shum{K_i^{\dd}}(G^{\dd}, X^{\dd})$ is the relative normalization of $\prod_{j \in J^i} \Shum{H_j^{\dd}}(G^{\dd}_j, X^{\dd}_j)$ in $\shu{K_i^{\dd}}(G^{\dd}, X^{\dd})$.
\end{lemma}
\begin{proof}
    From moduli interpretations, we see that the horizontal morphisms are proper and quasi-finite, thus are finite. The second statement comes from Zariski main theorem as usual.
\end{proof}

With well chosen cone decompositions, $\Shumc{K_i^{\dd}, \Sigma_i^{\dd}}(G^{\dd}, X^{\dd})$ is defined as the relative normalization of $\prod_{j \in J^i} \Shumc{H_j^{\dd}, \Sigma_i^{j, \dd}}(G^{\dd}_j, X^{\dd}_j)$ in $\shuc{K_i^{\dd}, \Sigma_i^{\dd}}(G^{\dd}, X^{\dd})$ (see \cite[\S 7]{lan2016compactifications}), and they fit into the commutative diagram which extends the diagram on the generic fiber:

\begin{equation}\label{eq: toroidal, intermediate}
\begin{tikzcd}
	{\Shumc{K_1^{\dd}, \Sigma_1^{\dd}}(G^{\dd}, X^{\dd})} & {\prod_{j \in J^1} \Shumc{H_j^{\dd}, \Sigma_i^{j, \dd}}(G^{\dd}_j, X^{\dd}_j)} \\
	{\Shumc{K_2^{\dd}, \Sigma_2^{\dd}}(G^{\dd}, X^{\dd})} & {\prod_{j \in J^2} \Shumc{H_j^{\dd}, \Sigma_i^{j, \dd}}(G^{\dd}_j, X^{\dd}_j)}
	\arrow[from=1-1, to=1-2]
	\arrow[from=1-1, to=2-1]
	\arrow[from=1-2, to=2-2]
	\arrow[from=2-1, to=2-2]
\end{tikzcd}    
\end{equation}
where the right vertical morphism is simply the projection, and the left vertical morphism comes from functoriality of relative normalizations.

Similarly, the minimal compactification $\Shumm{K_i^{\dd}}(G^{\dd}, X^{\dd})$ of $\Shum{K_i^{\dd}}(G^{\dd}, X^{\dd})$ is also defined as the relative normalization of $\prod_{j \in J^i} \Shumm{H_j^{\dd}}(G^{\dd}_j, X^{\dd}_j)$ in $\Shumm{K_i^{\dd}}(G^{\dd}, X^{\dd})$ (see \cite[\S 6]{lan2016compactifications}), and we have a commutative diagram extending the one on generic fibers:
\begin{equation}\label{eq: intermediate, minimal}
\begin{tikzcd}
	{\Shumm{K_1^{\dd}}(G^{\dd}, X^{\dd})} & {\prod_{j \in J^1} \Shumm{H_j^{\dd}}(G^{\dd}_j, X^{\dd}_j)} \\
	{\Shumm{K_2^{\dd}}(G^{\dd}, X^{\dd})} & {\prod_{j \in J^2} \Shumm{H_j^{\dd}}(G^{\dd}_j, X^{\dd}_j)}
	\arrow[from=1-1, to=1-2]
	\arrow[from=1-1, to=2-1]
	\arrow[from=1-2, to=2-2]
	\arrow[from=2-1, to=2-2]
\end{tikzcd}
\end{equation}

\subsubsection{}

On the other hand, by combining the lattices $\lrbracket{\Lambda_j}_{j \in J^i}$ as in \cite{kisin2010integral}, \cite{kisin2018integral}, \cite{pappas2024p}, and use Zarhin's trick, we have an adjusted Siegel embedding $(G^{\dd}, X^{\dd}) \hookrightarrow (H_i, T_i)$ which replace $y_i \in \Bui_{\ext}(G^{\dd}, \Q_p)$ with a hyperspecial point $z_i \in \Bui_{\ext}(H_i, \Q_p)$. We take $M_{i, p} = \GSP^{\prime}_{z_i}(\Z_p)$ as usual. We can apply the construction in \cite{pera2019toroidal} to $(G^{\dd}, X^{\dd}) \hookrightarrow (H_i, T_i)$ (we view $(G^{\dd}, X^{\dd})$ as a Hodge-type Shimura datum).

\begin{lemma}\label{lem: two constructions equal}
    Let $\Sigma_i^{\dd}$ be the collection of cone decompositions induced by some $\Sigma_{H_i}$, then $\Shumc{K_i^{\dd}, \Sigma_i^{\dd}}(G^{\dd}, X^{\dd})$ is the relative normalization of $\Shumc{M_i, \Sigma_{H_i}}(H_i, T_i)$ in $\shuc{K_i^{\dd}, \Sigma_i^{\dd}}(G^{\dd}, X^{\dd})$. Also, $\Shumm{K_i^{\dd}}(G^{\dd}, X^{\dd})$ is the relative normalization of $\Shumm{M_i}(H_i, T_i)$ in $\shum{K_i^{\dd}}(G^{\dd}, X^{\dd})$.
\end{lemma}
\begin{proof}
    We omit the upperscript $(\ast)^{\dd}$ everywhere only in this proof, to simplify the writing. We denote such relative normalizations by $\Shumc{K_i, \Sigma_i}(G, X)'$ and $\Shumm{K_i}(G, X)'$ respectively. We need to show there are natural isomorphisms in the sense of Definition \ref{def: natural isomorphism}:
    \[  \Shumc{K_i, \Sigma_i}(G, X) \cong \Shumc{K_i, \Sigma_i}(G, X)',\quad \Shumm{K_i}(G, X) = \Shumm{K_i}(G, X)'. \]
    This follows from the proof of Proposition \ref{prop: uniqueness of toroidal compactification}: the compactifications $\Shumc{K_i, \Sigma_i}(G, X)'$ and $\Shumm{K_i}(G, X)'$ come from adjusted Siegel embedding, and we need to show the compactifications $\Shumc{K_i, \Sigma_i}(G, X)$ and $\Shumm{K_i}(G, X)$ also satisfy the propositions listed in Theorem \ref{theorem: axiomatic descriptions of compactifications}, and satisfy the two conditions in Remark \ref{rmk: real conditions we need}, in order to run through the arguments in the proof of Proposition \ref{prop: uniqueness of toroidal compactification}:

    The main result of the paper \cite{lan2016compactifications} shows that such $\Shumc{K_i, \Sigma_i}(G, X)$ and $\Shumm{K_i}(G, X)$ satisfy the propositions listed in Theorem \ref{theorem: axiomatic descriptions of compactifications}, see \cite[Proposition 2.1.2]{lan2018compactifications} (\cite[Case (3) in Assumption 2.1]{lan2018compactifications}), and Remark \ref{rmk: cusp labels coincide}. 

    We explain those two conditions in Remark \ref{rmk: real conditions we need}:
    \begin{enumerate}
        \item As in Lemma \ref{lemma: interior, rela norm}, $\Shum{K_i}(G, X)$ is defined using relative normalization via the nice embedding $(G, X, K_i) \to \prod_{j \in J^i}(G_j, X_j, H_j)$, then due to \cite[Theorem 4.1.8, 4.1.12]{daniels2024conjecture}, $\Shum{K_i}(G, X)$ is an integral model of $\shu{K_i}(G, X)$ in the sense of \cite{pappas2024p}. Since $(G, X, K_i) \to \prod_{j \in J^i}(G_j, X_j, H_j)$ is a nice embedding, then
        \begin{equation}\label{eq: nice embedding on boundary}
            (G_{\Phi, h}, X_{\Phi, h}, K_{i, \Phi, h}) \to \prod_{j \in J^i}(G_{j, \Phi_j, h}, X_{j, \Phi_j, h}, H_{j, \Phi_j, h})
        \end{equation}
        are also nice embeddings due to Remark \ref{rmk: intermediate steps, G_h}, where $\Phi_j \in \Cusp_{H_j}(G_j, X_j)$ are induced by $\Phi \in \Cusp_{K_i}(G, X)$. By construction in \cite[\S 7]{lan2016compactifications}, the morphism between boundary strata $\Shum{K_{i, \Phi, h}}(G_h, X_{\Phi, h}) \to \prod_{j \in J^i} \Shum{H_{j, \Phi_j, h}}(G_{j, h}, X_{j, \Phi_j, h})$ are defined using the relative normalization via \ref{eq: nice embedding on boundary}, thus $\Shum{K_{i, \Phi, h}}(G_h, X_{\Phi, h})$ are also integral models of $\shu{K_{i, \Phi, h}}(G_h, X_{\Phi, h})$ for all $\Phi$ in the sense of \cite{pappas2024p}.
        \item We apply the same lines of the proofs of \cite[Proposition 4.2.2, Corollary 4.2.3]{pera2019toroidal}, with $A_{K^{\dd}}(\Phi^{\dd})$ replaced by $\prod_{j \in J^i} A_{H_j}(\Phi_j)$, where each 
        \[ A_{H_j}(\Phi_j) \to \Shum{H_{j, \Phi_j, h}}(G_{j, h}, X_{j, \Phi_j, h})\] 
        is an abelian scheme. The assumption in \cite[Proposition 4.2.2]{pera2019toroidal} is verified in the similar way as in subsection \ref{subsubsec: about unipotent groups}, where $G^{\dd}$ is replaced with $\prod_{j \in J^i}G_j$. Note that group theoretically this case is even easier, since we only need to consider step \ref{subsubsec: modifications} and moreover consider the easier case $\prod$ instead of the restricted product $\prod^{\prime}$. We also modify the arguments in the part \ref{subsec: A Y} accordingly.
    \end{enumerate}
    Note that we can let $K^p$ be principal level group, since $K_p$ is parahoric, then $K = K_pK^p$ satisfies the condition in Corollary \ref{cor: Lambda trivial action}. Then $\Lambda_{K_i, \Phi}$ acts trivially on $\Shum{K_{i, \Phi, h}}(G_h, X_{\Phi, h})$ for all $[\Phi] \in \Cusp_{K_i}(G, X)$.
\end{proof}
\begin{corollary}\label{cor: intermediate shimura datum}
    Let $\Sigma_i$ be induced by $\Sigma_i^{\dd}$ (induced by some $\Sigma_{H_i}$), then $\Shumc{K_i, \Sigma_i}(G, X)$, which is defined as the relative normalization of $\Shumc{M_i, \Sigma_{H_i}}(H_i, T_i)$ in $\shuc{K_i, \Sigma_i}(G, X)$ via the adjusted Siegel embedding $(G, X) \to (G^{\dd}, X^{\dd}) \to (H_i, X_i)$, is the relative normalization of $\Shumc{K_i^{\dd}, \Sigma_i^{\dd}}(G^{\dd}, X^{\dd})$ in $\shuc{K_i, \Sigma_i}(G, X)$. Similar statement holds for minimal compactifications.
\end{corollary}

Although there is no functorial morphism between $(H_1, T_1)$ and $(H_2, T_2)$, we prove Theorem \ref{Theorem: change of parahoric} with the help of the intermediate Shimura datum $(G^{\dd}, X^{\dd})$ using Corollary \ref{cor: intermediate shimura datum}: with the help of commutative diagrams \ref{eq: toroidal, intermediate} and \ref{eq: intermediate, minimal}, use functoriality of relative normalizations, the morphism between integral models of different quasi-parahoric levels \ref{eq: change of parahoric, interior} extends to
\begin{equation}\label{eq: change of parahoric, compactification}
\begin{tikzcd}
	{\Shumc{K_1, \Sigma_1}(G, X)} & {\Shumc{K_1^{\dd}, \Sigma_1^{\dd}}(G^{\dd}, X^{\dd})} & {\Shumm{K_1}(G, X)} & {\Shumm{K_1^{\dd}}(G^{\dd}, X^{\dd})} \\
	{\Shumc{K_2, \Sigma_2}(G, X)} & {\Shumc{K_2^{\dd}, \Sigma_2^{\dd}}(G^{\dd}, X^{\dd})} & {\Shumm{K_2}(G, X)} & {\Shumm{K_2^{\dd}}(G^{\dd}, X^{\dd})}
	\arrow[from=1-1, to=1-2]
	\arrow[from=1-1, to=2-1]
	\arrow[from=1-2, to=2-2]
	\arrow[from=1-3, to=1-4]
	\arrow[from=1-3, to=2-3]
	\arrow[from=1-4, to=2-4]
	\arrow[from=2-1, to=2-2]
	\arrow[from=2-3, to=2-4]
\end{tikzcd}
\end{equation}
combine with Remark \ref{rmk: remaining of the proof, change of parahoric}, we finish the proof of Theorem \ref{Theorem: change of parahoric}.

\subsection{Change-of-Parahoric: Propositions}

Keep notations from last subsection. We give a more detailed examination of the morphisms $\pi^{\tor}: \Shumc{K_1, \Sigma_1} \to \Shumc{K_2, \Sigma_2}$ and $\pi^{\min}: \Shumm{K_1} \to \Shumm{K_2}$. In particular, we show that $\pi^{\tor}$ and $\pi^{\min}$ satisfy similar propositions listed in \cite[Proposition 3.4]{lan2022closed}, we follow the list there and give the proof:

\begin{proposition}\label{Proposition: change of parahoric, prop}
    The morphisms $\pi^{\tor}$, $\pi^{\min}$ in Theorem \ref{Theorem: change of parahoric} satisfy the following properties:
    \begin{enumerate}
        \item For each stratum $\mathcal{Z}_1([\Phi_1])$ of $\Shumm{K_1}$, there exists a (unique) stratum $\mathcal{Z}_2([\Phi_2])$ of $\Shumm{K_2}$ such that $\pi^{\min}(\mathcal{Z}_1([\Phi_1])) \subset \mathcal{Z}_2([\Phi_2])$ (as subsets). Moreover, $\mathcal{Z}_1([\Phi_1])$ is open and closed in $\pi^{\min, -1}(\mathcal{Z}_2([\Phi_2]))$. $\pi^{\min}$ induces a proper surjective morphism $\mathcal{Z}_1([\Phi_1]) \to \mathcal{Z}_2([\Phi_2])$.
        \item We have proper surjective morphisms $\Zb^{\bigsur}_1(\Phi_1) \to \Zb^{\bigsur}_2(\Phi_2)$ that are compatible with $\mathcal{Z}_1([\Phi_1]) \to \mathcal{Z}_2([\Phi_2])$ under the canonical isomorphisms $\Zb_i(\Phi_i) \cong \mathcal{Z}_i([\Phi_i])$ and projections $\Zb^{\bigsur}_i(\Phi_i) \to \Zb_i(\Phi_i)$ for $i = 1, 2$. Moreover, $\Zb^{\bigsur}_1(\Phi_1) \to \Zb^{\bigsur}_2(\Phi_2)$ is the change-of-parahoric morphism bewteen Pappas-Rapoport integral models of Hodge-type Shimura varieties with quasi-parahoric level structures.
        \item Over any $\Zb^{\bigsur}_1(\Phi_1) \to \Zb^{\bigsur}_2(\Phi_2)$ as above, we have a proper surjective morphism $C_1(\Phi_1) \to C_2(\Phi_2)$ equivariant with a morphism $\ab_1(\Phi_1) \to \ab_2(\Phi_2)$, where for $i = 1, 2$, $\ab_i(\Phi_i)$ is an abelian scheme over $\Zb^{\bigsur}_i(\Phi_i)$, $C_i(\Phi_i) \to \Zb^{\bigsur}_i(\Phi_i)$ is an abelian scheme torsor under $\ab_i(\Phi_i)$. Moreover, $C_1(\Phi_1) \to C_2(\Phi_2)$ induces a finite morphism $C_1(\Phi_1) \to C_2(\Phi_2) \times_{\Zb^{\bigsur}_2(\Phi_2)} \Zb^{\bigsur}_1(\Phi_1)$.
        \item Over any $C_1(\Phi_1) \to C_2(\Phi_2)$ as above, we have a proper surjective morphism $\Xi_1(\Phi_1) \to \Xi_2(\Phi_2)$, which induces a finite morphism $\Xi_1(\Phi_1) \to \Xi_2(\Phi_2) \times_{C_2(\Phi_2)} C_1(\Phi_1)$ which is equivalent with the pullback of a group homomorphism of tori $\mathbf{E}_{K_1}(\Phi_1) \to \mathbf{E}_{K_2}(\Phi_2)$ with finite kernel over $\Spec(\Z)$ that is dual to a homomorphism $\mathbf{S}_{K_2}(\Phi_2) \to \mathbf{S}_{K_1}(\Phi_1)$ of character groups with finite cokernel. The $\R$-dual of this homomorphism is an injective homomorphism $\mathbf{B}_{K_1}(\Phi_1)_{\R} \hookrightarrow \mathbf{B}_{K_2}(\Phi_2)_{\R}$, inducing morphisms $\mathbf{P}_{K_1}(\Phi_1) \hookrightarrow \mathbf{P}_{K_2}(\Phi_2)$, and $\mathbf{P}_{K_1}(\Phi_1)^+ \hookrightarrow \mathbf{P}_{K_2}(\Phi_2)^+$, here $\mathbf{P}_{K_i}(\Phi_i)$ is a subset of $\mathbf{B}_{K_i}(\Phi_i)_{\R}$, $\mathbf{P}_{K_i}(\Phi_i)$ is the union of the interior $\mathbf{P}_{K_i}(\Phi_i)^+$ of a homogeneous self-adjoint cone, $\Sigma_i(\Phi_i)$ (resp. $\Sigma_i(\Phi_i)^+$) is a cone decomposition of $\mathbf{P}_{K_i}(\Phi_i)$ (resp. $\mathbf{P}_{K_i}(\Phi_i)^+$).
        \item If the image $\sigma_1 \in \Sigma_1(\Phi_1)$ under $\mathbf{P}_{K_1}(\Phi_1) \hookrightarrow \mathbf{P}_{K_2}(\Phi_2)$ is contained in some $\sigma_2 \in \Sigma_2(\Phi_2)$, then we have a canonical morphism $\Xi_1(\Phi_1)(\sigma_1) \to \Xi_2(\Phi_2)(\sigma_2)$ extending $\Xi_1(\Phi_1) \to \Xi_2(\Phi_2)$, and it induces a morphism $\Xi_{1, \sigma_1}(\Phi_1) \to \Xi_{2, \sigma_2}(\Phi_2)$.
        \item Suppose that $\Sigma_1$ is induced by $\Sigma_2$ from now on, if the image $\sigma_1 \in \Sigma_1(\Phi_1)$ is contained in some $\sigma_2 \in \Sigma_2(\Phi_2)$, then $\pi^{\tor}$ induces a morphism $\mathcal{Z}_1([\Phi_1, \sigma_1]) \to \mathcal{Z}_2([\Phi_2, \sigma_2])$ that can be identified with $\Xi_{1, \sigma_1}(\Phi_1) \to \Xi_{2, \sigma_2}(\Phi_2)$ above, under the canonical identifications $\mathcal{Z}_i([\Phi_i, \sigma_i]) \cong \Xi_{i, \sigma_i}(\Phi_i)$.
        \item There is a proper morphism $\ovl{\Xi}_1(\Phi_1) \to \ovl{\Xi}_2(\Phi_2)$, whose formal completions give a proper morphism $\mathfrak{X}_1(\Phi_1) \to \mathfrak{X}_2(\Phi_2)$. If the image $\sigma_1 \in \Sigma_1(\Phi_1)$ is contained in some $\sigma_2 \in \Sigma_2(\Phi_2)$, then we have induced morphism $\mathfrak{X}_1(\Phi_1)_{\sigma_1}^{\circ} \to \mathfrak{X}_2(\Phi_2)_{\sigma_2}^{\circ}$, which can be identified with 
        \[ (\Shumc{K_1, \Sigma_1})_{\underset{\tau_1\in \Sigma_1(\Phi_1)^{+},\ \bar{\tau}_1\subset\bar{\sigma}_1}{\cup} \mathcal{Z}_1([\Phi_1, \tau_1])}^{\wedge}  \to (\Shumc{K_2, \Sigma_2})_{\underset{\tau_2\in \Sigma_2(\Phi_2)^{+},\ \bar{\tau}_2\subset\bar{\sigma}_2}{\cup} \mathcal{Z}_2([\Phi_2, \tau_2])}^{\wedge}  \]
        induced by $\pi^{\tor}$ under the canonical identifications \ref{eq: mathfrak X, identification}.
        \item Moreover, for each affine open formal subscheme $\mathfrak{W}_2 = \Spf(R_2) \subset \mathfrak{X}_2(\Phi_2)_{\sigma_2}^{\circ}$, let $\mathfrak{W}_1 = \Spf(R_1) \subset \mathfrak{X}_1(\Phi_1)_{\sigma_1}^{\circ}$ be an affine open formal subscheme contained in the preimage of $\mathfrak{W}_2$. Under the induced morphism $W_2 = \Spec R_2 \to \Shumc{K_2, \Sigma_2}$, $W_2 \to \Xi_2(\Phi_2)(\sigma_2)$, $W_1 = \Spec R_1 \to \Shumc{K_1, \Sigma_1}$, $W_1 \to \Xi_1(\Phi_1)(\sigma_1)$, the preimages of $\Shum{K_i}$ and $\Xi_i(\Phi_i)$ coincides as an open subscheme $W^0_i \subset W_i$, and the preimage of $W^0_2$ is $W^0_1$.
        
    \end{enumerate}
\end{proposition}
\begin{proof}
    Over generic fiber $\Spec E$, all the statements are true due to \cite[Proposition 3.4]{lan2022closed}, here we use \cite[Assumption 2.1, case (1)]{lan2022closed}. In particular, we can find $\Sigma_1$ induced by $\Sigma_2$ due to \cite[Proposition 4.10]{lan2022closed}, note that this is a result on generic fiber. 
    
    Also, given noetherian reduced schemes $X, Y$ that are flat over a base $S$, here $S$ is noetherian and flat over $\OO_E$, if a proper morphism $f: X \to Y$ is surjective over the generic fiber $S_{\eta}$ of $S$, then $f$ is surjective and its scheme theoratical image is $Y$.

    $(1)$ Apply \cite[Lemma A.3.4]{pera2019toroidal}, then the morphism $\mathcal{Z}_1([\Phi_1])_{\eta} \to \mathcal{Z}_2([\Phi_2])_{\eta}$ extends to $\mathcal{Z}_1([\Phi_1]) \to \mathcal{Z}_2([\Phi_2])$, and $\mathcal{Z}_2([\Phi_2])$ is the unique stratum such that $\pi^{\min}(\mathcal{Z}_1([\Phi_1])) \subset \mathcal{Z}_2([\Phi_2])$. Since $\mathcal{Z}_1([\Phi_1])_{\eta}$ is open and closed in $\pi^{\min, -1}(\mathcal{Z}_2([\Phi_2])_{\eta})$, and the closure relations of various $\lrbracket{\mathcal{Z}_1([\Phi_1])}_{[\Phi_1] \in \Cusp_{K_1}}$ are determined by $\lrbracket{\mathcal{Z}_1([\Phi_1])_{\eta}}_{[\Phi_1] \in \Cusp_{K_1}}$, then $\mathcal{Z}_1([\Phi_1])$, as the closure of $\mathcal{Z}_1([\Phi_1])_{\eta}$ in $\pi^{\min, -1}(\mathcal{Z}_2([\Phi_2]))$, is also open and closed in $\pi^{\min, -1}(\mathcal{Z}_2([\Phi_2]))$. Since $\pi^{\min}$ is proper, then $\pi^{\min}: \mathcal{Z}_1([\Phi_1]) \to \mathcal{Z}_2([\Phi_2])$ is proper.

    $(2)$ Note that $\pi_{\eta}$ induces $\Zb_{1, \eta}^{\bigsur} \to \Zb_{2, \eta}^{\bigsur}$ which is the projection $\shu{K_{1, \Phi_1, h}}(G_{\Phi_1, h}, X_{\Phi_1, h}) \to \shu{K_{2, \Phi_2, h}}(G_{\Phi_2, h}, X_{\Phi_2, h})$, where by definition $G_{\Phi_1, h} = G_{\Phi_2, h}$, $X_{\Phi_1, h} = X_{\Phi_2, h}$, $K_{1, \Phi_1, h} \subset K_{2, \Phi_2, h}$ are quasi-parahoric subgroups due to part $(3)$ of Theorem \ref{thm: boundary of kp is again kp}, and $\Zb_{i}^{\bigsur}$ are the Pappas-Rapoport integral models due to part $(5)$ of Theorem \ref{thm: boundary of kp is again kp}. Due to \cite[Corollary 4.3.2]{pappas2024p} and \cite[Corollary 4.1.10]{daniels2024conjecture}, $\Zb_{1, \eta}^{\bigsur} \to \Zb_{2, \eta}^{\bigsur}$ extends to $\Zb_{1}^{\bigsur} \to \Zb_{2}^{\bigsur}$. Due to Remark \ref{rmk: remaining of the proof, change of parahoric}, such extension is proper, and is unique and compatible with $\mathcal{Z}_1([\Phi_1]) \to \mathcal{Z}_2([\Phi_2])$.

    $(3)$ Let $\ab_1(\Phi_1)'\to \Zb^{\bigsur}_1(\Phi_1)$ be the pullback of $\ab_2(\Phi_2) \to \Zb^{\bigsur}_2(\Phi_2)$ along $\Zb^{\bigsur}_1(\Phi_1) \to \Zb^{\bigsur}_2(\Phi_2)$, then $\ab_1(\Phi_1)'$ and $\ab_1(\Phi_1)$ are abelian schemes over the normal scheme $\Zb^{\bigsur}_1(\Phi_1)$, due to Corollary \ref{cor: pera, 4.1.5}. Since there exists a homomorphism $\ab_1(\Phi_1)_{\eta} \to \ab_1(\Phi_1)'_{\eta}$ over $\Zb^{\bigsur}_1(\Phi_1)_{\eta}$, and $\Zb^{\bigsur}_1(\Phi_1)$ is flat and finite type over $\Spec \OO_{E(v)}$, then due to Weil extension property (or use \cite[Chatper I, Proposition 2.7]{faltings2013degeneration}), such homomorphism extends uniquely to a morphism $\ab_1(\Phi_1) \to \ab_1(\Phi_1)'$ over $\Zb^{\bigsur}_1(\Phi_1)$, thus we have a morphism $\ab_1(\Phi_1) \to \ab_2(\Phi_2)$. Note that $\ab_1(\Phi_1)_{\eta} \to \ab_1(\Phi_1)'_{\eta}$ is an isogeny, thus the proper morphism $\ab_1(\Phi_1) \to \ab_1(\Phi_1)'$ is surjective and has scheme theoretical image $\ab_1(\Phi_1)'$, thus its kernel $\mathcal{A}$ is a flat group scheme over $\Zb^{\bigsur}_1(\Phi_1)$ due to \cite[Theorem(A)]{achter2024images}. Since $\mathcal{A}_{\eta}$ is quasi-finite over $\Zb^{\bigsur}_1(\Phi_1)_{\eta}$, thus $\mathcal{A}$ is quasi-finite over $\Zb^{\bigsur}_1(\Phi_1)$, then the proper morphism $\ab_1(\Phi_1) \to \ab_1(\Phi_1)'$ is quasi-finite, thus is a finite morphism.

    Recall that we have the following diagram coming from diagrams \ref{eq: toroidal, intermediate} and \ref{eq: change of parahoric, compactification}:
\[\begin{tikzcd}
	{C_1(\Phi_1)_{\eta}} & {\prod_{j \in J^1} C_j^{\dd}(\Phi_j^{\dd})_{\eta}} \\
	{C_2(\Phi_2)_{\eta}} & {\prod_{j \in J^2}C_j^{\dd}(\Phi_j^{\dd})_{\eta}}
	\arrow[from=1-1, to=1-2]
	\arrow[from=1-1, to=2-1]
	\arrow[from=1-2, to=2-2]
	\arrow[from=2-1, to=2-2]
\end{tikzcd}\]
     where $\Phi_j^{\dd} \in \Cusp_{H_j^{\dd}}(G^{\dd}_j, X^{\dd}_j)$ is induced by $\Phi_1^{\dd} \in \Cusp_{K_1^{\dd}}(G^{\dd}, X^{\dd})$ (by $\Phi_1 \in \Cusp_{K_1}(G, X)$). Since $\prod_{j \in J^1} C_j^{\dd}(\Phi_j^{\dd})_{\eta} \to \prod_{j \in J^2}C_j^{\dd}(\Phi_j^{\dd})_{\eta}$ naturally extends to $\prod_{j \in J^1} C_j^{\dd}(\Phi_j^{\dd}) \to \prod_{j \in J^2}C_j^{\dd}(\Phi_j^{\dd})$, and $C_i(\Phi_i)$ is the relative normalization of  $\prod_{j \in J^i} C_j^{\dd}(\Phi_j^{\dd})$ in $C_i(\Phi_i)_{\eta}$ due to \cite[Proposition 8.4]{lan2016compactifications} and Lemma \ref{lem: two constructions equal}, therefore we have a morphism $C_1(\Phi_1) \to C_2(\Phi_2)$ extending $C_1(\Phi_1)_{\eta} \to C_2(\Phi_2)_{\eta}$ by the functoriality of relative normalizations. Due to Remark \ref{rmk: remaining of the proof, change of parahoric}, such an extension is unique, and we have following commutative diagram:
\[\begin{tikzcd}
	{\ab_1(\Phi_1) \times C_1(\Phi_1)} & {C_1(\Phi_1)} \\
	{\ab_2(\Phi_2) \times C_2(\Phi_2)} & {C_2(\Phi_2)}
	\arrow[from=1-1, to=1-2]
	\arrow[from=1-1, to=2-1]
	\arrow[from=1-2, to=2-2]
	\arrow[from=2-1, to=2-2]
\end{tikzcd}\]

  The induced morphism $C_1(\Phi_1) \to C_2(\Phi_2) \times_{\Zb^{\bigsur}_2(\Phi_2)} \Zb^{\bigsur}_1(\Phi_1)$ is \'etale locally $\ab_1(\Phi_1) \to \ab_1(\Phi_1)'$ over $\Zb^{\bigsur}_1(\Phi_1)$, thus is finite.
  
    $(4)$ Similar as $(3)$, note that $\mathbf{E}_{K_1}(\Phi_1) \to \mathbf{E}_{K_2}(\Phi_2)$ is a morphism of constant split tori over, $\Z$ we have a proper morphism $\Xi_1(\Phi_1) \to \Xi_2(\Phi_2)$ over $C_1(\Phi_1) \to C_2(\Phi_2)$ which induces a finite morphism $\Xi_1(\Phi_1) \to \Xi_2(\Phi_2) \times_{C_2(\Phi_2)} C_1(\Phi_1)$ which is equivariant with the pullback of a homomorphism $\mathbf{E}_{K_1}(\Phi_1) \to \mathbf{E}_{K_2}(\Phi_2)$. The rest of the statement $(3)$ are results on generic fiber.

    $(5)$ This is a formal consequence of $(4)$ from the definitions of $\Xi(\Phi)(\sigma)$ and $\Xi(\Phi)_{\sigma}$, see \cite[Proposition 3.4(2)(3)]{lan2022closed} for details.

    $(6)$ This is true over $\Spec E$, we apply \cite[Lemma A.3.4]{pera2019toroidal} as in the proof of part $(1)$. The compatibility of $\mathcal{Z}_1([\Phi_1, \sigma_1]) \to \mathcal{Z}_2([\Phi_2, \sigma_2])$ and $\Xi_{1, \sigma_1}(\Phi_1) \to \Xi_{2, \sigma_2}(\Phi_2)$ is true over $\Spec E$ thus is true over $\Spec \OO_{E(v)}$ due to Remark \ref{rmk: remaining of the proof, change of parahoric}.

    $(7)$ Use $(5)$, $(6)$ and Remark \ref{rmk: remaining of the proof, change of parahoric}.

    $(8)$ Fix $i = 1, 2$, then the results concerning $W_i^0 \subset W_i$ are due to \cite[Proposition 3.1(8)]{lan2022closed}. The preimage of $W^0_2$ is $W^0_1$ since the preimage of $\Shum{K_2}$ under $\Shumc{K_1, \Sigma_1} \to \Shumc{K_2, \Sigma_2}$ is $\Shum{K_1}$.
\end{proof}


\bibliographystyle{abbrv}
\bibliography{boundary}

\end{document}